\documentclass[12pt,twoside,openright,titlepage]{report}
\usepackage{frontespizio}
\usepackage[utf8]{inputenc}
\usepackage{amsmath,amssymb, mathrsfs}
\usepackage{amsthm}
\usepackage{bbm}
\usepackage{lipsum}

\usepackage[absolute]{textpos}
\usepackage{mathtools,emptypage}
\usepackage[inner=3cm,outer=3cm]{geometry} 
\usepackage{graphicx}
\graphicspath{{images/}}
\usepackage{marvosym}
\usepackage{wasysym}%
\newcommand{\arian}{\textbf{\currency}}

\usepackage{bbold}
\usepackage{centernot}
\usepackage{cite}
\usepackage[dvipsnames]{xcolor}

\usepackage{tikz}
\usetikzlibrary{arrows}
\usetikzlibrary{patterns}
\usetikzlibrary{shapes}
\usepackage{mathrsfs}
\usepackage{bbm}
\usepackage{enumerate}
\usepackage[shortlabels]{enumitem}
\usepackage{tikz}
\usepackage{tikzcd}

\usepackage{amssymb}
\usepackage{tikz-cd}
\usetikzlibrary{cd}

\usepackage{appendix}
\usepackage{chngcntr}
\usepackage{etoolbox}
\AtBeginEnvironment{subappendices}{%
\chapter*{Appendices}
\addcontentsline{toc}{chapter}{Appendices}
\counterwithin{figure}{section}
\counterwithin{table}{section}
}

\AtEndEnvironment{subappendices}{%
\counterwithout{figure}{section}
\counterwithout{table}{section}
}

\newcommand{\R}{\mathbb{R}}
\newcommand{\N}{\mathbb{N}}
\newcommand{\Z}{\mathbb{Z}}
\newcommand{\C}{\mathbb{C}}
\newcommand{\de}{\partial}
\renewcommand{\-}{\smallsetminus}
\newcommand{\D}{\mathbb{D}}
\newcommand{\B}{\mathbb{B}}

\renewcommand{\a}{\alpha}

\newcommand{\f}{\varphi}
\newcommand{\e}{\varepsilon}
\newcommand{\w}{{\omega}}

\usepackage{pifont}

\newcommand{\transv}{\mathrel{\text{\tpitchfork}}}
\makeatletter
\newcommand{\tpitchfork}{%
  \raise-0.1ex\vbox{
    \baselineskip\z@skip
    \lineskip-.52ex
    \lineskiplimit\maxdimen
    \m@th
    \ialign{##\crcr\hidewidth\smash{$-$}\hidewidth\crcr$\pitchfork$\crcr}
  }%
}
\makeatother

\newcommand{\en}{\E\#_{X\in W}}
\newcommand{\vol}{\mathrm{vol}}
\newcommand{\ci}{\mathcal{C}^1}
\newcommand{\RP}{\mathbb{R}\mathrm{P}}
\newcommand{\PP}{\mathbb{P}}

\newcommand{\Prob}{(\Omega, \mathfrak{S},\P)}
\renewcommand{\P}{\mathbb{P}}
\newcommand{\E}{\mathbb{E}}
\newcommand{\nrw}{\Rightarrow}
\newcommand{\spt}{\text{supp}}

\newcommand{\g}[3]{\mathcal{G}^{#1}(#2,\R^{#3})}

\newcommand{\goo}[3]{\mathcal{G}^{#1}(#2|#3)}

\newcommand{\coo}[3]{\mathcal{C}^{#1}(#2,#3)}
\newcommand{\Cr}[3]{\mathcal{C}^{#1}(#2,\R^{#3})}
\newcommand{\mC}{\mathcal{C}}
\newcommand{\mS}{\mathcal{S}}

\newcommand{\Gam}{\Gamma(X,W)}
\newcommand{\nrm}[3]{\|{#1}\|_{#2,#3}}
\newcommand{\G}{\mathscr{G}}
\newcommand{\K}{\mathcal{K}}

\newcommand{\HPK}{Let $(O_p,N_0,\phi,g)$ be a KROK model for a KROK couple $(X,W)$}
\newcommand{\HP}{Let $(X,W)$ be KROK}
\newcommand{\MW}{\mathcal{M}}
\newcommand{\XW}{\mathcal{X}}
\newcommand{\NW}{\mathcal{N}}
\newcommand{\DElta}{\mathcal{W}}

\usepackage{hyperref}
\hypersetup{
	colorlinks=true,       
	linkcolor=BlueViolet,         
	citecolor=BlueViolet
}

 \usepackage[usenames,dvipsnames]{pstricks}
 \usepackage{epsfig}
 \usepackage{pst-grad} 
 \usepackage{pst-plot} 

\newtheorem{thm}{Theorem}
\newtheorem{lemma}[thm]{Lemma}
\newtheorem{cor}[thm]{Corollary}
\newtheorem{prop}[thm]{Proposition}
\newtheorem{claim}[thm]{Claim}
\newtheorem*{stella}{($*$)}
\newtheorem{innercustomthm}{Theorem}
\newenvironment{customthm}[1]
  {\renewcommand\theinnercustomthm{#1}\innercustomthm}
  {\endinnercustomthm}

\theoremstyle{definition}

\newtheorem{defi}[thm]{Definition}
\newtheorem{remark}[thm]{Remark}

\newtheorem{example}[thm]{Example}

\newcommand{\be}{\begin{equation}}
\newcommand{\ee}{\end{equation}}

\newcommand{\bega}{\begin{equation}\begin{aligned}}
\newcommand{\eega}{\end{aligned}\end{equation}}

\usepackage{mathtools} 
\mathtoolsset{showonlyrefs,showmanualtags} 

\numberwithin{equation}{section}

\title{
{Random Differential Topology}\\
{\large SISSA}\\
\vspace{50 pt}
{\includegraphics[scale=1]{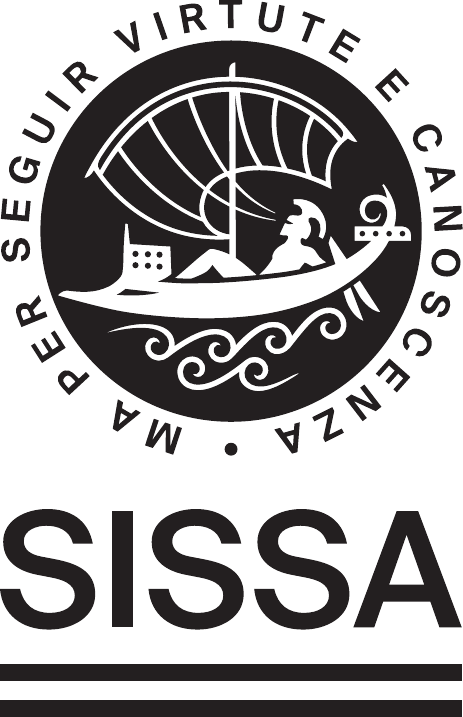}}
}
\author{Michele Stecconi}
\date{1/9/2020}

\begin{document}

\newenvironment{changemargin}[2]{
\begin{list}{}{
\setlength{\topsep}{0pt}
\setlength{\leftmargin}{#1}
\setlength{\rightmargin}{#2}
\setlength{\listparindent}{\parindent}
\setlength{\itemindent}{\parindent}
\setlength{\parsep}{\parskip}
}
\item[]}{\end{list}}


\begin{titlepage}
\pagestyle{empty}
\newgeometry{left=2cm, right=2cm, top=5cm, bottom=0cm}
\begin{textblock*}{1.58in}[0,0](0mm,0.8mm)
    \includegraphics[width=1.58in]{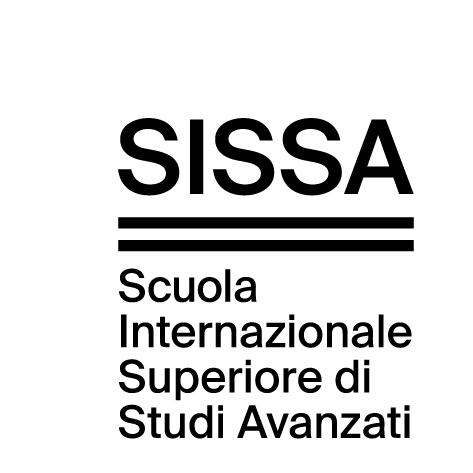}
\end{textblock*}

\centering
\LARGE
Mathematics Area - PhD course in\\
Geometry and Mathematical Physics

\vfill
\Huge
\textbf{Random\\ 
Differential Topology}
 
\vfill
\Large
\begin{changemargin}{0cm}{0cm}
\begin{minipage}{0.3\textwidth}
\begin{flushleft}
Candidate: \\
Michele Stecconi
\end{flushleft}
\hfill
\end{minipage}
\hfill
\begin{minipage}{0.6\textwidth}
\begin{flushright}
Advisors: \\
Antonio Lerario
\\
Andrei Agrachev
\\

\end{flushright}
\end{minipage}      
\end{changemargin}   
\vspace{1.5cm}
 
\Large
Academic Year 2019-20 
\vspace{0.3cm}
\begin{center}
  \includegraphics[width=1.18in]{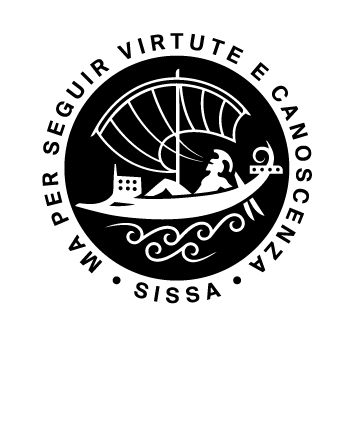}
\end{center}
\restoregeometry
\end{titlepage}


\chapter*{ }
To my Bimba\dots


\chapter*{Acknowledgements}
\vspace*{-\baselineskip}
I would like to start by expressing my deepest gratitude to my supervisors, not only for their mathematical guidance, but also for their inspiring enthusiasm and the undefinable wisdom that I learned to admire and hope to have partially absorbed.

Antonio Lerario: professor, coauthor,  psychotherapist and friend, he always supported and helped me in every aspect of my life at SISSA, going far beyond the institutional duties. Also, this thesis would probably be just a sequence of definitions, theorems and proofs if it weren't for his efforts to teach me how to properly write a paper. For all this, I will always be deeply indebted to him.

Andrei Agrachev, whose  padronance of Mathematics is a rare example that I will never forget.  I am grateful to him for believing in me in these years and for giving me the opportunity to work together on the ``Floer'' project which, I am sorry to say, has nothing to do with this thesis.
\par 
I also wish to thank the people who have shown interest in my work or have offered me their point of view, among them Mattia Sensi, Giulio Ruzza, Paul Breiding, Riccardo Ghiloni, Mikhail Sodin.

A special mention is due to Riccardo Tione, Nico Samà, Maria Strazzullo and Léo Mathis, who surely belong to the previous cathegory, but also to that of those adventure companions that remain in your life forever. Thanks to you, to Marino and to Slobodan.

I would like to extend my gratitude to SISSA, to the school itself for being what it is, to the infinite patience of the administrative staff like Riccardo Iancer, Marco Marin and Emanuele Tuillier Illingworth and to the kindness of the people who work in the canteen: Silvia, Eleonora, Patrizia, Lucia, Lucio, etc. Finally, a special thank to the people who has been part of my daily life in and out the school: Alessandro Carotenuto, Guido Mazzuca, Enrico Schlitzer, Emanuele Caputo, Stefano Baranzini, Mario De Marco, Luca Franzoi, Matteo Zancanaro, Federico Pichi and many others.

Last, but by no means least, thanks to my (extended) family: Mamma, Papi, Ceci, Cami, Zia Alma, Monica, Francesca, Tommi and my beloved Arian, my favourite person on earth. Because without them everything would be much harder. 

\tableofcontents

\bibliographystyle{plain}

\chapter{Introduction}
\section{Motivations}
This manuscript collects three independent works: \emph{Differential Topology of Gaussian Random Fields} (Antonio Lerario and M.S. \cite{dtgrf}), \emph{Kac Rice formula for Transversal Intersections} (M.S. \cite{krstec}) and  \emph{Maximal and Typical Topology of real polynomial Singularities} (Antonio Lerario and M.S. \cite{mttrp}), together with some additional results, observations, examples and comments, some of which were taken up in the subsequent work \emph{What is the degree of a smooth hypersurface?} (Antonio Lerario and M.S. \cite{witdoash}). 

The common thread of these works is the study of topological and geometric properties of random smooth maps, 
specifically we are interested in the asymptotical behaviour of things when the random map depends continuously, in some sense to be specified, on a parameter.

\subsection{Limit probabilities}\label{moti:limitprob}
The first situation of interest is the following. Let $X_d\colon M\to \mathbb{R}^k$ be a sequence of smooth Gaussian Random Fields (see Definition \ref{def:GRF}) defined on a smooth compact\footnote{We suppose that $M$ is compact here, in order to simplify the exposition, although we will hardly make compactness assumption in the next chapters.} manifold and assume we want to show that $X_d$ satisfies a given condition with positive probability, for all $d$ big enough. Technically this translates into proving that
\be\label{moti:P}
\liminf_{d\to +\infty}\P\{X_d\in U\}>0,
\ee
for some subset $U\subset \Cr rMk$.
In this direction, in \cite{GaWe1} Gayet and Welschinger proved \eqref{moti:P} in the case when $X_d$ is the Kostlan polynomial of degree $d$ on the $m$-sphere and $U=U_d$ is the set of all maps whose zero set has a connected component contained in a disk of radius $d^{-\frac12}$ that is isotopic to a given compact hypersurface $\Sigma\subset \R^m$. A similar result, due to Lerario and Lundberg \cite{Lerariolemniscate}, states that every nesting of the ovals composing a random leminiscate\footnote{A leminiscate of degree $d$ is an algebraic curve on the Riemann's sphere $\C\P^1$, defined by an equation of the form $\left|\frac{p(x,y)}{q(x,y)}\right|=1$, where $p,q$ are complex homogeneous polynomials of the same degree $d$.} is realized with positive probability in any given spherical disk of radius $d^{-\frac12}$, for large enough degree.

We propose a method to investigate  \eqref{moti:P} using the tools and the framework developed in \cite{dtgrf} that works well in the most common situations, such as the previous two examples, when there is a convergence in law $X_d\nrw X$. 
The method consists of three basic steps:
\begin{enumerate}[i), wide]
\item\label{itm:conv} Establish the convergence in law $X_d\nrw X$, corresponding to the narrow convergence of the probability measures induced on the space $\Cr \infty Mk$.
This problem is reduced to a simpler deterministic one, by means of the following theorem (Theorem \ref{thm:1} of Chapter \ref{chap:dtgrf}).
\begin{thm}\label{moti:thm:koconv}
A sequence $X_d\colon M\to \R^k$ of smooth gaussian random fields converges in law to $X$, if and only if the corresponding sequence of covariance functions converges in the $\mC^\infty$ topology.
    \end{thm}
    This notion of convergence suits perfectly our context since it means that for every Borel subset $U\subset \Cr \infty Mk$, we have that
    \be\label{moti:port}
    \liminf_{d\to +\infty}\P\{X_d\in U\}\ge \P\{X\in \text{int}(U)\}.
    \ee
    \item Show that $\P\{X\in \text{int}(U)\}> 0$ by studying the support\footnote{The support of a $\mC^r$ random field is the smallest closed subset in $\coo rM{\R^k}$ having probability one.} of the limit, for instance via the next theorem which is a combination of Theorem \ref{thm:dtgrf:3} and Theorem \ref{thm:sputo}. 
\begin{thm}\label{moti:thm:support}
Let $X=\sum_{n\in \N}\xi_n f_n$ be a $\mC^r$ gaussian function, where $f_n\colon M\to \R$ are $\mC^r$ functions and $\xi_n$ is a sequence of independent normal gaussians. Then
\be 
\spt(X)=\overline{\text{span}\{f_n\}_{n\in\N}}^{\mC^r}.
\ee
In particular, given a dense sequence of points $p_n$ in $M$, then $X$ admits a series expansion of the above kind, with $f_n=K_X(p_n,\cdot)$, where $K_X$ is the covariance function of $X$.
\end{thm}  
Notice that the only object involved in Theorem \ref{moti:thm:support} is the limit law $X$.
    \item In some cases \eqref{moti:P} can be improved to an equality
    \be\label{moti:Pexact}
    \lim_{d\to +\infty}\P\{X_d\in U\}=\P\{X\in U\}.
    \ee
    This, again, can be checked by looking only at the limit law, indeed \eqref{moti:Pexact} holds if $\P\{X\in\de U\}=0$. Most often the elements of the set $\de U$ are nongeneric functions having certain kind of singularities, thus \eqref{moti:Pexact} can be deduced by an application of the  result below (a simple consequence of Theorem \ref{transthm2} of Chapter \ref{chap:dtgrf}). It is a probabilistic version of Thom's jet Transversality Theorem and states that, for any given notion of singularity, the probability that a gaussian field $X$ has one or more such singular point is zero.
    \begin{thm}\label{intro:thm:transv}
Let $X\colon M\to \mathbb{R}^k$ be a smooth gaussian random field with full support\footnote{A random map has full support $\spt(X)=\coo rM{\R^k}$ when it satisfies any open (with respect to the weak Whitney topology) condition with positive probability.} and $W\subset J^r(M,\mathbb{R}^k)$ a submanifold of the space of $r-$jets. Then $X$ is transversal to $W$ almost surely.
\end{thm}
\end{enumerate}
 \subsection{Expectation of Betti numbers}\label{moti:expectabet}
Let us consider another object of investigation: let $Z_d\subset M$ be a random subset and we want to estimate the asymptotic behavior of the expected topological complexity of $Z_d$. Mathematically, we quantify the latter concept as the sum $b(Z_d)$ of all of its Betti numbers, hence the question is
\be\label{moti:Eb}
\exists \lim_{d\to +\infty}\E\{b(Z_d)\}=?
\ee
A much studied example is that of a random real projective variety $Z_d\subset \R\P^m$ defined by Kostlan polynomials of degree $d$ (see  \cite{EdelmanKostlan95}, \cite{shsm}, \cite{GaWe2,GaWe3}, \cite{buerg:07} for instance), followed by the nodal set of a gaussian combination of eigenfunctions of the Laplacian on a riemannian manifold (in this case the parameter $d$ is substituted by the set of the corresponding eigenvalues). 
Other examples are the level set of a random function, the set of its critical points, the preimage of a submanifold $W\subset \R^n$, etc.

More generally, we will consider $Z_d$ to be the set of points $p\in M$ such that the $r$-jet of a random field $X_d\colon M\to \R^k$, denoted by $j_p^rX_d$, belongs to some given subset of the jet space $W\subset J^r(M,\R^k)$. In other words, we assume that $Z_d$ is defined by a given set of conditions involving the random map $X_d\colon M\to \R^k$ and a finite number of its derivatives. This definition includes almost any random set of differential geometric nature: we will call it a \emph{singularity}.

Here, the following improved version of Theorem \ref{intro:thm:transv} (Theorem \ref{transthm2}, Chapter \ref{chap:dtgrf}) allows to make a preliminary observation to better understand the question in  \eqref{moti:Eb}, in that it implies that, unless $X_d$ is degenerate, the random set $Z_d=\{j^r_pX_d\in W\}$ is a submanifold of $M$ having the same codimension as $W$.
 \begin{thm}\label{intro:thm:transvjet}
Let $X\colon M\to \mathbb{R}^k$ be a smooth gaussian random field with nondegenerate $r$-jet, that is, for each $p\in M$, the gaussian random vector $j_p^rX$ is nondegenerate. Let $W\subset J^r(M,\mathbb{R}^k)$ be a submanifold of the space of $r$-jets. Then $X$ is transverse to $W$ almost surely.
\end{thm}

A first clue to answer \eqref{moti:Eb} is provided by the classical fact of Morse theory, stating that the Betti numbers of a manifold are bounded by the number of critical points of any Morse function. In Chapter \ref{chap:mttps} we will show that this criteria can be implemented to fit our setting, resulting in Theorem \ref{thm:strat}, of which we report a simplified version.
\begin{thm}\label{moti:thm:morse}
There exists another submanifold $\hat{W}\subset J^{r+1}(M,\R^k)$ such that
\be\label{moti:mors}
b\left(Z_d\right)\le \#(\hat{Z}_d),
\ee
where $\hat{Z}_d$ is the singularity corresponding to $\hat{W}$. Moreover $\hat{Z}_d$ is zero dimensional.
\end{thm}
The most remarkable aspect of Theorem \ref{thm:strat} is that it holds also if $W$ is not smooth, but semialgebraic, with the addition of a constant $C_W$ depending on $W$ in \eqref{moti:mors}. In virtue of this inequality the problem of estimating $b(Z_d)$ reduces to the case in which $Z_d$ is a set of points and thus the only positive Betti number is the cardinality $\#(Z_d)$. This is the topic of Chapter \ref{chap:kr}, where we prove the following.
\begin{thm}
If $X_d$ and $W$ are compatible then there exists a density $\delta_{j^rX_d\in W}$ on $M$ such that
\be \label{moti:KR}
\E\#(Z_d)=\int_{M}\delta_{j^rX_d\in W}.
\ee 
\end{thm}
The formula for $\delta_{j^rX_d\in W}$ is given in Theorem \ref{thm:introKRform}. It is a generalization of the celebreted Kac-Rice formula and it expresses the expected number of points of $Z_d$ as the integral on $M$ of a function depending only on the $r+2$ jet of the covariance function $K_{X_d}$.

The compatibility condition to which we refer is that the $r$-jet of $X$ and $W$ have to be a KROK (Kac Rice OK) couple (see Definition \ref{def:krok} in Chapter \ref{chap:kr}). This holds whenever a list of quite natural hypotheses is satisfied, however in the gaussian case, such list reduces significantly: it is required only that $W$ is \emph{projectable}, meaning that it should be transverse to the spaces $J^r_p(M,\R^k)=\{j_p^rf\colon f\in \coo rM{\R^k}\}$ for all $p\in M$, and that $X_d$ and $X$ have nondegenerate $r$-jet, in the same sense of Theorem \ref{intro:thm:transvjet}. 
Moreover, in presence of a convergence $X_d\nrw X$, as that of Theorem \ref{moti:thm:koconv}, we have the convergence of the expectations (Theorem \ref{thm:mainEgau}, Chapter \ref{chap:kr}).
\begin{thm}
Suppose that $W$ is a projectable submanifold of $J^r(M,\R^k)$ that is either compact or locally semialgebraic. Let $X_d$ and $X$ be gaussian random fields with nondegenerate $r$-jet. Assume that the covariance function of $X_d$ converges to that of $X$ in the $\mC^{\infty}$ topology, then
\be\label{moti:asureconvbettizero}
\lim_{d\to +\infty}\E\#(Z_d)=\E\#(Z).
\ee
\end{thm}
This, together with \eqref{moti:mors} and some notions of differential topology, allows to use the dominated convergence theorem\footnote{More precisely, one should use Fatou's lemma two times.} to deduce that the same property holds for the sum of the Betti numbers:
\be\label{moti:convbettitutti}
\lim_{d\to +\infty}\E b(Z_d)=\E b(Z).
\ee
\subsection{Application to Kostlan polynomials}
In Chapter \ref{chap:mttps} we apply all the above methods to study Kostlan random homogeneous polynomials, viewing them as Gaussian Random Fields $\psi_d\colon S^m\to\R^k$. 

It is well known that the restriction of $\psi_d$ to a spherical ball $B_d$ of radius $\sim d^{-\frac12}$ behaves like the Bargmann-Fock random field, as the degree $d$ grows. Thanks to the results of Chapter \ref{chap:dtgrf} we can give a detailed description of such convergence (see Theorem \ref{thm:Kostlan}). In this direction, Theorem \ref{moti:thm:koconv} makes the proof of the following fact really simple and easy to generalize. 

\begin{thm}\label{moti:thm:kostkoconv}
Let $U\subset \Cr \infty{\D^m}k$ be a Borel subset. Let $X_d=\psi_d|_{B_d}$ and let $X$ be the Bargmann-Fock random field, 
then
\be 
\P\{X \in \emph{int}(U)\}\le \liminf_{d\to +\infty}\P\{X_d\in U\}\le \limsup_{d\to +\infty}\P\{X_d\in U\}\le \P\left\{X\in \overline{U}\right\}.
\ee
\end{thm}

From this we deduce that  \eqref{moti:P} holds for every open subset $U\subset \Cr r{\D^m}k$, becoming the equality \eqref{moti:Pexact} in the most common cases. A consequence of this is that the results from \cite{GaWe1} and \cite{Lerariolemniscate} cited at the beginning of Section \ref{moti:limitprob} hold in the stronger form, with the equality \eqref{moti:Pexact}.

When dealing with well-behaved random maps such as Kostlan polynomials it is common practise to authomatically translate \emph{generic} properties of maps into \emph{almost sure} ones. This is correct in most cases and an application of Theorem \ref{intro:thm:transvjet} shows to what extent. In particular for any given \emph{singularity} $W$, the $r^{th}$ jet of the local limit field (on the standard disk $\D^m$) is transverse to $W$ almost surely. 
\begin{thm}
For any submanifold of $r$-jets $W$ we have
\be 
\P\{j^r X\transv W\}=1.
\ee 
This applies both when $W\subset J^r(S^m,\R^k)$ and $X=\psi_d$ with $d\ge r$, and when $W\subset J^r(\R^m,\R^k)$ and $X$ is the Bargmann-Fock random field.
\end{thm}

In the context of Kostlan polynomials, the generalized Kac Rice formula (Theorem \ref{thm:maingau}), developed in Chapter \ref{chap:kr}, is applicable to a large class of singularity conditions of interest (for example zeroes, critical points or degenerations of higher jets). The most remarkable aspect of this is that, for such singularities, we can argue as for \eqref{moti:convbettitutti} to prove convergence of all the expected Betti numbers:
\begin{thm}
Let $W\subset J^r(\D^m,\R^k)$ be a semialgebraic subset that is invariant under diffeomorphisms\footnote{Most natural differential geometric examples fall into this category. For example the set of zeroes, critical points or the degeneracy locus of higher derivatives.}
    \be 
    \exists \lim_{d\to +\infty}\E b_i\left(\{p\in B_d\colon j^r_p\psi_d\in W\}\right)=\E b_i\left(\{p\in \D^m\colon j^r_pX\in W\}\right),
    \ee
    where $X$ is the Bargmann-Fock random field.
\end{thm}   
We conclude this story with Theorem \ref{thm:bettiorderintro}, establishing a ``generalized square root law'' for the asymptotics of the expected Betti numbers of the whole singularity.
\begin{thm}[Generalized Square Root Law]
Let $W\subset J^{r}(S^m,\R^k)$ be semialgebraic and diffeomorphism invariant. Then
\be\label{moti:gensquarerootlaw}
\E b(j^r\psi_d^{-1}(W))=\Theta(d^{\frac{m}{2}}).
\ee
\end{thm}
Estimates of this kind reflect a phenomenon that has been observed by several authors (see \cite{EdelmanKostlan95, Kostlan:93, shsm, buerg:07, Po, GaWe1, GaWe2, GaWe3, FyLeLu}) in different contexts: random real algebraic geometry seems to behave as the ``square root'' of generic complex geometry. The estimates obtained in the latter context frequently provide a sharp upper bound for the real case (this is what happens, for instance, in Bézout Theorem). In this sense the name ``square root law'' is understood in relation with the following result, again from Chapter \ref{chap:mttps}. 
\begin{thm}
Let $W\subset J^r(S^m,\R^k)$ be a semialgebraic subset. Then, almost surely,
\be 
b(j^r\psi_d^{-1}(W))\le O(d^{m}).
\ee
\end{thm}
This theorem is interesting in that it provides a better estimate than the one that can be deduced from the classical Milnor-Thom bound \cite{MilnorBound}, which would be $O(d^{m+1})$.
\subsection{Deterministic estimates of Betti numbers.}
As we mentioned in Section \ref{moti:expectabet}, the sum of the Betti numbers $b(Z)$ can be considered a measure of the complexity of a topological space $Z$. 
Often the topological space under study is defined in an implicit manner and it is desirable to estimate its Betti numbers in terms of the data provided. 
%
This consideration motivates the search for estimates of the form:
\be \label{moti:eq:bettiexpression}
b(f^{-1}(W))\le \mathcal{E}(f), \qquad \text{(or $\ge$)}
\ee
where $\mathcal{E}(f)$ is intended as an expression involving data that can be inferred from the ``equation'' $f$ (and possibly from $W$) that defines $Z$.
Here we are focusing on the situation, similar to that of the previous sections, in which $Z=f^{-1}(W)$ is the preimage via a smooth map $f\colon M\to N$ of a submanifold $W\subset N$.

With this purpose in mind, we derive an original theorem of Differential Topology (Theorem \ref{semi:thm:main}, Chapter \ref{chap:semicontop}), saying that, under a $\mC^0$ small perturbation of a regular system of equations, the topology (meant as the sum of the Betti numbers) of the set of solutions cannot decrease. 

\begin{thm}[Semicontinuity of Betti numbers]\label{moti:thm:top}
Let $f\colon M\to N$ be a smooth map that is transverse to a given closed submanifold $W\subset N$. If $g\colon M\to N$ is another smooth map transverse to $W$ that is $\mathcal{C}^0$-close enough to $f$, then all the Betti numbers of $g^{-1}(W)$ are larger or equal than those of the initial manifold $f^{-1}(W)$.
\end{thm}
Theorem \ref{moti:thm:top} is a good option in those situations in which a good $\mC^1$ approximation is not available.
In fact, coupling it with the Holonomic Approximation Lemma  of \cite{eliash}, that is $\mC^0$ at most, we will use it to prove the lower bound implied in \eqref{moti:gensquarerootlaw} by showing the following.
\begin{thm}
For any \emph{formal} section $F\colon M\to J^r(M,\R^k)$, there exists a smooth function $\psi\colon M\to\R^k$ such that
\be
b\left(j^r\psi^{-1}(W)\right)\ge b\left(F^{-1}(W)\right).
\ee
\end{thm}

Moreover, Theorem \ref{moti:thm:top} can be of great use to produce quantitative bounds, in that it suggests that the expression $\mathcal{E}(f)$ appearing in \eqref{moti:eq:bettiexpression} should not be sensible to a change in $f$ due, for example, to the addition of an oscillating term, as long as the amplitude of the oscillation is small. This is an advantage with respect to more standard methods which require more rigidity, since they involve $\mC^1$ tools such as the Implicit Function Theorem or Thom's isotopy Lemma\footnote{Which says, roughly speaking, that if the perturbation of the equations is $\mC^1$ small, then the topology of the set of solutions doesn't change at all.}. 
We will use this aspect of the previous theorem in Chapter \ref{chap:semicontop} to produce an analogue of Milnor-Thom bound \cite{MilnorBound}, valid for any compact smooth hypersurface of $\R^n$ (Theorem \ref{semi:thm:witdash}). In this analogy the degree is replaced by an appropriate notion of condition number $\kappa_1(f)$.  
\begin{thm}\label{moti:thm:witdash}
Let $Z$ be a compact hypersurface of $\R^n$ defined by a smooth regular equation $f=0$.  Then there is a constant $c_1>0$ depending on the diameter of $Z$ and on $n$, such that
\be\label{moti:eq:bobound}
b\left(Z(f)\right)\le c_1\left(\kappa_1(f)\right)^n.
\ee
\end{thm}
A more geometric version of the inequality \eqref{moti:eq:bobound} is given in terms of the reach $\rho$ of the submanifold $Z$. This is due to the fact, establisehd by Theorem \ref{semi:thm:reach}, that it is possible to construct a function $f$ such that $f=0$ is a regular equation for $Z$, with $\kappa_1(f)\le 2(1+\rho^{-1})$ .
\begin{thm}\label{moti:thm:reach}
Let $Z$ be a compact hypersurface of $\R^n$ having reach $\rho$.  Then there is a constant $c_2>0$ depending only on $n$ and on the diameter of $Z$, such that
\be\label{moti:eq:reachbound}
b\left(Z\right)\le c_2\left(1+\frac 1\rho\right)^n.
\ee
\end{thm}
\section{Main results}
\subsection{Differential Topology of Gaussian Random Fields}\label{sec:gp}
Chapter \ref{chap:dtgrf} contains the paper \cite{dtgrf}, coauthored by Antonio Lerario, which is devoted to the set up of an efficient framework to treat problems of differential geometric and topological nature that arise frequently, when dealing with sequences of random maps and transversality. 

We propose a perspective focused on measure theory to work with  sequences of gaussian random fields, whose key point is to view them as random elements of the topological space $\mathcal{C}^r(M,\mathbb{R}^k)$ inducing a gaussian measure, in order to use the full power of the general theories of metric measure spaces and gaussian measures.
The paper contains the following results.

The first characterizes the topology of the space $\G(E^r)$ of gaussian measures on $\Cr rMk$, i.e. gaussian random fields up to equivalence, in terms of the covariance functions.
\begin{customthm}{\ref{thm:1}}[Measure-Covariance]\label{moti:thm:1}The natural map assigning to each gaussian measure its covariance function
\be\label{intro:eq:covariance} \mathcal{K}^r:\G(E^r)\to \mathcal{C}^{r}(M\times M, \R^{k\times k}),\ee
 is injective and continuous for all $r\in \N\cup\{ \infty\}$; when $r=\infty$ this map is also a closed topological embedding\footnote{A continuous injective map that is an homeomorphism onto its image.}.
\end{customthm}

We will show in Example \ref{escempio} that the function $\mathcal{K}^r$ is never a topological embedding if $r<\infty$. In this case we have the following partial result, formulated in terms of the convergence in law of a sequence of gaussian random fields.

\begin{customthm}{\ref{thm:2}}[Limit probabilities]
Let $\{X_d\}_{d\in \mathbb{N}}$ be a sequence of $\mC^{r}$ gaussian random fields such that the sequence $\{K_d\}_{d\in \mathbb{N}}$ of the associated covariance functions converges to $K$ in $\mathcal{C}^{r+2}(M\times M, \R^{k\times k})$. Then there exists a $\mC^r$ gaussian random field $X$ with $K_X=K$ such that for every Borel set $A\subset \Cr rMk$ we have
\be\label{intro:eq:limitprob} \PP(X\in \mathrm{int}(A))\leq \liminf_{d\to \infty}\PP(X_d\in A)\leq \limsup_{d\to \infty}\PP(X_d\in A)\leq \PP(X\in \overline{A}).\ee
\end{customthm}

A big class of problems in random algebraic geometry\footnote{We pick the letter ``$d$'' for the parameter, because $X_d$ is often a random polynomial of degree $d$.} can be reduced to the analysis of a deterministic object (see Section \ref{moti:limitprob}) by means of Theorem \ref{thm:2}.

The second result is a criterium establishing when a family of random maps depending smoothly on a countable quantity of gaussian parameters is almost surely transversal to a given submanifold. It is a generalization of the classical parametric transversality theorem to an \emph{infinite} dimensional and \emph{gaussian} situation. This is  nontrivial considering that in general the latter, being based on Sard's theorem, fails in the infinite dimensional setting.
\begin{customthm}{\ref{thm:transthm}}[Probabilistic transversality]
Let $X$ be a $\mC^r$ gaussian random field on a manifold $M$ and denote $F=\spt(X)$. Let $P,N$ be smooth manifolds and $W\subset N$ a submanifold. Assume that $\Phi\colon P\times F\to N$ is a smooth map  such that $\Phi\transv W$. Then
\be 
\P\{\phi(X)\transv W\}=1,
\ee
where $\phi(f)$ is the map that sends $p$ to $\Phi(p,f)$.
\end{customthm}
As a corollary one gets the following probabilistic version of Thom's Theorem, from differential topology, see \cite{Hirsch}.

\begin{customthm}{\ref{transthm2}}
Let $X\colon M\to \R^k$ a smooth gaussian random field and denote $F=\spt(X)$. Let $r\in \N$. Assume that for every $p\in M$ we have
\be \label{intro:eq:trans}
\mathrm{supp}(j^r_pX)=J^r_p(M,\R^k).\ee
Then for any submanifold $W\subset J^r(M,\R^k)$, we have $\P(j^rX \pitchfork W)=1$.
\end{customthm}

 These theorems help to move nimbly between the realms of Differential Topology and that of Stochastic Geometry, where the ``generic'' events are replaced by the ``almost sure'' ones.

\subsection{Kac Rice Formula for Transversal Intersections}\label{sec:kr}
Chapter \ref{chap:kr} is devoted to developing a generalization to the famous Theorem of Kac and Rice, which provides a formula for calculating the expected number of zeroes of a random function $X\colon M \to \mathbb{R}^{\dim M}$ in terms of the pointwise joint distributions of $X(p)$ and $X'(p)$. We prove a more general formula\footnote{It is rather easy to reduce to the standard case, by taking $W=\varphi^{-1}(0)$ for some map $\varphi\colon N\to \mathbb{R}^m$, and considering the random map $\varphi\circ X$. However then it is difficult to reinterpret the everything in terms of $X$.} that, in the same spirit, computes the expectation of the number of preimages via a random map $X\colon M\to N$ of a smooth submanifold $W\subset N$ having codimension equal to $\dim M$. 
\begin{customthm}{\ref{thm:main}}\label{thm:introKRform}
Let $X\colon M \to N$ be a random $\mC^1$ map between two riemannian manifolds and let $W\subset N$ such that $(X,W)$ is a \emph{good} (see below) couple. Then for every Borel subset $A\subset M$ we have
\begin{multline}\label{eq:introformrho}
\E\#\left(X^{-1}(W)\cap A\right)=\\=\int_{A}\left(\int_{ S_p\cap W}\E\left\{J_pX\frac{\sigma_x(X,W)}{\sigma_x(S_p,W)}\bigg|X(p)=x\right\}\rho_{X(p)}(x)d(S_p\cap W)(x)\right)dM(p),
\end{multline}
where $\rho_{X(p)}$ denotes the density of the random variable $X(p)$, $d(S_p\cap W)$ and $dM$ denote the volume densities of the corresponding riemannian manifolds and $J_pX$ is the jacobian of $X$; 
besides, $\sigma_x(X,W)$ and $\sigma_x(S_p,W)$ denote the ``angles'' \footnote{In the sense of Definition \ref{defi:angle}.}
made by $T_xW$ with, respectively, $d_pX(T_pM)$ and $T_xS_p$. 
\end{customthm}
The hypothesis of validity of this formula, namely the requirement that $(X,W)$ is a KROK couple, is satisfied in particular when $X$ is a smooth gaussian random section of a smooth vector bundle $\pi:E\to M$. In this case we show that the expectation depends continuously on the covariance of $X$.
\begin{customthm}{\ref{cor:mainjgau}}\label{thm:intromainEjgau}
Let $r\in \N$. Let $\pi\colon M\to E$ and let $X_d, X_\infty$ be a sequence of  smooth Gaussian random sections with non-degenerate $r^{th}$ jet and assume that 
\be K_{X_d}\to_{d\to +\infty} K_{X_\infty}\ee in the $\mC^{2r+2}$ topology (weak Whitney). Let $W=\sqcup_{i\in I}W_i\subset J^rE$ be a smooth closed Whitney stratified subset of codimension $m$ such that $W_i\transv E_p$ for every $p\in M$ and for every stratum $W_i$. Assume that $W$ is either \textbf{compact} or \textbf{locally semialgebraic}. Then 
\be 
\lim_{d\to +\infty}\E\#\{p\in A\colon j^r_pX_d\in W\}=\E\#\{p\in A\colon j^r_pX_{\infty}\in W\}
\ee
for every relatively compact Borel subset $A\subset M$.
\end{customthm}
\subsection{Maximal and Typical Topology of Real Polynomial Singularities}\label{sec:rag}
In Chapter \ref{chap:mttps}, containing the paper \cite{mttrp}, we investigate the topological structure of the singularities of real polynomial maps, unveiling a remarkable interplay between their maximal and typical complexity (meant as the sum of their Betti numbers).
In particular we prove the following theorem.

\begin{customthm}{\ref{thm:bettibound} and \ref{thm:bettiorderintro}}\label{thm:kost}
Let $\psi\colon S^m\to \mathbb{R}^k$ be a (homogenous) polynomial map of degree $d$ and $Z\subset S^m$ be the set of points in the sphere where the map $\psi$, together with its derivatives satisfies some given semialgebraic conditions. Then 
\begin{enumerate}[$\bullet$]
    \item $b(Z)\le O(d^m)$ for the generic choice of $\psi$\footnote{Here the classical bound, proved by Milnor \cite{MilnorBound}, would give $O(d^{m+1})$.}.
    \item $\mathbb{E}\{b(Z)\}=\Theta( d^\frac{m}{2})$\footnote{$a_d=\Theta(b_d)$ if there are two constants $0<c\le C$ and a big enough $d_0$, such that $cb_d\le a_d\le C b_d$ for all $d\ge d_0$.} for the typical choice of $\psi$.
\end{enumerate}
\end{customthm}
The set $Z$ is interpreted as a \emph{singularity} of $\psi$ of a certain kind. It is the preimage via the jet map $j^r\psi$ of a given diffeomorphism-invariant semialgebraic subset $W\subset J^r(S^m,\mathbb{R}^k)$.
The word ``typical'' and the expectation ``$\mathbb{E}$'' in the above Theorem refer to a natural probability measure on the space of polynomials, called Kostlan distribution. This is the most studied model because it's the only (up to homoteties) gaussian measure whose complexification is invariant under the action of the unitary group.


The problem of giving a good upper bound on the complexity of $Z=j^r\psi^{-1}(W)$ will require us to produce a version of Morse's inequalities for stratified maps.
\begin{customthm}{\ref{thm:strat}}
Let $W\subset J$ be a semialgebraic subset of a real algebraic smooth manifold $J$, with a given semialgebraic Whitney stratification $W=\sqcup_{S\in \mathscr{S}}S$ and let $M$ be a real algebraic smooth manifold. 
There exists a semialgebraic subset $\hat{W}\subset J^{1}(M,J\times \R)$ having codimension larger or equal than $\dim M$, equipped with a semialgebraic Whitney stratification that satisfies the following properties with respect to any couple of smooth maps $\psi\colon M\to J$ and $g\colon M\to \R$.
\begin{enumerate}
\item 
If $\psi\transv W$and $j^1(\psi,g)\transv \hat{W}$, then $g|_{\psi^{-1}(W)}$ is a Morse function with respect to the stratification $\psi^{-1}\mathscr{S}$ and
\be 
\text{Crit}(g|_{\psi^{-1}(W)})=\left(j^{1}(\psi,g)\right)^{-1}(\hat{W}).
\ee

\item
There is a constant $N_W>0$ depending only on $W$ and $\mathscr{S}$, such that if $\psi^{-1}(W)$ is compact, $\psi\transv W$ and $j^1(\psi,g)\transv \hat{W}$, then
\be 
b_i(\psi^{-1}(W))\le N_W\#\text{Crit}(g|_{\psi^{-1}(W)}),
\ee
for all $i=0,1,2\dots$
\end{enumerate}
\end{customthm}
For what regards the lower bound, we will use the Semicontinuity Theorem, described in the next subsection.

\subsection{Semicontinuity of Topology}\label{sec:top}
A consequence of Thom's Isotopy Lemma is that the set of solutions of a regular smooth equation is stable under $\mathcal{C}^1$-small perturbation, but what happens if the perturbation is just $\mathcal{C}^0$-small? Well, it turns out that the topology can only increase. In Chapter \ref{chap:semicontop} we will prove the following.
\begin{customthm}{\ref{semi:cor:mainnotco}}[Semicontinuity Theorem]
Let $M,N$ be smooth manifolds, let $Y\subset N$ be a closed smooth submanifold. Let $f\colon M\to N$ be a $\mC^1$ map such that $f\transv Y$. 
There is a neighborhood $\mathcal{U}\subset \coo0MN$ of $f$, open with respect to Whitney's strong topology, with the property that for any $\Tilde{f}\in\mathcal{U}$ there exist abelian groups $G_i$, for each $i\in \N$, such that
    \be
     \check{H}^i(\Tilde{f}^{-1}(Y))\cong H^i(f^{-1}(Y))\oplus G_i.
    \ee
\end{customthm}
This result, while enlightening a very general topological fact, finds nice applications when coupled with holonomic approximation arguments \cite{eliash}, which generally is only $\mathcal{C}^0-$small, or in contexts where a quantitative bounds is sought. 

We will devote the second part of Chapter \ref{chap:semicontop} to demonstrate how this theorem\footnote{Actually we will need to use Theorem \ref{semi:thm:main} wich allows in addition to control the size of the $\mC^0$ neighborhood $\mathcal{U}$.} can be used to produce a quantitative estimate on the sum of the Betti numbers $b(Z)$ of a compact hypersurface $Z$ of $\R^n$, in analogy with the classical bound proved by Milnor \cite{MilnorBound} in the algebraic case. Milnor's result applies when $Z$ is the zero set of a polynomial, providing a bound in terms of its degree: $b(Z)\le O(d^m)$. 

We will consider a smooth function $f$ with regular zero set $Z(f)$ contained in the interior of a disk $D$ and replace the degree with 
\be 
\kappa_1(f, D)=\frac{\|f\|_{\coo 1D{\R}}}{\min\left\{\min_{x\in\de D}|f(x)|,\min_{x\in D}\left(|f(x)|^2+\|\nabla f (x)\|^2\right)^{1/2}\right\}}.\ee
This parameter is essentially the $\mC^1$ distance of $f$ from the set of all non-regular equations (inside the projectivization of the space of smooth functions), for this reason it can be thought of as a sort of condition number. Moreover, it can be estimated in term of the reach of the submanifold $Z(f)$ itself (see Theorem \ref{semi:cor:reachwitdash}) .
\begin{customthm}{\ref{semi:thm:witdash}}\label{intro:thm:witdash}
Let $D \subset \R^n$ be a disk. Let $f\in \coo 1{D}\R$ have regular zero set $Z(f)\subset \emph{int}(D)$. Then there is a constant $c_1(D)>0$ such that
\be\label{intro:eq:bobound}
b\left(Z(f)\right)\le {c_1(D)\cdot}\left({\kappa_1(f,D)}\right)^n
\ee
\end{customthm}
Similar estimates can be obtained by means of Thom's isotopy lemma, but depend on the $\mC^2$ norm of $f$, while, using Theorem \ref{semi:cor:mainnotco}, only the $\mC^1$ information on $f$ is needed. In fact, we will also show that this bound is ``asymptotically sharp'', in the sense of the following statement.
\begin{customthm}{\ref{semi:thm:sharpdeg}}
There exists a sequence $\{f_k\}_{k\in \mathbb{N}}$ of maps in $C^1(D, \R)$, such that the sequence $\kappa_1(f_k,D)$ converges to $+\infty$, and a constant $c_0(D)>0$ such that for every $k\in \mathbb{N}$ the zero set $Z(f_k)\subset \mathrm{int}(D)$ is regular and
\be b(Z(f_k))\geq c_0(D)\cdot(\kappa_1(f_k, D))^{n}.
\ee 
\end{customthm}

\chapter{DIFFERENTIAL TOPOLOGY OF GAUSSIAN RANDOM FIELDS}\label{chap:dtgrf}
Motivated by numerous questions in random geometry, given a smooth manifold $M$, we approach a systematic study of the differential topology of Gaussian random fields (GRF) $X:M\to \R^k$, that we interpret as random variables with values in $\mathcal{C}^r(M, \R^k)$, inducing on it a Gaussian measure.

When the latter is given the weak Whitney topology, the convergence in law of $X$ allows to compute the limit probability of certain events in terms of the probability distribution of the limit. This is true, in particular, for the events of a geometric or topological nature, like: ``$X$ is transverse to $W$'' or ``$X^{-1}(0)$ is homeomorphic to $Z$''.

We relate the convergence in law of a sequence of GRFs with that of their covariance structures, proving that in the smooth case ($r=\infty$), the two conditions coincide, in analogy with what happens for finite dimensional Gaussian measures. We also show that this is false in the case of finite regularity ($r\in\N$), although the convergence of the covariance structures in the $\mathcal{C}^{r+2}$ sense is a sufficient condition for the convergence in law of the corresponding GRFs in the $\mathcal{C}^r$ sense.

We complement this study by proving an important technical tools: an infinite dimensional, probabilistic version of the Thom transversality theorem, which ensures that, under some conditions on the support, the jet of a GRF is almost surely transverse to a given submanifold of the jet space.

\section{Introduction}
\subsection{Overview}The subject of \emph{Gaussian random fields} is classical and largely developed (see for instance\footnote{This list is by no means complete!} \cite{AdlerTaylor, NazarovSodin2, EdelmanKostlan95, bogachev}). Motivated by problems in differential topology, in this paper we adopt a point of view which complements the classical one and we view Gaussian random fields as random variables in the space of smooth maps. Inside this space there is a rich structure coming from the geometric conditions that we can impose on the maps we are studying. There are some natural events, described by differential properties of the maps under consideration (e.g. being transverse to a given submanifold; having a certain number of critical points;  having a fixed homotopy type for the set of points satisfying some regular equation written in term of the field...), which are of specific interest to differential topology and it is desirable to have a verifiable notion of convergence of Gaussian random fields which ensures the convergence of the probability of these natural events. At the same time, once the space of functions is endowed with a probability distribution, it is natural to investigate the stability of these properties using the probabilistic language (replacing the notion of ``generic'' from differential topology with the notion of ``probability one''). 

The purpose of this paper is precisely to produce a general framework for investigating this type of questions. Specifically, Theorem \ref{thm:2} below allows to study the limit probabilities of these natural events for a family of Gaussian random fields (the needed notion of convergence is ``verifiable'' because it is written in terms of the convergence of the covariance functions of these fields). Theorem \ref{thm:1} relates this notion to the convergence of the fields in an appropriate topology: we achieve this by proving that there is a topological embedding of the set of smooth Gaussian random fields into the space of covariance functions. The switch from ``generic'' to ``probability one'' happens with Theorem \ref{transthm2}, which gives a probabilistic version of the Thom Transversality Theorem (again the needed conditions for this to hold can be checked using the covariance function of the field). This is actually a corollary of the more general Theorem \ref{thm:transthm}, which provides an infinite dimensional, probabilistic version, of the Parametric Transversality Theorem.

\subsection{Topology of random maps} Let $M$ be a smooth $m$-dimensional manifold (possibly with boundary). We denote by $E^r=\mathcal{C}^r(M, \R^k)$ the space of differentiable maps endowed with the weak Whitney topology, where $r\in \N\cup\{\infty\}$, and we call $\mathscr{P}(E^r)$ the set of probability measures on $\mathcal{C}^r(M, \R^k)$, endowed with the narrow topology (i.e. the weak* topology of $\mathcal{C}_b(E^r)^*$, see Definition \ref{defi:narrowtop}). 

In this paper we are interested in a special subset of $\mathscr{P}(E^r)$, namely the set $\G(E^r)$ of \emph{Gaussian measures}: these are probability measures with the property that for every finite set of points $p_1, \ldots, p_j\in M$ the evaluation map $\varphi:C^r(M, \R^k)\to \R^{jk}$ at these points induces (by pushforward) a Gaussian measure on $\R^{jk}.$ \footnote{In remark \ref{rem:bridge} we explain how this definition is equivalent to that of a Gaussian measure on the topological vector space $E^r$.} We denote by $\mathcal{G}^r(M, \R^k)$ the set of $\mathcal{C}^r$ \emph{Gaussian random fields} (GRF) i.e. random variables with values in $E^r$ that induce a Gaussian measure (see Definition \ref{def:GRF} below).
\begin{example}
The easiest example of GRF is that of a random function of the type $X=\xi_1 f_1 +\dots +\xi_n f_n$, where $\xi_1,\dots,\xi_n$ are independent Gaussian variables and $f_i\in \Cr rMk$. A slightly more general example is an almost surely convergent series
\be\label{eq:KarLo}
X=\sum_{n=0}^\infty \xi_n f_n.
\ee 
In fact, a standard result in the general theory of Gaussian measures (see \cite{bogachev}) is that every GRF admits such a representation, which is called the Karhunen-Loève expansion. We give a proof of this in Appendix \ref{app:rgrf} (see Theorem \ref{thm:rep}), adapted to the language of the present paper.
\end{example}
\begin{remark}One can define a Gaussian random section of a vector bundle $E\to M$ in an analogous way (the evaluation map here takes values in the finite dimensional vector space $E_{p_1}\oplus\dots\oplus E_{p_j}$). We choose to discuss only the case of  trivial vector bundles to avoid a complicated notation, besides, any vector bundle can be linearly embedded in a trivial one, so that any Gaussian random section can be viewed as a Gaussian random field. For this reason, the results we are going to present regarding GRFs are true, mutatis mutandis, for Gaussian random sections of vector bundles.
\end{remark}

We have the following sequence of continuous injections:
\be \G(E^{\infty})\subset \cdots \subset \G(E^r)\subset \cdots \subset \G(E^0)\subset \mathscr{P}(E^0),\ee
with the topologies induced by the inclusion $\G(E^r)\subset  \mathscr{P}(E^r)$ as a closed subset.

By definition, a Gaussian random field $X$ induces a Gaussian measure on $\mathcal{C}^r(M, \R^k)$, measure that we denote by $[X]$. Two fields are called \emph{equivalent} if they induce the same measures. A Gaussian measure $\mu=[X]\in \G(E^r)$ gives rise to a differentiable function $K_\mu\in \mathcal{C}^{r}(M\times M, \R^{k\times k})$ called the \emph{covariance function} and defined for $p,q\in M$ by:
\be K_\mu(p,q)=\E \left\{X(p) X(q)^T\right\}=\int_{E^r} f(p)f(q)^Td\mu(f).\ee
Equivalent fields give rise to the same covariance function, and to every covariance function there corresponds a unique (up to equivalence) Gaussian field. 

\begin{remark}In this paper we are interested in random fields up to equivalence, this is the reason why we choose to focus on the narrow topology. Indeed the notion of narrow convergence of a family $\{X_d\}_{d\in \mathbb{N}}$ of GRFs corresponds to the notion of \emph{convergence in law} of random elements in a topological space and it regards only the probability measures $[X_d]$. By the Skorohod theorem (see \cite[Theorem 6.7]{Billingsley}) this notion corresponds to almost sure convergence, up to equivalence of GRFs. In case one is interested in the almost sure convergence or in the convergence in probability of a particular sequence of GRFs one should be aware that these two notions take into account also the joint probabilities. For example, convergence in probability is equivalent to narrow convergence of the couple $(X_d,X)\nrw (X,X)$ (see Theorem \ref{lem:convprob}).
\end{remark}
Our first theorem translates convergence in $\G(E^r)$ of Gaussian measures
with respect to the narrow topology in terms of the corresponding sequence of covariance functions in the space $\mathcal{C}^{r}(M\times M, \R^{k\times k})$, endowed with the weak Whitney topology. 

\begin{thm}[Measure-Covariance]\label{thm:1}The natural map
\be\label{eq:covariance} \mathcal{K}^r:\G(E^r)\to \mathcal{C}^{r}(M\times M, \R^{k\times k}),\ee
given by $\mathcal{K}^r:\mu\mapsto K_\mu,$ is injective and continuous for all $r\in \N\cup\{ \infty\}$; when $r=\infty$ this map is also a closed topological embedding\footnote{A continuous injective map that is an homeomorphism onto its image.}.
\end{thm}

We observe at this point that the condition $r=\infty$ in the second part of the statement of Theorem \ref{thm:1} is necessary: as Example \ref{escempio} and Theorem \ref{thm:counter} show, it is possible to build a family of $\mathcal{C}^r$ ($r\neq \infty$) GRFs with covariance structures which are $\mathcal{C}^r$ converging but such that the family of GRFs does not converge narrowly to the GRF corresponding to the limit covariance.

Theorem \ref{thm:1} is especially useful when one has to deal with a family of Gaussian fields depending on some parameters, as it allows to infer asymptotic properties of probabilities on $\mathcal{C}^\infty(M, \R^k)$ from the convergence of the covariance functions (notice that this ``implication'' goes the opposite way of the arrow in \eqref{eq:covariance}). 

\begin{thm}[Limit probabilities]\label{thm:2}Let $\{X_d\}_{d\in \mathbb{N}}\subset \mathcal{G}^{r}(M,\R^k)$ be a sequence of Gaussian random fields such that the sequence $\{K_d\}_{d\in \mathbb{N}}$ of the associated covariance functions converges to $K$ in $\mathcal{C}^{r+2,r+2}(M\times M, \R^{k\times k})$\footnote{It is the space of functions $K(x,y)$ having continuous partial derivatives of order at least $r$ in both variables $x$ and $y$, see Section \ref{sec:spasmo}.}. Then there exists $X\in \g rMk$ with $K_X=K$ such that for every Borel set $A\subset E^r$ we have
\be\label{eq:limitprob} \PP(X\in \mathrm{int}(A))\leq \liminf_{d\to \infty}\PP(X_d\in A)\leq \limsup_{d\to \infty}\PP(X_d\in A)\leq \PP(X\in \overline{A}).\ee
In particular, if $\P(X\in\de A)=0$, then the limit exists:
\be\label{eq:limitprob2}
\lim_{d\to \infty}\P(X_d\in A)=\P(X\in A).
\ee
\end{thm}


\subsection{The support of a Gaussian random field}
The previous Theorem \ref{thm:2} raises two natural questions: 
\begin{enumerate}[$(1)$]
\item\label{itm:Q1} when is the leftmost probability in \eqref{eq:limitprob} strictly positive?
\item\label{itm:Q2} For which sets $A\subset E^\infty$ do we have  \eqref{eq:limitprob2}?
\end{enumerate}
Answering question \ref{itm:Q1} for a given Gaussian random field $X\in \mathcal{G}^r(M, \R^k)$, amounts to determine its topological support:
\be\label{eq:defsupport} \textrm{supp}(X)=\left\{f\in E^r\colon \PP(X\in U)>0\textrm{ for every neighborhood $U$ of $f$}\right\}.
\ee

We provide a description of the support of a Gaussian field $X=(X^1, \ldots, X^k)\in \mathcal{G}^{r}(M, \R^k)$ in terms of its covariance function $K_X$. 
 \begin{thm}\label{thm:dtgrf:3}
Let $X\in \g rMk$, consider all functions $h_p^j\in E^r$ of the form 
\be \label{suppbasiseq}
\begin{aligned}
\left. h^j_p(q) \right. = \begin{pmatrix}
K_X(q,p)^{1j} \\
\vdots \\
K_X(q,p)^{kj}
\end{pmatrix},\quad \text{for $p\in M$ and $j\in\{1,\dots,k\}$}. 
\end{aligned}
\ee
then
\be 
\spt(X)=\overline{\textrm{span}\{h_p^j\colon p\in M,j=1\dots k\}}^{E^r}.
\ee
\end{thm}
In particular, note that the support of a GRF is always a vector space, thus any neighborhood of $0$ has positive probability. 

Theorem \ref{thm:dtgrf:3} is just a general property of Gaussian measures, translated into the language of the present paper. In Section \ref{sec:PT3}, we prove it as a consequence of \cite[Theorem 3.6.1]{bogachev} together with a description of the Cameron-Martin space of $X$. 
In Appendix \ref{app:rgrf}, we present a direct proof of such result, adapted to our language (see Corollary \ref{cor:sputoCMapp}). We do this by generalizing the proof given in \cite[Section A.3-A.6]{NazarovSodin2} for the case in which $r=0$ and $k=1$. 

\subsection{Differential topology from the random point of view} Addressing question \ref{itm:Q2} above, let us observe that the probabilities in \eqref{eq:limitprob} are equal if and only if $\PP(X\in\partial A)=0$, and the study of this condition naturally leads us to the world of Differential Topology.

 When studying smooth maps, most relevant sets are given imposing some conditions on their jets (this is what happens, for instance, when studying a given singularity class). For example, let us take for $A\subset E^\infty$ in Theorem \ref{thm:2} a set defined by a condition on the $r$-th jet of $X$:
\be A=\{\textrm{$f\in E^\infty$ such that $j_x^r f\in V\subseteq J^{r}(M, \R^k)$ for all $x\in M$}\}.\ee
One can show that if $V$ is an open set with smooth boundary $\de V$, then there is no map $f\in \de A$ satisfying $j^rf \pitchfork \de V$. This is a frequent situation, indeed in most cases, the boundary of $A$ consists of functions whoose jet \emph{is not} transverse to a given submanifold $W\subset J^r(M, \R^k)$, and then the problem of proving the existence of the limit \eqref{eq:limitprob2} reduces to show that $\P(j^r X\pitchfork W)=1$.  Motivated by this, we prove the following.
\begin{thm}\label{transthm2}
Let $X\in \g\infty Mk$ and denote $F=\spt(X)$. Let $r\in \N$. Assume that for every $p\in M$ we have
\be \label{eq:trans}
\mathrm{supp}(j^r_pX)=J^r_p(M,\R^k).\ee
Then for any submanifold $W\subset J^r(M,\R^k)$, we have $\P(j^rX \pitchfork W)=1$.
\end{thm}
Let us explain condition \eqref{eq:trans} better. Given $X\in \mathcal{G}^r(M, \R^k)$ and $p\in M$ one can consider the random vector $j_p^rX\in J_p^r(M, \R^k)$: this is a Gaussian variable and \eqref{eq:trans} is the condition that the support of this Gaussian variable is the whole $J_p^r(M, \R^k)$.
For example, if the support of a $\mathcal{C}^r$-Gaussian field $X$ equals the whole $E^r$, then for every $W\subset J^{r}(M, \R^k)$ we have $X\pitchfork W$ with probability one. 

We will actually prove Theorem \ref{transthm2} as a corollary of the following more general theorem, that is an infinite dimensional version of the Parametric Transversality Theorem \ref{parametric transversality}.
\begin{thm}[Probabilistic transversality]\label{thm:transthm}
Let $X\in \g r Mk$, for $r\in\N\cup \{\infty\}$, and denote $F=\spt(X)$. Let $P,N$ be smooth manifolds and $W\subset N$ a submanifold. Assume that $\Phi\colon P\times F\to N$ is a smooth map  such that $\Phi\transv W$. Then
\be 
\P\{\phi(X)\transv W\}=1,
\ee
where $\phi(f)$ is the map $p\mapsto\Phi(p,f)$.
\end{thm}
\begin{remark}
We stress the fact that the space $F$ in Theorem \ref{thm:transthm} might be infinite dimensional. This is remarkable in view of the fact that the proof of the finite dimensional analogue of Theorem \ref{thm:transthm} makes use of Sard's theorem, which is essentially a finite dimensional tool. In fact, such result is false in general for smooth maps defined on an infinite dimensional space (see \cite{counterSardKupka}). In this context, the only alternative tool is the Sard-Smale theorem (see \cite{smalesard}), which says that the set of critical values of a smooth Fredholm map between Banach spaces is meagre (it is contained in a countable union of closed sets with empty interior). However, this is not enough to say something about the evaluation of a Gaussian measure on such set, not even when the dimension is finite.

Moreover, in both the proof of the finite dimensional transversality theorem and of Sard-Smale theorem an essential instrument is the Implicit Function theorem. Although this result, in its generalized version developed by Nash and Moser, is still at our disposal in the setting of Theorem \ref{thm:transthm} (at least when $M$ is compact), it fails to hold in the context of more general spaces.

 That said, the proof of theorem \ref{thm:transthm} relies on finite dimensional arguments and on the Cameron-Martin theorem (see \cite[Theorem 2.4.5]{bogachev}), a result that is specific to Gaussian measures on locally convex spaces. In fact, the careful reader can observe that the only property of $X$ that we use is that $[X]$ is a nondegenerate Radon Gaussian measure (in the sense of \cite[Def. 3.6.2]{bogachev}) on the second-countable, locally convex vector space $F$.
\end{remark}

\section{Preliminaries}\label{sec:preliminaries}
\subsection{Space of smooth functions}\label{sec:spasmo}
Let $M$ be a smooth manifold of dimension $m$. We will always implicitely assume that $M$ is Hausdorff and second countable, possibly with boundary. Let $k\in \N$ and $r\in\N\cup\{+\infty\}$. We will consider the set of $\mathcal{C}^r$ functions 
\be E^r=\Cr{r}{M}{k}\ee as a topological space with the weak Whitney topology as in ~\cite{Hirsch,NazarovSodin2}. Let $Q\colon D\hookrightarrow M$ be an embedding of a compact set $D\subset \R^n$, we define for any $f\in \Cr rMk$, the seminorm
\be\label{eq:seminormQR}
\|f\|_{Q,r}\doteq \sup\Big\{\big|\de_\a \big(f\circ Q\big) (x)\big| \colon \a\in \N^{m},\ |\a|\le r,\ x\in \text{int}(D)\Big\}.
\ee
Then, for $r\in \N$ finite, the weak topology on $\Cr rMk$ is defined by the family of seminorms $\{\nrm \cdot Qr\}_Q$, while the topology on $\Cr\infty Mk$ is defined by the whole family $\{\nrm \cdot Qr\}_{Q,r}$. 
We recall that for any $r\in\N\cup\{\infty\}$, the topological space $\Cr rMk$ is a Polish space: it is separable and homeomorphic to a complete metric space (indeed it is a Fréchet space). 
We will also need to consider the space $\Cr {r,r}{M\times M}k$ consisting of those functions that have continuous partial derivatives of order at least $r$ with respect to both the product variables. The topology on this space is defined by the seminorms
\[
\|f\|_{Q,(r,r)}\doteq \sup\Big\{\big|\de_{(\a,\beta)} \big(f\circ Q\big) (x,y)\big| \colon \a,\beta \in \N^{m},\ |\a|, |\beta|\le r,\ x,y\in \text{int}(D)\Big\},
\]
where now $Q$ varies among all product embeddings: $Q(x,y)=(Q_1(x), Q_2(y))\in M\times M$ and $Q_1, Q_2$ are embeddings of two compact sets $D_1, D_2$.
\begin{lemma}\label{Lemma:14}
Let $f,f_n\in \Cr rMk$. $f_n\to f$ in $\Cr rMk$ if and only if for any convergent sequence $p_n\to p$ in $M$, 
\be
j^r_{p_n}f_n\to j^r_{p}f \quad \text{in $J^r(M,\R^k)$}.
\ee
\end{lemma}
\begin{proof}See \cite[Chapter 2, Section 4]{Hirsch}.
\end{proof}
Given an open cover $\{U_\ell\}_{\ell\in L}$ of $M$, the restriction maps define a topological embedding $
\Cr rMk \hookrightarrow \prod_{\ell\in L} \Cr r{U_\ell}k,
$
indeed any converging sequence $p_n\to p$ belongs to some $U_\ell$ eventually. In particular, suppose that $Q_\ell\colon\D^m\hookrightarrow M$ are a countable family of embeddings of the unit $m-$disk $\D^m$ such that int$(Q_\ell(\D^m))=U_\ell$ is a covering of $M$\footnote{This is always possible in a smooth manifold without boundary, by definition, and it is still true if the manifold has boundary: if $p\in \de M$, take an embedding of the unit  disk $Q\colon\D^m  \hookrightarrow M$ such that $Q(\de \D^m)$ intersects $\de M$ in an open neighborhood of $p$, then the interior of $Q(\D^m)$, viewed as a subset of $M$, contains $p$.}. Then the maps $Q_\ell^*\colon f\mapsto f\circ Q_\ell$ define a topological embedding
\be\label{prodembballeq}
\{Q_\ell^*\}_\ell \colon \Cr rMk \hookrightarrow \left(\Cr r{\D^m}k\right)^L.
\ee

We refer to the book \cite{Hirsch} for the details about topologies on spaces of differentiable functions.
\subsection{Gaussian random fields}
Most of the material in this section, can be found in the book \cite{AdlerTaylor} and in the paper \cite{NazarovSodin2}; we develop the language in a slightly different way so that it suits our point of view focused on measure theory.

A real random variable $\gamma$ on a probability space $\Prob$ is said to be Gaussian if there are real numbers $\mu\in \R$ and $\sigma\geq 0$, such that $\gamma\sim N(\mu,\sigma^2)$, meaning that it induces the $N(\mu,\sigma^2)$ measure on the real numbers, which is $\delta_{\mu}$ if $\sigma=0$, and for $\sigma>0$ it has density 
\[
\rho(t)=\frac1{\sqrt{2\pi\sigma^2}}e^{-\frac{(t-\mu)^2}{2\sigma^2}}.
\]
In this paper, unless otherwise specified, all Gaussian variables and vectors are meant to be \emph{centered}, namely with $\mu=0$.

A (centered) Gaussian random vector $\xi$ in $\R^k$ is a random variable on $\R^k$ s.t. for any covector $\lambda\in(\R^k)^*$, the real random variable $\lambda \xi$ is (centered) Gaussian. In this case we write $\xi\sim N(0,K)$ where $K=\E\{\xi \xi^T\}$ is the so called \emph{covariance matrix}.
If $\xi$ is a Gaussian random vector in $\R^k$, there is a random vector $\gamma\sim N(0,\mathbbm{1}_j)$ in $\R^j$ and an injective $k\times j$ matrix $A$ s.t.
\[
\xi=A\gamma.
\]
In this case $K=AA^T$ and the support of $\xi$ is the image of $A$, which concides with the image of the matrix $K$, that is
\[
\textrm{supp}(\xi)=\{p\in\R^k \colon \P\{U_p\}>0 \text{  for all neighborhoods $U_p\ni p$ }\}=\text{Im}K\footnote{This is the finite dimensional version of Theorem \ref{thm:dtgrf:3}.},
\]
indeed $\xi\in \text{Im}K=\text{Im}A$ with $\P=1$. If $A$ is invertible, $\xi$ is said to be nondegenerate, this happens if and only if $\det K\neq 0$, if and only if $\textrm{supp}(\xi)=\R^n$, if and only if the probability induced by $\xi$ admits a density, which is given by the formula
\be\label{gaussianlaweq}
\P\{\xi\in U\}=\frac1{(2\pi)^\frac{n}{2}\det K^\frac12}\int_{U}e^{-\frac12 W^TK^{-1}W}dW^n.
\ee

\begin{defi}[Gaussian random field]\label{def:GRF}
Let $M$ be a smooth manifold.
Let $\Prob$ be a probability space. An $\R^k$-valued \emph{random field (RF)} on $M$ is a measurable map
\[
X:\Omega\to {(\R^{k})}^M,
\]
with respect to the product $\sigma-$algebra on the codomain.
An $\R$-valued \emph{RF} is called a \emph{random function}.

Let $r\in\N\cup \{\infty\}$. We say that $X$ is a $\mathcal{C}^r$ random field, if $X_{\w} \in \Cr{r}{M}{k}$ for $\P$-almost every $\w\in \Omega$. 
We say that $X$ is a \emph{Gaussian random field (GRF)}, or just \emph{Gaussian field}, if for any finite collection of points $p_1,\dots, p_j\in M$, the random vector in $\R^{jk}$ defined by $(X(p_1),\dots,X(p_j))$ is Gaussian.
We denote by $\mathcal{G}^r(M, \R^k)$ the set of $\mathcal{C}^r$ Gaussian fields.
\end{defi}
When dealing with random fields $X:\Omega\to (\R^k)^M$, we will most often use the shortened notation of omitting the dependence from the variable $\w$. In this way $X:M\to\R^k$ is a  \emph{random} map, i.e. a random element\footnote{We recall that, given a measurable space $(S,\mathfrak{A})$, a measurable map from a probability space $\Prob$ to $S$ is also called a \emph{Random Element} of $S$ (see \cite{Billingsley}). Random variables and random vectors are random elements of $\R$ and $\R^k$, respectively.} of $(\R^k)^M$.

\begin{remark}
In the above definition, the sentence:
\be \textrm{``$X_{\w} \in \Cr{r}{M}{k}$ for $\P$-almost every $\w\in \Omega$''}\ee 
means that the set $\{\w\in\Omega\colon X_\w\in \Cr{r}{M}{k}\}$ contains a \emph{measurable} set $\Omega_0$ which has probability one. We make this remark because the subset $\Cr rMk$ doesn't belong to the product $\sigma-$algebra of $(\R^{k})^M$. 
\end{remark}

\begin{lemma}\label{borelemm1}
For all $r \in \N\cup\{+\infty\}$ the Borel $\sigma$-algebra $\mathcal{B}\Big(\Cr rMk\Big)$ is generated by the sets 
\[
\{f\in\Cr rMk \ :\ f(p)\in A\},
\]
with $p\in M$ and $A\subset \R^k$ open. 
Moreover $\Cr rMk$ is a Borel subset of $\Cr 0Mk$, for all $r \in \N\cup\{+\infty\}$.\end{lemma}
\begin{proof}See \cite[p. 43,44]{NazarovSodin2} or \cite[p. 374]{bogachev}.
\end{proof}

As a consequence we have that the Borel $\sigma$-algebra $\mathcal{B}(\Cr rMk)$ is  the restriction to $\Cr rMk$ of the product $\sigma$-algebra of $(\R^k)^M$. It follows that $X$ is a $\mathcal{C}^r$ RF on $M$ if and only if it is $\P-$almost surely equal to a random element of $\Cr rMk$.

A second consequence is that if $X$ is a $\mathcal{C}^r$ RF, then the associated map $\tilde{X}\colon \Omega\times M\to \R^k$ is measurable, being the composition $e\circ (X\times \text{id})$, where $e\colon \Cr rMk \times M \to \R^k$ is the continous map defined by $e(f,p)=f(p)$.




If $X$ is a $\mathcal{C}^r$ RF, then it induces a probability measure $X_*\P$ on $\Cr rMk$, or equivalently (because of Lemma \ref{borelemm1}) a probability measure on $\Cr 0Mk$ that is supported on $\Cr rMk$. We say that two RFs are \emph{equivalent} if they induce the same measure; note that this can happen even if they are defined on different probability spaces.

It is easy to see that every probability measure $\mu$ on $\Cr rMk$ is induced by some RF (just take $\Omega=\Cr rMk$, $\mu=\P$ and define $X$ to be the identity, then clearly $\mu=X_*\P$). This means that the study of $\mathcal{C}^r$ random fields up to equivalence corresponds to the study of Borel probability measures on $\Cr rMk$. 

Note that, as a consequence of Lemma \ref{borelemm1}, a Borel measure $\mu$ on $\Cr rMk$ is uniquely determined by its \emph{finite dimensional distributions}, which are the measures induced on $\R^{kj}$ by evaluation on $j$ points.

We will write $\mu=[X]$ to say that the probability measure $\mu$ is induced by a random field $X$.
In particular we define the \emph{Gaussian measures} on $\Cr rMk$ to be those measures that are induced by a $\mathcal{C}^r$ GRF, equivalently we give the following measure-theoretic definition.
\begin{defi}[Gaussian measure]\label{def:gauss}
Let $M$ be a smooth manifold and let $r\in \N\cup \{\infty\}$, $k\in \N$.
A \emph{Gaussian measure} on $\Cr rMk$ is a probability measure on the topological space $\Cr rMk$, with the property that for any finite set of points $p_1,\dots p_j\in M$, the measure induced on $\R^{jk}$ by the map $f\mapsto (f(p_1),\dots,f(p_j))$ is Gaussian (centered and possibly degenerate).
We denote by $\G(E^r)$ the set of Gaussian probability measures on $E^r=\mathcal{C}^r(M, \R^k).$
\end{defi}
\begin{remark}\label{rem:bridge}
In general a Gaussian measure on a topological vector space $W$ is defined as a Borel measure on $W$ such that all the elements in $W^*$ are Gaussian random variables (see \cite{bogachev}). In the case $W=\Cr rMk$, this is equivalent to Definition \ref{def:gauss}, because the set of functionals $f\mapsto a_1f(p_1)+\dots+a_jf(p_j))$ is dense in the topological dual $W^*$ (Theorem \ref{thm:deltadense} of Appendix \ref{app:dual}), therefore every continuous linear functional $\lambda\in W^*$ can be obtained as the almost sure limit of a sequence of Gaussian variables and thus it is Gaussian itself.
\end{remark}

We prove now a simple Lemma that will be needed in the following. Given a differentiable map $f\in \mathcal{C}^r(M, \R^k)$ with $r\geq 1$, and a smooth vector field $v$ on $M$, we denote by $vf$ the derivative of $f$ in the direction of $v$.

\begin{lemma}\label{vXisGRF}
Let $X\in\g rMk$ and let $v$ be a smooth vector field on $M$. Then $vX\in\g {r-1}Mk$.
(Notice that, as a consequence, the $r$-jet of a $\mathcal{C}^r$ GRF is a $\mathcal{C}^0$ GRF.)
\end{lemma}
\begin{proof}
Since $ X\in \Cr rMk$ almost surely, then $vX\in \Cr {r-1}Mk$ almost surely, thus $vX$ defines a probability measure supported on $\Cr {r-1}Mk$. To prove that it is a Gaussian measure, note that $vX(p)$ is either a $N(0,0)$ Gaussian, if $v_p=0$, or an almost sure limit of Gaussian vectors, indeed passing to a coordinate chart $x^1,\dots,x^m$ centered at $p$ s.t. $v_p=\frac{\de}{\de x^1}$, we have 
\[
vX(p)=\lim_{t\to 0}\frac{X(t,0,\dots, 0)-X(0,0,\dots, 0)}{t}\quad\text{a.s.}
\]
therefore it is Gaussian. The analogous argument can be applied when we consider a finite number of points in $M$.
\end{proof}

\subsection{The topology of random fields.} We denote by $\mathscr{P}(E^r)$, the set of all  Borel probability measures on $E^r$. We shall endow the space $\mathscr{P}(E^r)$ with the \emph{narrow topology}, defined as follows. Let $\mathcal{C}_b(E^r)$ be the Banach space of all bounded continuous functions from $E^r$ to $\R$.
\begin{defi}[Narrow topology]\label{defi:narrowtop}
The narrow topology on $\mathscr{P}(E^r)$ is defined as the coarsest topology such that for every $\varphi\in \mathcal{C}_b(E^r)$ the map $\textrm{ev}_{\varphi}:\mathcal{P}(E^r)\to \R$ given by:
\be\textrm{ev}_{\varphi}:\P \mapsto \int_{E^r} \varphi \,d\P\ee
is continuous.
\end{defi}
In other words, the narrow topology is the topology induced by the weak-$*$ topology of $\mathcal{C}_b(E^r)^*$, via the inclusion
\[
\mathscr{P}(E^r)\hookrightarrow \mathcal{C}_b(E^r)^*
\]
\[
\P\mapsto \E\{\cdot\}
\]
\begin{remark} The narrow topology is also classically refered to as the \emph{weak topology} (see \cite{Parth}, \cite{bogachev} or \cite{Billingsley}). We avoid the latter terminology to prevent confusion with the topology induced by the weak topology of $\mathcal{C}_b(E^r)^*$, which is strictly finer. Indeed if a sequence of probability measures $\mu_n$ converges to a probability measure $\mu$ in the weak topology of $\mathcal{C}_b(E^r)^*$, then $\lim_{n\to\infty}\mu_n(A)=\mu(A)$ for any measurable set $A\in E^r$. This is a strictly stronger condition than narrow convergence, see Portmanteau's theorem \cite{Billingsley}.\end{remark}


Convergence of a sequence of probability measures $\mu_n\in \mathscr{P}(E^r)$ in the narrow topology is denoted as $\mu_n\nrw \mu$.
From the point of view of random fields, $[X_n]\nrw[X]$ in $\mathscr{P}(E^r)$, if and only if 
\[
\lim_{n\to \infty}\E\{\varphi(X_n)\}=\E\{\varphi(X)\} \qquad \forall \varphi\in \mathcal{C}_b(E^r)
\]
and in this case we will simply write $X_n\nrw X$.
This notion of convergence of random variables is also called convergence \emph{in law} or \emph{in distribution}. 

To understand the notion of narrow convergence it is important to recall Skorohod's theorem (see \cite[Theorem 6.7]{Billingsley}), which states that $\mu_n\nrw \mu_0$ in $\mathscr{P}(E^r)$ if and only if there is a sequence $X_n$ of random elements of $E^r$, such that $\mu_n=[X_n]$ and $X_n\to X_0$ almost surely. In other words, narrow convergence is equivalent to almost sure convergence from the point of view of the measures $\mu_n$.

However, for a given sequence of random fields $X_n$, the notion of narrow convergence is even weaker than that of convergence in probability. The subtle difference, as showed in Lemma \ref{lem:convprob} below, is that the latter takes into account the joint distributions. 
\begin{lemma}\label{lem:convprob} 
Let $X_d, X \in \g r Mk$. The sequence $X_d$ convergese to $X$ in probability if and only if $(X_d,X)\nrw (X,X)$.
\end{lemma}
\begin{proof}
First, note that if $X_d\to X$ in probability, then $(X_d,X)\to (X,X)$ in probability and therefore $(X_d,X)\nrw (X,X)$. For the converse,
let $d$ be any metric on $\Cr rMk$. Since $d$ is a continuous function, if $(X_d,X)\nrw (X,X)$ then $d(X_d,X)\nrw 0$, which is equivalent to convergence in probability, by definition.
\end{proof}

We recall the following useful fact relating properties of the topology of $E$ to properties of the narrow topology on $\mathscr{P}(E);$ for the proof the reader is referred to \cite[p. 42-46]{Parth}.
\begin{prop}\label{prop:equiva}The following properties are true:
\begin{enumerate}
\item $E$ is separable and metrizable if and only if $\mathscr{P}(E)$ is separable and metrizable. In this case, the map $E\hookrightarrow \mathscr{P}(E)$, defined by $f\mapsto \delta_f$, is a closed topological embedding and the convex hull of its image is dense in $\mathscr{P}(E)$.

\item $E$ is compact and metrizable if and only if $\mathscr{P}(E)$ is compact and metrizable.

\item $E$ is Polish if and only if $\mathscr{P}(E)$ is Polish.

\end{enumerate}
%
\end{prop}
\begin{cor}\label{inducedmap}
Let $E_1$ and $E_2$ be two separable metric spaces. Let $\pi\colon E_1\to E_2$ be continuous. Then the induced map $\pi_*\colon \mathscr{P}(E_1)\to \mathscr{P}(E_2)$ is continuous. If moreover $\pi$ is a topological embedding, then $\pi_*$ is a topological embedding as well.
\end{cor}
\begin{proof}
If $\pi$ is continuous, then for any bounded and continuous real function $\varphi\in\mathcal{C}_b(E_2)$, the composition $\varphi\circ \pi$ is in $\mathcal{C}_b(E_1)$. Hence, the function $\int_{E_1}(\varphi\circ \pi) \colon \mathscr{P}(E_1) \to \R$ defined as $\P\mapsto \int_{E_1}(\varphi\circ \pi)d\P$ is continuous. Observe that for any $\P\in \mathscr{P}(E_1)$
\be 
\int_{E_1}(\varphi\circ \pi)d\P =\int_{E_2}\varphi \,d(\pi_*\P)=\left(\int_{E_2}\varphi\right)\circ \pi_*(\P),
\ee
thus the composition $(\int_{E_2}\varphi)\circ \pi_*\colon \mathscr{P}(E_1)\to \R$ is continuous for any $\varphi\in\mathcal{C}_b(E_2)$. From the definition of the topology on $\mathscr{P}(E_2)$, it follows that $\pi_*$ is continuous.


Assume now that $\pi$ is a topological embedding. This is equivalent to say that $\pi$ is injective and any open set $U\subset E_1$ is of the form $\pi^{-1}(V)$ for some open subset $V\subset E_2$, and the same for Borel sets. It follows that $\pi_*$ is injective, indeed if two probability measures $\P_1, \P_2\in \mathscr{P}(E_1)$, have equal induced measures $\pi_*\P_1=\pi_*\P_2$, then 
\be 
\P_1\{\pi^{-1}(V)\}=\P_2\{\pi^{-1}(V)\}
\ee
for any Borel subset $V\subset E_2$, thus $\P_1\{U\}=\P_2\{U\}$ for any Borel subset $U\subset E_1$ and $\P_1=\P_2$.

It remains to prove that $\pi_*^{-1}$ is continuous on the image of $\pi_*$. Let $\P_n\in \mathscr{P}(E_1)$ be such that $\pi_*\P_n\nrw \pi_*\P_0$. Let $U\subset E_1$ be open, then there is some open subset $V\subset E_2$ such that $\pi^{-1}(V)=U$ and, by Portmanteau's theorem (see \cite[p. 40]{Parth}), we get
\be 
\liminf_n \P_n\{U\}=\liminf_n \pi_*\P_n\{V\}\ge \pi_*\P_0\{V\}=\P_0\{U\}.
\ee
This implies that $\P_n\nrw \P_0$. We conclude using point (1) of Proposition \ref{prop:equiva},and the fact that on metric spaces, sequential continuity is equivalent to continuity.
\end{proof}
\begin{example}
Let $\phi \colon M \to N$ be a $\mathcal{C}^r$ maps between smooth manifolds, then the map $\phi^*\colon \mathcal{C}^r(N,W)\to \mathcal{C}^r(M,W)$ defined as $\phi^*(f)=f\circ \phi$ is continuous, therefore the induced map between the spaces of probabilities, which we still denote as $\phi^*$, is continuous. The same holds for the map $\phi_*\colon \mathcal{C}^r(W,M)\to \mathcal{C}^r(W,N)$, such that $\phi_*(f)=\phi\circ f$. 
\end{example}

Note that $\mathcal{C}^r$ narrow convergence implies $\mathcal{C}^s$ narrow convergence, for every $s\le r$, but not vice versa. Indeed there are continuous injections
\be \G(E^{\infty})\subset \cdots \subset \G(E^r)\subset \cdots \subset \G(E^0)\subset \mathscr{P}(E^0).\ee

\begin{prop}
$\G(E^r)$ is closed in $\mathscr{P}(E^r)$.
\end{prop}
\begin{proof}
Let $X_n\in \g rMk$ s.t. $X_n\nrw X\in \mathscr{P}(E^r)$. Then for any $p_1\dots, p_j\in M$ we have
\[
\left(X_n(p_1),\dots ,X_n(p_j)\right)\nrw \left(X(p_1),\dots ,X(p_j)\right) \]
in $\mathscr{P}(\R^{jk})$. Therefore the latter is a Gaussian random vector and thus $[X]\in \G(E^r).$
\end{proof}
\subsection{The covariance function.}
Given a Gaussian random vector $\xi$, it is clear by equation \eqref{gaussianlaweq} that the corresponding measure $\mu=[\xi]$ on $\R^m$ is determined by the covariance matrix $K=\E\{\xi\xi^T\}$. Similarly, if $X\in \g rMk$, then $[X]$ is a measure on $\Cr rMk$ and it is uniquely determined by its finite dimensional distributions, which are the Gaussian measures induced on $\R^{kj}$ by evaluation on $j$ points. It follows that $[X]$ is uniquely determined by the collection of all the covariances of the evaluations at couples of points in $M$, which we call \emph{covariance function}.
\begin{defi}[covariance function]
Given $X\in \g rMk$, we define its \emph{covariance function} as:
\be 
K_X\colon M\times M \to \R^{k\times k} 
\ee
\be 
K_X(p,q)=\E\{X(p)X(q)^T\}.
\ee
\end{defi}
The function $K_X$ is symmetric: $K_X(p,q)^T=K_X(q,p)$ and non-negative definite, which means that for any $p_1,\dots, p_j\in M$ and $\lambda_1,\dots, \lambda_j \in \R^k$, $
\sum_{i=1}^j \lambda_j^TK_X(p_i,p_i)\lambda_j\ge 0.$

The covariance function of a $\mathcal{C}^r$ random field is of class $\mathcal{C}^{r,r}$, see Section \ref{sec:spasmo}. This is better understood by introducing the following object.
Suppose that $X$ is a Gaussian random field on $M$, defined on a probability space $\Prob$, then it defines a map
\be\label{eq:gammaX}
 \gamma_X :M\to L^2\Prob^k\ee
such that $\gamma_X(p)=X(p)$.

To say that $X$ is Gaussian is equivalent to say that span$\{\gamma_X(M)\}$ is a Gaussian subspace of $L^2\Prob^k$, namely a vector subspace whose elements are Gaussian random vectors. 
Next proposition from \cite{NazarovSodin2} will be instrumental for us.

\begin{prop}[Lemma A.3 from \cite{NazarovSodin2}]\label{Ltwocont}
Let $X\in \g rMk$, then the map $\gamma_X\colon M\to L^2\Prob^k$ is $\mathcal{C}^r$. Moreover if $x,y$ are any two coordinate charts on $M$, then
\be 
\E\left\{\de_\a X(x)\left(\de_\beta X(y)\right)^T\right\}=
\langle \de_\a \gamma_X (x), \de_\beta \gamma_X(y)^T\rangle_{L^2\Prob}=
\de_{(\a,\beta)} K_X(x,y),
\ee
for any multi-indices $|\a|, |\beta|\le r$.
\end{prop}
\begin{cor}\label{cor:covcrr}
Let $X\in \g rMk$, then $K_X\in \Cr {r,r}{M\times M}{k\times k}$.
\end{cor}

\subsection{A Gaussian inequality}
The scope of this section is to prove Theorem \ref{Gaussinequality}, which contains a key technical inequality. Although such inequality can be seen as a consequence of Kolmogorov's theorem for $\mathcal{C}^{k,k}$ kernels, as discussed in \cite[Sec. A.9]{NazarovSodin2}, we report here a simpler proof. 
In fact, the result follows from a general inequality valid for all GRFs, not necessarily continuous. 

Given a GRF $X\colon M\to \R$, we define for all $\e>0$ the quantity $N(\e)$, to be the minimum number of $L^2$-balls of radius $\e$ needed to cover $\gamma_X(M)$. This number is always finite if $\gamma_X(M)$ is relatively compact in $L^2$. We will need the following Theorem from \cite{AdlerTaylor}.
\begin{thm}[Theorem 1.3.3 from \cite{AdlerTaylor}]\label{Adlerineq}
Let $\gamma_X(M)$ be compact in $L^2\Prob$. Let $\Delta_X=\text{diam}(\gamma_X(M))$. There exists a universal constant $C>0$ such that
\be 
\E\left\{\sup_{x\in M}X(t)\right\}\le C\int_0^{\Delta_X}\sqrt{\ln N(\e)}d\e
\ee

\end{thm}
As a corollary, in our setting we can derive the following.
\begin{lemma}\label{zeroineqlemma}
Let $X\in \g {1}M{}$ and consider an embedding $Q: D \hookrightarrow M$ of a compact disk $D\subset \R^m$. There is a constant $C_Q>0$ such that
\be 
\E\{\|X\|_{Q,0}\}\le C_Q \sqrt{\|K_X\|_{Q\times Q,1 }}
\ee
\end{lemma}
\begin{proof}
It is not restrictive to assume that $M=D$ and $Q=\textrm{id}$.
Notice that since the map $\gamma_X$ is continuous, by Proposition \ref{Ltwocont}, it follows that $\gamma_X(D)$ is compact in $L^2\Prob$, so that we can apply Theorem \ref{Adlerineq} to get that
\be
\E\{\|X\|_{D, 0}\}\le 2C\int_0^{\Delta_X}\sqrt{\ln N(\e)}d\e
\ee
Moreover, for any $q,p\in D$, we have that 
\be 
\begin{aligned}
\|X(p)-X(q)\|_{L^2}^2 &=K(p,p)+K(q,q)-2K(p,q) \\
&\le \left|K(p,p)-K(q,p)\right|+\left|K(q,q)-K(p,q)\right|\\
&\le
2 \sup_{x,y\in D} \left|\frac{\de K}{\de x}(x,y)\right||p-q|,
\end{aligned}
\ee
where $K=K_X$. Thus, denoting $\Lambda^2=2\|K\|_{Q\times Q, 1}$, we obtain that 
\be\label{ineq}
\|X(p)-X(q)\|_{L^2}\le \Lambda |q-p|^{\frac12}.
\ee
Let now $\tilde{N}(\rho)$ be the minimum number of standard balls in $\R^m$ with radius $\rho$, required to cover $D$. A consequence of \eqref{ineq} is that every ball of radius $\rho$ in $D$ is contained in the preimage via $\gamma_X$ of a ball of radius $\Lambda\rho^{\frac12}$ in $L^2$, therefore $N(\e)\le \tilde{N}(\frac{\e^2}{\Lambda^2})$. Besides, $\Delta_X\le \Lambda\sqrt{R}$, where $R$ is the diameter of $D$, so that
\be 
\E\{\|X\|_{D,0}\} \le 2C\int_0^{\Lambda \sqrt{R}}\sqrt{\ln \tilde{N}\left(\frac{\e^2}{\Lambda^2}\right)}d\e 
= 
2C\Lambda\int_0^{\sqrt{R}}\sqrt{\ln \tilde{N}\left(s^2\right)}ds.
\ee
Now, since $D\subset \R^m$, there is a constant $c_m$ such that $\tilde{N}(\rho)\le c_m\left(\frac{R}{\rho}\right)^m$, therefore
\[
I(R)=\int_0^{\sqrt{R}}\sqrt{\ln \tilde{N}\left(s^2\right)}ds\le \int_0^{\sqrt{R}}\sqrt{\ln c_m\left(\frac{R^2}{s^2}\right)^{m} }ds<\infty.
\]
We conclude that $\E\{\|X\|_{D,0}\}\le 2\sqrt{2}\cdot C\cdot I(R)\sqrt{\|K_X\|_{Q\times Q,1 }}$.
\end{proof}
We are now able to prove the required Gaussian inequality.
\begin{thm}\label{Gaussinequality}
Let $X\in \g {r}Mk$ and consider an embedding $Q: D \hookrightarrow M$ of a compact disk $D\subset \R^m$. Then
\[
\E\{\|X\|_{Q, {r-1}}\}\le C \sqrt{\|K_X\|_{Q\times  Q, (r,r)}},
\]
Where $C$ is a constant depending only on $Q$, $r$ and $k$.
\end{thm}
\begin{proof}
A repeated application of Lemma \ref{vXisGRF} proves that $\de_\a X^i$ is Gaussian, so that we can use Lemma \ref{zeroineqlemma} as follows.
\be 
\begin{aligned}
\E\{\|X\|_{Q,r-1}\} &\le \sum_{|\a|< r, i\le k}\E\{\|\de_\a X^i\|_{Q,0}\} \\
&\le \sum_{|\a|< r, i\le k}C_Q \sqrt{\|K_{\de_\a X^i}\|_{Q\times Q,1}} \\
& =\sum_{|\a|< r, i\le k}C_Q \sqrt{\|\de_{(\a,\a)} K^{i,i}_{X}\|_{Q\times Q,1}} \\
&\le C(Q,r,k)\sqrt{\|K_{X}\|_{Q\times Q,(r,r)}}\ .
\end{aligned}
\ee
\end{proof}

%
%
%
%
%
%
%
%
%
%
%
%
%
%
%
\section{Proof of Theorem \ref{thm:1} and Theorem \ref{thm:2}}\label{sec:thm12}

\subsection{Proof that $\K^r$ is injective and continuous}We already noted that $K_X$ determines $[X]$, and this is equivalent to say that $\K^0$ is injective. It follows that $\K^r$ is injective for every $r$, since $\K^r$ is just the restriction of $\K^0$ to $\G(E^r)$. 

Let us prove continuity. Since both the domain and the codomain are metrizable topological spaces, it will be sufficient to prove sequential continuity. 
Let $\mu_n\nrw\mu\in \G(E^r)$. Let $X\in \mathcal{G}^r(M, \R^k)$ be a GRF such that $\mu=[X]$ and for every $n\in \N$ let $X_n\in \mathcal{G}^r(M, \R^k)$ be such that $\mu_n=[X_n].$  By Skorohod's representation theorem (see \cite[Theorem 6.7]{Billingsley}) we can assume that the $X_n$ are GRFs defined on a common probability space $\Prob$ and that $X_n\to X$ almost surely in the topological space $\Cr rMk$.

To prove $\mathcal{C}^{r,r}$ convergence of $K_n=K_{X_n}$ to $K=K_X$, it is sufficient (and necessary) to show that given coordinate charts $(x,y)$ on $M\times M$, a sequence $(x_n,y_n)\to (x_0, y_0)$, a couple of indices $|\a|,|\beta|\le r$ and two indices $i,j \in\{1,\dots,k\}$, then
\be\label{Kconvergeseq}
\de_{(\a,\beta)} K_n^{i,j}(x_n,y_n)\to \de_{(\a,\beta)} K^{i,j}(x_0,y_0).
\ee
Let $\gamma_n=\de_\a X_n^i(x_n)$ and $\xi_n=\de_\beta X_n^j(y_n)$. By Lemma \ref{vXisGRF}, these two random vectors are Gaussian; moreover $\gamma_n\to \gamma$ and $\xi_n\to \xi$ almost surely. It follows that the convergence holds also in $L^2\Prob$, so that
\be 
\E\{\gamma_n\xi_n\}\to \E\{\gamma\xi\},
\ee
which is exactly \eqref{Kconvergeseq}.

\subsection{Relative compactness}
As we will see with Theorem \ref{thm:counter}, the map $\K^r$ is not proper when $r$ is finite. However, we have the following partial result.
\begin{thm}\label{rproper}
Let $r\in \N$ and consider $[X_n]\in \G(E^{r+2})$ and let $\{Q_\ell\}_{\ell\in \N}$ be a countable family of embeddings $Q_\ell:\D^m\hookrightarrow M$, such that the family of open sets $\textrm{int}(Q_\ell(\D))$ is a covering of $M$ (so that condition \eqref{prodembballeq} holds). Then the following conditions are related by the implications: \ref{itm:Ir} $\implies$ \ref{itm:IIr} $\implies$ \ref{itm:IIIr}.
\begin{enumerate}[$(1)$]
\item\label{itm:Ir} $ \sup_{n}\|K_{X_n}\|_{Q_\ell\times Q_\ell, (r+2,r+2)} <\infty$, for very $\ell\in \N$.
\item \label{itm:IIr} 
$\sup_{n}\E\left\{\|X_n\|_{Q, r+1}\right\} <\infty$, for any embedding $Q\colon D\hookrightarrow M$ of a compact $D\subset \R^m$.
\item \label{itm:IIIr}
The sequence $\{[X_n]\}_{n\in \N}$ is relatively compact in $\G(E^r)$.
\end{enumerate}
\end{thm}
Since there is a continuous inclusion $\G (E^\infty)\subset \G (E^r)$ and the function $\K^r$ is continuous, for any $r\in\N$, we immediately obtain the following corollary for the case $r=\infty$.
\begin{thm}\label{infproper}
Let $[X_n]\in \G(E^{\infty})$ and $\{Q_\ell\}_{\ell\in \N}$ be a countable family of embeddings $Q_\ell:\D^m\hookrightarrow M$, such that the family of open sets $\textrm{int}(Q_\ell(\D))$ is a covering of $M$ (so that condition \eqref{prodembballeq} holds). Then the following conditions are equivalent.
\begin{enumerate}[$(1)$]
\item\label{itm:Iinf} $ \sup_{n}\|K_{X_n}\|_{Q_\ell\times Q_\ell, (r,r)} <\infty$, for very $r,\ell\in \N$.
\item \label{itm:IIinf} 
$\sup_{n}\E\left\{\|X_n\|_{Q, r}\right\} <\infty$, for any embedding $Q\colon D\hookrightarrow M$ of a compact $D\subset \R^m$ and every $r\in\N$.
\item \label{itm:IIIinf}
The sequence $\{[X_n]\}_{n\in \N}$ is relatively compact in $\G(E^\infty)$.
\end{enumerate}
\end{thm}

Before proving this Theorem, recall that $\Cr rMk$ has the product topology with respect to the countable family of maps $\{Q_\ell^*\}_{\ell\in \N}$, defined as in \eqref{prodembballeq}.  It follows that a subset $\mathscr{A}\subset \Cr rMk$ is relatively compact if and only if $Q_{\ell}^*\mathscr{A}\subset \Cr r{\D^m}k$ is relatively compact for all $\ell$. In particular, if $r<\infty$, given constants $A_\ell >0$, the set
\be 
\mathscr{A}^r=\left\{f\in \Cr rMk \colon \|f\|_{Q_\ell,r+1}\le A_\ell \ \forall \ell\right\}
\ee
is compact in $\Cr rMk$. Similarly, given $A_\ell^r>0$ for all $r, \ell\in \N$, the set 
\be 
\mathscr{A^\infty}=\left\{f\in \Cr \infty Mk \colon \|f\|_{Q_\ell,r}\le A_\ell^r \ \forall r,\ell\right\}
\ee
is compact in $\Cr \infty Mk$. An important thing to note here is that every compact set in $\Cr \infty Mk$ is contained in a set of the form $\mathscr{A}^\infty$, while the analogous fact is not true when $r$ is finite.
\begin{proof}[Theorem \ref{rproper}]
Let $Q\colon D\hookrightarrow M$ be the embedding of a compact subset $D\subset \R^m$. Then we can cover $D$ with a finite family of disks $D_1,\dots,D_N$ such that $Q(D_i)\subset \text{int}(Q_{\ell_i}(\D^m))$ for some $\ell_1,\dots,\ell_N$. It follows that there exists a constant $c>0$ such that $\|f\|_{Q,r+1}\le c\sum_{i=1}^N\|f\|_{Q_{\ell_i},r+1}$ for all $i=1,\dots, N$ and for all $f\in \Cr {r+1}Mk$.

 By applying Theorem \ref{Gaussinequality} to each $Q_\ell$, we get the inequality 
\be 
\begin{aligned}
\sup_n\E\left\{ \|X_n\|_{Q,r+1} \right\}
&\le c\sum_{i=1}^N\sup_n\E\left\{ \|X_n\|_{Q_{\ell_i},r+1} \right\}
\\
&\le c \sum_{i=0}^N C(Q_{\ell_i},r+1,k) \sqrt{\sup_{n}\|K_{X_n}\|_{Q_\ell\times Q_\ell,(r+2,r+2)}}
\\
&<+\infty.
\end{aligned}
\ee
This proves the implication \ref{itm:Ir} $\implies$ \ref{itm:IIr}.

By Prohorov's Theorem (see \cite[Theorem 5.2]{Billingsley}), to prove the second implication, it is sufficient to show that $\{[X_n]\}_n$ is tight in $\G(E^{r})$, i.e. that for every $\e>0$ there is a compact set $\mathscr{A}\subset E^{r}$, such that $\P(
X_n\in\mathscr{A})> 1-\e$ for any $n\in \N$.

Fix $\e>0$. By \ref{itm:IIr}, the number $A_\ell= {\sup_{n}\E\left\{\|X_n\|_{Q_\ell,r+1}\right\}}$  is finite for each $\ell\in\N$.
Thus, we can consider the compact subset $\mathscr{A}\subset \Cr rMk$ defined as follows,
\be 
\mathscr{A}=\left\{f\in \Cr r Mk \colon \|f\|_{Q_\ell,r+1}\le \frac{2^{\ell+1}}{\e}A_\ell, \ \forall \ell\in \N\right\}.
\ee
By subadditivity and Markov's inequality we have that for all $n\in \N$:
\be 
\begin{aligned}
\P\{X_n\notin \mathscr{A}\} & \le \sum_{\ell\in \N}\P\left\{\|X_n\|_{Q_\ell,r+1}>\frac{2^{\ell+1}}{\e}A_\ell\right\} \\
& \le \sum_{\ell\in \N} \frac{\e}{2^{\ell+1}}\cdot \frac{\E\{\|X_n\|_{Q_\ell,r+1}\}}{A_\ell} \\
& \le \sum_{\ell\in\N}2^{-(\ell+1)}\e=\e.
\end{aligned}
\ee
We conclude that $\{[X_n]\}_n$ is tight.
\end{proof}

%

\subsection{Proof that $\K^\infty$ is a closed topological embedding} We already know that $\K^\infty$ is injective and continuous. To prove that it is a closed topological embedding it is sufficient to show that $\K^\infty$ is proper: both $\G(E^{\infty})\subset \mathscr{P}(E^\infty)$ and $E^\infty$ are metrizable spaces, and a proper map between metrizable spaces is closed. 

Let $\mathscr{A}\subset \Cr \infty {M\times M}{k\times k}$ be a compact set; then for any $Q: D\hookrightarrow M$ embedding of a compact subset $D\subset \R^m$ and for every $r\in \N$, it holds
\be 
\sup_{K\in \mathscr{A}}\|K\|_{Q\times Q,r }<\infty .
\ee
Therefore Theorem \ref{infproper} implies that the closed subset $(\K^\infty)^{-1}(\mathscr{A})\subset\G(E^\infty)$ is also relatively compact, hence it is compact.
\subsection{Proof of Theorem \ref{thm:2}}By Theorem \ref{thm:1}, if $K_d\xrightarrow{\mathcal{C}^\infty}K$ then $\mu_d\nrw\mu$. Observe also that, by definition for every $A\subset E^\infty$:
\be \P(X\in A)=\mu(A)\quad \textrm{and}\quad \P(X_d\in A)=\mu_d(A).\ee
Consequently \eqref{eq:limitprob} follows from Portmanteau's theorem (see \cite[Theorem 2.1]{Billingsley}).
\subsection{Addendum: a ``counter-theorem''}
It is possible to improve Theorem \ref{Gaussinequality} in order to control $\E\{\|X\|_{Q,r}\}$ with a $(r+\a,r+\a)$ Holder norm of the covariance function, if the latter is finite for some $\a\in (0,1)$ (see \cite[Sec. A.9]{NazarovSodin2}). But there is no way to get such an estimate with $\a=0$, as the following example shows.
\begin{example}\label{escempio}
Let $D\subset \R^{m}$ compact with non empty interior. We now construct a sequence of smooth GRFs $X_n\in \g 0 D{}$, with $\|K_{X_n}\|_{D,0}\to 0$, such that
\be
\liminf_{n\to \infty}\E\{\|X_n\|_{D, 0}\}\ge 1 .
\ee

Let $I^{(n)}_1,\dots, I^{(n)}_{n^2}$ be disjoint open sets in $D$ (their size doesn't matter), containing points $x^{(n)}_1,\dots, x^{(n)}_{n^2}$. Let $\f^{(n)}_1,\dots, \f^{(n)}_{n^2}$ be smooth functions $\varphi_i^{(n)}:D\to [0,1]$ such that $\f^{(n)}_i$ is supported in $I^{(n)}_i$ and $\f^{(n)}_i(x^{(n)}_i)=1$ (see Figure \ref{fig:compact}).
Let $\gamma_i$ be a countable family of independent standard Gaussian random variables.
Let $a_n\in \R$ be the real number such that $\P\{|\gamma|>a_n\}=\frac{1}{n}$, for any $\gamma\sim N(0,1)$, hence $a_n \to +\infty$.
Define
\be 
X_n=\frac{1}{a_n}\sum_{i=1}^{n^2}\gamma_i \f^{(n)}_i \in \g 0D{}.
\ee
Then $K_{X_n}(x,y)=\frac{1}{a_n^2}\f^{(n)}_i(x)\f^{(n)}_j(y)$ for some $i=i_x,j=j_x$, thus $\|K_{X_n}\|_{D,0}\to 0$.

We can now estimate the probability that the $\mathcal{C}^0$-norm of $X_n$ is small by
\be\label{Pmisfacteq}
\begin{aligned}
\P\{\|X_n\|_{D,0}< 1\} &\le \P\left\{\max_{i=1,\dots, n^2}|X_n(x^{(n)}_i)|< 1\right\} \\
&= \P\{|\gamma|< a_n\}^{n^2} \\
&= \left(1-\frac 1n\right)^{n^2} \xrightarrow[{n\to \infty}]{} 0.
\end{aligned}
\ee
Consequently, by Markov's inequality
\be \label{Emisfacteq}
\liminf_{n\to \infty}\E\left\{\|X_n\|_{D,0}\right\}\ge \liminf_{n\to \infty}\P\left\{\|X_n\|_{D,0}\ge 1\right\}=1.
\ee
\begin{figure}\begin{center}
\includegraphics[width=0.7\textwidth]{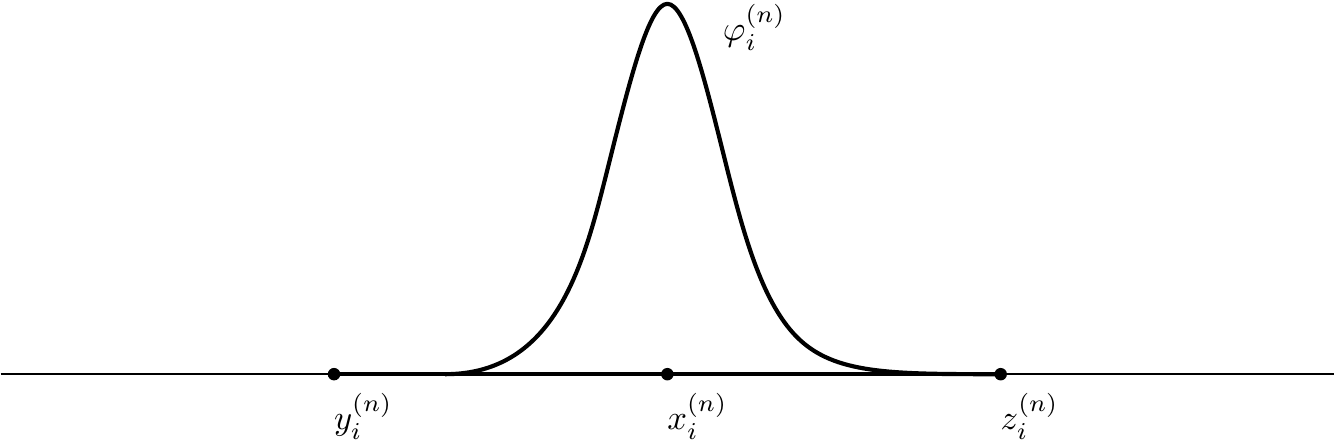}
	\caption{The function $\varphi_i^{(n)}$ from Example \ref{escempio} is supported on the interval $I_i^{(n)}=(y_i^{(n)}, z_i^{(n)})$ and takes value $1$ at $x_{i}^{(n)}.$}\label{fig:compact}
\end{center}
	\end{figure}
\end{example}

Note that the function $K(x,y)=0$ is the covariance function of the GRF $X_0$, which corresponds to the probability measure $\delta_{0}\in \G(E^0)$ concentrated on the zero function $0\in \Cr 0D{}$. Since $\P\left\{\|X_0\|_{D,0} \ge 1\right\}=0$, equation \eqref{Emisfacteq} proves also that $[X_n]$ does not converge to $[X_0]$ in $\G(E^0)$, even if $K_{X_n}\to K_{X_0}$ in $\Cr 0D{}$.

The previous Example \ref{escempio} can be generalized to prove the following result, which shows that the condition $r=\infty$ in the second part of the statement of Theorem \ref{thm:1} is necessary.

\begin{thm}\label{thm:counter}
If $r$ is finite, the map $(\mathcal{K}^r)^{-1}$ is not continuous.
\end{thm}
\begin{proof}
Construct $X_n\in \g 0D{}$ as in Example $\ref{escempio}$, with $D=[0,1]\subset \R$. Since $X_n$ is a sum of functions with compact support, we can as well consider $X_n$ as a random element of $\Cr 0\R{}$. So that $K_{X_n}\to 0$ in $\Cr 0{\R\times \R}{}$, because their support is contained in $D\times D$, but $X_n \centernot{\nrw} 0$. 

Let $c\notin D$ and let $Y_n$ be the GRF defined as
\be 
Y_n(\cdot)=\int_c^{(\cdot)}\int_c^{s_r}\dots \int_c^{s_2}X_n(s_1)ds_1 \dots ds_{r}.
\ee
Then $Y_n\in \g r\R{}$ (indeed $Y_n$ is a smooth GRF), and $\frac{d^r}{dx^r}Y_n=X_n$. Moreover,
\be 
\frac{d^{2r}}{dx^rdy^r}K_{Y_n}=K_{X_n}\to 0
\ee
in $\Cr 0{\R\times\R}{}$ and $K_{Y_n}=0$ in a neighborhood of $(c,c)$, therefore $K_{Y_n}\to 0$ in $\Cr {r,r}{\R\times\R}{}$, but $Y_n\not\nrw 0$.

Let $M$ be a smooth manifold of  dimension $m$. Denoting $(t,x)\in \R\times \R^{m-1}=\R^m$, we define a smooth function $\rho\colon \R^{m}\to [0,1]$ with compact support and such that $\rho=1$ in a neighborhood of the set $[0,1]\times \{0\}$. 
Let $j\colon \R^m\to M$ be any embedding and fix $v\in \R^k$. Define the transformation $T\colon \mathcal{C}^r(\R,\R)\to \Cr r Mk$ such that $f\mapsto g=Tf$, where
\be 
\begin{aligned}
g(j(t,x))=\rho(t,x)f(t)v, && \\
g(p)=0 \quad \text{if }p\notin j(\R^m).
\end{aligned}
\ee
Since $\rho$ has compact support, $T$ is continuous for all $r\in\N\cup\{\infty\}$
, so that $Z_n=TY_n$ is a well defined smooth GRF with compact support on $M$. Thanks to the continuity of $T$, we have that $K_{Z_n}\to 0$ in $\Cr {r,r}M{k\times k}$, but $Z_n \centernot{\nrw} 0$ in $\g rMk$ because $
\frac{d^r}{dt^r}Z_n(j(t,0))=X_n(t)
$
for every $t\in[0,1]$ and $X_n|_{D}\not\nrw 0$.
\end{proof}

%
%
%
%
%
%
%
%

\section{Proof of Theorem \ref{thm:dtgrf:3}}\label{sec:PT3}

Given a Gaussian field $X=(X^1, \ldots, X^k)\in \mathcal{G}^{r}(M, \R^k)$ defined on a probability space $(\Omega, \mathfrak{S}, \PP)$, we consider the Hilbert space $\Gamma_X$ defined by:
\be \Gamma_X=\overline{\textrm{span}\{X^j(p),\,p\in M,\, j=1, \ldots, k\}}^{L^2(\Omega, \mathfrak{S}, \PP)}.\ee
Since $X$ is a Gaussian field, all the elements of $\Gamma_X$ are gaussian random variables (and viceversa).
By Lemma \ref{Ltwocont}, we know that the function $\gamma_X\colon M\to L^2\Prob$ defined as in equation \eqref{eq:gammaX} is of class $\mathcal{C}^r$, therefore we can define a linear map $\rho_X:\Gamma_X\to E^r$ by:
\be \rho_X(\gamma)=\E\left(X(\cdot)\gamma\right)=\left(\langle \gamma_X^1(\cdot), \gamma\rangle_{L^2(\Omega, \mathfrak{S}, \PP)},\dots,\langle \gamma_X^k(\cdot), \gamma\rangle_{L^2(\Omega, \mathfrak{S}, \PP)}\right).\ee
\begin{prop}\label{propo:injection}The map $\rho_X\colon \Gamma_X\to E^r$ is a linear, continuous injection.\end{prop}
\begin{proof}Let $\gamma\in\Gamma_X$ and assume that $\rho_X(\gamma)=0$. Then $\langle\gamma, X^j(p)\rangle_{L^2}=0$ for all $p\in M$ and $j=1,\ldots,k$, so that $\gamma \in \Gamma_X^\perp$, thus $\gamma=0$. This proves that $\rho_X$ is injective.

By linearity, it is sufficient to check continuity at $\gamma=0$.
Let $Q\colon D\hookrightarrow M$ be the embedding of a compact set $D\subset \R^m$. If $r$ is finite, we have
\begin{align}
\|\rho_X(\gamma)\|_{Q, r} &=\sup_{|\a|\le r, x\in D} |\E\left\{ \de_\a(X\circ Q)(x)\gamma\right\}|  \\
&\le \sqrt{k} \sup_{|\a|\le r, x\in D} \E\left\{|\de_\a(X\circ Q)(x)|^2\right\}^{\frac12}\|\gamma\|_{L^2} \\
&=\sqrt{k}\sup_{|\a|\le r, x\in D} \left(\sum_{j=1}^k \de_{(\a,\a)}(K_X^{j,j}\circ Q\times Q) (x,x)\right)^{\frac12}\|\gamma\|_{L^2}\\
&\le 
k\left(\|K_X\|_{Q\times Q, (r,r)}
\right)^{\frac12} \|\gamma\|_{L^2}.
\end{align}
Therefore $\lim_{\gamma \to 0}\|\rho_X(\gamma)\|_{Q, r}=0$ for every $Q$, hence $\rho_X$ is continuous. For the case $r=\infty$, it is sufficient to note that continuity with respect to $E^r$ for every $r$, implies continuity with respect to $E^\infty$.
\end{proof}
 
 \begin{prop}\label{prop:CamMa} The image of $\rho_X$ coincides with the Cameron-Martin space (see \cite[p. 44, 59]{bogachev}) of the measure $[X]$ and we denote it by $\mathcal{H}_X$. 
 \end{prop}
 \begin{proof}
 According to \cite[Lemma 2.4.1]{bogachev} $\mathcal{H}_X$ is the set of those $h\in E^r$ for which there exists a $T\in (E^r)^*$ such that
 \be\label{eq:CM}
 L(h)=\E\{T(X)L(X)\}, \quad \text{for all $L\in (E^r)^*$}.
 \ee
 Observe that the map $T\mapsto T(X)$ defines a surjection $(E^r)^*\to\Gamma_X$, because every continuous linear functional $T\in (E^r)^*$ can be approximated by linear combinations of functionals of the form $\delta_p^j:f\mapsto f^j(p)$ (see Theorem \ref{thm:deltadense} in Appendix \ref{app:dual}). For the same reason, condition \eqref{eq:CM} is equivalent to the existence of $\gamma\in \Gamma_X$ such that
 \be 
 h^j(p)=\E\{\gamma X^j(p)\}\}, \quad \text{for all $p\in M$ and $j=1,\dots, k$},
 \ee
 that is, by definition, $h=\rho_X(\gamma)$. Thus $\mathcal{H}_X=\rho_X(\Gamma_X)$.
 \end{proof}

Observe that $\mathcal{H}_X$ contains all the functions $h_p^j=\rho_X(X^j(p))$ satisfying equation\eqref{suppbasiseq} in Theorem \ref{thm:dtgrf:3}.
Moreover, it carries the Hilbert structure induced by the map $\rho_X$, which makes it isometric to $\Gamma_X$. It follows that $\mathcal{H}_X$ is the Hilbert completion of the vector space span$\{h_p^j\colon p\in M,\ j=1,\dots,k\}$, endowed with the scalar product
\be 
\langle h_p^j,h_q^\ell\rangle_{\mathcal{H}_X}\doteq \left\langle X^j(p),X^\ell(q)\right\rangle_{L^2}=K_X^{j,\ell}(p,q).
\ee
Now, Theorem \ref{thm:dtgrf:3} follows from \cite[Theorem 3.6.1]{bogachev}:
\be\label{eq:sputoCM}
\text{supp}(X)=\overline{\mathcal{H}_X}^{\Cr rMk}.
\ee
In Appendix \ref{app:rgrf} (equation \eqref{eq:sputoCMapp}) the reader can find a proof of \eqref{eq:sputoCM} adapted to our language.
\begin{remark}
Note that the Hilbert space $\mathcal{H}_X$ depends only on $K_X$, thus it depends only on the measure $[X]$.
\end{remark}

%
%
%
%
%
%
%
%

\section{Proof of Theorems \ref{transthm2} and \ref{thm:transthm}}\label{sec:transversality}
\subsection{Transversality}
We want to prove some results analogous to Thom's Transversality Theorem (see \cite[Section 3, Theorem 2.8]{Hirsch}) in our probabilistic setting. We first recall the definition of transversality. Let $f\colon M\to N$ be a smooth map, $W\subset N$ a submanifold and $K\subset M$ be any subset. Then we say that $f$ is transverse to $W$ on $K$ and write $f\transv_K W$, if and only if for every $x\in K\cap f^{-1}(W)$ we have:
\be 
df_x(T_xM)+T_{f(x)}W=T_{f(x)}N.
\ee
We will simply write $f\transv W$ if $K=W$.
We recall the following classical tool, usually called the \emph{Parametric Transversality Theorem}.
\begin{thm}[Theorem 2.7 from \cite{Hirsch}, Chapter 3]
\label{parametric transversality}
Let $g\colon P\times F\to N$ be a smooth map between smooth manifolds of finite dimension.
Let $W\subset N$ be a smooth submanifold and $K\subset P$ be any subset. If $g \transv_{K\times F} W$, then $g(\cdot, f)\transv_K W$ for almost every $f\in F$.
\end{thm}
In our context we prove the following infinite-dimensional, probablistic version of Theorem \ref{parametric transversality}.
\begin{thm}\label{transthm}
Let $F\subset E^r$ such that $F=\spt(X)$ for some $X\in \g r Mk$, with $r\in\N\cup\{\infty\}$. Let $P,N$ be smooth manifolds and $W\subset N$ a submanifold. Assume that $\Phi\colon P\times F\to N$ is a ``smooth''\footnote{Here by ``smooth'' we mean that:
\begin{enumerate}
    \item the map $\Phi$ is smooth when restricted to finite dimensional subspaces;
    \item the linear map $(p, f, v)\mapsto D_{(p, f)}\Phi v= D_{(p, f)}\left(\Phi|_{\textrm{span}\{f,v\}}\right) v$ is continuous in all its arguments.
    \end{enumerate}} map such that $\Phi\transv W$. Then
\be 
\P\{\Phi(\cdot, X)\transv W\}=1.
\ee

\end{thm}
A particular case in which we can apply Theorem \ref{transthm} is when $P=M$,  $N=J^r=J^r(M,\R^k)$, $r=\infty$ and $\Phi$ is the jet-evaluation map
\be\label{eq:jr}
j^r\colon M \times E^\infty  \to J^r, \qquad (p,f)\mapsto j^r_pf.
\ee
It is straightforward to see that this map is ``smooth'' in the sense of the statement of Theorem \ref{transthm}.
\begin{proof}
(In order to simplify the notations, we denote by $\phi(X)$ the map $p\mapsto \Phi(p, X)$.)

First we show that we can assume $W$ to be compact (possibly with boundary). Indeed let $W=\cup_{k\in \N} W_k$, such that $W_k$ is compact. Then $\Phi\transv W_k$ for any k, and 
\be 
\P\{\phi(X)\transv W\}\ge 1- \sum_{k\in\N}\left(1-\P\{\phi(X)\transv W_k\}\right).
\ee
\begin{claim}\label{claim 55} Moreover, it is sufficient to prove the following weaker statement.
\end{claim}
\begin{stella}\label{claim:HeyStella}
 For all $p\in P$ and $x\in F$ there are neighborhoods $Q_p$ of $p$ in $P$ and $U_x$ of $x$ in $E^r$ such that: \be \P
\big\{\phi(X)\transv_{\overline{Q}_p}W\big|X\in U_x\big\}:=\frac{\P\left(\{\phi(X)\transv_{\overline{Q}_p}W\}\cap\{X\in U_x\}\right)}{\P\left(\{X\in U_x\}\right)}=1.\ee
\end{stella}

Assume that $(*)$ is true, then there exists a countable open cover of $P\times F$ of the form $Q_k\times U_l$ such that $\P\{\phi(X)\transv_{\overline{Q}_k}W|X\in U_l\}=1$, i.e. the probability that $X\in U_l$ and  $\phi(X)$ is not transverse to $W$ at some point $p\in \overline{Q}_k$ is zero. 
Thus
\be 
\P\{\phi(X)\not\transv W\}\le\sum_{l,k}\P\left\{\phi(X)\not\transv_{\overline{Q}_k}W,X\in U_l\right\}=0,
\ee
hence Claim \ref{claim 55} is true.

Let us prove $(*)$. Let $p\in P$ and $x\in F$. Since $W\subset N$ is closed, if $\Phi(p,x)\notin W$, then $\Phi(q,\tilde{x})\notin W$ for all $q$ in a compact neighborhood $Q$ of $p$ and $\tilde{x}$ in some neighborhood $N_x$ of $x$ in $E^r$, so that, in particular $\P\{\phi(X)\transv_{Q}W|X\in N_x\}=1$. 

Assume now that $\Phi{(p,x)}=\theta\in W$, then by hypothesis we have that 
\be
D_{(p,x)}\Phi(T_pP+F) + T_\theta W=T_\theta N,
\ee
hence there is a finite dimensional space $F_0=\textrm{span}\{f_1,\dots,f_a\}\subset F$ such that 
\be
D_{(p,x)}\Phi\left(T_pP+F_0\right) + T_\theta W=T_\theta N.
\ee
Note that $F_0=T_xF_x$, where $F_x=x+\text{span}\{f_1,\dots,f_a\}$. Therefore 
$
\Phi|_{P\times F_x}\transv_{(p,x)} W
$
(here we are in a finite dimensional setting). Moreover, there is a compact neighborhood  $p\in Q\subset P$ and a $\e>0$ such that
\be \label{difftrasveq}
\Phi|_{P\times F_x}\transv_{Q\times D_\e} W.
\ee
where $D_\e=D_\e(x,f)=\{x+f_1u^1+\dots +f_au^a\colon u\in \R^a,\ |u|\le\e\}$.
Observe that the set of $(a+1)-$tuples $(x,f)=(x,f_1,\dots,f_a)\in F\times F^a$ for which \eqref{difftrasveq} holds (with fixed $\e$), form an open set, indeed the map 
\be 
\tau\colon F\times F^a\to \mathcal{C}^\infty(P\times \R^a,N),\qquad \tau(x,f):(p,u)\mapsto \Phi(p,x+fu)
\ee
is continuous and the set $\Theta=\{T\in \mathcal{C}^\infty(P\times \R^a,N)\colon T\transv_{Q\times D_\e}W\}$ is open in the codomain because $Q\times D_\e $ is compact and $W$ is closed (see \cite[p. 74]{Hirsch}); therefore 
\be 
\tau^{-1}(\Theta)=\{(x,f)\in F\times F^a\colon \text{\eqref{difftrasveq} holds}\}
\ee
is open. It follows that there is an open neighborhood $V_x$ of $x$ and an $h\in (\mathcal{H}_X)^a$ such that \eqref{difftrasveq} holds with $(\tilde{x},h)$ for any $\tilde{x}\in V_x$, indeed the Cameron-Martin space $\mathcal{H}_X$ is dense in $F$ (see equation \eqref{eq:sputoCM}). 

Define $\Lambda=\{e\in E^r \colon \phi(e)\transv_{Q}W\}$. 
By Theorem \ref{parametric transversality} we get that if $\tilde{x}\in V_x$, then $\phi(\tilde{x}+hu)\transv_{Q}W$, equivalently $(\tilde{x}+hu)\in \Lambda$, for almost every $|u|\le\e$. Denote by $1_\Lambda\colon E^r\to \{0,1\}$ the characteristic function of the (open) set $\Lambda$.
Using the Fubini-Tonelli theorem, we have
\be 
0=\int_{V_x}\left(\int_{\D^m_\e} 1-1_\Lambda(\tilde{x}+hu)du\right)d[X](\tilde{x})=\int_{\D^m_\e}\P\{X+hu\notin \Lambda, X\in V_x\}du.
\ee
hence $\P\{X+hu\in (V_x+hu) \- \Lambda\}=0$ for almost every $|u|\le\e$. Let $u$ be also so small that $x\in V_x+hu$. Then, taking $U_x=V_x+hu$, we have that $\P\{X+hu\in U_x\- \Lambda\}=0$. Since $hu\in \mathcal{H}_X$, the Cameron-Martin theorem (see \cite[Theorem 2.4.5]{bogachev}) implies that 
$[X]$ is absolutely continuous with respect to $[X+hu]$ and consequently $\P\{X\in U_x\- \Lambda\}=0$. In other words, $\P\{\phi(X)\transv_{Q}W|X\in U_x\}=1$, that proves $(*)$.
\end{proof}

We give know a criteria to check the validity of the hypothesis of Theorem \ref{transthm}, without necessarily knowing the support of $X$. Before that, let's observe that the canonical map $J^r:=J^r(M,\R^k)\to M$ is a smooth vector bundle over $M$ with fiber $J^r_p$, so that $T_{\theta}J_p^r$ is canonically identified with $J^r_p$ itself, for all $\theta\in J^r_p$. 
\begin{prop}
Let $X\in\g rMk$ and $F=\spt(X)$. Let $W\subset J^r$ be a smooth submanifold and fix a point $p\in M$. The next conditions are related by the following chain of implications: $
\text{\ref{itm:transjetA}}\Longleftarrow \text{\ref{itm:transjetB}}\Longleftarrow \text{\ref{itm:transjetC}}\iff \text{\ref{itm:transjetD}}
$.
\begin{enumerate}[$(1)$]
\item\label{itm:transjetA}  $j^r|_{M\times F}\transv_{\{p\}\times F} W$, where $j^r\colon M\times E^\infty\to J^r$ is the map defined in \eqref{eq:jr}; 
\item\label{itm:transjetB}  the vector space $\spt(j_p^rX)$ is transverse to  $(T_\theta W\cap T_\theta J^r_p)$ in $J^r_p$, for all $\theta \in j^r_p(F)\cap W$;
\item\label{itm:transjetC}  $\spt(j^r_pX)=J_p^r$;
\item\label{itm:transjetD}  given a chart of $M$ around $p$, the matrix below has maximal rank.
\be\label{covjeteq}
\left(\de_{(\a,\beta)}K_X(p,p)\right)_{|\a|,|\beta|\le r}.
\ee
\end{enumerate}
\end{prop}
\begin{proof}
\ref{itm:transjetB}$\implies$ \ref{itm:transjetA}. Let $f\in F$ such that $\theta=j^r_pf\in W$. Under the identification $T_\theta J^r_p\cong J^r_p$, mentioned above, we have
\be 
\left(D_{(p,f)}j^r\right)(0,g)= \frac{d}{dt}\Big|_0 j^r_p(f+tg)\cong j^r_pg,
\ee
so that, for any $x\in F$, $D_{(p,f)}j^r(T_xF)=j^r_p(F)=\spt(j^r_pX)$. Then for all $(p,f)\in (j^r)^{-1}(W)\cap M\times F$, we have
\be 
\begin{aligned}
D_{(p,f)}j^r\left(T_{(p,f)}(M\times F)\right)+T_\theta W &\supset 
D_{(p,f)}j^r(T_pM)+ \spt(j^r_pX)+T_\theta W\cap J^r_p=\\
& =D_{(p,f)}j^r(T_pM)+ J_p^r=\\ &=T_\theta J^r.
 \end{aligned}
\ee
The last equality follows from the fact that the map $j^rf$ is a section of the bundle $J^r\to M$. 

\ref{itm:transjetC}$\implies$ \ref{itm:transjetB}.  Obvious.

\ref{itm:transjetC}$\iff$ \ref{itm:transjetD}. Any chart around $p$ defines a linear isomorphism 
\be 
J_p^r\to \R^{\{\a\colon |\a|\le r\}}, \quad j^r_pf\mapsto \left(\de_af(p)\right)_\a.
\ee
With this coordiante system, the covariance matrix of the Gaussian random vector $j^r_pX$, is exactly the one in $\eqref{covjeteq}$, hence the result follows from the fact that the random Gaussian vector $j^rf$ has full support if and only if its covariance matrix is nondegenerate.
\end{proof}
Given $X\in \g \infty Mk$, we can also consider it as an element of $ \g rMk$ such that $\P\{X\in \Cr\infty Mk\}=1$. We use the notation $\spt_{\mathcal{C}^r}(X)\subset E^r$ to denote the support of the latter, namely 
\be 
\spt_{\mathcal{C}^r}(X)=\overline{\mathcal{H}_X}^{\mathcal{C}^r}.
\ee
\begin{cor}\label{transcor}
Let $X\in \g \infty Mk$, such that $\mathrm{supp}_{\mathcal{C}^{r}}(X)=\Cr {r} Mk$. Then for every submanifold $W\subset J^r(M,\R^k)$, one has
\be 
\P\{j^rX\transv W\}=1.
\ee
\end{cor}
\begin{proof}
Clearly $X$ satisfies for every $p\in M$ condition \ref{itm:transjetC} of the proposition above, hence the hypotheses of Theorem \ref{transthm} are satisfied for every $W\subset J^r(M,\R^k)$.
\end{proof}
\begin{subappendices}
\section{The dual of $E^r$}\label{app:dual}
The purpose of this section is to fill the gap between Definition \ref{def:gauss} of a Gaussian measure on the space $E^r=\mathcal{C}^r(M,\R^k)$ and the abstract definition of Gaussian measures on topological vector spaces, for which we refer to the book \cite{bogachev}. As we already mentioned in Remark \ref{rem:bridge}, the two definitions coincide. In order to see this clearly, one simply have to  understand the topological dual of $E^r$, that is the space $(E^r)^*$ defined as follows.

Let $(E^r)^*$ be the set of all linear and continuous functions $T\colon E^r \to \R$, endowed with the weak-$*$ topology, namely the topology induced by the inclusion $(E^r)^*\subset \R^{E^r}$, when the latter is given the product topology.
\begin{lemma}\label{duallemma}
Let $T\in (E^r)^*$, with $r\in\N\cup\{\infty\}$. There exists a finite set $\mathscr{Q}$ of embeddings $Q\colon \mathbb{D}^m\hookrightarrow M$, a constant $C>0$ and a finite natural number $s\le r$, such that
\be\label{eq:dualiT}
|T(f)|\le C\max_{Q\in\mathscr{Q}}\|f\|_{Q,s},\footnote{The seminorm $\|\cdot\|_{Q,r}$ is defined as in equation \eqref{eq:seminormQR}.}
\ee
for all $f\in E^r$.
As a consequence, denoting $K=\cup_{Q\in\mathscr{Q}}Q(\mathbb{D})$, there is a unique $\hat{T}\in (\Cr sK k)^*$ such that $T(f)=\hat{T}(f|_{K})$ for all $f\in E^r$. 
\end{lemma}
Let $K\subset M$ be as in Lemma \ref{duallemma}. The vector space $\Cr sKk$ is, by definition, the image of the restriction map 
\be 
\Cr sMk \to \Cr 0K k, \quad f\mapsto f|_K.
\ee
Denote by $\Omega=\text{int}(K)\subset M$. Notice that the derivatives, of order less than $s$, of a function $f\in\Cr sKk$ are well defined and continuous at points of $\overline{\Omega}$, thus when $K=\overline{\Omega}$, we have a well defined continuous function \be j^sf\colon K\to J^s(K,\R^k)=\{j^s_pf\in J^s(M,\R^k)\colon p\in K\}.\ee In this case ($\overline{\text{int}(K)}=K$) we endow the space $\Cr sKk$ with the topology that makes Lemma \ref{Lemma:14} true with $M=K$. 
Such topology is equivalent to the one defined by the norm $\|\cdot\|_{K,s}$ below (it depends on $\mathscr{Q}$), with which $\Cr sKk$ becomes a Banach space:
\be\label{normKeq}
\|f\|_{K,s}=\max_{Q\in\mathscr{Q}}\|f\|_{Q,s}.
\ee
(Note that, if $K=\overline{\Omega}$, then $\|f\|_{K,s}$ depends only on $f|_{\Omega}$.)
\begin{remark}
When $M$ is an open subset $M\subset\R^m$ and $k=1$, the elements of $(E^\infty)^*$ are exactly the distributions with compact support (in the sense of \cite{schwartz1957}).
\end{remark}
\begin{proof}[Proof of Lemma \ref{duallemma}]
Let $Q_\ell\colon \mathbb{D}\hookrightarrow M$ be a countable family of embeddings such that $g_n\to 0$ in $E^r$ if and only if $\|g_n\|_{Q_\ell,s}\to 0$ for all $\ell\in\N$ and $s\le r$ (it can be constructed as in \eqref{prodembballeq}).
Assume that for all $n\in\N$ there is a function $f_n\in E^r$, such that
\be \label{dualabsurdeq}
|T(f_n)|> n\max_{\ell\le n}\|f_n\|_{Q_\ell,s_n},
\ee
where $s_n:=n$ if $r=\infty$, otherwise $s_n:=r$.
Then the sequence 
\be 
g_n=\frac{f_n}{n\max_{\ell\le N}\|f_n\|_{Q_\ell,s_n}}
\ee
converges to $0$ in $E^r$, indeed $\|g_n\|_{Q_\ell,s}\le \frac1N$ for any fixed $\ell\in \N$ and $s\le r$. Therefore, by the continuity of $T$, we get that $T(g_n)\to 0$. But $|T(g_n)|>1$ according to \eqref{dualabsurdeq}, so we get a contradiction. It follows that there exists $N$ such that for all $f\in E^r$ we have
\be
|T(f)|\le N\max_{\ell\le N}\|f\|_{Q_\ell,s_N}.
\ee
This proves the first part of the Lemma, with $\mathscr{Q}=\{Q_0,\dots,Q_N\}$, $C=N$ and $s=s_N$.

Define $\Omega=\text{int}(K)$. Note that, since $Q(\text{int}(\D^m))\subset \Omega$ for all $Q\in\mathscr{Q}$, if $p\in K\-\Omega$, then $p\in Q(\de\mathbb{D}^m)$ for some $Q\in \mathscr{Q}$ and therefore $p\in \overline{Q(\text{int}(\mathbb{D}^m))}\subset \overline{\Omega}$. This proves that $K=\overline{\Omega}$, hence $\Cr sKk$ is a Banach space with the norm \eqref{normKeq}. 

Let $f,g\in E^r$ be such that $f|_K=g|_K$, then \be 
\begin{aligned}
|T(f)-T(g)|=|T(f-g)|\le C\max_{Q\in\mathscr{Q}}\|f-g\|_{Q,s} =C\|f|_K-g|_K\|_{K,s}=0.
\end{aligned}
\ee
It follows that the function $L\colon \Cr rKk\to \R$ such that $L(f|_\Omega)=T(f)$ for all $f\in E^r$, is well defined and continuous with respect to the norm $\|\cdot\|_{K,s}$. Since $\Cr rKk$ is dense in $\Cr sKk$, there is a unique way to extend $L$ to a bounded linear functional on $\Cr sKk$, that we call $\hat{T}.$
\end{proof}
We recall the following classical theorem from functional analysis (see \cite[Theorem 1.54]{ambrofuscopalla}), which we can use to give a more explicit description of $(E^r)^*$.
\begin{thm}[Riesz's representation theorem]\label{thm:riesz}
Let $K$ be a compact metrizable space. Let $\mathcal{M}(K)$ be the Banach space of Radon measures on $K$ (on a compact set it is the set of finite Borel signed measures), endowed with the total variation norm.
Then the map 
\be 
\mathcal{M}(K)\to \left(\mathcal{C}(K)\right)^*, \quad \mu\mapsto \int_K(\cdot)d\mu 
\ee
is a linear isometry of Banach spaces.
\end{thm}
\begin{thm}\label{thm:Mloc}
Let $\mathcal{M}^r_{loc}$ be the set of all $T\in (E^r)^*$ of the form 
\be 
T(f)=\int_{\mathbb{D}^m}\de_\a(f^j\circ Q) d\mu,
\ee
for some embedding $Q\colon\mathbb{D}^m\hookrightarrow M$, some finite multi-index $\a\in\N^m$ such that $|\a|\le r$, some $j\in \{1,\dots,k\}$ and some $\mu\in \mathcal{M}(\mathbb{D}^m)$.
Then $(E^r)^*=\text{span}\{\mathcal{M}^r_{loc}\}$.
\end{thm}
\begin{proof}
Let $T\in (E^r)^*$ and let $\mathscr{Q}$, $s$, $K$,  $C$ and $\hat{T}$ defined as in lemma \ref{duallemma}. Consider the topological space
\be 
D=\mathbb{D}^m\times \mathscr{Q}\times \{\a\in \N^{m}\colon |\a|\le s\}\times\{1,\dots , k\}.
\ee
$D$ is homeomorphic to a finite union of disjoint copies of the closed disk, therefore it is compact and metrizable. There is a continuous linear embedding with closed image 
\be 
\begin{aligned}
\mathcal{J}^s\colon \Cr sKk \hookrightarrow  \mathcal{C}(D), \quad
\mathcal{J}^sf(u,Q,\a,j)=\de_\a (f^j\circ Q)(u).
\end{aligned}
\ee
Indeed $\|\mathcal{J}^sf\|_{\mathcal{C}(D)}\le \|f\|_{K,s}\le \sqrt{k}\|\mathcal{J}^sf\|_{\mathcal{C}(D)}$, if $\|\cdot\|_{K,s}$ is defined as in \eqref{normKeq}. By identifying $\Cr sKk$ with its image under $\mathcal{J}^s$, we can extend $\hat{T}$ to the whole $\mathcal{C}(D)$, using the Hahn-Banach theorem and the extension can then be represented by a Radon measure $\mu\in \mathcal{M}(D)$, by the Riesz theorem \ref{thm:riesz}. 

Denote by $\mu_{Q,\a,j}\in \mathcal{M}(\mathbb{D}^m)$ the restriction of $\mu$ to the connected component $\mathbb{D}^m\times\{Q\}\times \{\a\}\times\{j\}$. Let $T_{Q,\a,j}$ be the element of $\mathcal{M}^r_{loc}$ associated with $Q$, $\a$, $j$ and $\mu_{Q,\a,j}$.
Then, for all $f\in E^r$, we have 
\be \begin{aligned}
T(f)&=\hat{T}(f|_{K})\\
&=\int_D \mathcal{J}^sfd\mu \\
&=\sum_{Q\in\mathscr{Q}, |\a|\le s, j=1,\dots,k}\int_{\mathbb{D}^m\times\{Q\}\times \{\a\}\times\{j\}}\mathcal{J}^sf d\mu \\
&=\sum_{Q\in\mathscr{Q}, |\a|\le s, j=1,\dots,k}\int_{\mathbb{D}^m}\de_\a(f^{j}\circ Q)d\mu_{Q,\a,j}\\
&=\sum_{Q\in\mathscr{Q}, |\a|\le s, j=1,\dots,k}T_{Q,\a,j}(f).
\end{aligned}
\ee
Therefore $T$ is the sum of all the $T_{Q,\a,j}$, thus $T\in \text{span}\{\mathcal{M}^r_{loc}\}$.
\end{proof}
We are now in the position of justifying Remark \ref{rem:bridge}. First, observe that the manifold $M$ is topologically embedded in $(E^r)^*$, via the natural association $p\mapsto \delta_p$. We denote by $\delta_M\subset (E^r)^*$ the image of the latter map (it is a closed subset). From this it follows that any abstract Gaussian measure on $E^r$ is also a Gaussian measure in the sense of Definition \ref{def:gauss}. The opposite implication is a consequence of the following Lemma, combined with the fact that the pointwise limit of a sequence of Gaussian random variable is Gaussian.
\begin{cor}\label{thm:deltadense}
$(E^r)^*=\overline{\text{span}\{\delta_M\}}$.
\end{cor}
\begin{proof}
By Theorem \ref{thm:Mloc}, it is sufficient to prove that $\mathcal{M}^r_{loc}\subset \overline{\text{span}\{\delta_M\}}$. To do this, we can restrict to the case $M=\mathbb{D}^m$, $Q=$id and $k=1$.

Observe  that any functional of the type $\delta_p\circ\de_\a$ belongs to $\overline{\text{span}\{\delta_M\}}$. This can be proved by induction on the order of differentiation $|\a|$: if $|\a|=0$ there is nothing to prove, otherwise we have 
\be 
\delta_u\circ\frac{\de}{\de u^j}\circ\de_{\a}= \lim_{n\to \infty}n\left(\delta_{u+\frac{1}{n}e_j}\circ\de_\a-\delta_u\circ\de_\a\right)\in \overline{\text{span}\{\delta_M\}}.
\ee
 Note also that any $T^r\in\mathcal{M}^r_{loc}$ is of the form $T^0\circ \de_\a$ for some $T^0\in \mathcal{M}^0_{loc}$ and $|\a|\le r$ and, together with the previous consideration, this implies that it is sufficient to prove the theorem in the case $r=0$, so that we can conclude with the following lemma.
 \begin{lemma}
Let $K$ be a compact metric space. The subspace $\text{span}\{\delta_K\}$ is sequentially dense (and therefore dense) in $\mathcal{M}(K)$, with respect to the weak-$*$ topology on $\mathcal{M}(K)=\mathcal{C}(K)^*$.
\end{lemma}
Let $\mu$ be a Radon measure on $K$. Define for any $n\in\N$ a partition $\{A^n_i\}_{i\in I_n}$ of $K$ in Borel subsets of diameter smaller that $\frac1n$ and let $a^n_i\in A^n_i$. Define $t_n\in \text{span}\{\delta_K\}$ as
\be 
t_n=\sum_{i\in I_n}\mu(A^n_i)\delta_{a^n_i}.
\ee
Given $f\in \mathcal{C}(K)$, we have 
\be 
\begin{aligned}\label{heinecantoreq}
\left|\int_{K}f d\mu-t_n(f)\right|&\le \sum_{i\in I_n}\left|\int_{A^n_i}f-f(a^n_i)d\mu\right| \\
&\le |\mu|(K)\sup_{|x-y|\le \frac1n}|f(x)-f(y)|.
\end{aligned}
\ee
By the Heine-Cantor theorem, $f$ is uniformly continuous on $K$, hence the last term in \eqref{heinecantoreq} goes to zero as $n\to \infty$. Therefore, for every $f\in \mathcal{C}(K)$ we have that $t_n(f)\to\int_Kfd\mu$, equivalently: $t_n\to \mu$ in the weak-$*$ topology.
\end{proof}
We conclude this Appendix with an observation on the case $r=\infty$.
\begin{prop}
Let $T\in \mathcal{M}^\infty_{loc}$. Then the associated measure $\mu$ can be assumed to be of the form $\rho du$ for some $\rho\in L^{\infty}(\mathbb{D})$.
\end{prop}
\begin{proof}
Let $T\in\mathcal{M}^\infty_{loc}$ be associated with $Q$, $\a$, $\mu$. 
It is not restrictive to assume $M=\mathbb{D}^m$, $Q=id$ and $k=1$.
Let us extend $T$ to the space $\mathcal{C}^\infty_c(\R^m)$ by declaring
\be 
T(\f)=T(\f|_{\mathbb{D}^m})=\int_{\mathbb{D}^m}\de_\a\f d\mu.
\ee
Let $e\in\N^m$ be the multi-index $e=(1,\dots,1)$. Note that for all $\f\in\mathcal{C}^\infty_c(\R^m)$
\be 
\max_{u\in \R^m} \|\de_\a \f(u)\|\le \int_{\R^m}\left|\de_{\a+e}(u) \f\right|du.
\ee
Define $V\subset L^1(\R^m)$ as $ 
V=\{\de_{\a+e}\f\colon \f\in \mathcal{C}^\infty_c(\R^m)\}$, and let $\lambda\colon V\to\R$ be the liner function defined by $\lambda(\de_{\a+e}\f)=T(\f)$. Then $\lambda$ is a (well defined) linear and bounded functional on $(V,\|\cdot\|_{L^1})$, since 
\be 
\begin{aligned}
|\lambda(\de_{\a+e}\f)|&=|T(\f)|
=|\int_{\mathbb{D}^m} \de_\a \f d\mu|
\le |\mu|(\mathbb{D}^m)\max_{\R^m} \|\de_\a \f\|
\le |\mu|(\mathbb{D}^m)\|\de_{\a+e}\f\|_{L^1}.
\end{aligned}
\ee
The Hahn-Banach theorem, implies that $\lambda$ can be extended to a continuous linear functional $\Lambda$ on the whole space $L^1(\R^m)$ and hence it coincides, as a distribution, with a function $\rho\in L^\infty(\D^m)=L^1(\D^m)^*$.
In particular, for all $\f\in \mathcal{C}^\infty_c(\R^m)$, we have that 
\be 
\begin{aligned}
T(\f)&= \lambda(\de_{\a+e} \f)
&= \int_{\R^n}\de_{\a+e}\f(u)\rho(u) du.\\
\end{aligned}
\ee
\end{proof}
%
%
%
%
%
%
\section{The representation of Gaussian Random Fields}\label{app:rgrf}
The purpose of this section is to prove that every GRF $X\in \g rMk$ (where $r$ may be infinite) is of the form \eqref{eq:KarLo}. This fact is well-known in the general theory of Gaussian measures on Fréchet spaces (see \cite{bogachev}) as the Karhunen-Loève expansion. Here we present and give a proof of such result using our language, with the scope of making the exposition more complete and self-contained. Our presentation is analogous to that in the Appendix of \cite{NazarovSodin2} which, however, treats only the case of continuous Gaussian random functions, namely GRFs in $\g 0M{}$.

Given a Gaussian random field $X\in \mathcal{G}^r(M, \R^k)$, we define its Cameron-Martin Hilbert space $\mathcal{H}_X\subset E^r$ to be the image of the map $\rho_X$, as we did in Section \ref{sec:PT3}.
\be \rho_X\colon \Gamma_X=\overline{\textrm{span}\{X^j(p),\,p\in M\}}^{L^2(\Omega, \mathfrak{S}, \PP)} \to E^r
\ee
\be \rho_X(\gamma)=\E\left(X(\cdot)\gamma\right)=\left(\langle X^1(\cdot), \gamma\rangle_{L^2(\Omega, \mathfrak{S}, \PP)},\dots,\langle X^k(\cdot), \gamma\rangle_{L^2(\Omega, \mathfrak{S}, \PP)}\right).\ee
(This is consistent with the more abstract theory from \cite{bogachev} because of Proposition \ref{prop:CamMa}.)
Note that $\mathcal{H}_X$ is separable, since $M$ is, hence it has a countable Hilbert-orthonormal basis $\{h_n\}_{n\in \N}$, corresponding via $\rho_X$ to a Hilbert-orthonormal basis $\{\xi_n\}_{n\in \N}$ in $\Gamma_X$. This means that for any $p$ and $j$, one has $h_n^j(p)=\langle X^j(p),\xi_n\rangle$, namely that $h_n^j(p)$ is precisely the $n^{th}$ coordinate of $X^j(p)$ with respect to the basis $\{\xi_n\}_{n\in\N}$. In other words: 
\be \label{repreq}
X(p)=\lim_{n\to \infty} \sum_{m\leq n}\xi_m h_m(p),
\ee
where the limit is taken in $L^2\Prob^k$. 
In particular, since the $L^2$ convergence of random variables implies their convergence in probability, we have that
\be\label{eq:rep}
\lim_{n\to \infty}\P\left\{\left|\sum_{m>n}\xi_m h_m(p)\right|>\e\right\}=0.
\ee 
\begin{thm}[Representation theorem]\label{thm:rep} Let $X\in \mathcal{G}^r(M, \R^k)$ be a GRF, with $r\in\N\cup\{\infty\}$. For every Hilbert-orthonormal basis $\{h_n\}_{n\in \N}$ of $\mathcal{H}_X$, there exists a sequence $\{\xi_n\}_{n\in N}$ of independent, standard Gaussians such that the series $\sum_{n\in \N}\xi_nh_n$ converges\footnote{Given a sequence $\{x_n\}_{n\in \N}\subset E$, the sentence ``the series $\sum_{n\in \N}x_n$ converges in $E$ to $x$'' means that $s_N=\sum_{n\leq N}x_n$ converges in $E$ to $x$ as $N\to \infty$.} in $E^r$ to $X$ almost surely.
\end{thm}
To prove the above result, we will need a convergence criterion for a random series (see Theorem \ref{converepr}). It essentially follows from the Ito-Nisio theorem, which we recall for the reader's convenience.
\begin{thm}[Ito-Nisio]\label{itonisiothm}
Let $E$ be a separable real Banach space. Let $M\subset E^*$ be such that the family of sets of the form $\{f\in E \ | \ \langle p,f \rangle \in A\}$, with $A\in \mathcal{B}(\R)$ and $p\in M$, generates the Borel $\sigma$-algebra of $E$. Let $\{x_n\}_{n\in \N}$ be independent symmetric random elements of $E$, define
\be 
X_n=\sum_{m\leq n} x_m.
\ee
Then the following statements are equivalent:
\begin{enumerate}
\item $X_n$ converges almost surely;
\item $\{X_n\}_{n\in \N}$ is tight in $\mathscr{P}(E)$;
\item There is a random variable $X$ with values in $E$ such that $\langle p, X_n \rangle \to \langle p, X \rangle$ in probability for all $p\in M$.
\end{enumerate}
 
\end{thm}
\begin{remark}
In the original paper \cite{ItoNisio}, the theorem is stated with the hypothesis that $M=E^*$, but the same proof works in the slightly weaker assumptions of Theorem \ref{itonisiothm}.
\end{remark}

\begin{thm}\label{converepr}
Let $x_n\in \g rMk$, for all $n\in \N$. Assume that $x_n$ are independent and consider the sequence $X_n$ of GRFs defined as
\be X_n=\sum_{j\leq n} x_j.\ee 
The following conditions are equivalent.
\begin{enumerate}
\item $X_n$ converges in $\Cr rMk$ almost surely.
\item Denoting by $\mu_n$ the measure associated to $X_n$, we have that $\{\mu_{n}\}_{n\in \N}$ is relatively compact in $\G(E^r)$. 
\item There is a random field $X$ such that, for all $p\in M$, the sequence $\{X_n(p)\}_{n\in \N}$ converges in probability to $X(p)$.
\end{enumerate}
\end{thm} 

\begin{proof}We prove both that $(1)\iff (2)$ and $(3)\iff(1)$. We repeatedly use the fact that a.s. convergence implies convergence in probability, which in turn implies convergence in distribution (narrow convergence).

$(1)\Rightarrow (2)$ This descends directly from the fact that almost sure convergence implies narrow convergence.

$(1)\Rightarrow (3)$ This step is also clear, since the almost sure convergence of $X_n$ to some random field $X$ implies that for any $p\in M$ the sequence of random vectors $X_n(p)$ converges to $X(p)$ almost surely and hence also in probability.

$(3)\Rightarrow (1)$ and $(2)\Rightarrow (1)$
Let $Q_\ell\colon D\hookrightarrow M$, be a countable family of embeddings of the compact disk, as in \eqref{prodembballeq}. Note that if $\{X_n\}_n$ is tight in $\g rMk$ (i.e. $\mu_n$ is tight in $\G(E^r)$), then $\{X_n\circ Q_\ell\}_n$ is tight in $\g rDk$. Moreover,
if $X_n\circ Q_\ell \to X\circ Q_\ell$ almost surely in $\Cr rDk$, for every $\ell\in \N$, then $X_n\to X$ almost surely in $\Cr rMk$. Therefore it is sufficient to prove the theorem in the case $M=D$. For analogous reasons, we can assume that $r$ is finite.

The topological vector space $E=\Cr rDk$ has the topology of a separable real Banach space, with norm 
\be 
\|\cdot\|_E=\|\cdot \|_{\text{id}_D, r}.
\ee
Since the $\sigma$-algebra $\mathcal{B}(\mathcal C^{r}(M, \R^k))$ is generated by sets of the form $\{f: f(p)\in A\}$, where $p\in M$ and $A\subset \R^k$ is open and since Gaussian variables are symmetric, we can conclude applying the Ito-Nisio Theorem \ref{itonisiothm} to the sequence $X_n$ of random elements of $E^r$.
\end{proof}
\begin{proof}[Theorem \ref{thm:rep}] Let $\{h_{n}\}_{n\in \N}$ be a Hilbert orthonormal basis for $\mathcal{H}_X$ and set $\xi_n=\rho_X^{-1}(h_n)$ (it is a family of independent, real Gaussian variables). From equation \eqref{eq:rep} we get that for every $p\in M$ and $j=1,\dots,k$ we have convergence in probability for the series:
\be
\label{eq:ppp}X(p)=\lim_{n\to \infty} \sum_{m\leq n}h_n(p) \xi_n.
\ee
Then, the a.s. convergence of the above series in $\mathcal{C}^r(M, \R^k)$ follows from point (1) of Theorem \ref{converepr}.
\end{proof}
\subsection{The support of a Gaussian Random Field}
By definition (see equation \eqref{eq:defsupport}), the support of $X\in\g rMk$ has the property that if it intersects an open set $U\subset E^r$, then $\P\{X\in U\}>0$. The following proposition guarantees that the converse is also true, namely that if $\P\{X\in U\}>0$, then $U\cap \spt(X)\neq \emptyset$.

\begin{prop}
The support of $X\in \g rMk$ is the smallest closed set $C\subset E^r$ such that $\P\{X\in C\}=1$. 
\end{prop}

\begin{proof}
By definition we can write the complement of $\textrm{supp}(X)$ as
\be
\left(\text{supp}(X)\right)^c=\bigcup\{U\subset \mathcal{C}^{r}(M, \R^k) \text{ open such that }\P\{X\in U\}=0\}.
\ee

Consequently $\textrm{supp}(X)$ equals the intersection of all closed sets $C\subset E^r$ such that $\P\{X\in C\}=1$, hence it is closed. Since $E^r$ is second countable the union above and the resulting intersection of closed sets can be taken over a countable family, so that $\P\{X\in \text{supp}(X)\}=1$.
\end{proof}
\begin{remark}
Assume that $X_n\nrw X \in \g rMk$ and recall Portmanteau's theorem (see \cite[Theorem 3.1]{Billingsley}). Then, for any open set $U\subset \Cr rMk$ such that $U\cap \text{supp}(X)\neq \emptyset$, there is a constant $p_U=\frac12 \P\{X\in U\}>0$ such that for $n$ big enough, one has
\be 
\P\{X_n\in U\}\ge p_U.
\ee
In particular, this implies that
\be 
\text{supp}(X)\subset \bigcap_{n_0}\overline{\left(\bigcup_{n\ge n_0}\text{supp}(X_n)\right)}=\limsup_{n\to \infty}\text{supp}(X
_N).
\ee
\end{remark}

\begin{thm}[The support of a Gaussian random map]\label{thm:sputo}Let $X\in \mathcal{G}^{r}(M, \R^k)$. Let $\{f_n\}_{n\in \N}\subset E^r$ and consider a sequence $\{\xi_n\}_{n\in\N}$ of independent, standard Gaussians. Assume that the series $\sum_{n\in \N}\xi_n f_n$ converges in $E^r$ to $X$ almost surely. Then
\be \mathrm{supp}(X)=\overline{\mathrm{span}\{f_n\}_{n\in \N}}^{E^r}.\ee
\end{thm}

\begin{proof}
We start by observing that  $X\in \overline{\text{span}\{f_n\}_n}$ with $\P=1$, thus the first inclusion ``$\subset$'' is proved.
Let now $c=\sum_{n=0}^{N_0}a_nf_n$ and let $U_c\subset \Cr rMk$ be an open neighborhood of $c$ of the form
\be 
U_c=\left\{ f\in \Cr rMk \colon \|f-c\|_{Q, r}<\e\right\}.
\ee
for some embedding $Q$.
Denote by $S_N=\sum_{n\leq N}\xi_n f_n$. Observe that if $N\ge N_0$, then $S_{N}-c\in \text{span}\{f_1\dots f_{N}\}$, which is a finite dimensional vector space, hence there is a constant $A_{N}>0$ such that $\|\sum_{n=0}^N a_nf_n\|_{Q,r}\le A_N\max\{|a_0|\dots |a_N|\}$. 
By the convergence in probability of $S_N$ to $X$, there exists $N>N_0$ so big that $\P\left\{\|X-S_{N}\|_{Q, r}\geq \frac{\epsilon}2\right\}<\frac12$.
Thus, setting $a_n=0$ for all $n>N_0$, we have:
\be
\begin{aligned}
\P\{X\in U_c\}&\ge \P\left\{\|X-S_{N}\|_{Q, r}<\frac{\e}2, \|S_{N}-c\|_{Q,r}<\frac\e 2\right\}\\
&\ge \P\left\{\|S_{N}-c\|_{Q,r}<\frac\e 2\right\}\frac12  \\
&\ge \left(\prod_{n=0}^{N}\P\left\{|\xi_n-a_n|<\frac{\e}{2A_{N}}\right\} \right)\frac12\\
&>0.
\end{aligned}
\ee

Every open neighborhood of $c$ in $\Cr rMk$ contains a subset of the form of $U_c$, therefore $c\in \text{supp}(X)$. Since $\text{supp}(X)$ is closed in $E^r$, we conclude.
\end{proof}
\begin{cor}\label{cor:sputoCMapp} 
Let $X\in\g rMk$ and let $\mathcal{H}_X\subset E^r$ be its Cameron-Martin space.
\be\label{eq:sputoCMapp}
\spt(X)=\overline{\mathcal{H}_X}^{E^r}.
\ee
\end{cor}
\begin{proof}
It is a consequence of Theorems \ref{thm:rep} and \ref{thm:sputo}.
\end{proof}
%
%
%
%
%
\end{subappendices}

%
%


\chapter{KAC RICE FORMULA FOR TRANSVERSAL INTERSECTIONS}\label{chap:kr}

We prove a generalized Kac-Rice formula that, in a well defined regular setting, computes the expected cardinality of the preimage of a submanifold via a random map, by expressing it as the integral of a density. Our proof starts from scratch and although it follows the guidelines of the standard proofs of Kac-Rice formula, it contains some new ideas coming from the point of view of measure theory. Generalizing further, we extend this formula to any other type of counting measure, such as the intersection degree. 

We discuss in depth the specialization to smooth Gaussian random sections of a vector bundle. Here, the formula computes the expected number of points where the section meets a given submanifold of the total space, it holds under natural non-degeneracy conditions and can be simplified by using appropriate connections. Moreover, we point out a class of submanifolds, that we call sub-Gaussian, for which the formula is locally finite and depends continuously  with respect to the covariance of the first jet. In particular, this applies to any notion of singularity of sections that can be defined as the set of points where the jet prolongation meets a given semialgebraic submanifold of the jet space.

Various examples of applications and special cases are discussed. In particular, we report a new proof of the Poincaré kinematic formula for homogeneous spaces and we observe how the formula simplifies for isotropic Gaussian fields on the sphere.

%
%
%


\subsection{Overview}

What motivates this work is the interest in studying the expected number of realizations of a geometric condition. This topic is at the heart of stochastic geometry and geometric probability and, in recent years, it has gained a role also in subjects of more deterministic nature like enumerative geometry (see \cite{Kostlan:93}, \cite{shsm}, \cite{GaWe1,GaWe2,GaWe3}, \cite{SarnakWigman}, \cite{NazarovSodin1,NazarovSodin2}, \cite{FyLeLu,Letwo,Lerarioshsp,Lerariolemniscate}) and physics (see for instance \cite{marinucci_peccati_2011}, \cite{Wig,wigEBNRF},\cite{park2013betti}).

As a first example, let us consider a random $\mC^1$ function $X\colon M\to \R^n$, where $M\subset \R^m$ is an open subset and let $t\in \R^n$. Then, under appropriate assumptions on $X$, we have the so called \emph{Kac-Rice formula} for the expected cardinality of the set of solutions of the equation $X=t$.
\be\label{eq:KR}
\E\#X^{-1}(t)=\int_M\E\left\{\left|\det(d_uX)\right|\Big|X(u)=t\right\}\rho_{X(u)}(t)du,
\ee
where $\rho_{X(u)}$ is the density of the random variable $X(u)$, meaning that for every event $A\subset \R^n$ we have $\P\{X(u)\in A\}=\int_{A}\rho_{X(u)}(t)dt$. 
First appeared in the independent works by M. Kac \cite[1943]{kac43} and S. O. Rice \cite[1944]{rice44}, this formula is today one of the most important tool in the application of smooth stochastic processes (here called ``random maps'' or ``random fields'') both in pure and applied maths. In fact, the ubiquity of this formula is suggested by its name and birth in that Kac's paper is about random algebraic geometry, while Rice's one deals with the analysis of random noise in electronic engineering. We refer to the book \cite{AdlerTaylor} for a detailed treatment of Kac-Rice formula. We were also inspired by the book \cite{Wschebor} containing a simpler proof in the Gaussian case.

The main result of this paper is the generalization of Kac-Rice formula to one that computes the expected number of points in the preimage of a submanifold $W$, namely the number of solutions of $X\in W$,  rather than just $X=t$. 
\begin{multline}\label{eq:mainformula}
\E\#(X^{-1}(W))= \\=\int_{M}\left(\int_{ S_p\cap W}\E\left\{J_pX\frac{\sigma_q(X,W)}{\sigma_q(S_p,W)}\bigg|X(p)=q\right\}\rho_{X(p)}(q)d(S_p\cap W)(q)\right)dM(p),
\end{multline}
This is the content of Theorems \ref{thm:main}, \ref{thm:maindens} and  \ref{thm:maingraph} reported in Section \ref{sec:main}, after a brief introduction to the problem in Section \ref{sec:intro}. Such theorems are essentially equivalent alternative formulations of the same result. In presenting them, we  pay a special attention to their hypotheses, aiming to propose a setting, that we call KROK\footnote{Stands for ``Kac Rice OK''.} hypotheses (Definition \ref{def:krok}), that appears frequently in random geometry and that is easy to recognize, especially in the Gaussian case.
\begin{remark}[A comment on the proof]
The first idea that comes to mind is to write, locally, the submanifold $W$ as the preimage $\f^{-1}(0)$ of a smooth function $\f\colon N\to \R^m$ and then apply the standard Kac-Rice formula to the random map $\f\circ X$. After that, however, one wants to get rid of $\f$ since it is desirable to have an intrinsic statement, independent from the arbitrary choice of this auxiliary function. In fact, this is the key issue, but it ends up being ugly. So, we chose to reprove everything from the beginning, instead. In doing this, we aim also at proposing an alternative reference for the proof of the standard case. 

Specializing the proof of Theorem \ref{thm:main} to the case in which $W$ is a point, one obtains a proof for the standard Kac-Rice formula, in the KROK setting. Although this setting is very general, the complexity of the proof is comparable to that of Azais and Wschebor \cite{Wschebor} for the Gaussian case and quite simple if compared to the one reported in Adler and Taylor's book \cite{AdlerTaylor}. Moreover, we use an argument that is new in this context: instead of dividing the domain in many little pieces, we interpret the expectation as a measure and use Lebesgue Differentiation Theorem. This makes the hard step of the proof (the ``$\ge$'' part) a little more elegant.
\end{remark}
In Section \ref{sec:gausscase} we focus on the case of Gaussian random sections of a vector bundle. Here, the formula specializes to Theorem \ref{thm:maingau}, where the hypotheses reduce to simple non-degeneracy conditions, thanks also to the Probabilistic Transversality theorem from Chapter \ref{chap:dtgrf}. We also provide alternative ways of writing the formula \eqref{eq:mainformula}: as a measure on the submanifold $W$ (Corollary \ref{cor:maingauW}), or using the canonical connection defined by the Gaussian field (Corollary  \ref{cor:connectedform}), see \cite{Nicolaescu2016, AdlerTaylor}. Moreover, in this case we establish a continuity property of the expected number of singular points of a Gaussian random section, with respect to the corresponding covariance tensor (Theorem \ref{thm:mainEgau}). 
This last result has a nice application in the study of semialgebraic singularities of Gaussian random fields (Corollary \ref{cor:mainjgau}).

We also discuss, in Section \ref{sec:weights}, the problem of counting solutions with ``weights'', for instance the intersection degree of $X$ and $W$. Here we show (Theorem \ref{thm:mainalph}) that, under KROK hypotheses, the formula can be directly generalized to hold for any counting measure with measurable weights.

Finally, in section \ref{sec:ex}, we test our formula in two prominent instances of random geometry. First, we show that it can be used to obtain a new quite elementary proof of Poincaré kinematic formula for homogeneous spaces (Theorem \ref{thm:howard}), in the case of zero dimensional intersection; then, we deduce a simple but general formula for isotropic Gaussian random fields on the sphere (Theorem \ref{thm:mainIsotrop}).

\subsection{Structure of the Chapter}
Sections \ref{sec:intro}-\ref{sec:weights}  contain the presentation of the results of the paper, without proofs. All of their proofs are contained in the Sections \ref{sec:proflemmameas}-\ref{sec:proofgauss}. In particular, Section \ref{sec:gen} is devoted to deduce from the coarea formula that the identity \eqref{eq:mainformula} holds for almost every $W$, under very general assumptions. This essentially allows to prove the ``$\le$'' part of \eqref{eq:mainformula} in the KROK setting, while the opposite inequality is proved in Section \ref{sec:mainproof}. 
Section \ref{sec:ex} contains minor results (and their proof) obtained from applications of the main formula. In the appendix we report some details regarding a few notions of which we make extensive use throughout this paper.


\begin{remark}
 The reader who wants to grasp the meaning of the generalized formula, without going into its more abstract aspects, may just skip Section \ref{sec:main} and go directly to the Gaussian case, discussed in Section \ref{sec:gausscase}. Enough references are provided so that this is a safe practice. However, we recommend that you take a look at Section \ref{sec:intro} first.
\end{remark}
\section{Introduction}\label{sec:intro}
\subsection{Notations}

\begin{enumerate}[$\bullet$]
 \item We write $\#(S)$ for the cardinality of the set $S$.
\item We use the symbol $A\transv B$ to say that objects $A$ and $B$ are in transverse position, in the usual sense of differential topology (as in \cite{Hirsch}).
    \item The space of $\mC^r$ functions between two manifolds $M$ and $N$ is denoted by $\coo rMN$.
    If $E\to M$ is a vector bundle, we denote the space of its $\mC^r$ sections by $\mathcal{C}^r(M|E)$. In both cases, we consider it to be a topological space endowed with the weak Whitney's topology (see \cite{Hirsch}).
    \item We call $X$ a \emph{random element} (see \cite{Billingsley}) of the topological space $T$ if $X$ is a measurable map $X\colon \Omega\to T$, defined on some probability space $\Prob$ and we denote by $[X]=\P X^{-1}$ the Borel probability measure on $T$ induced by pushforward. We will alternatively use the following equivalent notations:
\be 
\P\{X\in U\}:=[X](U)=\P\left(X^{-1}(U)\right)=\int_{U}d[X],
\ee 
to denote the probability that $X\in U$, for some measurable subset $U\subset T$, and
\be 
\E\{f(X)\}:=\int_{T}f(t)d[X](t),
\ee
to denote the expectation of a measurable vector-valued function $f\colon T\to \R^k$.
Two random elements $X_1,X_2$ are said to be \emph{equivalent} and treated as if they were equal if $[X_1]=[X_2]$. We might call $X$ a \emph{random variable}, \emph{random vector} or \emph{random map} if $T$ is the real line, a vector space or a space of functions $\coo rMN$, respectively.
    \item The sentence: ``$X$ has the property $\mathcal{P}$ almost surely'' (abbreviated ``a.s.'') means that the set  $S=\{f\in \coo1MN| f \text{ has the property }\mathcal{P}\}$ contains a Borel set of $[X]$-measure $1$. It follows, in particular, that the set $S$ is $[X]$-measurable, i.e. it belongs to the $\sigma$-algebra obtained from the completion of the measure space $(\coo1MN,\mathcal{B},[X])$.
    
    \item The bundle of densities of a manifold $M$ is denoted by $\Delta M$, see Appendix \ref{app:densities} for details. If $M$ is a Riemannian manifold, we denote its volume density by $dM$. The subset of positive density elements is denoted by $\Delta^+M$. We denote by $B^+(M)$ the set of positive Borel functions $M\to [0,+\infty]$ and by $L^+(M)$ the set of positive densities, i.e. densities of the form $\rho dM$, where $\rho\in B^+(M)$ and $dM$ is the volume density of some Riemannian metric on $M$. In other words 
\be 
L^+(M)=\left\{ \text{measurable functions }M\ni p\mapsto \delta(p)\in\Delta^+_pM \cup \{+\infty\}\right\}
\ee 
is the set of all non negative, non necessarily finite Borel measurable densities.
The integral of a density $\delta\in L^+(M)$ is written as $\int_M \delta=\int_M \delta(p) dp$.
\item The Jacobian of a $\mC^1$ map $f\colon M\to N$ between Riemannian manifolds (see Definition \ref{def:jacob}), evaluated at a point $p\in M$ is denoted by $J_pf\in [0,+\infty)$. The Jacobian density is then $\delta_pf=J_pf dM(p)$ (see Appendix \ref{app:densities}). If moreover $M$ and $N$ have the same dimension then we may stress this fact by writing $J_pf=|\det (d_pf)|$. In case $f$ is a linear map between Euclidean spaces, then we will just write $Jf:=J_0f=J_pf$. 
    \item
Given a finite dimensional Hilbert space $E$, the expression $\sigma_E(V,W)$ denotes the ``angle''\footnote{Actually a better analogy is with the sine of the angle.} between two vector subspaces $V,W\subset E$, see Appendix \ref{app:angle}. If $f\colon M\to N$ is a $\ci$ map between Riemannian manifolds and $W\subset N$ is a submanifold, we will write shortly
\be 
\sigma_x(W,f):=\sigma_x(W,d_pf):=\sigma_{T_xN}(T_xW, d_pf(T_pM)),
\ee
whenever $f(p)=x\in W$. If $S\subset N$ is another submanifold and $x\in S\cap W$, then
\be 
\sigma_x(W,S)=\sigma_{T_xN}(T_xW, T_xS)).
\ee
\end{enumerate}
\subsection{The expected counting measure}
Let us start by considering the following setting. 
\begin{enumerate}[i.]\label{setting}
    \item $M, N$ smooth manifolds ($\mC^\infty$ and without boundary) of dimension $m,n$.
    \item $W\subset N$ smooth submanifold (image of a smooth embedding) of codimension $m$.
    \item $X\colon M\to N$ random $\mC^1$ map, i.e. it represents a Borel probability measure $[X]$ on the topological space $\coo1MN$ endowed with the (weak) Whitney $\mC^1$ topology.
    \item $X\transv W$ almost surely.
\end{enumerate}
If moreover $W$ is closed (this assumption can and will be removed with Lemma \ref{lemma:meas}), then the random set $X^{-1}(W)$ is almost surely a discrete  subset of $M$, so that for every $U\subset M$ relatively compact open set, the number 
\be 
\#_{X\in W}(U):=\#(X^{-1}(W)\cap U)
\ee
is almost surely finite (if $W$ is not closed, this number can be $+\infty$) and it is a continuous function with respect to $X\in \{f\in \coo1MN \colon f\transv W\}$, thus it defines an integer valued random variable. Now, Lemma \ref{lemma:meas} below guarantees that its mean value
\be 
\E\#_{X\in W}(U)=\E\{\#(X^{-1}(W)\cap U)\}
\ee
can be extended to a Borel measure on $M$.

\begin{lemma}\label{lemma:meas}
Let $X\colon M\to N\supset W$ satisfy \emph{i-iv}. For any $A\in \mathcal{B}(M)$, the number $\#_{X\in W}(A)=\#(X^{-1}(W)\cap A)$ is a measurable random variable 
and the set function
\be 
\E\#_{X\in W}\colon\mathcal{B}(M)\to [0,+\infty],\qquad A\mapsto \E\#_{X\in W}(A)
\ee
is a Borel (not necessarily finite) measure on $M$.
\end{lemma}

The proof of this Lemma is postponed to Section \ref{sec:proflemmameas}. At this point, a couple of curiosities about this measure naturally arise: 
when is it a Radon measure? When is it absolutely continuous (in the sense of Definition \ref{def:abscont})? 
In this paper we are going to address these questions giving sufficient conditions for $\E\#_{X\in W}$ to be an absolutely continuous Radon measure and a formula to compute it in this case.

\section{KROK hypotheses and the main result.}\label{sec:main}

By considering the following particularly simple examples, that we should always bear in mind, we can observe that the setting i-iv described above is far too general to allow to give a yes/no answer to the previous questions.
\begin{enumerate}[$\bullet$]
\item Let $X$ be \emph{deterministic}, in the sense that it is constantly equal to a function $f\transv W$. Then $\E\#_{X\in W}$ is the counting measure of the set $f^{-1}(W)$. This measure is Radon if and only if the set has no accumulation point, which is a consequence of transversality when $W$ is closed. In this situation the only case where $\E\#_{X\in W}$ is absolutely continuous is if $f^{-1}(W)=\emptyset$.
\item Let $\dim M=0$, i.e. $M=\{p_i\}_{i\in\N}$, and let $W\subset N$ be an open subset. In this case $X$ is a random element in the product space $X\in N^M$ and it can be easily checked that
\be\label{eq:stupid} 
\E\#\{p\in M\colon X_p\in W\}=\sum_{p\in M}\P\{X_p\in W\}.
\ee
This is a rather stupid case, however, the above formula \eqref{eq:stupid} is close in spirit to the one we are going to prove (in fact it is a special case of Theorem \ref{thm:main}), in that the right hand side depends only on the marginal probabilities of the random variables $\{X_p\}_{p\in N}$.
\item Let $X\colon M\to M\times M$ be the map $p\mapsto (p,\xi)$, for some random element $\xi\in M$ and let $W=\Delta$ be the diagonal. Then the measure $\E\#_{X\in W}$ is the law of $\xi$. Since the hypotheses i-iv are satisfied for every random variable $\xi$, this example shows that certainly every Borel probability measure on $M$ can be realized in this way (it is more difficult to realize an arbitrary measure with total mass greater than $1$).
\item Let $G$ be a group. Let $X\colon M\to E$ be a random section of a $G$-equivariant\footnote{Meaning that $G$ acts on the left on both $E$ and $M$ and the action commutes with the projection: $\pi(g\cdot x)=g\cdot \pi(x)$, for any $x\in E$ and $g\in G$. Thus, the function $gXg^{-1}$ such that $p\mapsto g\cdot X(g^{-1}\cdot p)$ is a section.} bundle $E\to M$, such that $gXg^{-1}$ has the same law (on the space of $\mC^1$ sections of $E$) as $X$, then the measure $\E\#_{X\in W}$ is $G$-invariant. This condition, in many situations, implies that $\E\#_{X\in W}$ is a constant multiple of the volume measure of some Riemannian metric on $M$ and therefore it is absolutely continuous. 
\item Let $M\subset \R^m$ be an open subset, $W=\{t\}$ a point of $N=\R^n$ and $X\colon M\to \R^n$  a $\mC^1$ Gaussian random field such that $X(p)$ is non degenerate for every $p\in M$. Then Kac-Rice formula \eqref{eq:KR} holds (see \cite{AdlerTaylor}) 
\end{enumerate}
To be able to say something meaningful we need to restrict our field of investigation. We will now make a series of assumption on the random map $X$ and on the submanifold $W$ under which the measure $ \E\#_{X\in W}\colon \mathcal{B}(M)\to [0,+\infty]$ is absolutely continuous and we can write a formula for its density. 
In doing so, one of our aim is to propose a setting that is easy to recognize in contexts involving differential topology and smooth random maps.
Although such hypotheses do not reach the highest level of generality in which Kac-Rice formula holds (see \cite{AdlerTaylor}), they describe a much more general setting than that in which the random map is assumed to be Gaussian and at the same time allow to give a proof whose simplicity is comparable to those for the Gaussian case.

\begin{defi}\label{def:krok}
 Let $M$ and $N$ be two smooth manifolds ($\mC^\infty$ and without boundary) of dimension $m$ and $n$. Let $W\subset N$ be a smooth submanifold (without boundary) and $X\colon M\to N$ a random map. We will say that $(X,W)$ is a \emph{KROK}
  couple if the following hypotheses are satisfied.
\begin{enumerate}[ , wide, left=0pt]
    \item Properties of $X$:
\begin{enumerate}[(i)]
\item\label{itm:krok:2} \boldmath{$X\in \mC^1$: }\unboldmath

$X\colon M\to N$ is a random $\mC^1$ map, i.e. it represents a Borel probability measure $[X]$ on the topological space $\coo1MN$ endowed with the weak Whitney $\mC^1$ topology (see \cite{Hirsch}).
\item\label{itm:krok:3} \boldmath{$d[X(p)]=\rho_{X(p)}dS_p$: }\unboldmath 

Let $N$ be endowed with a Riemannian metric. Assume that for each $p\in M$, the probability measure $[X(p)]$ is absolutely continuous with respect to the Riemannian volume density of a certain smooth submanifold $S_p$. 
In other words, there exists a measurable function $\rho_{X(p)}\colon S_p\to [0,+\infty]$ such that
    \be 
    \E\{F(X(p))\}=\int_{S_p}F(q)\rho_{X(p)}(q)dS_p(q).
    \ee
    for every Borel function $F\colon N\to\R$. ($\rho_{X(p)}$ is allowed to vanish on $S_p$.)
\end{enumerate}
\item Properties of $W$ (transversality):
\begin{enumerate}[(i)]
\setcounter{enumii}{2}
\item\label{itm:kroktransvas} \boldmath{$X\transv W$ }\unboldmath

 almost surely.
\item\label{itm:krok:6} \boldmath{$S_p\transv W$ }\unboldmath  

for every $p\in M$. 
\item\label{itm:krok:4} \boldmath{$\dim X^{-1}(W)=0$: }\unboldmath 

The codimension of $W$ is $m=\dim M$.
\end{enumerate}
\item Continuity properties:
\begin{enumerate}[(i)]
\setcounter{enumii}{5}
\item\label{itm:krok:7} \boldmath{$S_{(\cdot)}\in\mC^\infty$: }\unboldmath

The set $\mathcal{S}=\{(p,q)\colon q\in S_p\}\subset M\times N$ is a 
closed smooth submanifold. 
This, together with \ref{itm:krok:6}, implies that the set defined as $\MW  =\{(p,q)\in M\times N\colon q\in S_p\cap W\}$ is a smooth manifold. 
\item \boldmath{$\rho_{X(\cdot)}(\cdot)\in\mC^0$: }\unboldmath

The function $M\times N\ni (p,q)\mapsto \rho_{X(p)}(q)\in \R_+$ is continuous at all points of $\MW  $.
\item{\label{itm:krokEcont}}\boldmath{$[X|X(\cdot)=\cdot]\in\mC^0$: }\unboldmath 

There exists a regular conditional probability\footnote{See \cite{dudley}, or wait for the next subsection \ref{intro:krok9}.} ${[X|X(p)=\cdot]}$ such that the expectation
$M\times N\ni (p,q)\mapsto \E\big\{\a(X,p)J_pX\big|X(p)=q\big\}$ is continuous \emph{at}\footnote{i.e. continuous at every point of $\MW  $.
} $\MW  $, for any bounded and continuous function $\a\colon \coo1MN\times M\to \R$. 
\end{enumerate}
\end{enumerate}
In the following, we will refer to this hypotheses as KROK.\ref{itm:krok:2}, KROK.\ref{itm:krok:3}, etc\dots
\end{defi}

\begin{thm}[Generalized Kac-Rice formula]\label{thm:main}
Let $X\colon M \to N$ be a random $\mC^1$ map between two Riemannian manifolds and let $W\subset N$ such that $(X,W)$ is a KROK couple. Then for every Borel subset $A\subset M$ we have
\begin{multline}\label{eq:formrho}
\E\#_{X\in W}(A)=\int_A \rho_{X\in W}(p)dM(p)= \\=\int_{A}\left(\int_{ S_p\cap W}\E\left\{J_pX\frac{\sigma_q(X,W)}{\sigma_q(S_p,W)}\bigg|X(p)=q\right\}\rho_{X(p)}(q)d(S_p\cap W)(q)\right)dM(p),
\end{multline}
where $d(S_p\cap W)$ and $dM$ denote the volume densities of the corresponding Riemannian manifolds and $J_pX$ is the jacobian of $X$ (see Definition \ref{def:jacob}); besides, $\sigma_q(X,W)$ and $\sigma_x(S_p,W)$ denote the ``angles'' (in the sense of Definition \ref{defi:angle}) made by $T_qW$ with, respectively, $d_pX(T_pM)$ and $T_qS_p$.
\end{thm}
\begin{remark}[Special cases]
 The standard Kac-Rice formula corresponds to the situation when $S_p=N$ and $W=\{q\}$. Here, the term $\frac{\sigma_q(X,W)}{\sigma_q(S_p,W)}$ disappears, since both angles are equal to $1$. 

 When $T_qS_p=T_qN$, then $\sigma(S_p,W)=1$.

 When $S_p\cap W=\{q\}$, there is no integration $\int_{S_p\cap W}$.

 When $W$ is an open subset, then $\sigma_q({X,W})=1$ and $\sigma_q(S_p,W)=1$, but the dimension hypothesis KROK.\ref{itm:krok:4} falls, unless $m=0$. In such case, the above formula reduces to equation \eqref{eq:stupid}.

 If $S_p=\{f(p)\}$ it means that $X=f$ is deterministic. Unless $m=0$, the couple $(f,W)$ is KROK only if  $f^{-1}(W)=\emptyset$ because of KROK.\ref{itm:krok:3}. Indeed, as we previously observed in the first of the examples above, in the deterministic case the measure $\E\#_{f\in W}=\#_{f\in W}$ is not absolutely continuous, for obvious reasons, unless it is zero.
This is one of the reason why we can't change KROK.\ref{itm:krok:6} into ``$S_p\transv W$ for a.e. $p$''.
\end{remark}
In formula \eqref{eq:formrho} above, $\rho_{X\in W}\in B^+(M)$ is a not necessarily finite Borel measurable function $\rho_{X\in W}\colon M\to [0,+\infty]$. It is precisely the Radon-Nykodim derivative of $\E\#_{X\in W}$ with respect to the Riemannian volume measure of $M$.

We can write the above formula in another equivalent way, using the jacobian density
\be\label{eq:intro:indensity}
\delta_pX=J_pXdM(p)\in \Delta_pM
\ee
defined in \eqref{eq:indensity}, which is a more natural object in that it doesn't depend on the Riemannian  structure of $M$. 

By using the notion of \emph{density} we can write a more intrinsic formula, without involving a Riemannian metric on $M$. 
A density is a section of the vector bundle $\Delta M=\wedge^m(T^*M)\otimes L$ obtained by twisting the bundle of top degree forms with the orientation bundle $L$ (see \cite[Section 7]{botttu}). 
  The peculiarity of densities is that they can be integrated over $M$ in a canonical way. In particular the volume density of a Riemannian manifold $M$ is a density in all respects and we denote it by $p\mapsto dM(p)\in\Delta_pM$. We collected some details and notations regarding densities in Appendix \ref{app:densities}.
Although the function $\rho_{X\in W}$ appearing in \eqref{eq:formrho} depends on the Riemannian structure of $M$, the expression $\rho_{X\in W}(p)dM(p)=:\delta_{X\in W}(p)$ defines an $L^+(M)$ density that is independent from the metric. This is clarified in the subsection  \ref{sec:indepmetr}, but it is actually a consequence of Theorem \ref{thm:main}, since the left hand side of \eqref{eq:formrho} depends merely on the ``set theoretic nature'' of the objects in play.

\begin{cor}[Main Theorem/Definition]\label{thm:maindens}
Let $X\colon M \to N$ be a random $\mC^1$ map between two Riemannian manifolds and let $W\subset N$ such that $(X,W)$ is a KROK couple. Then the measure $\en$ is absolutely continuous on $M$ with density $\delta_{X\in W}\in L^+(M)$ defined as follows.
\begin{multline}\label{eq:formdel}
\delta_{X\in W}(p):=\left(\int_{ S_p\cap W}\E\left\{\delta_pX\frac{\sigma_q(X,W)}{\sigma_q(S_p,W)}\bigg|X(p)=q\right\}\rho_{X(p)}(q)d(S_p\cap W)(q)\right)\in\Delta_pM,
\end{multline}
where $d(S_p\cap W)$ denotes the volume densities of the corresponding Riemannian manifold and $\delta_pX$ is the Jacobian density of $X$ defined as in \eqref{eq:intro:indensity}; besides, $\sigma_q(X,W)$ and $\sigma_q(S_p,W)$ denote the ``angles'' (in the sense of Definition \ref{defi:angle}) made by $T_qW$ with, respectively, $d_pX(T_pM)$ and $T_qS_p$. Therefore, for every Borel subset $A\subset M$,
\be 
\en(A)=\int_A \delta_{X\in W}.
\ee
\end{cor}
\begin{remark}\label{rem:altforms}
Other alternative forms of the above formula can be obtained from the identities: 
\be \delta_p X\sigma(X,W)=J(\Pi_{TW^\perp}\circ dX)dM(p)=|X^*\nu|,
\ee where $\nu=\nu^1\wedge \dots \wedge \nu^r$ for some orthonormal basis on $TW^\perp$. The first identity follows from Proposition \ref{prop:chi}, while in the second we are representing the density element as the modulus of a differential form, via the function $|\cdot |\colon \wedge^mT^*_pM\to \Delta_pM$, defined in Appendix \ref{app:densities}.
\end{remark}
\begin{remark}
Theorem \ref{thm:main} does not guarantee that $\E\#_{X\in W}$ is a Radon measure. This condition corresponds to the local integrability of the density: $\delta_{X\in W}\in L^1_{loc}(M)$, while in general $\rho_{X\in W}\in B^+(M)$ and $\delta_{X\in W}\in L^+(M)$. This issue, as it will be clear from Theorem \ref{thm:maingraph} below, comes from the integration over $W\cap S_p$, which can be a non-compact submanifold.
\end{remark}

The strength of this formula is that, exactly as in the standard Kac-Rice case (when $W$ is a point), although the left hand side depends, a priori, on the whole probability $[X]$ on $\coo 1MN)$, the right hand side depends only on the pointwise law of the first jet $j^1_pX=(X(p),d_pX)$. This is a significant simplification in that the former is the joint probability of all the random variables $\{X(p)\}_{p\in M}$, while the latter is the collection of the marginal probabilities $\{[j^1_pX]\}_{p\in M}$, which is a simpler data. 

\subsection{Explanation of condition KROK.\ref{itm:krokEcont}}\label{intro:krok9}
Given a random element $X\in \mC^1(M,N)$ as in KROK.\ref{itm:krok:3} and a point $p \in M$,
a \emph{regular conditional probability}\footnote{
See \cite{dudley} or \cite{Erhan}. In the latter the same object is called a \emph{regular version of the conditional probability}.
} for $X$ given $X(p)$ is a function 
\be
\P\{X\in \cdot\ |X(p)=\cdot\}\colon \mathcal{B}(\mC^1(M,N))\times N\to [0,1],
\ee
that satisfies the following two properties.
\begin{enumerate}[(a)]
\item For every $B\in\mathcal{B}(\mC^1(M,N))$, the function $\P(X\in B|X(p)=\cdot)\colon N\to [0,1]$ is Borel and for every $V\in \mathcal{B}(N)$, we have
\be 
\P\{X\in B ; X(p)\in V\}=\int_V\P(X\in B|X(p)=q\}d[X(p)](q).
\ee
\item For all $q\in N$, $\P\{X\in \cdot\ |X(p)=q\}$ is a probability measure on $\mC^1(M,N)$.
\end{enumerate} 
The notation that we use is what we believe to be the most intuitive one and consistent with the other used in this paper. Given $p\in M$ and $q\in N$, we write ${[X|X(p)=q]}:={\P\{X\in \cdot\ |X(p)=q\}}$ for the probability measure and ${\E\{\a(X)|X(p)=q\}}$ for the expectation/integral of a function $\a\colon \mC^1(M,N)\to \R$ with respect to the probability measure ${[X|X(p)=q]}$.

The fact that the space $\mC^1(M,N)$ is Polish ensures that a regular conditional probability measures $[X|X(p)=\cdot]$ exists (see \cite[Theorem 10.2.2]{dudley}) and it is unique up to $[X(p)]$-a.e. equivalence on $N$. However, strictly speaking, it is not a well defined function, although the notation can mislead to think that. 

In our case such ambiguity may be traumatic, since we are interested in evaluating $\E\{\dots|X(p)=q\}$ for $q\in W$ which, under KROK.\ref{itm:krok:3} and KROK.\ref{itm:krok:6}, is negligible for the measure $[X(p)]$, i.e. $\P\{X(p)\in W\}=0$. Therefore it is essential to choose a regular conditional probability that has some continuity property at $W$, otherwise formula \eqref{eq:formrho} doesn't make sense, as well as all of its siblings. This is the motivation for the hypothesys KROK.\ref{itm:krokEcont}.
\begin{remark}
In the common statements of Kac-Rice formula ($W$ is a point), one finds the analogous hypothesis that there exists a density for the measure $[J_pX|X(p)=q]$ that is continuous at $q\in W$ (see for instance \cite{AdlerTaylor}). This is different than KROK.\ref{itm:krokEcont}, in that we don't need to assume that $[J_pX|X(p)=q]$ has a density.
\end{remark}

To have a complete perspective, let us rewrite the hypothesis KROK.\ref{itm:krokEcont} in a more suggestive way.
Let $\mathcal{F}=\mC^1(M,N)$. 
Consider the space $\mathscr{P}(\mathcal{F})$ of all Borel probability measures on $\mathcal{F}$, endowed with the narrow topology (also called weak topology: see Chapter \ref{chap:dtgrf}), namely the one induced by the inclusion $\mathscr{P}(\mathcal{F})\subset \mC_b(\mathcal{F})^*$. A sequence of measures $[X_n]$ converges in this topology: $X_n\nrw X$, if and only if $\E\{\a(X_n)\}\to\E\{\a(X)\}$ for every $\a\in \mC_b(\mathcal{F})$, see \cite{Parth,Billingsley}.

Let $[X|X(p)=\cdot]$ be a regular conditional probability. Consider, for each $p\in M$ and $q\in N$, the probability $\mu(p,q)\in\mathscr{P}(\mathcal{F})$, given by
\be 
\mu(p,q)(B)=\E\{1_B(X)J_pX|X(p)=q\}=\int_B (J_pf) d[X|X(p)=q](f).
\ee
In other words, $\mu(p,q)=J_p\cdot [X|X(p)=q]$ is the multiplication of the measure ${[X|X(p)=q]}$ by the positive function $J_p\colon \mathcal{F}\to \R$, such that $J_p(f)=J_pf$.
This defines a function $\mu\colon M\times N\to \mathscr{P}(\mathcal{F})$.
\begin{prop}\label{prop:CazzoDuro}
KROK.\ref{itm:krokEcont} holds if and only if $\mu$ is continuous at $\MW  $.
\end{prop}
\begin{proof}
If the function $\a$ in KROK.\ref{itm:krokEcont} was not allowed to depend on $p\in M$, this fact would be obvious from the definition of the topology on $\mathscr{P}(\mathcal{F})$.  This, in particular, implies the \emph{only if} part of the statement.

Let us show the converse.
Fix $\a\in \mC_b(\mathcal{F}\times M)$ and let $(p_n,q_n)\in M\times N$ be any  sequence of points  such that $(p_n,q_n)\to (p,q)\in \MW  $. Then $\mu_n:=\mu(p_n,q_n)\nrw \mu=\mu(p,q)$ in $\mathscr{P}(\mathcal{F})$. By Skorohod theorem (see \cite{Billingsley, Parth}) there exists a representation $\mu_n=[Y_n]$, $\mu=[Y]$ for some sequence of random elements $Y_n,Y\in \mathcal{F}$ such that $Y_n\to Y$ in $\mathcal{F}$ almost surely. It follows that $\a(Y_n,p_n)\to \a(Y,p)$ almost surely and since $\a$ is bounded, we conlcude by dominated convergence that $\E\{\a(Y_n,p_n)\}\to \E\{\a(Y,p)\}$. This concludes the proof, since for all $n\in\N\cup\{\infty\}$:
\be 
\E\{\a(X,p_n)J_{p_n}X|X(p_n)=q_n\} =\int_{\mathcal{F}}\a(f,p_n)d\mu(p_n,q_n)(f)=\E\{\a(Y_n,p_n)\}.
\ee 
\end{proof}

When dealing with a KROK couple $(X,W)$, we will always implicitely assume that the function $(p,q)\mapsto [X|X(p)=q]$ is chosen among those for which $\mu$ is continuous at $\MW$. This arbitrary choice does not influence the final result, in that formula \eqref{eq:mainformula} depends only on $\mu|_{\MW}$.
\subsection{A closer look to the density}
In order to have a better understanding of the density $\delta_{X\in W}$, it is convenient to adopt a more general point of view. Let us consider, for any $V\subset M\times N$, the random number $\#_{\Gam}(V)=\#(\Gam\cap V)$, where $\Gam$ is the graph of the map $X|_{X^{-1}(W)}$, that is:
\be 
\Gamma(X,W):=\left\{(p,q)\in M\times W\colon X(p)=q\right\}, \quad \#_{\Gam}(V):=\#(\Gam\cap V).
\ee
The expectation $\E\#_{\Gam}$ of such random variable can be proven\footnote{The argument is exactly the same as that used to prove Lemma \ref{lemma:meas}} to be a Borel measure on $M\times N$ and by viewing it as an extension of the measure $\E\#_{X\in W}$, we can deduce Theorem \ref{thm:main} from the following slightly more general result.
\begin{thm}\label{thm:maingraph}
Let $(X,W)$, be a KROK couple, then the measure $\E\#_{\Gam}$ is supported on $\MW  =\{(p,q)\in M\times W\colon q\in S_p\}$ and it is an absolutely continuous measure on it, with a continuous density
\be\label{eq:gamdens}
\delta_{\Gam}(p,q)=\E\left\{J_pX{\sigma_q(X,W)}\bigg|X(p)=q\right\}\rho_{X(p)}(q)\delta \MW  (p,q),
\ee
where $\delta \MW  $ is the density on $\MW  $ defined below\footnote{We are implicitely making the identification $\Delta_{(p,q)}\MW  \cong \Delta_q(S_p\cap W)\otimes \Delta_pM$. By the KROK hypotheses \ref{def:krok}, $\MW  $ is a smooth submanifold of $M\times N$. However, $\delta {\MW  }$ is not the volume density of the metric induced by inclusion in the product Riemannian manifold $M\times N$.}.
\be 
\delta \MW  (p,q)=\frac{1}{\sigma_q(S_p,W)}d(S_p\cap W)(q)dM(p),
\ee  
\end{thm}
Precisely, this means that $ \E\#_{\Gam}(V)= \int_{V\cap \MW  }\delta_{\Gam}$,
for any Borel subset $V\subset M\times N$. In particular, if $V=A\times B$ we get
\be\label{eq:AperB}
\E\#_{X\in W\cap B}(A)=\E\#_{\Gam}(A\times B)=\int_A\left(\int_{B\cap W\cap S_p}\delta_{\Gam} (p,q)dq\right) dp.
\ee
for every $A\in \mathcal{B}(M)$ and $B\in \mathcal{B}(N)$.
\begin{remark}
The density $\delta_{X\in W}\in L^+(M)$ of the measure $\E\#_{X\in W}$ is obtained from the continuous density $\delta_{\Gam}\in \mathscr{D}^0(\MW  )$, by integration over the fibers of the projection map $\MW  \to M$.
\be\label{eq:twodens}
\delta_{X\in W}(p,q)=\int_{S_p\cap W}\delta_{\Gam}(p,q)dq.
\ee
This has to be intended as follows. The splitting $T_{(p,q)}\MW  \cong T_q(S_p\cap W)\oplus T_pM$ yields a natural identification $\Delta_{(p,q)} \MW  \cong \Delta_q (S_p\cap W)\otimes \Delta_p M$, allowing to define the \emph{partial integral} $\int_{S_p\cap W}\colon \mathscr{D}(\MW  )\to \Delta_p(M)$. 
\end{remark}
\begin{remark}
If $\mu=\E\#_{\Gamma(X,W)}$ on $M\times N$, then the integral of a measurable function $f\colon M\times N\to \R$ is given by the formula
\be 
\int fd\mu=\int_{\MW  }f\cdot \delta_{\Gamma(X,W)}=\E\left\{\sum_{p\in X^{-1}(q),\ q\in W}f(p,q)\right\}.
\ee
The proof of this fact, by monotone convergence, can be reduced to the case of characteristic functions $f=1_{A\times B}$, case in which  the formula is equivalent to equation \eqref{eq:AperB}.
\end{remark}

\subsection{The case of fiber bundles and the meaning of $\delta \MW$}

Let us consider the situation in which $\pi: N\to M$ is a smooth fiber bundle with fiber $S_p=\pi^{-1}(p)$
and let $W\subset N$ be a smooth submanifold such that $W\transv S_p$ for every $p\in M$.
Assume that $X\colon M\to N$ is a $\mC^1$ random section of $\pi$ 
and that $(X,W)$ is a KROK couple. 

In this case, the projection $\MW  =\{(p,q)\in M\times W\colon q\in S_p\cap W\}\to W$ is a bijection and we can identify the two spaces $\MW  \cong W$. Assume that both manifolds are endowed with Riemannian metrics in such a way that $\pi$ is a Riemannian submersion, meaning that the next map is an isometry\footnote{Such pair of metrics, always exists. To construct them, first define any metrics on $M$ and $N$. Then consider the subbundle $H\subset TN$ given by the orthogonal complement of the vertical one, namely $H_q=\ker(d_q\pi)^\perp$ (alternatively, take $H$ to be any Ehresmann connection). Now modifiy the metric on $H_q$ by declaring the map $d_q\pi\colon H_q\to T_pM$ a linear isometry.},
\be 
d_q\pi \colon \ker(d_q\pi)^\perp\to T_pM .
\ee
Then the formula \eqref{eq:gamdens} for the density, given in Theorem \ref{thm:maingraph} becomes easier and more meaningful.
\begin{thm}\label{thm:megafica}
Let $\pi\colon N\to M$ be a fiber bundle and Riemannian submersion. Let $(X,W)$ be a KROK couple such that $S_p=\pi^{-1}(p)$ for each $p\in M$. Then $\delta_{\Gam}$ is the continuous density on $W$ defined by the formula
\be\label{eq:megafica} 
\delta_{\Gam}(q)=\E\left\{J_{\pi(q)}X\sigma_q(X,W)\bigg|X(\pi(q))=q\right\}\rho_{X(\pi(q))}(q)dW(q).
\ee
\end{thm}
This is due to the fact that in this case we have $\delta \MW  =dW$.

\subsection{Independence on the metric}\label{sec:indepmetr}
It is important to note that the Riemannian structure on $N$ is just an auxiliary object that allows to write the formulas \eqref{eq:formdel}, \eqref{eq:gamdens}. In fact, the densities $\delta_{\Gam}$ and $\delta_{X\in W}$ must be independent of the chosen metric on $N$, since the corresponding measures have nothing to do with the Riemannian structure. 
Indeed, let us define a notation for the following expression:
\be 
d_pX\ \underset{T_qS_p}{\overset{T_qW}{\llcorner}}\ \delta_{X(p)}(q):=J_pX\frac{\sigma_q(X,W)}{\sigma_q(S_p,W)}\rho_{X(p)}(q)d(S_p\cap W)(q)dM(p).
\ee
This defines a density element in $\Delta_{(p,q)}\MW  $ depending only on the transverse vector subspaces $T_qW, T_qS_p\subset T_qN$, on the linear map $d_pX\colon T_pM\to T_qN$, and on the density element $\delta_{X(p)}(q)=\rho_{X(p)}(q)dS_p(q)\in \Delta_{q}S_p$. 
With this notation we can give a totally intrinsic version of formula \eqref{eq:gamdens}:
\be 
\delta_{\Gam}(p,q)=\E\left\{d_pX\ \underset{T_qS_p}{\overset{{T_qW}}{\llcorner}}\ \delta_{X(p)}(q)\ \Bigg|X(p)=q\right\}.
\ee

\section{The Gaussian case}\label{sec:gausscase}
\subsection{Smooth Gaussian random sections}\label{sec:gausscaseintro}
The first type of random maps that one encounters in random geometry are, with a very high probability, \emph{Gaussian random fields}, which are random maps $X\colon M\to \R^n$, whose evaluations at points are Gaussian (compare with Chapter \ref{chap:dtgrf}). In this section we are going to deal with a slight generalization of this concept, namely Gaussian random sections of a vector bundle. 

Precisely, let $\pi\colon E\to M$ be a smooth vector bundle of rank $s$ over a smooth $m$-dimensional manifold $M$ and let $X\colon M\to E$ be a random $\mC^r$ section of $\pi$. The random section $X$ is said to be \emph{Gaussian} if for any finite set of points $p_1,\dots, p_j$ the random vector 
\be 
\left(X(p_1), \dots , X(p_j)\right)\in E_{p_1}\oplus \dots \oplus E_{p_j}
\ee
is Gaussian. For simplicity in this paper we will assume all Gaussian variables to be centered, although this assumption is not necessary. Taking up the notation of Chapter \ref{chap:dtgrf}, we will denote as $\goo rME$ the set of $\mC^r$ Gaussian Random Sections (GRS) of a vector bundle $E$ over $M$. As for every Gaussian stochastic process, a GRS $X\in \goo rME$ is completely determined by its covariance tensor, which is the section $K_X\colon M\times M\to  \text{pr}_1^*E\otimes \text{pr}_2^*E$, defined by the following identity holding for every $\lambda_1,\lambda_2 \in E^*$:
\be 
K_X(p,q)\langle\lambda_1,\lambda_2\rangle=\E\left\{\lambda_1(X(p))\lambda_2(X(q))\right\}.
\ee
In particular, $K_X(p,p)=K_{X(p)}$ is a symmetric, semipositive, bilinear form on $E_p^*$.
\begin{defi}
 If $K_{X(p)}$ is positive definite (equivalently, $\spt[X(p)]=E_p$) for every $p\in M$, then we say that $X$ is \emph{non-degenerate}. 
 \end{defi}
In this case, if moreover $E$ is endowed with a bundle metric $g\colon E\cong E^*$, 
one can define the inverse covariance tensor, which is a bilinear form on $E_p$ that we denote by $K_{X(p)}^{-1}\left\langle \cdot,\cdot\right\rangle$. Then
we have $d[X(p)]=\rho_{X(p)}dE_p$ (in the sense of point \ref{itm:krok:3} of Definition \ref{def:krok}), where $dE_p$ is the Riemannian volume density of the fiber $E_p$ and
\be\label{eq:dgauvec}
\rho_{X(p)}(x)= \frac{\exp\left(-\frac{1}{2}K_{X(p)}^{-1}\left\langle x,x\right\rangle\right)}{\pi^\frac s2\sqrt{\det\left(K_{X(p)}\right)}}.
\ee
(The same formula is true in coordinates, if $K_{X(p)}$ denotes the covariance matrix.)

We want to apply Theorem \ref{thm:main} to compute the average number of points $p\in M$ such that $X(p)$ belongs to a given smooth submanifold $W\subset E$ of the total space $E$, having codimension $m$. In the Gaussian case it is particularly easy to verify the hypotheses of the theorem, indeed with the help of the (Gaussian) Probabilistic Transversality theorem, Theorem \ref{thm:transthm} from Chapter \ref{chap:dtgrf}, one deduces immediately that the couple $(X,W)$ is KROK (Definition \ref{def:krok}) if $X$ is non-degenerate and $W\transv E_p$ for every $p\in M$. It is easy to see that Theorem \ref{transthm2} from Chapter \ref{chap:dtgrf} can be extended to Gaussian sections, giving the following result. 
\begin{thm}\label{thm:probtransv}
Let $X\in\goo \infty ME$. Assume that for every $p\in M$
\be 
\spt[X(p)]=E_p.
\ee
Then for any smooth submanifold $W\subset E$, we have that $\P\{X\transv W\}=1$.
\end{thm}
From this, one deduces easily that the couple $(X,W)$ is KROK (Definition \ref{def:krok}) if $X$ is non-degenerate and $W\transv E_p$ for every $p\in M$. The only non obvious condition to check is KROK.\ref{itm:krokEcont}, which turns out to be a consequence of the Gaussian regression formula. This argument is used also in the proof of the standard Kac-Rice formula given in \cite{Wschebor}. Here, it is proved in Lemma \ref{lem:gioiellino}.
\begin{remark}
If $X\in \goo 1 ME$ is non-degenerate and $W\subset E$ is a submanifold such that $W\transv E_p$ for every $p\in M$, then $(X,W)$ is a KROK couple provided that $X\transv W$ almost surely. In the smooth case, the last hypothesis is redundant, due to Theorem \ref{thm:probtransv}.
This result holds only for sufficiently smooth fields, as well as its finite dimensional analogue, because it relies on Sard's Theorem. For this reason, we chose to focus on smooth GRS.
\end{remark}
The following theorem is the translation of the main Theorem \ref{thm:main} in the Gaussian setting. Although it is stated in a simpler way, it actually holds whenever the couple $(X,W)$ is KROK.
\begin{thm}\label{thm:maingau}
Let $\pi\colon E\to M$ and let $X\in\goo \infty ME$ be a non-degenerate $\mC^\infty$ Gaussian random section. Let $W\subset E$ be a smooth submanifold of codimension $m$ such that $W\transv E_p$ for every $p\in M$ and let $W_p=W\cap E_p$. Let the total space of $E$ be endowed with a Riemannian metric that is euclidean on fibers. 
Then for any Borel subset $A\subset M$
\begin{multline}\label{eq:formgaus}
\E\#_{X\in W}(A)= \int_{A}\delta_{X\in W}=\\
=\int_{A}\int_{W_p}\E\left\{J_pX\frac{\sigma_x(X,W)}{\sigma_x(E_p,W)}\bigg|X(p)=x\right\}\frac{e^{\left(-\frac{1}{2}K_p^{-1}
\left\langle x,x\right\rangle\right)}}{\pi^\frac s2\sqrt{\det(K_p)}}dW_p(x)dM(p).
\end{multline}
Here $K_p=K_X(p,p)$
; $dW_p$ and $dM$ denote the volume densities of the corresponding Riemannian manifolds; $J_pX$ is the Jacobian of $d_pX:T_pM\to T_xE$ (see Definition \ref{def:jacob}); besides, $\sigma_x(X,W)$ and $\sigma_x(E_p,W)$ denote the ``angles'' (in the sense of Definition \ref{defi:angle}) made by $T_xW$ with, respectively, $d_pX(T_pM)$ and $T_xE_p$.
\end{thm} 
We say that a Riemannian metric on the vector bundle $\pi\colon E\to M$ is Euclidean on fibers when the metric induced on each fiber $E_p\subset E$ is a vector space metric, meaning that $E_p$ is linearly isometric to the Euclidean space $\R^s$, as Riemannian manifolds. 

Such metric always exists on any vector bundle. The natural way to construct one is by defining a metric on $M$, a vector bundle metric on $E$ and an Ehressmann connection $H$ for the bundle $\pi$, that is: a vector subbundle of $TE$, the \emph{horizontal} bundle, such that $d_q\pi|_{H_q}\colon H_q\to T_{\pi(q)}M$ is a bijection. Then, the metric on $E$ is defined by declaring the implied isomorphism $T_qE\cong T_{\pi(p)}M\oplus_{\perp} E_q$ to be an isometry. A metric defined with this procedure is Euclidean on fibers, but it also make $\pi\colon E\to M$ a Riemannian submersion.
\begin{defi}
Let $\pi\colon E\to M$ be a vector bundle, such that $E$ is endowed with a metric constructed via a connection, with the above procedure. Then, we say that $\pi$ is a \emph{connected Riemannian bundle} or that it has a \emph{connected Riemannian metric}. We will say \emph{linearly connected} if the connection is linear.\footnote{A connection $H\subset E$ is linear if $H_{\lambda q}=d_qL_{\lambda}(H_q)$ for every $\lambda \in\R$, where $L_\lambda$ denotes the scalar multiplication. in this case, the operator $\nabla\colon \mC^\infty(M|E)\to \mC^{\infty}(M|E\otimes T^*M)$ satisfies the Leibnitz rule and thus it defines a covariant derivative.}
\end{defi}

Notice that in the case of Theorem \ref{thm:maingau} it is easy to see that $\MW $ is diffeomorphic to $W$ and $\delta_{\Gam}$ is a continuous density on it, although $\delta_{X\in W}\in L^+(M)$. Thus, by endowing $E$ with a connected Riemannian metric, Theorem \ref{thm:megafica} implies the following more natural formula. 
\begin{cor}\label{cor:maingauW}
In the same setting of Theorem \ref{thm:maingau}, assume that $\pi\colon E\to M$ is endowed with a connected Riemannian metric. Let $V\subset E$ be any Borel subset, then there exists a smooth density $\delta_{\Gam}\in \mathscr{D}^\infty(W)$ such that
\bega\label{eq:gaudelgamma}
\E\#\left(W \cap V\right)&=\int_{W\cap V}\delta_{\Gam}\\
&=\int_{W\cap V}
\E\left\{J_pX\sigma_x(X,W)\bigg|X(p)=x\right\}\frac{e^{\left(-\frac{1}{2}K_p^{-1}\left\langle x,x\right\rangle\right)}}{\pi^\frac s2\sqrt{\det(K_p)}}dW(x).
\eega
\end{cor}
\begin{remark}\label{rem:connectionForm}
If moreover $W$ is parallel, for the given connection, that is: $T_qW^{\perp}\subset E_q$, then at a point $q=X(p)$ we have
\be 
J_pX\sigma_q(X,W)=J\left(\Pi_{T_qW^\perp}\circ d_pX\right)=|\det\left(\Pi_{T_qW_p^\perp}\circ\left(\nabla X\right)_p\right) |,
\ee
where $W_p=W\cap E_p$ and $(\nabla X)_p\colon T_pM\to E_p$ is the vertical projection of $d_pX$. 
\end{remark}
We are (ab)using the symbol $\nabla$, since this notion coincides with that of \emph{covariant derivative}, in the case in which the connection is linear. 
Given that $W$ is transverse to the fibers of $\pi$,  one can always define a horizontal space $h_q=T_qW\cap (T_qW\cap T_qE_{\pi(q)})^\perp$ for each point $q\in W$. Then, if $h$ can be extended to  the whole $E$, it defines a connection (non linear, in general) $H$ for which $W$ is parallel. This construction is possible whenever $W\subset E$ is closed, by Tietze's extension theorem, but in general, there can be problems at $\overline{W}\- W$.

A particularly special case is when the connection is $\nabla=\nabla^X$ and the bundle metric on $E$ are the ones naturally defined by $X$ (see \cite{Nicolaescu2016}), namely $K_p$ is the dual metric and
\be \label{eq:nablaX}
\nabla^Xs:=Ds-\E\{DX|X=s\},
\ee
for any other connection $D$.
Since $\nabla^X$ is a metric connection in this case, it follows that for any Riemannian metric on $M$, a non-degenerate Gaussian random section defines a connected Riemannian structure on $E$. Moreover, since $(\nabla^XX)_p$ and $X(p)$ are independent, the formula in this case becomes much simpler.
\begin{cor}\label{cor:connectedform}
In the same setting of Theorem \ref{thm:maingau}, assume that $\pi\colon E\to M$ is endowed with the connected Riemannian structure defined by $X$. Let $W\subset E$ be parallel for this structure. Then
\bega\label{eq:connectedgaudelgamma}
\E\#\left(W \cap V\right)&=\int_{W\cap V}
\E\left\{|\det\left(\Pi_{T_xW_p^\perp}\circ\left(\nabla^X X\right)_p\right) |\right\}\frac{e^{-\frac{|x|^2}{2}}}{(2\pi)^{\frac s2}}dW(x).
\eega
\end{cor}

Thanks again to the Probabilistic Transversality theorem, the above result immediately generalizes to the case of a Whitney stratified submanifold (see \cite{GoreskyMacPherson}) $W\subset E$ of codimension $m$, simply because the probability that $X(M)$ intersects the lower dimensional strata is zero, therefore one can replace $W$ with its smooth locus $W^{(s)}$, namely the stratum of dimension $s$. In this case we still write $\delta_{X\in W}$ for the density, in place of $\delta_{X\in W^{(s)}}$.

\subsection{Finiteness and continuity}

\begin{enumerate}[\textbf{Q.}$1$:]
\item \emph{Is $\E\#_{X\in W}$ a Radon measure?} This question has positive answer precisely when the density is locally integrable, that is $\delta_{X\in W}\in L^1_{loc}(M)$. Theorem \ref{thm:maingau} leaves open the possibility that the density $\delta_{X\in W}(p)=\int_{W_p}\delta_{\Gam}(p,x)dx$ is even infinite. 
\item  \emph{Is the function $X\mapsto \E\#_{X\in W}(A)$ continuous?}  Understanding this is really useful in those situations where one is interested in the asymptotic behavior of things, for instance when dealing with Kostlan polynomials (see \cite{mttrp}). 
\end{enumerate}


From Corollary \ref{cor:maingauW} it is clear that $\delta_{X\in W}$ is finite at least in the case in which $W$ has finite volume. However, this would not be satisfying, since in many possible applications, $W$ has infinite volume. For instance, when $W\subset E$ is a vector subbundle of $E$, in fact, we will see that the density is finite in this case.  On the other hand, it should be clear that, due to the natural additivity of the formula:
$
\delta_{X\in \cup_n W_n}=\sum_n \delta_{X\in W_n},$
there are cases in which $\delta_{X\in W}(p)=+\infty$.
The intuition behind this is that if $W$ is too much ``concentrated'' over the fiber over a point $p_0$, then the probability that $X(p)\in W$ for some point near $p$ is too big, resulting in having $\E\#_{X\in W}(O)=+\infty$ for some neighborhood $O\subset M$ of $p_0$.

To express such concept, we introduce the notion of \emph{sub-Gaussian concentration}. This will allow us to compare the magnitude of $W$ with that of Gaussian sections, by passing through the linear structure of the bundle.
\begin{defi}\label{def:subexpvol}
Let $\pi\colon E\to M$ be a linearly connected Riemannian vector bundle. Let $B_R\subset E$ be the subset of vectors $e\in E$ with length at most $R$ for the given bundle metric.
We say that a smooth submanifold $W\subset E$ \emph{has sub-Gaussian concentration} if: for every compact $D\subset M$, the volume of $\pi^{-1}(D)\cap W\cap B_R$ (in the Riemannian manifold $W$) grows less than any Gaussian density, that is: $\forall \e>0$ $\exists C(\e)\ge 0$ such that $\forall R>0$
\be 
\text{Vol}_W\left(
\pi^{-1}(D)\cap W\cap B_R
\right)
=\int_{\pi^{-1}(D)\cap W\cap B_R}dW\le C(\e)e^{\e R^2}.
\ee
If $W$ is a Whitney stratified submanifold, we say that it has sub-Gaussian concentration if its smooth locus has sub-Gaussian concentration.

\end{defi}
It turns out that the property of having sub-Gaussian concentration is local and it is independent from the choice of a metric. In fact, this condition can be checked by proving that $W$ has sub-Gaussian concentration in each chart $E|_U\cong\D^m\times \R^s$ of a trivialization atlas for the bundle $E\to M$, and with respect to the standard metric. This is proved in Lemma \ref{lem:subGwell}. 
For this reason, in the following results we won't need to mention the Riemannian structure at all.

\begin{thm}\label{thm:mainEgau}
Let $\pi\colon M\to E$.
Let $W=\sqcup_{i\in I}W_i\subset E$ be a smooth Whitney stratified subset of codimension $m$ such that $W^{(s)}\transv E_p$ for every $p\in M$, where $W^{(s)}$ is the union of the higher dimensional strata.
Assume that $W$ has sub-Gaussian concentration. 
\begin{enumerate}[]
\item  Let $X\in\goo \infty ME$ be a non-degenerate $\mC^\infty$ Gaussian random section. Then $\delta_{X\in W}$ is locally integrable, hence $\E\#_{X\in W}$ is an absolutely continous Radon measure on $M$. 
\item  Let $X_d,X_\infty\in\goo \infty ME$ be a sequence of non-degenerate $\mC^\infty$ Gaussian random sections such that $K_{X_d}\to K_{X_\infty}$ in the $\mC^2$ topology (weak Whitney), as $d\to +\infty$. Assume that the limit $X_\infty$ is also non-degenerate. Then
\be 
\lim_{d\to +\infty}\E\#_{X_d\in W}(A)=\E\#_{X_\infty \in W}(A)
\ee
for every relatively compact Borel subset $A\subset M$.
\end{enumerate}
\end{thm}
\subsection{Semialgebraic singularities}
Clearly, if $W$ is compact, or a linear subbundle, or a cylinder over a compact, then it has sub-Gaussian concentration. The example that we are most interested in, though, is the case in which $W\subset E$ is locally semialgebraic. By this we mean that every $p\in M$ has a neighborhood $U\subset M$ such that there is a trivialization of the bundle $E|_{U}\cong \R^s\times \R^m$ such that $W\cap E|_{U}$ is a semialgebraic subset of $\R^s\times \R^m$. In this case the volume of $W\cap E|_{U} \cap B_R$ evidently grows in a polynomial way and thus...
\begin{remark}
...if $W\subset E$ is locally semialgebraic, it has sub-Gaussian concentration.
\end{remark}
The reason why we put the accent on the semialgebraic case is that Theorems \ref{thm:maingau} and \ref{thm:mainEgau} can be used to study the expected number of singular points of a GRS. The meaning of ``singular point'' depends on the situation, but in general it is a point $p\in M$ where the section satisfies some condition involving its derivatives. A general model for that  (the same proposed in \cite{mttrp} and \cite{HaniehPaulAnto}) is to consider a subset $W\subset J^rE$ of the bundle of $r$ jets (if the derivatives involved are of order less than $r$) of sections of $E$ (see \cite{Hirsch} for a definition of the space of jets) and call \emph{singular points of class $W$} those points $p\in M$ such that the $r^{th}$ jet of $X$ at $p$ belongs to $W$. 
Examples are:
\begin{enumerate}[$\bullet$]
\item ``Zeroes'', when $W=M\subset E=J^0(E)$.
\item ``Critical points'', when $W\subset J^1(M,\R^k)$, is  such that $j^1_pX\in W$ if and only if $d_pX$ is not surjective.
\item 
Combining the two previous examples, one can consider
$W\subset J^1(E\times \R)$, such that given a function $g\colon M\to \R$ and a section $s\colon M\to E$, then
$(j_p^1(s,g))^{-1}(W)=\text{Crit}(g|_{s^{-1}(0)})$. This is useful in that it provides an upper bound for the total Betti number of the set of zeroes of $s$. Indeed, generically, by Morse theory the latter must be smaller than the \emph{number} of  singular points of class $W$.
\item The Boardman singularity classes: $W=\Sigma^{(i_1,\dots,i_r)}\subset J^rE$, see \cite{arnold2012singularities}.
\end{enumerate}
In all of the above examples, and in most natural situations, the singularity class is given by a locally semialgebraic subset $W\subset J^rE$.

Considered this, we rewrite the statements of Theorems \ref{thm:maingau} and \ref{thm:mainEgau} in the case when the vector bundle is $J^rE\to M$ and the Gaussian random section is holonomic, namely it is of the form $j^rX$.

\begin{cor}\label{cor:mainjgau}
Let  $\pi\colon E\to M$ a smooth vector bundle. Let $W\subset J^rE$, with $r\in \N$, be a smooth submanifold of codimension $m$ such that $W\transv J^r_pE$ for every $p\in M$. 
Let $A\subset M$ be any relatively compact Borel subset.
\begin{enumerate}[$1$]
\item Let $X\in\goo \infty ME$ be a $\mC^{\infty}$ Gaussian random section with non-degenerate $r^{th}$ jet. Then there exists a smooth density $\delta_{\Gamma(j^rX,W)}\in \mathscr{D}^\infty(W)$ such that, 
\be\label{eq:formjgaus}
\E\#\{p\in A\colon j^r_pX_d\in W\}= \int_{W\cap \pi^{-1}(A)}\delta_{\Gamma(j^rX,W)}.
\ee
Moreover, if $W$ has sub-Gaussian concentration, then the above quantity is finite.
\item Assume that $W$ has sub-Gaussian concentration. Let $X_d, X_\infty\in\goo {\infty}ME$ be a sequence of  $\mC^{\infty}$ Gaussian random sections with non-degenerate $r^{th}$ jet and assume that 
$ 
K_{X_d}\to K_{X_\infty}
$ in the $\mC^{2r+2}$ topology (weak Whitney), as $d\to+\infty$. Assume that the limit $X_\infty$ also has non-degenerate $r^{th}$ jet. Then
\be 
\lim_{d\to +\infty}\E\#\{p\in A\colon j^r_pX_d\in W\}=\E\#\{p\in A\colon j^r_pX_{\infty}\in W\}.
\ee
\end{enumerate}
\end{cor} 
 
\begin{remark}
Here the formula for $\delta_{\Gamma(j^rX,W)}$ is obtained from formula \ref{eq:gaudelgamma} by replacing $E$ with $J^rE$, $X$ with $j^rX$ and $K_X$ with the covariance tensor of $j^rX$. Notice that the latter can be derived from the jet of order $2r$ of $K_X$.
\end{remark}

%
%
%

\section{Expectation of other counting measures}\label{sec:weights}
Let $X\colon M\to N \supset W$ be a random $\mC^1$ map. In this manuscript we chose to focus on the expectation of the actual number of points of intersection of $X$ and $W$. However, all the discussion can be generalized with minimal effort to the case in which a different weight is assigned to each point, in the following way.

For any $\a\colon \mC^1(M,N)\times M\to \R^k$ Borel measurable and $A\subset M$, define
\be\label{eq:Na}
\#^\a_{X\in W}(A)=\sum_{p\in A\cap X^{-1}(W)}\a(X,p) \in \R^k,
\ee 
and a density $\delta^\a_{X\in W}\colon M \to \Delta M$, such that
\be\label{eq:deltaa}
\delta_{X\in W}^\a(p)=\int_{ S_p\cap W}\E\left\{\a(X,p)\delta_pX\frac{\sigma_q(X,W)}{\sigma_q(S_p,W)}\bigg|X(p)=q\right\}\rho_{X(p)}(q)d(S_p\cap W)(q).
\ee
The following result extends Theorem \ref{thm:maindens} (Compare with \cite[Theorem 12.4.4]{AdlerTaylor} and \cite[Proposition 6.5]{Wschebor}, in the standard case.).
\begin{thm}\label{thm:mainalph}
Let $(X,W)$ be a KROK couple. Then Theorem \ref{thm:maindens} holds for $\#^\a_{X\in W}$ and $\delta^\a_{X\in W}$: for any Borel subset $A\subset M$ we have
\be\label{eq:mainalph}
\E\#^\a_{X\in W}(A)=\int_A\delta^\a_{X\in W} 
\ee
When not finite, both sides take the same infinite value among $+\infty,-\infty$ or $\infty-\infty$.
\end{thm}
\subsection{The intersection degree}
\newcommand{\dg}{\text{deg}}
\newcommand{\lk}{\text{link}}
Let $M$ be oriented and let $W\subset N$ be a closed cooriented submanifold. Then, given $q\in W$ and a linear bijection $L\colon T_pM\to T_qN/T_qW$, the sign of $\det(L)$ is well defined, since both vector spaces are oriented. If $f\colon M\to N$ is a $\mC^1$ map such that $f\transv W$, then the \emph{intersection degree} of $f$ and $W$ is defined as $\text{deg}(f,W)=\#^\a_{f\in W}(M)\in \Z$, where
\be 
\a(f,p)=\text{sgn}\left(\frac{d_pf}{T_{f(p)}W}\right).
\ee
In this situation, we can incorporate the sign in the definition of the angle, by defining $\text{sgn}\sigma_q(f,W)=\a(f,p)\sigma_q(f,W)$, so that the formula \eqref{eq:mainalph} for the expected intersection degree of a KROK couple $(X,W)$ and an open subset $U\subset M$ becomes
\be\label{eq:Edeg}
\E\{\dg(f|_U,W)\}=\int_U \int_{ S_p\cap W}\E\left\{\delta_pX\frac{\text{sgn}\sigma_q(X,W)}{\sigma_q(S_p,W)}\bigg|X(p)=q\right\}\rho_{X(p)}(q)d(S_p\cap W)(q).
\ee
\begin{remark}\label{rem:anto}
This confirms a general idea suggested to the author by Antonio Lerario, according to which the general philosophy to deal with these kind of formulas should be:
\emph{To get the formula for the signed count, add the sign to both members of the formula for the normal count.} Theorem \ref{thm:mainalph} can be thought as an extension of this philosophy from the \emph{sign} to any choice of \emph{weight} $\a$.
\end{remark}
\begin{remark} If $U\subset M$ be an open set whose closure is a compact topological submanifold with boundary, such that $f(\de U)\subset N\- W$. It can be seen that $\dg (f|_U,W)$ actually depends only on homological data and thus it is defined and locally constant on the space of continuous functions $f\colon \overline{U}\to N$, such that $f(\de U)\subset N\- W$. Moreover, if the Poincaré dual of $W$ in $N$ vanishes, then one can define the \emph{linking number}\footnote{
Since $W\subset N$ is cooriented, there exists a Thom class $\tau\in H^m(N,N\smallsetminus W)$. By definition, $\dg(f|_U,W)=\int_Uf^*\tau$ in the case $f\transv W$, but now this identity remains true for any continuous $f$ such that $f(\de U)\subset N\smallsetminus W$.
Looking at the long exact sequence for the pair $(N,N\smallsetminus W)$, we see that $\tau\mapsto e$ maps to the Poincaré dual of $W$, so that if $e=0$ then there exists an element $\a\in H^{m-1}(N\smallsetminus W)$ that maps to $\tau$ (i.e. $d^*\a=\tau$), therefore $\int_Uf^*\tau=\int_{\de U} f^*\a$ by naturality (or Stokes theorem from De Rham's point of view). In such case, $\lk(f|_U,W)=\int_{\de U}f^*\a$. If $N$ is a tubular neighborhood of $W$, then $e,\tau,\a$ are  respectively the Euler class, the Thom class and the class of a closed global angular form. Compare with \cite{botttu}.
} 
$\lk(g,W)\in \Z$  for any $g\colon B\to N\-W$, where $B$ is a closed manifold of dimension $m-1$. In this case $\dg(f|_U,W)=\lk(f|_{\de U},W)$ depends only on the restriction to the boundary. For this reason, $\dg(X|_U,W)$ can be thought to be \emph{less random} than $\#_{X\in W}(U)$; in fact, often it ends up being deterministic.
\end{remark}
\subsection{Multiplicativity and currents}\label{sub:multcurr}
The formula \eqref{eq:Edeg} can be written as the integral of a differential form $\E\{\dg(f|_U,W)\}=\int_{M_+} \Omega_{X\in W}$ over the oriented manifold $M_+$, defined as
\be\label{eq:Omega}
\Omega_{X\in W}(p)=\int_{ S_p\cap W}\E\left\{{(X^*\nu)_p}\big|X(p)=q\right\}\rho_{X(p)}(q)\frac{d(S_p\cap W)(q)}{\sigma_q(S_p,W)} \in \wedge^m T^*_pM,
\ee
where $\nu\in \wedge^m TW^\perp$ is the volume form of the oriented metric bundle $TW^\perp$. This follows simply from remark \ref{rem:altforms} and the fact that if $M$ and $N$ are Riemannian, then
\be 
\text{sgn}\left(\frac{d_pf}{T_{f(p)}W}\right)=\text{sgn}\det\left(\Pi_{ T_qW^\perp}\circ d_pX\right).
\ee
What is interesting about this is that the form $\Omega_{X\in W}\in \Omega^r(M)$ is still defined if the codimension of $W$ is $r< m$. 

If the couple $(X,W)$ satisfies the KROK hypotheses except for the requirement on the codimension of $W$ in KROK.\ref{itm:krok:4}, then let us say that $(X,W)$ is semi-KROK.
Let us consider semi-KROK couples $(X_i,W_i)$ on the manifolds $M, N_i$ for $i=1,\dots,k$ where the codimensions $r_i$ of $W_i$ are such that $r_1+\dots+r_k=m$ and all the $W_i$s are cooriented. Then it is easy to see that the product map $X:=X_1\times \dots \times X_k\colon M\to N:=N_1\times \dots \times N_k$ and the submanifold $W=W_1\times \dots \times W_k$ form a KROK couple $(X,W)$, and formula \eqref{eq:Omega} gives
\be 
\Omega_{X\in W}= \Omega_{X_1\in W_1}\wedge \dots \wedge \Omega_{X_k\in W_k}.
\ee
This can be interpreted in the language of currents, in the same spirit of \cite{Nicolaescu2016}.

\begin{claim}\label{claimcurr}
The $r$-form $\Omega_{X_i\in W_i}\in \Omega^r(M)$ is equal, as a current, to the expectation of the random current $[[X_i^{-1}(W)]]$ defined by the integration over the oriented submanifold $X_i^{-1}(W)$. 
\end{claim}

This follows from to the fact that the intersection degree $\dg(X|_U,W)$ may be viewed as the evaluation of the $0$-dimensional current $[[X^{-1}(W)]]=[[X_1^{-1}(W_1)\cap\dots\cap X^{-1}_k(W_k)]]$ on the function $1_U$. And since the intersection of currents is linear, we have 
\bega \label{eq:doubt}
\Omega_{X_1\in W_1}\wedge &\dots \wedge \Omega_{X_k\in W_k}
\\
=\Omega_{X\in W}&
=\E[[X^{-1}(W)]]
\\
=\E[[X_1^{-1}(W_1)]]\cap &\dots \cap \E[[X_k^{-1}(W_k)]].
\eega
\begin{remark}
Claim \ref{claimcurr} is not a theorem yet, for it has to be proved that the validity of the equation \eqref{eq:doubt} for every family of semi-KROK couples, implies that $\Omega_{X_i\in W_i}=\E[[X_i^{-1}(W_i)]]$. We plan to elaborate more on this subject in a subsequent paper.
\end{remark}

\subsection{Euler characteristic}
A special case of intersection degree is when $N=E\to M$ is the total space of an oriented vector bundle over an oriented manifold $M$ and $W=0_E\subset E$ is the zero section. Then $\dg(f,0_E)=\chi(E)$ is the Euler characteristic. In this case we can present the formula \eqref{eq:mainalph} in the form of Remark \ref{rem:connectionForm}.
\bega
\E\#_{X\in W}(U) &=\int_U \E\left\{|\det(\nabla X)_p|\bigg|X(p)=0\right\}\rho_{X(p)}(0)dM(p),
\\
\E\dg(X|_U,0_E) &=\int_U \E\left\{\ \det(\nabla X)_p\ \bigg|X(p)=0\right\}\rho_{X(p)}(0)dM(p)
\\
(&= \chi(E) ,\quad \text{If $U=M$ and is compact}),
\eega
where $E$ is endowed with some connected Riemannian metric\footnote{
Actually here the connection is not needed, since if $X(p)=0$ then $(\nabla X)_p$ is independent from $\nabla$.
 }.
  Notice that here, the density $\det(\nabla X)\rho_{X(\cdot)}(0)dM\colon M\to \Delta M$ is actually an intrinsic object, independent from the chosen Riemannian structures. 
  
  Following the discussion in the previous subsection \ref{sub:multcurr}, we now view the intersection degree $\dg(X|_{(\cdot)},0_E)$ as a random current $[[X^{-1}(0_E)]]\in\Omega_c^{0}(M)^*$. Its expectation is thus given by $\E[[X=0]] = \Omega_{X\in 0_E}$, defined in \eqref{eq:Omega}.

  
  Suppose now that $X\in\goo \infty M E$ is a nondegenerate smooth Gaussian random section and $s\colon M\to E$ is any smooth section. Let $E$ be endowed with the bundle metric defined by $X$ and let $\nabla=\nabla^X$ be the metric connection naturally associated with $X$ (see \cite{Nicolaescu2016}). Let $\nabla^2\in \bigwedge^2T^*M\otimes \bigwedge^2 E$ be its curvature. 
  Assume that $m$ is even, then a formula for $\Omega_{X\in 0_E}$ was computed in \cite{Nicolaescu2016}.
 \bega
 \Omega_{X\in 0_E}=
\E\left\{ \det(\nabla X)\right\}\frac{\w_{M}}{(2\pi)^{\frac m2}}=\frac{\text{Pf}(\nabla^2)}{(2\pi)^{-\frac{m}{2}}
}  \eega
 where $\w_{M}\in\Omega^m(M)$ is the volume form of the oriented Riemannian manifold $M$.  
\begin{remark}  
  The result of Nicolaescu \cite{Nicolaescu2016}, extends this to vector bundles with arbitrary even rank $r$. He proves that the expectation of the random current $[[X=0]]\in \Omega_c^{m-r}(M)^*$ is the current defined by the $r$-form $e(E,\nabla)=(2\pi)^{-\frac{r}{2}}Pf(\nabla^2)$. (Our sign convention in the definition of the Pfaffian $\text{Pf}$ is the same as in \cite{milnor-stasheff}, which differs to that of \cite{Nicolaescu2016} by a factor $(-1)^\frac {m}{2}$.)
 \end{remark}
 \subsection{Absolutely continuous measures}
 \HP\  on $M$ and $N$. Let $\mu$ be a finite Borel signed measure on $\mC^1(M,N)$ that is absolutely continuous with respect to $X$. By the Radon-Nikodym theorem, this is equivalent to the existence of an integrable Borel function $\a\colon \mC^1(M,N)\to \R$, with $\E\{|\a(X)|\}<\infty$, such that 
 \be 
 \int_{\mC^1(M,N)} \mathcal{F}(f)d\mu(f)=\E\{\mathcal{F}(X)\a(X)\}.
 \ee
Considering the case in which the function $\a$ in Theorem \ref{thm:mainalph} does not depends on the point, we deduce that the generalized Kac-Rice formula holds for every such measure $\mu$.
\be 
\int_{\mC^1(M,N)}\#(f^{-1}(W)\cap A)d\mu(f)=\E\#^\a_{X\in W}(A)=\int_{A}\delta^\a_{X\in W}.
\ee

\section{Examples}\label{sec:ex}

\subsection{Poincaré Formula for Homogenous Spaces}
In this section we will use Theorem \ref{thm:main} to give a new proof of the following Theorem. It is a special case of the Poincaré Formula for homogenous spaces, see \cite[Th. 3.8]{Howard}. 
\begin{thm}\label{thm:howard}
Let $G$ be a Lie group and let $K\subset G$ be a compact subgroup. Assume that $G$ is endowed with a left-invariant Riemannian metric that is also right-invariant for elements of $K$ and define the metric on $G/K$ to be the one that makes the projection map $\pi\colon G\to G/K$ a Riemannian submersion. Let $M,W\subset G/K$ be two smooth submanifolds (possibly with boundary) of complementary dimensions. Then
\be\label{eq:howardpoincare}
\int_G\#(gM\cap W)dG(g)=\int_M\int_W \bar{\sigma}_K\left(T_xM^\perp,T_yW^\perp\right)\Delta_K(x)dM(x)dW(y).
\ee
\end{thm}
Here, $\Delta_K(y)$ is defined as follows. Let $y=\nu K=\pi(\nu)$, then $\Delta_K(y)$ is defined to be the jacobian of the right multiplication by $\nu^{-1}$ in $G$,
\be 
\Delta_K(y):=\Delta(\nu)=J(R_{\nu^{-1}}).
\ee
To see the this definition is well posed, observe that since $\Delta\colon G\to \R$ is a group homomorphism and $K$ is compact, the set $\Delta(K)$ must be a compact subgroup of $\R$, thus $\Delta(K)=\{1\}$. It follows that $\Delta$ factorizes to a well defined function $\Delta_K\colon G/K\to \R$.

Notice that the angle $\bar{\sigma}_K\left(T_xM^\perp,T_yW^\perp\right)$ takes two subspaces that do not belong to the same tangent space $T_x(G/K)\neq T_y(G/K)$, in general. In fact $\bar{\sigma}_K$ is not the usual angle of Definition \ref{defi:angle}, but it is defined as follows. 
If $x=\mu K$ and $y=\nu K$ and $U,V$ are vector subspaces of, respectively, $T_x(G/K)\subset T_\mu G$ and $T_y(G/K)\subset T_\nu G$, then $\mu^{-1}_*U$ and $\mu^{-1}_*V$ are both subspaces of $T_1G$ and we can compute the angle $\sigma(\mu^{-1}_*U,\nu^{-1}_*V)$ using Definition \ref{defi:angle}. Then $\bar{\sigma}_K\left(U,V\right)$ is obtained by taking the average among all choices of $\mu$ and $\nu$:
\be 
\bar{\sigma}_K(U,V)=\int_K \sigma_{T_1 G}(\mu^{-1}_*U,k^{-1}_*\nu^{-1}_*V)dK(k).
\ee
Observe that, by Proposition \ref{prop:angleperp}, we have that $\bar{\sigma}_K(U,V)=\bar{\sigma}_K(U^\perp,V^\perp)$.
\begin{proof}[Proof of Theorem \ref{thm:howard}]
Let $\Omega\subset G$ be an open subset with compact closure and with $\vol_G(\de \Omega)=0$. Define the smooth random map $X\colon M\to G/K$ such that
\be
\P\{X=(L_g)|_{M}, \text{ for some $g\in A$}\}=\frac{\vol_G(A\cap \Omega)}{\vol_G(\Omega)}.
\ee
In other words, $X$ is uniformly distributed on the set of left multiplications by elements of $\Omega\subset G$, that means that $X\colon M\to G/K$ is the map $X(p)=g\cdot p$, where $g$ is a uniform random element of $\Omega$. We want to apply Theorem \ref{thm:main} to the random map $X$, where $W\subset N=G/K$. To this end, let us show that the couple $(X,W)$ is KROK (see Definition \ref{def:krok}).

\begin{enumerate}[(i)]
\item Ok.
\item \label{itm:how3} Let us fix $p=\mu K\in M$. The support of $[X(p)]$ is the open set $\Omega \cdot p\subset G/K$. Let $A\subset G/K$
\bega
\P\{X(p)\in A\}&=\frac{1}{\vol_G(\Omega)}\vol_G\left(\{g\in \Omega\colon g\mu K\in A\}\right)\\
&=\frac{1}{\vol_G(\Omega)}\vol_G\left(\Omega\cap\left(\pi^{-1}(A)\cdot\mu^{-1}\right) \right)\\
&=\frac{\Delta(\mu)}{\vol_G(\Omega)}\vol_G\left( \Omega\cdot\mu\cap\pi^{-1}(A)\right)\\
&=\frac{\Delta_K(p)}{\vol_G(\Omega)}\int_A\left(\int_{\Omega\cdot\mu\cap \pi^{-1}(q)}d\left(\pi^{-1}(q)\right)\right)d\left(G/K\right)(q) 
\\
&= \frac{\Delta_K(p)}{\vol_G(\Omega)}\int_A\vol_K(K\cap \nu^{-1}\Omega \mu)\ d\left(G/K\right)(\nu K),
\eega
where in the fourth step, we used the Coarea Formula (Theorem \ref{thm:coarea}). 

At this point, we can define a continuos function $f_\Omega\colon G/K\times G/K\to \R$ such that 
\be 
f_\Omega(p,q)=\vol_K(K\cap q^{-1}\Omega p):=\vol_K(K\cap \nu^{-1}\Omega\mu)\ee
 for any $\nu\in\pi^{-1}(q)$ and $\mu\in\pi^{-1}(p)$. This definition is well posed since the metric, and hence the volume form, on $G$ is invariant with respect to elements of $K$, both on the left and on the right. To see that $f_\Omega$ is continuous it is enough to show that if $\nu_n\to \nu$ and $\mu_n\to \mu$ in $G$, then $\vol_K(K\cap \nu_n^{-1}\Omega\mu_n)\to \vol_K(\nu^{-1}\Omega\mu)$.
 We prove this via the dominated convergence theorem (here we use the fact that $\vol_G(\de\Omega)=0$), since 
 \be 
 \vol_K(K\cap \nu_n^{-1}\Omega\mu_n)=\int_K 1_{\Omega}(\nu_n a \mu_n^{-1})dK(a),
 \ee
where $1_{\Omega}(\nu_n a \mu_n^{-1})\to 1_{\Omega}(\nu a \mu^{-1})$ for every $a\notin  \nu^{-1}\de\Omega \mu$ and, eventually,
\be 
1_{\Omega}(\nu_n a \mu_n^{-1})\le 1_{B_1(\nu^{-1})\cdot\Omega\cdot B_1(\mu)}(a)\in L^1(G).
\ee
It follows that $[X(p)]$ is absolutely continuous on $S_p:=G/K$, with a continuous density 
\be \label{eq:howarddensity}
\delta_{X(p)}(q)=\rho_{X(p)}(q)\cdot d(G/K)(q)= \frac{\Delta_K(p)}{\vol_G(\Omega)}f_\Omega(p,q)\cdot d(G/K)(q).
\ee
\begin{remark}
Notice that $\rho_{X(p)}(q)\neq 0$ if and only if $q\in \Omega\cdot p$. However, it would be a bad idea to define $S_p=\Omega\cdot p$. Indeed, with that choice the set $\mathcal{S}$ would not be closed and thus KROK.\ref{itm:krok:7} would not hold.
\end{remark}
\item Let us consider the smooth map $F\colon \Omega\times M\to G/K$, given by $F(g,p)=g\cdot p$. Clearly $F(\cdot,p)$ is a submersion, for every $p\in M$, therefore $F\transv W$. This, by a standard argument (Parametric Transversality, see \cite[Theorem 2.7]{Hirsch}) implies that $(L_g)|_M=F(g,\cdot)\transv W$ for almost every $g\in \Omega$. We conclude that
\be 
\P\{X\transv W\}=1.
\ee
\item Since $S_p=G/K$, the transversality assumption $S_p \transv W$ is certainly satisfied for every $p\in M$.
\item Ok.
\item In this case we have that $\mathcal{S}=M\times (G/K)$ is, without a doubt, a closed submanifold of $M\times (G/K)$. Moreover,
$\MW =M\times W$.
\item By equation \eqref{eq:howarddensity}, we have that $\rho_{X(p)}(q)$ is continuous with respect to all $(p,q)\in M\times (G/K)$.
\item This is the most complicated condition to check. Let $\a\colon \mC^1(M,G/K)\to \R$ be a continuous function (the case in which $\a$ depends on $p$ follows automatically, because of Proposition \ref{prop:CazzoDuro}). We have to show that the function $e\colon M\times G/K\to \R$, defined as $e(p,q)=\E\left\{\a(X)\big|X(p)=q\right\}$, is continuous at all points of $M\times W$. Let $p,q\in G/K$ and let $\nu,\mu\in G$  be such that $p=\mu K$ and $q=\nu K$. Notice that if $X=(L_g)|_{M}$, then $X(p)=q$ if and only $g\in \nu K\mu^{-1}$. Therefore if $\w$ denotes a uniformly distributed random element of $\Omega$ and $\xi$ denotes a uniformly distributed random element of $K$, then
\bega\label{eq:howcondi}
\E\left\{\a(X)\big|X(p)=q\right\}
&= \E\left\{\a\left((L_\w)|_M\right)\big|\w\in \Omega \cap \nu K \mu^{-1} \right\}
\\
&= \E\left\{\a\left((L_{\nu \xi\mu^{-1}})|_M\right)\big|k\in K\cap \nu^{-1} \Omega \mu \right\}
\\
&= \int_{K\cap \nu^{-1} \Omega \mu}\a\left((L_{\nu k\mu^{-1}})|_M\right) \frac{dK(k)}{\vol_K\left(K\cap \nu^{-1} \Omega \mu\right)}
\\
&=.\frac{1}{f_\Omega(p,q)}\int_{K\cap \nu^{-1} \Omega \mu}f_\a(\nu k\mu^{-1}) dK(k),
\eega
where $f_\a\colon G\to \R$ is a continuous function, defined as $f_\a(g)=\a\left(L_g|_M\right)$. Notice, that the last integral depends only on $p,q$, for if $\mu'=\mu a$ and $\nu'=\nu b$ for some $a,b\in K$, then the change in the integral corresponds to the change of coordinates $k'= b^{-1}ka$. 
To see that $e$ is continuous it is enough to check the continuity of the composed function $e(\pi(\cdot),\pi(\cdot))$. Let $\mu_n\to \mu$ and $\nu_n\to \nu$ be converging sequences in $G$ and define, for any $n_0\in\N$, the number
\bega
s_{n_0}=\sup_{n\ge n_0}\sup_{k\in K}|f_\a(\nu_nk\mu_n^{-1})-f_\a(\nu k\mu^{-1})|.
\eega
Since $f_\a$ is continuous and $K$ is compact, we have that $s_n\to 0$. As a consequence we get that
\bega
\lefteqn{\limsup_{n\to \infty}|e(\mu_nK,\nu_n K)-e(\mu K, \nu K)|\le}
\\
& &\le\lim_{n_0\to \infty}\sup_{n\ge n_0}\ \int_{K\cap \nu^{-1} \Omega \mu}\frac{\left|f_\a(\nu_n k\mu_n^{-1})-f_\a(\nu k\mu^{-1})\right|}{\vol_K\left(K\cap \nu^{-1} \Omega \mu\right)}{dK(k)}
\\
& &\le \lim_{n_0\to \infty}s_{n_0}=0.
\eega
This proves that $e(\pi,\pi)\colon G\times G\to\R$ is continuous.
\end{enumerate}
At this point, we know that the couple $(X,W)$ is KROK, therefore the expected number of intersections of the submanifolds $X(M)$ and $W$ is given by the generalized Kac-Rice formula of Theorem \ref{thm:main}, where $S_p=G/K$.
\begin{multline} 
\E\{\#_{X\in W}(M)\}=\int_{M}\rho_{X\in W}dM(p)
\\
=\int_M\int_W\E\left\{J_pX\sigma_q\left(d_pX(T_pM),T_qW\right)\big|X(p)=q\right\}\rho_{X(p)}(q)dW(q)dM(p).
\end{multline}
Recall that $X$ ranges among left translations, which are isometries, thus $JX=1$ with probability one. Moreover, we already computed $\rho_{X(p)}(q)$ (see \eqref{eq:howarddensity}) and we already understood how to compute the conditional expectation $\E\{\cdot|X(p)=q\}$ (see \eqref{eq:howcondi}). 
\begin{multline}
\E\{\#_{X\in W}(M)\}= \\
\int_M\int_W\frac{\Delta_K(\mu K)}{\vol_G(\Omega)}\left(\int_{K\cap \nu^{-1} \Omega \mu}\sigma_{\nu K}\left(\nu k \mu^{-1}T_{\mu K}M,T_{\mu K}W\right) dK(k)\right)dW(\nu K)dM(\mu K)=
\\
\int_M\int_W\frac{\Delta_K(\mu K)}{\vol_G(\Omega)}\left(\int_{K\cap \nu^{-1} \Omega \mu}\sigma_{T_1G}\left(\mu^{-1}T_{\mu K}M,k^{-1}\nu^{-1}T_{\nu K}W\right) dK(k)\right)dW(\nu K)dM(\mu K).
\end{multline}
Let $\Omega_n\subset \Omega_{n+1}\subset G$ be a sequence of relatively compact subsets such that $\cup_{n}\Omega_n =G$ and with $\vol_G(\de \Omega_n)=0$, and let $X_n\colon M\to G/K$ be the random map defined as above, with $\Omega=\Omega_n$. We obtain the thesis \eqref{eq:howardpoincare} by monotone convergence:
\bega
\int_{G}\#(gM\cap W)dG(g) = 
\sup_{n\in \N}\int_{\Omega_n}\#(gM\cap W)dG(g)
= \sup_{n\in\N}\E\{\#_{X_n\in W}(M)\}\vol_G(\Omega_n)
\\
=\int_M\int_W\Delta_K(\mu K)\left(\int_{K}\sigma_{T_1G}\left(\mu^{-1}T_{\mu K}M,k^{-1}\nu^{-1}T_{\nu K}W\right) dK(k)\right)dW(\nu K)dM(\mu K)
\\
=\int_M\int_W\Delta_K(p)\bar{\sigma}_{K}\left(T_{p}M,T_{q}W\right)dW(q)dM(p).
\eega
\end{proof}
\subsection{Isotropic Gaussian fields on the Sphere}
Let $X\colon S^m\to \R^k$ be a $\mathcal{C}^1$ Gaussian random field on the sphere. Using the notation of section \ref{sec:gausscaseintro} (consistently with Chapter \ref{chap:dtgrf}), we say that $X\in \g 1{S^m}k=\goo 1{S^m}E$, meaning that $X$ is a Gaussian random section of the trivial bundle $E=S^m\times \R^k$. We say that $X$ is \emph{isotropic} if $[X\circ R]=[X]$ for any rotation $R\in O(m+1)$. In particular, the covariance of $X(x)$ and $X(y)$ depends only on the angular distance $\a(x,y)=\arccos(\langle x,y\rangle)$ (because $S^m\times S^m/O(m+1)\cong [0,\pi]$):
\be 
K_X(x,y)=\E\left\{X(x)\cdot X(y)^T\right\}=K(\langle x,y\rangle)=F(\a(x,y)),
\ee
where $K\colon [-1,1]\to \R^{k\times k}$ and $F\colon \R\to \R^{k\times k}$ are $\mathcal{C}^2$ functions, such that $F$ is even and $2\pi-$periodic and
\be \label{eq:stropFK}
F(\a)=K(\cos\a).
\ee
As a consequence, the covariance structure of the first jet of $X$ at a given point $p\in S^m$, namely the couple $j^1_pX=(X(p),d_pX)$, is understood as follows. Define \be
\Sigma_0:=K(1)=F(0) \quad \text{and} \quad \Sigma_1=K'(1)=-F''(0).
\ee
Then given any orthonormal basis $\de_1,\dots, \de_m$ of $T_p S^m$, we have the following identities
\be\label{eq:istropId}
\E\left\{X(p)X(p)^T\right\}=\Sigma_0 ,
\quad 
\E\left\{\de_i X\cdot X(p)^T\right\}=0, \quad
\E\left\{\de_iX \cdot \de_jX^T\right\}=\Sigma_1\delta_{i,j}.
\ee
\begin{example}[Kostlan Polynomials]
Let $\psi_d\in\g 1{S^m}{}$, be defined as the restriction to $S^m$ of the random homogeneous Kostlan polynomial of degree $d$:
\be 
\psi_d(x)=\sum_{|\a|=d}{d\choose\a}^\frac12 \gamma_\a x^\a,
\ee
where $\gamma_\a\sim N(0,1)$ are independent normal Gaussian. Then $\psi_d$ is a smooth isotropic Gaussian field $X\in\g 1{S^m}k$ with $K(t)=t^d$. In fact, any isotropic Gaussian field for which the function $K(t)$ has the form
\be
K(t)=K_0+K_1t+\dots +K_dt^d,
\ee 
for some positive definite simmetric matrices $K_\ell$, is a linear combination of Kostlan fields. To see this, let $A_\ell$ be a $k\times k$ matrix such that $A_\ell A_\ell^T=K_\ell$ and define 
\be \label{eq:stropKostlanmix}
\tilde{X}=A_0\begin{pmatrix}
\psi_0^1 \\ \vdots \\ \psi_0^k
\end{pmatrix}+A_1\begin{pmatrix}
\psi_1^1 \\ \vdots \\ \psi_1^k
\end{pmatrix}
+\dots +
A_d\begin{pmatrix}
\psi_d^1 \\ \vdots \\ \psi_d^k
\end{pmatrix},
\ee
where $\{\psi_\ell^i\}_{i,\ell}$ are independent copies of Kostlan polynomials of degree $\ell=1,\dots,d$. Then $\tilde{X}$ is equivalent to $X$, since they have the same covariance function. In the particular case in which $X=(\psi_{d_1},\dots,\psi_{d_k})^T$ where $\psi_{d_\ell}$ are independent Kostlan polynomials of degree $d_\ell$, then $\Sigma_0=\mathbb{1}_k$ and $\Sigma_1=\text{diag}(d_1,\dots,d_k)$.
\end{example}
We recall that the density function of a Gaussian random vector $\xi\sim N(0,\Sigma)$ in $\R^k$ with nondegenerate covariance matrix $\Sigma$ is given by
\be\label{eq:isodensga}
\rho_{\Sigma}(y)=\frac{e^{-\frac12y^T\Sigma^{-1}y}}{(2\pi)^\frac{k}2(\det \Sigma)^\frac12}.
\ee
\begin{lemma}\label{lem:MWRn}
Let $M$ be a Riemannian manifold. Let $X\in \g \infty Mk$ and assume that $X(p)$ has nondegenerate covariance matrix  $\Sigma_0(p)$. Let $W\subset \R^k$ be any submanifold (possibly stratified) of codimension $m$. Then 
\begin{multline}
\E\{\#X^{-1}(W)\}=\int_M\delta_{X\in W}(p)
\\
=\int_M\int_W\E\left\{\left|\det\left(\Pi_{T_yW^\perp}\circ d_p X\right)\right|\Bigg| X(p)=y\right\}\rho_{\Sigma_0}(y) dW(y)dM(p).
\end{multline}
\end{lemma}
\begin{proof}
Let $\Gamma_X$ be the graph of $X$ and $\MW :=M\times W$. Then $\Gamma_X$ is a non-degenerate Gaussian random section of the trivial bundle $E=M\times \R^k$ and $\MW $ is clearly transverse to all fibers of the bundle: $M\times W\transv \{p\}\times \R^k$. Therefore we can apply Theorem \ref{thm:maingau} to obtain a density $\delta_{X\in W}=\delta_{\Gamma(X)\in \MW }\in L^+(M)$. Moreover, the trivial connection on $E=M\times \R^k$ makes it into a linearly connected Riemannian bundle, for which $\MW $ is a parallel submanifold, so that we can present the density $\delta_{X\in W}$ with the formula of Remark \ref{rem:connectionForm}, knowing that $\nabla \Gamma_X=dX$ and $\MW \cap E_p=W$.
\be 
\delta_{X\in W}(p)=\int_W\E\left\{\left|\det\left(\Pi_{T_yW^\perp}\circ d_p X\right)\right|\Bigg| X(p)=y\right\}\rho_{\Sigma_0}(y) dW(y)dM(p).
\ee
\end{proof}
Given a smooth submanifold $W\subset \R^k$ (possibly stratified) of codimension $m$, we say that a measurable map $\nu\colon W\to \R^{m\times k}$ is a \emph{measurable normal framing} for $W$, if for almost every $y\in W$ the columns of the matrix $\nu(y)$ form an orthonormal basis of $T_yW^\perp$ (if $W$ is stratified, then this has to hold only for almost every $y$ in the top dimensional stratum of $W$).
\begin{thm}\label{thm:mainIsotrop}
Let $X\colon S^m\to \R^k$ be an isotropic $\mathcal{C}^1$ Gaussian random field. Let $\Sigma_0$ and $\Sigma_1$ be the $k\times k$ matrices defined as above and assume that $\Sigma_0$ is nondegenerate. Let $W\subset \R^k$ be any submanifold (possibly stratified) of codimension $m$ and let $\nu \colon W\to \R^{m\times k}$ be a measurable normal framing. Then
\be
\E\left\{\# X^{-1}(W)\right\}=2\cdot \int_W\sqrt{\det(\nu(y)^T\Sigma_1\nu(y))}\cdot\frac{e^{- \frac12y^T\Sigma_0^{-1}y}}{(2\pi)^\frac{k-m}2(\det \Sigma_0)^\frac12}dW(y)
.
\ee
\end{thm}
\begin{proof}
By Lemma \ref{lem:MWRn}, we have that 
\be 
\E\{\#X^{-1}(W)\}=\int_W\E\left\{\left|\det\left(\Pi_{T_yW^\perp}\circ d_p X\right)\right|\right\}\rho_{\Sigma_0}(y) dW(y)dS^2(p).
\ee
We can omit the conditioning $X(p)=y$, since in this case $X(p)$ and $d_pX$ are independent. The fact that the field is isotropic implies that the measure $\E\#_{X\in W}$ is an invariant measure on $S^m$, so that 
\be 
\E\{\#X^{-1}(W)\}=\vol(S^m)\int_W\E\left\{\left|\det\left(\Pi_{T_yW^\perp}\circ d_p X\right)\right|\right\}\rho_{\Sigma_0}(y) dW(y),
\ee
where $p\in S^2$ is any point. Let $\de_1,\dots,\de_m$ be an orthonormal basis of $T_pM$. It remains only to compute the expectation of the determinant of the random matrix $A(y)$ with columns $A_1(y),\dots, A_m(y)$ defined as
\be 
A_i(y):=\nu(y)^T\de_i X=\left(\Pi_{T_yW^\perp}\circ d_p X\right)(\de_i).
\ee
By the third set of the identities in  \eqref{eq:istropId}, we deduce that the columns of $A(y)$ are independent Gaussian random vectors, each of them having covariance matrix $K(y)=\nu(y)^T\Sigma_1\nu(y)$. Therefore
\bega
\E\left\{\left|\det\left(\Pi_{T_yW^\perp}\circ d_p X\right)\right|\right\}=\E\left\{\left|\det A(y)\right|\right\}=\sqrt{\det \nu(y)^T\Sigma_1\nu(y)} \frac{m!\vol(\B^m)}{(2\pi)^{\frac{m}{2}}}.
\eega

To conclude, let us observe that $\vol(S^m)\vol(\B^m) m!=2(2\pi)^m$, because of a special property of the Gamma function: $\Gamma(z+\frac12)\Gamma(z)=2^{1-2z}\sqrt{\pi}\Gamma(2z)$.
\bega 
\vol(S^m)\vol(\B^m) m!&= \frac{\pi^\frac{m+1}2(m+1)}{\Gamma\left(\frac{m+1}2 +1\right)} \frac{\pi^\frac{m}2}{\Gamma\left(\frac{m}2 +1\right)}m!
=\frac{
\pi^{m+\frac12}(m+1)!
}{
2^{- (m+1)}\sqrt{\pi}\Gamma(m+2)
}
=2^{m+1}\pi^m.
\eega
\end{proof}
\begin{remark} In the case $\Sigma_0=\sigma^2\mathbb{1}_k$, we obtain a particularly nice formula
\be 
\E\{\#X^{-1}(W)\}=\int_W \E\{\#X^{-1}(T_yW)\}\rho_{\sigma^2\mathbb{1}_{k-m}}(y)
dW(y).\ee
This allows to reduce to the case $k=m$ and to the standard version of Kac-Rice formula. Indeed $X\in (T_yW)$ if and only if $\nu(y)^TX=0$. For completeness, we report two results that can be proved by applying the standard Kac-Rice formula (Corollary \ref{thm:stropzero} and \ref{thm:stropshubsmale} are not new results).
\end{remark}
\begin{cor}(Gaussian Isotropic Kac-Rice formula)\label{thm:stropzero}
Let $X\colon S^m\to \R^m$ be an isotropic $\mathcal{C}^1$ Gaussian random field. 
Let $\Sigma_0$ and $\Sigma_1$ be the $m\times m$ matrices defined by the identities \eqref{eq:istropId} and assume that $\Sigma_0$ is nondegenerate. Then, for any $y\in\R^m$,
\be 
\E\left\{\# X^{-1}(y)\right\}=2\cdot \sqrt{\frac{\det(\Sigma_1)}{\det(\Sigma_0)}}e^{-\frac12 y^T\Sigma_0^{-1}y}.
\ee
\end{cor}
\begin{cor}(Shub-Smale Theorem \cite{shsm})\label{thm:stropshubsmale}
 Let $\psi_{1},\dots ,\psi_{m}$ be independent Kostlan homogenous polynomials of degrees $d_1,\dots, d_m$ and denote by $Z\subset\RP^m$ the random subset defined by the equations $\psi_i=0$. Then $\E\{\#Z\}=\sqrt{d_1 \cdot\cdot\cdot  d_m}$.
\end{cor}
Notice that Corollary \ref{thm:stropzero} covers also the case in which the equations in the Shub-Smale Theorem \ref{thm:stropshubsmale} are dependant, as long as they are jointly orthogonally invariant, in particular in the case of a mixed Kostlan polynomial defined as in equation \eqref{eq:stropKostlanmix}, we have 
\be
\E\{\#X^{-1}(0)\}=2\sqrt{\frac{\det\left(A_0A_0^T+A_1A_1^T+\dots +A_dA_d^T\right)}{\det\left(A_1A_1^T+2A_2A_2^T+\dots +dA_dA_d^T\right)}}.
\ee
%


\section{Proof of Lemma \ref{lemma:meas}}\label{sec:proflemmameas}

\begin{proof}[Proof of Lemma \ref{lemma:meas}]
We can assume that $W$ is closed, by replacing $X\colon M\to N$, with the random map $\XW\colon M\times W\to N\times N$, such that $(p,w)\mapsto (X(p),w)$ and $W\subset N$ with the diagonal $\DElta\subset N\times N$, which is surely a closed submanifold. It is straightforward to see that $X\transv W$ if and only if $\XW\transv \DElta$ and that $\#_{X\in W}(A)=\#_{\XW\in \DElta}(A\times W)$, for every $A\subset M$, moreover if $A\in \mathcal{B}(M)$, then $A\times W\in \mathcal{B}(M\times W)$\footnote{A consequence of this trick is that $\mu(A\times B):=\E\#_{X\in W\cap B}(A)$ defines a Borel measure on the product space $M\times W$. We are going to develop this idea properly later, in Section \ref{sec:constru}.}.

Let $\mathcal{D}$ be the family of subsets $A\subset M$ such that the function $\#_{X\in W}(A)$ is measurable. The class $\mathcal{D}$ contains the subfamily $\mathcal{P}$ of all relatively compact open sets in $M$, which is closed under intersection, hence the idea is to prove that it is also a Dynkin class\footnote{Let $M$ be a nonempty set; a \emph{Dynkin class} $\mathcal{D}$ is a collection of subsets of $M$ such that:
\begin{enumerate}
    \item $M\in \mathcal{D}$;
    \item if $A, B\in \mathcal{D}$ and $A\subset B$, then $B\backslash A\subset \mathcal{D};$
    \item given a family of sets $\{A_k\}_{k\in \N}$ with $A_k\in \mathcal{D}$ and $A_k\subset A_{k+1}$, then $\bigcup_{k}A_k\in \mathcal{D}$.
\end{enumerate}
The Monotone Class Theorem (see \cite[p. 3]{Erhan}) says that if a family $\mathcal{D}$ is a Dynkin class which contains a family $\mathcal{P}$ closed by intersection, then it contains also the  $\sigma$-algebra generated by $\mathcal{P}$.} to conclude, by the Monotone Class Theorem (see \cite[p. 3]{Erhan}), that $\mathcal{D}$ contains the $\sigma$-algebra generated by $\mathcal{P}$, which is precisely the Borel $\sigma$-algebra $\mathcal{B}(M)$.  

Actually, to prevent $\#_{X\in W}(A)$ from taking infinite values, it is more convenient to consider a countable increasing family of relatively compact open subsets $M_i$ such that $\cup_i M_i=M$ and work with the class $\mathcal{D}_i=\{A\in\mathcal{D}\colon A\subset M_i\}$, since $\#_{X\in W}(M_i)$ is almost surely finite.

By previous considerations, $M_i\in \mathcal{D}_i$. If $A,B\in\mathcal{D}_i$ and $A\subset B$, then  $B\- A \in \mathcal{D}_i$, because since $\#_{X\in W}(B)$ is almost surely finite, we can write $\#_{X\in W}(B\- A)=\#_{X\in W}(B)-\#_{X\in W}(A)$ . Suppose that $A_k\in \mathcal{D}$ is increasing, then 
\be \label{eq:contbelow}
\#_{X\in W}(\cup_k A_k)=\lim_{k}\#_{X\in W}(A_k),
\ee
thus $A=U_k A_k\in \mathcal{D}$ because $\#_{X\in W}(A)$ is the pointwise limit of measurable functions, thus in particular if $A_k\in\mathcal{D}_i$, then $A\in \mathcal{D}_i$. It follows that $\mathcal{D}_i$ is indeed a Dynkin class, hence $\mathcal{D}\supset\mathcal{D}_i\supset \mathcal{B}(M_i)$. Now let $A\in\mathcal{B}(M)$, then $A$ is the union of the increasing sequence $A\cap M_i$ and since $A\cap M_i\in\mathcal{B}(M_i)\subset \mathcal{D}$, we can use again the formula in \eqref{eq:contbelow}, to conclude that $A\in\mathcal{D}$.

Clearly $\en$ is finitely additive and $\en(\emptyset)=0$, therefore to prove that $\en$ is a measure it is enough to show that it is continuous from below. This can be seen just by taking the mean value in \eqref{eq:contbelow} and using the Monotone Convergence theorem, since $\#_{X\in W}(A_k)$ is increasing for any increasing sequence $A_k\in\mathcal{B}(M)$.
\end{proof}

%
%

\section{General formula}\label{sec:gen}

This section is devoted to the proof of the following theorem. It is a more general result than the main Theorem \ref{thm:main}, but in fact it is too abstract to be useful on its own. Its role is to create a solid first step for the proof of the main theorem and to better understand its hypotheses.

\begin{figure}\begin{center}
\includegraphics[scale=0.2]{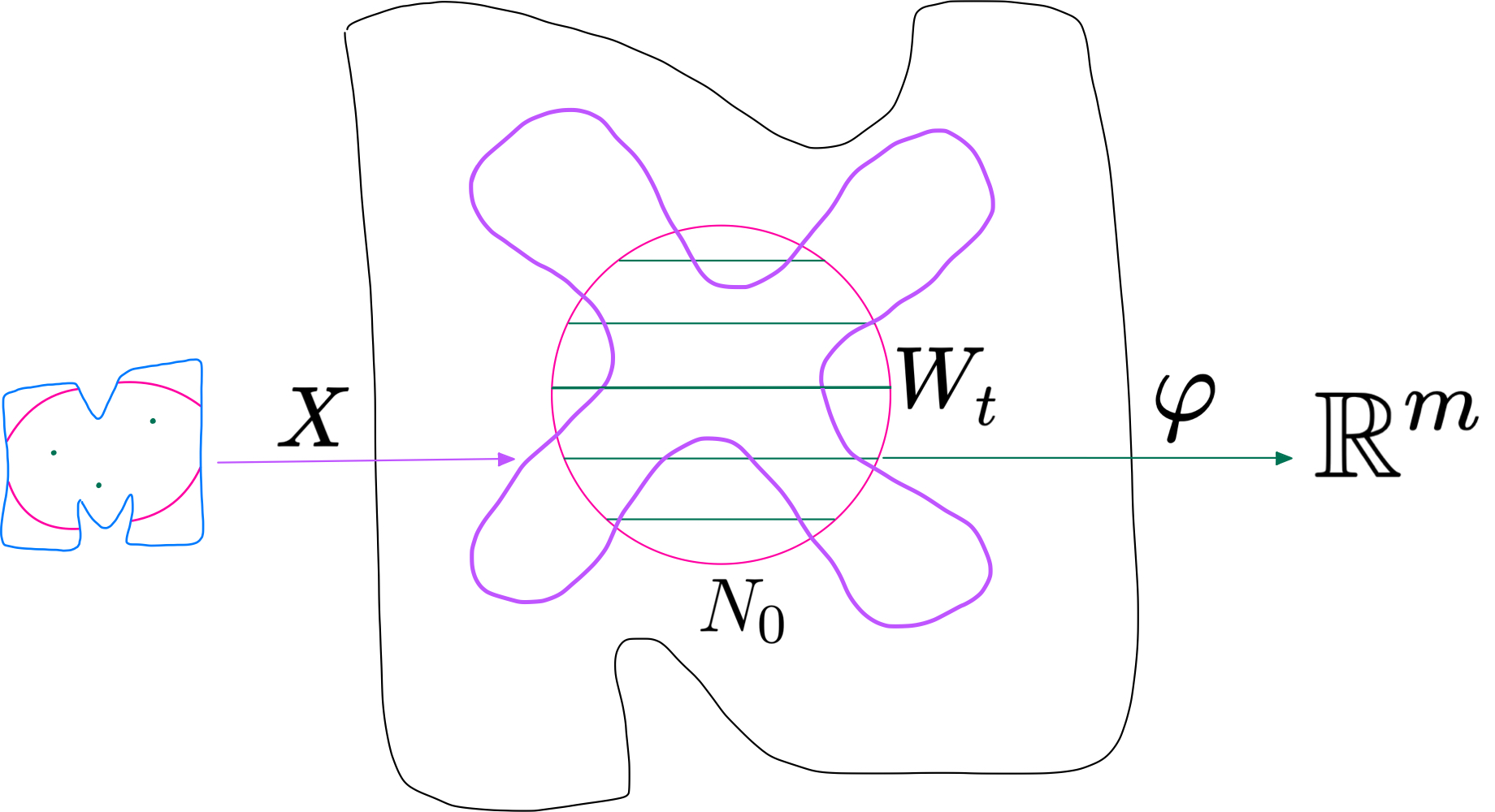}\caption{This figure is meant to give an idea of the set theoretic positions of the objects involved in Theorem \ref{thm:general}.}\label{fig:MXNR}
\end{center}
\end{figure}
\begin{thm}\label{thm:general}
Let $X\colon M\to N$ be a $\mC^1$ random map, such that $d[X(p)]=\rho_{X(p)}dS_p$\footnote{See point \ref{itm:krok:3} of Definition \ref{def:krok}.}, for some Riemannian submanifold $S_p\subset N$ and measurable function $\rho_{X(p)}\colon S_p\to [0,+\infty]$. Let $\{W_t\}_{t\in \R^m}$ be a smooth foliation of an open set $N_0\subset N$, defined by a submersion $\f\colon N_0\to \R^m$ via $W_t=\f^{-1}(t)$,  such that $W_t\transv S_p$ for all $p$ and $t$. Consider the density $\delta_{X\in W_t}(p)$ defined by the same formula given in \eqref{eq:formdel}:
\be
\delta_{X\in W_t}(p)=\int_{ S_p\cap W_t}\E\left\{J_pX\frac{\sigma_x(X,W_t)}{\sigma_x(S_p,W_t)}\bigg|X(p)=x\right\}\rho_{X(p)}(x)d(S_p\cap W_t)(x)dM(p),
\ee
and assume it to be a measurable function with respect to the couple $(p,t)$. Let $A\in \mathcal{B}(M)$ be any Borel subset of $M$.
\begin{enumerate}[$(i)$]
    \item\label{itm:general:1} If $A\in\mathcal{B}(M)$, 
then for almost every $t\in \R^m$
    \be\label{eq:askr}
    \E\#[X^{-1}(W_t)\cap A]=\int_A \delta_{X\in W_t}.
    \ee 
Equivalently, there is a full measure set $T\subset \R^m$, such that for all $t\in T$, the set function $A\mapsto \E\#_{X\in W_t}(A)=\E\#[X^{-1}(W_t)\cap A]$ is an absolutely continuous Borel measure on $M$ with density $\delta_{X\in W_t}$.
    \item\label{itm:general:3} Let $X\transv W_{t_0}$ almost surely and assume that there exists a density $\delta_{top}\in L^1_{loc}(M)$ such that $\limsup_{t\to t_0}\int_{K}\delta_{X\in W_t}\le \int_K{\delta_{top}}$ for every compact set $K\subset M$. Then the measure $\E\#_{X\in W_{t_0}}$ is an absolutely continuous Radon measure on $M$. 
    In this case the corresponding density $\delta\in L^1_{loc}(M)$ satisfies $\delta\le \delta_{top}$.
    
    In particular, if $\delta_{X\in W_t}\to \delta_{X\in W_{t_0}}$ in $L^1_{loc}(M)$, then 
    \be 
    \E\#_{X\in W_{t_0}}(A)=\int_{A}\delta\le\int_{A}\delta_{X\in W_{t_0}}.
   \ee
\end{enumerate}
\end{thm}
\begin{remark}\label{rem:toni}
The left hand side of equation \eqref{eq:askr} is well defined for almost every $t\in\R^k$. This can be seen as follows. Using Tonelli's theorem, we can prove that $\P\{X\transv W_t\}=1$ for almost every $t\in \R^m$:
\be 
\begin{aligned}
\int_{\R^m}\P\{X\not\transv W_t\}dt
&=d[X]\otimes \mathscr{L}^m(\left\{(f,t)\in \ci(M,N)\times \R^m:(\f\circ f)|_{f^{-1}(N_0)}\transv \{t\} \right\}) =\\
&= \E\{\mathscr{L}^m(\{\text{critical values of }(\f\circ X)|_{X^{-1}(N_0)}\})\}=0,
\end{aligned}
\ee
The last equality following from Sard's theorem: critical values of a $\ci$ map between two manifolds of the same dimension form a set of zero Lebesgue measure. Combining this fact with Lemma \ref{lemma:meas}, we deduce that there is a full measure set $T\subset \R^m$ such that the set function
\be 
\mathcal{B}(M)\ni \quad A\mapsto \E\#_{X\in W_t}(A):=\E\#[X^{-1}(W_t)\cap A]\quad \in [0,+\infty]
\ee
is well defined and is a Borel measure.
\end{remark}
The density $\delta_{X\in W_t}$ appearing above is the same that appears in the statement of Theorem \ref{thm:maindens},
where $M$ and $N$ are endowed with auxiliary Riemannian metrics. Since the conditional expectation  
\be 
x\mapsto \E\left\{J_pX\frac{\sigma_x(X,W_t)}{\sigma_x(S_p,W_t)}\bigg|X(p)=x\right\}
\ee
is defined (for every $p$) up to almost everywhere equivalence and $S_p\cap W_t$ has Lebesgue measure zero, it follows that the value of $\delta_{X\in W_t}(p)$ is not uniquely determined. By saying that  $\delta_{X\in W_t}(p)$ should be measurable in $p$ and $t$, we mean that we assume to have chosen a representative of the conditional expectation above in such a way that the function $(p,t)\mapsto \delta_{X\in W_t}(p)$ is measurable. In the KROK case (Definition \ref{def:krok}), however, there is no such ambiguity (see Subsection KROK.\ref{itm:krokEcont}).

Moreover, notice that the value of $\delta_{X\in W_t}(p)$ depends on the choices of $S_p$ and $\rho_{X(p)}$ rather than just $[X(p)]$. 
Indeed, the choice of the submanifold $S_p$ such that $d[X(p)]=\rho_{X(p)}dS_p$ is not unique in general. In fact, $S_p$ can even be replace with $S_p\-W_0$, so that $\delta_{X\in W_0}=0$. This is not in contradiction with the theorem, because the identity 
\eqref{eq:askr} is valid for all $t$ out of a measure zero set. Again, in the KROK case we don't have to worry about that, because, by point \ref{itm:krok:7} of \ref{def:krok}, $S_p\subset N$ is required to be closed.
\begin{remark}\label{rem:sorry}
There is a treacherous measurability issues, that the author wasn't able to solve.
 \emph{Given a random $\mC^1$ map $X\colon M\to N$ and a map $\f\colon N_0\subset N\to \R^m$ satisfying the hypotheses of Theorem \ref{thm:general}, is it true that there exists a version of the density $\delta_{X\in W_t}$ that is measurable with respect to $(p,t)$?} Measurability is crucial to use Fubini-Tonelli's Theorem, thus in Theorem \ref{thm:general} we assume that we are in a situation where the previous question has a positive answer, though such hypothesis may be redundant.
\end{remark}

Before going into the proof of Theorem \ref{thm:general}, let us prove two important preliminary results. The first (see Corollary \ref{cor:measurable}) ensures that the number $\#[X^{-1}(W_t)\cap A]$ is measurable in $(X,t)$ for every Borel set $A$. The second (see Lemma \ref{lem:sette}) gives an alternative expression for $\delta_{X\in W_t}$, that is convenient to use in the proof of the Theorem.
\begin{lemma}\label{lem:fsigmeas}
Let $M,N$ and $T$ be metrizable topological spaces and let $\f\colon K\to T$ be a continuous function,where $K\subset N$ is a closed subset. Then, for every compact set $A\subset M$, the function below is Borel
\be 
\coo{0}{M}{N}\times T \ni \quad (f,t)\mapsto \#(f^{-1}(\f^{-1}(t))\cap A)\quad\in \N\cup \{\infty\}.
\ee
\end{lemma}
\begin{proof}

Fix $\e>0$ and define, for each subset $S\subset M$ the number $\#_{\e}(S)$ to be the minimum number of open subsets of diameter smaller than $\e$ that are needed to cover $S$\footnote{It corresponds to the set function $\mathcal{H}_\e^0(S)$ used in the construction of the Hausdorff measure. }. Observe that $\#(S)=\sup_{\e>0}\#_\e(S)$, therefore we can conclude the proof by showing that the function $(f,t)\mapsto \#_\e(f^{-1}(\f^{-1}(t)))$ is Borel for every $\e$.

Let us consider two convergent sequences $f_n\to f$ in $\coo0MT$ and $t_n\to t$. Assume that $B_1,\dots, B_k$ are open balls in $M$ of diameter smaller than $\e$ such that $\cup_{i=1}^kB_i\supset f^{-1}(t)$. We claim that for $n$ big enough, we have an inclusion
\be 
\bigcup_{i=1}^kB_i\supset f_n^{-1}(\f^{-1}(t_n)).
\ee 
If not, there would be a sequence $p_n\in f_n^{-1}(\f^{-1}(t_n))$ such that $p_n\notin \cup_{i=1}^kB_i$; now, by the compactness of $M$, we can assume that $p_n\to p \notin \cup_{i=1}^kB_i$, but then we find a contradiction as follows. We have $f_n(p_n)\in K$ for all $n$, so that $f(p)\in K$, thus
 $f(p)\in\f^{-1}(t)$ because
$t_n=\f(f_n(p_n))\to \f(f(p))$. 
Hence $p\in f^{-1}(\f^{-1}(t))\subset \cup_{i=1}^kB_i$, which contradicts a previous statement.

It follows that $\#_\e(f^{-1}(\f^{-1}(t)))$ is an upper semicontinuous function:
\be 
\limsup_{(f_n,t_n)\to (f,t)}\#_\e(f_n^{-1}(\f^{-1}(t_n)))\le \#_\e(f^{-1}(\f^{-1}(t)))
\ee
and therefore it is measurable.
\end{proof}
\begin{cor}\label{cor:measurable}
Let $X\colon M\to N\supset W_t\subset  N_0\xrightarrow{\f}\R^m$ satisfy the hypotheses of Theorem \ref{thm:general}. Let $A\subset \mathcal{B}(M)$ be any Borel set. Then the function $u_A(f,t):= \#[f^{-1}(W_t)\cap A]$ is a measurable function on the completion of the measure space $\left(\mC^1(M,N)\times \R^m,\mathcal{B},[X]\otimes \mathcal{L}^m\right)$.
\end{cor}
\begin{proof}
Let us consider the Borel set $\mathcal{F}=\{(f,t)\colon f\transv W_t\}\subset \mC^1(M,N)\times \R^m$. Remark \ref{rem:toni} implies that $\mathcal{F}$ is a full measure subset of $\left(\mC^1(M,N)\times \R^m,\mathcal{B},[X]\otimes \mathcal{L}^m\right)$, therefore we can conclude by showing that the restriction $(u_A)|_{\mathcal{F}}$ is Borel for every $A\in\mathcal{B}(M)$. 

Let $K_i\subset N_0$, for $i\in \N$, be a sequence of increasing closed subsets (in $N$) whose union is equal to $N_0$ and define
\be 
u_A^{(i)}(f,t)=\#[f^{-1}(W_t\cap K_i)\cap A].
\ee
Then $u_A(f,t)=\sup_{i\in \N} u_A^{(i)}(f,t)$, and we already know that $u_A^{(i)}$ is Borel for any compact subset $A\subset M$ because of Lemma \ref{lem:fsigmeas}. Moreover, if $(f,t)\in\mathcal{F}$, we have that $f^{-1}(W_t\cap K_i)$ is a closed discrete subset of $M$, 
 because $f_t\transv W_t$ and $W_t\cap K_i=(\f|_{K_i})^{-1}(t)$ is closed. Therefore $(u_A^{(i)})\big|_{\mathcal{F}}$ is finite whenever $A$ is contained in a compact set.

Thanks to this observation, we can argue exactly as in the proof of Lemma \ref{lemma:meas}: write $M=\cup_j M_j$ as an increasing union of compact subsets; show that the family $\mathcal{D}_j$ of subsets $A\subset M_j$ such that $U_A^{(i)}|_{\mathcal{F}}$ is measurable forms a Dynkin class; conclude with the Monotone Class Theorem (see \cite{Erhan}) that $(u^{(i)}_A)|_{\mathcal{F}}$ is Borel for any $A\in\mathcal{B}(M)$.
\end{proof}
The following Lemma helps to rewrite the candidate formula \eqref{eq:formdel} for the Kac-Rice density into something that is more directly comparable with the Coarea formula (see Theorem \ref{thm:coarea}), which will be one of the main ingredient in the proof of Theorem \ref{thm:general}.
\begin{lemma}\label{lem:sette}
Let $X\colon M\to N\supset W_t\subset  N_0\xrightarrow{\f}\R^m$ satisfy the hypotheses of Theorem \ref{thm:general}. Then 
\be 
\delta_{X\in W_t}(p)=\int_{S_p\cap W_t}\E\left\{\frac{\delta_p(\f\circ X)}{J_x(\f|_{S_p})}\bigg| X(p)=x\right\}\rho_{X(p)}(x)d(S_p\cap W_t)(x).
\ee
\end{lemma}
\begin{proof}
 It is sufficient to prove the following identity (the Definitions of the object in play are in the Appendices \ref{app:densities} and \ref{app:angle}):
\be 
\frac{\delta_p(\f\circ X)}{J_x(\f|_{S})}= \delta_pX \frac{\sigma_x(W,d_pX)}{\sigma_x(W,S)},
\ee
where $X(p)=x$, $\f(x)=t$ and $W=\f^{-1}(t)$. Let us take an orthonormal basis $(\tau w \eta)$ of $T_xN$ such that the first vectors $\tau$ form a basis of $T_xW\cap T_xS$ and the $\eta$ form a basis of $T_xW^\perp$. The matrix of $d_x\f$ in such basis has the form $(00A)$ for some $m\times m$ invertible matrix $A$ and the space $T_xS$ (written in terms of the basis $(\tau w \eta)$) is spanned by the columns of a matrix of the form
\be 
\begin{pmatrix}
\mathbb{1} & 0\\
0 & B \\
0 & C
\end{pmatrix}.
\ee
Without loss of generality we can choose $\begin{pmatrix} B^T & C^T\end{pmatrix}^T$ to be an orthonormal frame. 
Notice that $\det C\neq 0$ because $T_xS\transv T_xW$ and that $\ker (d_u\f|_{S})$ is spanned by the first $m$ columns, since they correspond to the basis $\tau$, hence $J_x(\f|_{S})=|\det A C|$.

For any $u\in M=\D^m$ we have, by definition, that $\delta_{u}(\f\circ X)=J_u(\f\circ X)du$.
\be 
\begin{aligned}
\frac{J_p(\f\circ X)}{J_x(\f|_{S})}&=
\left|\frac{\det\left(d_p\f\frac{\de X}{\de u}\right)}{\det(A)\det(C)} \right|=\\
&=
\frac{|\det(A)| \vol\left(\Pi_{T_xW^\perp}\left(\frac{\de X}{\de u}\right)\right)}{|\det(A)|\vol\left(\Pi_{T_xW^\perp}\begin{pmatrix}
0\\ B \\ C
\end{pmatrix}\right)}  =\\
&= \frac{\sigma\left(T_xW, \text{span}\left(\frac{\de X}{\de u}\right)\right)\cdot \vol\left(\frac{\de X}{\de u}\right)}{\sigma\left(T_xW, T_xS\right)} .
\end{aligned}
\ee
\end{proof}
\subsection{Proof of the general formula: Theorem \ref{thm:general}}
\begin{proof}[Proof of Theorem \ref{thm:general}]
\ref{itm:general:1} Let $A\in\mathcal{B}(M)$. In what follows, let us keep in mind that, by Corollary \ref{cor:measurable}, the function $u_A(f,t)=\#(f^{-1}(W_t)\cap A)$ is measurable on the completion of the measure space $\left(\mC^1(M,N)\times \R^m,\mathcal{B},[X]\otimes \mathcal{L}^m\right)$. This, together with the hypothesis that $(p,t)\mapsto \delta_{X\in W_t}(p)$ is measurable, means that we don't need to worry about measurability issues.

Let $a\in B^+(\R^m)$ and apply the Area formula (see Theorem \ref{thm:area}) to the (deterministic) map $(\f\circ X)|_{X^{-1}(N_0)}\colon X^{-1}(N_0)\to \R^k$, with 
\be 
g(p)=a\left(\f\left(X(p)\right)\right)1_A(p)1_{N_0}\left(X(p)\right).
\ee
We obtain an identity, valid for all $X\in \mathcal{C}^1(M,N)$:
\be \label{eq:detform}
\int_{\R^m}a(t)\#[X^{-1}(W_t)\cap A] dt=\int_{A} \left( (a\circ\f\circ X)\cdot (1_{N_0}\circ X)\right) \delta(\f\circ X).
\ee
Using the Coarea formula \ref{thm:coarea}, we deduce a second identity, as follows:
\be\label{eq:appcoarea}
\E\left\{\left( (a\circ\f\circ X)\cdot (1_{N_0}\circ X)\right) \delta(\f\circ X)\right\}|_p= \ee
\be
\begin{aligned}
&= \E\left\{ a(\f(X(p))1_{N_0}(X(p))\delta_p(\f\circ X)\right\}=\\
&= \int_{N_0} a(\f(x))\E\{\delta_p(\f\circ X)|X(p)=x\}d[X(p)](x) =\\
&= \int_{S_p\cap N_0}a(\f(x))\E\{\delta_p(\f\circ X)|X(p)=x\}\rho_{X(p)}(x)dS_p =\\
&= \int_{\R^m}a(t)\int_{S_p\cap \f^{-1}(t)}\E\{\delta_p(\f\circ X)|X(p)=x\}\rho_{X(p)}(x)\frac{d\left((\f|_{S_p})^{-1}(t)\right)}{J_x(\f|_{S_p})}dt=\\
&=\int_{\R^m}a(t)\left(\int_{S_p\cap W_t}\E\left\{\frac{\delta_p(\f\circ X)}{J_x(\f|_{S_p})}\bigg| X(p)=x\right\}\rho_{X(p)}(x)d(S_p\cap W_t)(x)\right)dt\\
&=\int_{\R^m}a(t)\delta_{X\in W_t}(p) dt.
\end{aligned}
\ee
The Coarea formula was applied to the function $\f|_{S_p\cap N_0}\colon S_p\cap N_0\to \R^k$ in the fourth line.

Taking the expectation on both sides of \eqref{eq:detform} and repeatedly interchanging the order of integration via Tonelli's theorem\footnote{It is possible because the functions involved are positive and measurable.}, we obtain
\be 
\begin{aligned}
\int_{\R^m}a(t)\E\left\{\#[X^{-1}(W_t)\cap A]\right\}dt
&=
\E\left\{\int_{\R^m}a(t)\#[X^{-1}(W_t)\cap A] dt\right\}=\\
&=\E\left\{\int_{A} \left( (a\circ\f\circ X)\cdot (1_{N_0}\circ X)\right) \delta(\f\circ X)\right\}\\
&=\int_A \int_{\R^m}a(t)\delta_{X\in W_t}(p)dt\\
&= \int_{\R^m}a(t)\left(\int_A \delta_{X\in W_t}(p)\right)dt.
\end{aligned}
\ee
Given the arbitrariness of $a$ this proves \ref{itm:general:1}.


\ref{itm:general:3} Lemma \ref{lemma:meas} guarantees that in this case the set function $\mu(A)=\E\#[X^{-1}(W_0)\cap A]$ is a well defined Borel measure on $M$. 

Let $K\subset M$ be a compact set and let $O$ be its interior. If $X\transv W_0$, then $0$ is a regular value for the map $(\f\circ X)\colon X^{-1}(N_0)\to \R^{m}$, hence there is an $\e(X)>0$ such that, for any $0<\e<\e(X)$, the set $(\f\circ X)^{-1}(B_\e)\cap O$ is a disjoint union of balls $U_i$ with the property that $(\f\circ X)|_{U_i}\colon U_i\to B_\e$ is a $\ci$ diffeomorphism, therefore the number $\#[X^{-1}(W_t)\cap O]$ is constant for all $|t|<\e(X)$. 

Let $a_\e\in \mathcal{C}^\infty_c(\R^m)$ be a non negative function supported in the ball $B_\e$ of radius $\e>0$ and with $\int_{\R^m} a_\e =1$. In this case, formula \eqref{eq:detform} implies that
\be 
\#[X^{-1}(W_0)\cap O]= \lim_{\e\to 0^+}\int_{O} \left( (a_\e\circ\f\circ X)\cdot (1_{N_0}\circ X)\right) \delta(\f\circ X).
\ee

Taking expectation on both sides we have
\be
\begin{aligned}
\mu(O)
&=\E\left\{\lim_{\e\to 0^+}\int_{O} \left( (a_\e\circ\f\circ X)\cdot (1_{N_0}\circ X)\right) \delta(\f\circ X)\right\} \\
&\le \liminf_{\e\to 0}\E\left\{\int_{O} \left( (a_\e\circ\f\circ X)\cdot (1_{N_0}\circ X)\right) \delta(\f\circ X)\right\}\quad \textrm{(Fatou)}\\
&= \liminf_{\e\to 0}\int_O\int_{\R^m}a_\e(t)\delta_{X\in W_t} dt\quad\textrm{(Tonelli and \eqref{eq:appcoarea})}\\
&=\liminf_{\e\to 0}\int_{B_\e}a_\e(t)\left(\int_O\delta_{X\in W_t} \right)dt\quad\textrm{(Tonelli again)}\\
&\le \int_K \delta_{top}.
\end{aligned}
\ee


In particular, if $O\subset M$ belongs to the class $\mathcal{S}$ of relatively compact open subsets such that $\de O$ has measure zero, then $\nu(O)\le \int_O\delta_{top}$.
It follows that the measure $\mu$ is absolutely continuous with respect to the measure $\nu=\int \delta_{top}$. Indeed if $A\in \mathcal{B}(M)$ is such that $\nu(A)=0$, then
\be \label{eq:munu}
\mu(A)\le \inf_{B_i \in\mathcal{S},\ A\subset \cup_{i}B_i}\sum_i\mu(B_i)\le \inf_{B_i \in\mathcal{S},\ A\subset \cup_{i}B_i}\sum_i\int_{B_i}\delta_{top}=\nu(A)=0.
\ee
By the Radon-Nikodym theorem (see \ref{thm:RadoNiko}), this implies the existence of a measurable density $\delta\in L^1_{loc}(M)$ such that $\mu(A)=\int_A \delta$ for every $A\in\mathcal{B}(M)$. Moreover, since $\mu(A)\le \nu(A)$ by equation \eqref{eq:munu}, $\delta$ satisfies almost everywhere the inequality: $\delta\le \delta_{top}$.
\end{proof}



\section{Proof of the main Theorem}\label{sec:mainproof}
The goal of this section is to specialize the abstract result of Theorem \ref{thm:general} to the, more friendly, KROK situation. Let us consider a $\mC^1$ random map $X\colon M\to N$ and a submanifold $W\subset N$ such that the couple $(X,W)$ is KROK, that is, it satisfies the conditions i-ix of Definition \ref{def:krok}. 

To facilitate the next proofs, we will show that, without loss of generality we can make some further assumptions on $X$ and $W$.
\begin{figure}\begin{center}
\includegraphics[scale=0.15]{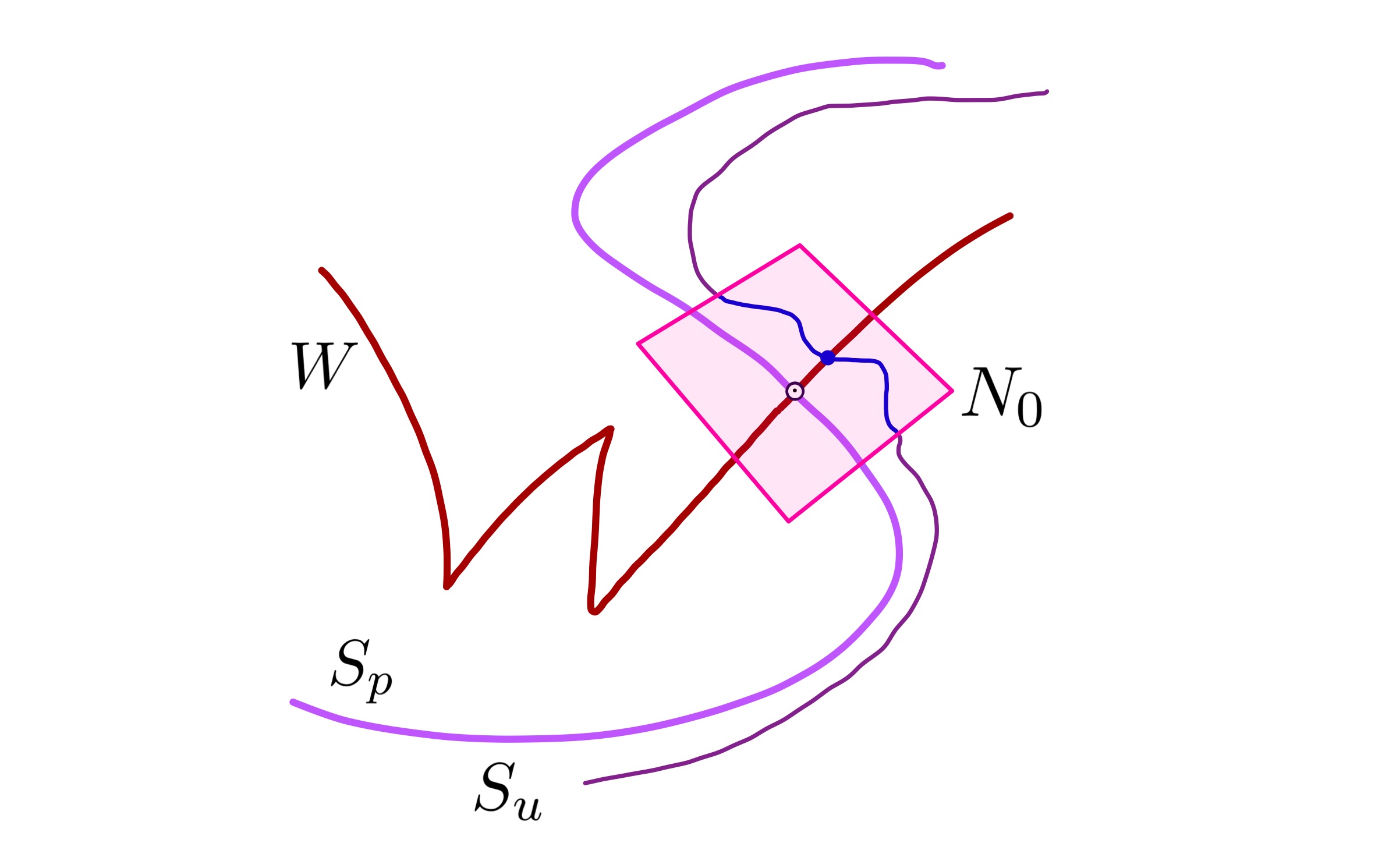}\caption{An illustration of Lemma \ref{lem:dieci}.}\label{fig:WSP}
\end{center}
\end{figure}

\begin{defi}\label{defi:krokcoords}
Let $(X,W)$ be a KROK couple. Let
\be 
    \phi\colon \overline{N_0}\xrightarrow{\sim} \D^m\times\D^{s-m}\times \D^{n-s}=\{ (t,y,z)\}
    \ee 
be a chart of $N$. We say that $\phi$ is a KROK chart at $(p,q)\in \MW $ if the following assumptions are satisfied.
\begin{enumerate}[(i)]
    \item $N_0\subset N$ is a relatively compact open subset such that $W\cap N_0= \phi^{-1}\{t=0\}$.
    \item  There is an open neighborhood $O_p\subset M$ of $p$ and a smooth map $g\in \coo\infty {O_p\times \D^m\times\D^{s-m}}{\D^{n-s}}$, such that
        \be 
        \phi(S_u\cap \overline{N_0})= \text{graph}(g(u,\cdot))=\{(t,y,g(u,t,y))\colon (t,y)\in \D^m\times\D^{s-m}\}.
        \ee
\end{enumerate}
In this case we say that the tuple $ (O_p,N_0,\phi, g)$ is a KROK model for $(X,W)$.
\end{defi}
\begin{lemma}\label{lem:dieci}
\HP. Then for all $(p,q)\in \MW $, there are open sets $p\in O_p\subset M$ and $q\in N_0\subset N$ and a KROK model $(O_p,N_0,\phi,g)$ for $(X,W)$.
\end{lemma}
\begin{proof}
Consider the point $q\in S_p\cap W$. Since $S_p\transv W$, there is a chart $\phi\colon \overline{N_0}\to\D^m\times\D^{s-m}\times \D^{n-s}=\{ (t,y,z)\}$ centered at $q$, such that $W\cap N_0=\{(0,y,z)\colon y,z\}$ and $S_p\cap N_0=\{(t,y,0)\colon t,y\}$. 
Now, let us consider the set $\mathcal{S}=\{(u,v)\in M\times N\colon u\in S_v\}$ (see Definition \ref{def:krok}) and observe that $T_{v}S_u=T_{(u,v)}\mS\cap T_v N$. Since 
$\mathcal{S}$ is a submanifold of $M\times N$ by Definition \ref{def:krok}, we deduce that for any $u,v$ close to $p,q$ in $M\times N$, the tangent spaces $T_vS_u$ and $T_qS_p$ are close to each other. This implies that for all $u$ in a neighborhood $O_p$ of $p$, de submanifold $S_u\cap N_0$ can be parametrized as the graph of a function $g_u(t,y)$ and that this function depends smoothly on $u$.
\end{proof}
\begin{lemma}\label{lem:undici}
\HPK. Let $(X^1,X^2,X^3):=\phi\circ X$ and define 
\be 
q\colon M\times \D^m\times \D^{s-m}\to N_0,
\ee
\be q(u,t,y):=\phi^{-1}(t,y,(g(u,t,y)))\in S_u\cap W_t\cap N_0.\ee Then the expression of $\delta_{X\in W_t}=\rho_{X\in W_t}du$ in the coordinates $t,y,z$ is the following.
\begin{multline}\label{eq:densieasy}
\rho_{X\in W_t}(u) =\\ 
=\int_{ \D^{s-m}}\E\left\{\left|\det\left(\frac{\de X^1(u)}{\de u}\right)\right|\Bigg| X(p)=q(u,t,y)\right\}\rho_{X(u)}(q(u,t,y))\sqrt{\det G(u,t,y)} dy,
\end{multline}
where $G(u,t,y)=\langle  \begin{pmatrix}\frac{\de q}{\de t} & \frac{\de q}{\de y}
\end{pmatrix}^T,\begin{pmatrix}\frac{\de q}{\de t} & \frac{\de q}{\de y}
\end{pmatrix} \rangle$ is the Gram matrix of $S_u$ with respect to the coordinates $t,y$.
\end{lemma}
\begin{proof}
It is convenient to write $\rho_{X\in W_t}$ with the formula of Lemma \ref{lem:sette}, where $\f\colon N_0\to \R^m$ is the function defined by $\f(\phi^{-1}(t,y,z))=t$. Thus, it is sufficient to show that 
\be 
\frac{d(S_u\cap W_t)}{J_v\left(\f|_{S_u}\right)}=\sqrt{\det G}\ dy.
\ee
Let us start by looking at the most tedious piece, namely the jacobian. The matrix of $d_u\f$ in the coordinates $t,y$ is $\begin{pmatrix}
\mathbb{1} & 0
\end{pmatrix}$, so that formula \eqref{eq:jacobidet} in Appendix, yelds
\be 
\begin{aligned}
J_u\left(\f|_{S_u}\right)&=\det\left(\begin{pmatrix}
\mathbb{1}& 0
\end{pmatrix}G^{-1}\begin{pmatrix}
\mathbb{1}\\ 0
\end{pmatrix}\right)^{\frac12}
=\det\left( S^{-1}\right)^{\frac12}= \left( \frac{\det G_{2,2}}{\det G}\right)^{\frac12},
\end{aligned}
\ee
where $S=G_{11}-G_{12}G_{22}^{-1}G_{12}$\footnote{$S$ is called the Schur complement of the block $G_{22}$ in $G=\begin{pmatrix}
G_{11} & G_{12} \\ G_{21} & G_{22}
\end{pmatrix}$.}. The last two equalities can be deduced from the identity: 
\be 
\begin{pmatrix}
G_{11} & G_{12} \\ G_{21} & G_{22}
\end{pmatrix}\cdot \begin{pmatrix}
\mathbb{1} & 0 \\ -G_{22}^{-1}G_{21} & \mathbb{1}
\end{pmatrix}=\begin{pmatrix}
S & G_{12} \\ 0 & G_{22}
\end{pmatrix}.
\ee
Now, since $\phi(N_0\cap S_u\cap W_t)=\{(t,y)\colon y\in \D^{s-m}\}$, the Gram matrix of $S_u\cap W_t$ in the coordinates $y$ is exactly $G_{22}$. Therefore, we conclude
\be 
\begin{aligned}
\frac{d(S_u\cap W_t)}{J_v\left(\f|_{S_u}\right)}
&=
\frac{\sqrt{\det G_{22}}dy}{\sqrt{ \frac{\det G_{2,2}}{\det G}}}
=\sqrt{\det G}\ dy.
\end{aligned}
\ee
\end{proof}
\begin{remark}\label{rem:independencemetrix}
Formula \eqref{eq:densieasy} can be rewritten in a form that doesn't involve the metrics on $M$ and $N$. Define a random element $X^1(u)\in \R^m\cup\{\arian\}$ such that $X^1(u)=\phi^1(X(u))$ if $X(u)\in N_0$ and $X^1(u)=\arian$ whenever $X(u)\notin N_0$. Then, if $A\subset \R^m$,
\be
\begin{aligned}
\P\{X^1(u)\in A\}&= \P\{X(u)\in N_0\cap \phi^{-1}(A\times \D^{s-m})\}\\
&=\int_{A}\left(\int_{\D^{s-m}}\rho_{X(u)}(q(u,t,y))\sqrt{\det G(u,t,y)} dy\right)dt,
\end{aligned}
\ee
hence the restriction of the measure $[X^1(u)]$ to $\R^m$ is absolutely continuous,
so that if we denote its density by $\rho_{X^1(u)}(t)dt$ we obtain an equivalent expression to that in formula \eqref{eq:densieasy}.
\be\label{eq:indepensity} 
\delta_{X\in W_t}(u)=\E\left\{\left|\det\left(\frac{\de X^1(u)}{\de u}\right)\right|\Bigg| \begin{aligned}&X(u)\in N_0 \\&X^1(u)=t\end{aligned}\right\}\rho_{X^1(u)}(t) du.
\ee
Notice that $\delta_{X\in W_t}=\delta_{X^1\in \{t\}}$. The above formula is completely independent from the Riemannian structures of $M$ and $N$.
\end{remark}
\subsection{A construction}\label{sec:constru}

The purpose of this section is to show that Theorems \ref{thm:main} and \ref{thm:maingraph} are actually equivalent. In fact, although the latter is evidently a more general result, it can be proved with a particularly simple application of the former. To understand this,
let us define a new $\mC^1$ random map $\XW\colon \MW \to \NW$, such that 
 \bega\label{eq:XW}
 \XW\colon \MW =\bigsqcup_{p\in M}\{p\}&\times (S_p\cap W) \longrightarrow \NW=\bigsqcup_{(p,q)\in \MW }\{(p,q)\}\times S_p,
 \\
 &(p,q) \mapsto (p,q,X(p)).
 \eega
 \begin{remark}
The fact that $\mathcal{S}$ is closed guarantees that $X(p)\in S_p$ for all $p\in M$, with probability one, hence this definition is well posed.  Indeed given a dense countable subset $D\subset M$, then 
\be 
\P\{X(p)\in S_p\ \forall p\in D\}=\P\left(\bigcap_{p\in D}\{X(p)\in S_p\}\right)=1
\ee 
and if the $X(p)\in S_p$ for all $p\in D$, then $(p,X(p))\in \overline{\mathcal{S}}=\mathcal{S}$ for all $p\in M$, by density and continuity.
\end{remark}
 Define also the diagonal submanifold $\DElta=\{(p,q,q)\colon (p,q)\in \MW \}\subset \NW$.
 Now, the random set $\Gam$ from Theorem \ref{thm:maingraph}, defined for a KROK couple $(X,W)$, can be interpreted as 
 \be 
\Gam=\left\{(p,q)\in M\times W\colon q=X(p)\right\}=\XW^{-1}(\DElta) ,
\ee
 therefore
 $\#_{\Gam}(V)=\#_{\XW\in \DElta}(V\cap \MW )$ for any $V\subset M\times N$.


Observe that if $\dim S_p=s$, then the dimensions of $\MW , \NW, \DElta$ are $s,2s,s$.

\begin{claim}
The couple $(\XW,\DElta)$ is KROK and
\be
d[\XW(p,q)](x,y,y')=\delta_p(x)\delta_q(y)\rho_{X(p)}(y')dZ_{p,q}\footnote{Here $\delta_p$ denotes the $\delta-$measure on the point $p$.}
\ee
where $Z_{p,q}=\{p,q\}\times S_p$. Moreover $ Z_p\cap \DElta=\{(p,q,q)\}$ is always finite.
\end{claim}
\begin{proof}
It is enough to observe that if $X\transv W$, then $\XW\transv \DElta$, and that if $S_p\transv W$, then $Z_{p,q}\transv \DElta$.
\end{proof}
It follows that the expression $V\mapsto \E\#_{\XW\in \DElta}(V)$  defines a Borel measure on $\MW $ and by putting $V=(A\times B)$, where $A\in \mathcal{B}(M)$ and $B\in\mathcal{B}(N)$, we have
\be 
\#_{X\in B\cap W}(A)=\#_{\XW\in \DElta}(A\times B\cap \MW ).
\ee
In other words, we can consider the measure $\E\#_{X\in W}(\cdot)$ as the section $\mu(W\times(\cdot))$ of a Borel measure $\mu=\E\#_{\XW\in \DElta}(\MW \cap (\cdot))$ on $M\times N$, with $\text{supp}(\mu)\subset \MW $\footnote{It is a strict inclusion, when $\rho_{X(p)}(q)= 0$ on a not negligible set of points $p,q$.}.

Moreover, the fact that $Z_{p,q}\cap \DElta$ is always equal to the point $(p,q,q)\in \NW$, permits to get rid of the integral in the formula for $\delta_{\XW\in \DElta}$:
\be
\delta_{\XW\in \DElta}(p,q)=\E\left\{\delta_{(p,q)}\XW\frac{\sigma_{(p,q,q)}(\XW,\DElta)}{\sigma_{(p,q,q)}(Z_{p,q},\DElta)}\bigg|X(p)=q\right\}\rho_{X(p)}(q).
\ee

\begin{lemma} Let $\XW\colon \MW  \to \NW$ be the map defined above. Then the density element $\delta_{\XW\in \DElta}(p,q)\in \Delta_{(p,q)}\MW$ can be written as follows.
\bega
\delta_{\XW\in \DElta}(p,q)&=
\E\left\{J_pX\frac{\sigma_x(X,W)}{\sigma_x(S_p,W)}\bigg|X(p)=x\right\}\rho_{X(p)}(x)d\MW (p,q)
\eega
where $\MW $ is endowed with the Riemannian metric induced by the isomorphism $T_{(p,q)}\MW \cong T_pM\oplus T_q(W\cap S_p)$, by declaring it to be an orthogonal splitting.
\end{lemma}
\begin{proof}
\HPK . Define $q(u,t,y)=\phi^{-1}(t,y,g(u,t,y))$. In the rest of the proof we will identify $O_p\cong \D^{m}$, so that in particular $p=0$ and $q=q(0,0,0)$. 
 Let $N_0^{(W)}$ be the open neighborhood of $(p,q)$ in $\NW$, defined as the image of the map
\be 
\D^s\times \D^s\ni(t^{(W)},z^{(W)})=\left((t',y'),(u,y)\right)\mapsto \left(u,q(u,0,y),q(u,t',y+y')\right).\footnote{After a rescaling of the coordinate $y$, we can assume that $\D^s\times\D^s\subset \phi(N_0)$. }
\ee
We call $\phi^{(W)}\colon \overline{N^{(W)}_0}\to\D^s\times \D^s$ the inverse of the above map, which provides a coordinate chart for $N_0^{(W)}$. 
Let us consider the following small open subset in $\MW $ 
\be O^{(W)}_{(p,q)}=\{(u,q(u,0,y))\colon (u,y)\in \D^{s}\}\cong \D^s\ee
and let us define the coordinate $u^{(W)}=(u,y)$ on it.

We can see now that $(O_{(p,q)}, N_0^{(W)},\phi^{(W)},g^{(W)})$ is a KROK model for $(\XW,\DElta)$, where $ g^{(W)}(u^{(W)},t^{(W)})=u^{(W)}$. Indeed $\DElta\cap N_0^{(W)}=\{t^{(W)}=0\}$ and $\phi^{(W)}(Z_{u^{(W)}}\cap N_0^{(W)})=\{(t^{(W)}, z^{(W)})\colon z^{(W)}=u^{(W)}\}$ is equal to the graph of the map $ g^{(W)}$.
It follows that $\delta_{\XW\in \DElta}=\rho_{\XW\in \DElta}du^{(W)}$ can be represented with the formula of Lemma \ref{lem:sette}, where $\f(p,q,q')=t^{(W)}$.
\bega 
\rho_{\XW\in \DElta}(p,q)
&=\E\left\{\left|\det\left(\frac{\de \XW^1}{\de u^{(W)}}\right)\right|\Bigg| \XW(p,q)=(p,q,q)\right\}\rho_{\XW(p,q)}(p,q,q)\sqrt{\det G^{(W)}}\\
&= 
\E\left\{\bigg|\det
\begin{pmatrix}
\frac{\de X^1}{\de u} & 0 \\ \frac{\de X^2}{\de u} & -\mathbb{1}
\end{pmatrix}
\bigg|\Bigg|X(p)=q\right\}\rho_{X(p)}(q) \vol\left(\frac{\de (\phi^{(W)})^{-1}}{\de t^{(W)}}\right)\\
&=
\E\left\{\left|\det\left(\frac{\de X^1}{\de u}\right)\right|\Bigg|X(p)=q\right\}\rho_{X(p)}(q)\vol\left(\frac{\de }{\de (t',y')}q(u,t',y+y')\right)\\
&=
\E\left\{\left|\det\left(\frac{\de X^1}{\de u}\right)\right|\Bigg|X(p)=q\right\}\rho_{X(p)}(q)\sqrt{\det G}.
\eega 
Here we used that, according to our definition, $\XW^1(u,y)=(X^1(u),X^2(u)-y)$.

In terms of the density,  recalling Lemma \ref{lem:undici} we just showed that
\bega 
\delta_{\XW\in \DElta}&=\E\left\{\left|\det\left(\frac{\de X^1}{\de u}\right)\right|\Bigg|X(p)=q\right\}\rho_{X(p)}(q)\sqrt{\det G}dydu\\
&=
\E\left\{\delta_pX\frac{\sigma_q(X,W)}{\sigma_q(S_p,W)}\bigg|X(p)=q\right\}\rho_{X(p)}(q)d(S_p\cap W)du.
\eega
and since, by definition, $d(S_p\cap W)du=d\MW $, we conclude.
\end{proof}
In other words, the relation between the absolutely continuous Borel measures with densities $\delta_{X\in W}=\rho_{X\in W}dM\in L^+(M)$ and $\delta_{\XW\in \DElta}=\rho_{\XW\in \DElta}d\MW \in L^+(\MW )$ mirrors the relation between the Borel measures $\E\#_{X\in W}$ on $M$ and $\E\#_{\XW\in \DElta}$ on $\MW$. Indeed, observing that $d\MW (p,q)=dM(p)\otimes d(S_p\cap W)(q)$, we have
\be 
\delta_{X\in W}(p)=\int_{W\cap S_p}\delta_{\XW\in \DElta}(p,q)dq.
\ee
Moreover, the hypotheses II and III, ensure that $\delta_{\XW\in \DElta}$ is a continous density on $\MW $.

\subsection{Proof of Theorems \ref{thm:main} and \ref{thm:maingraph}}
\begin{thm}
\HP. Then for any $A\in\mathcal{B}(M)$
\be 
\E\#_{X\in W}(A)=\int_{A}\delta_{X\in W}.
\ee
\end{thm}
\begin{proof}
Because of the construction of Section \ref{sec:constru}, we can assume that $S_p\cap W=\{q(p)\}$ is a point, for every $p\in M$.

 By Lemma \ref{lem:dieci} we can cover $M$ with a countable collection of embedded $m$-disks $D_{i}$, such that for each $i$, there is an open set $N_i\subset N$ containing $q(D_i)$ and such that there exists a KROK model $(D_i,N_i,\phi_i,g_i)$ for each $i$ (see Definition \ref{defi:krokcoords}). in particular $X(O_i)\subset N_i$. If we assume that the theorem holds for each $X_i=X|_{O_i}$, then for every $A\subset O_i$ we have:
\be 
\E\#_{X\in W}(A)=\E\#_{X_i\in W_i}(A)=\int_{A}\delta_{X_i\in W_i}=\int_{A}\delta_{X\in W}.
\ee 
The last equality is due to the fact that for every $p\in A$ the point $q(p)$ is already in $N_i$.
This implies that the two Borel measures $\E\#_{X\in W}$ and $\int \delta_{X\in W}$ coincide for every $A$ small enough ($A\subset O_i$ for some $i$), thus they are equal.

For this reason, we can assume that there is a global KROK model $(M,N,\phi,g)$. In this case the variable $y$ is not needed because $\dim(S_p)=m$ and we have, from Lemma \ref{lem:undici}, that $\delta_{X\in W_t}=\rho_t du$ with
\be 
\rho_t(u):=\E\left\{\left|\det\left(\frac{\de X^1(u)}{\de u}\right)\right|\Bigg|X(u)=q(u,t)\right\}\rho_{X(u)}(q(u,t))\sqrt{\det G(u,t)},
\ee
where $q(u,t)=\phi^{-1}(t,g(u,t))$ and $q(u,0)=q(u)$.
The KROK assumptions ensure that the function $(u,t)\mapsto \rho_{t}(u)$ is continuous at $M\times 0$, thus $\delta_{X\in W_t}\to \delta_{X\in W_0}$ in $L^1_{loc}$, so
 that from point \ref{itm:general:3} of Theorem \ref{thm:general}, applied with $\f(t,z)=t$, it follows that there exists a measurable function $g\colon M\to [0,1]$, such that $\E\#_{X\in W} =\int g du$ and $g\le \rho_0$.
To end the proof it is sufficient to show that $g(p)\ge \rho_0(p)$ for almost every $p\in M$. 

Let us consider the subset $P\subset M$ (recall that we are assuming $M=\D^{m}$) of all Lebesgue points for $g$, so that $p\in P$ if and only if
\be 
g(p)=\lim_{r\to 0}\frac{\E\#_{X\in W}(B_r)}{\mathscr{L}^m(B_r)}.
\ee
 We will prove that if $p\in P$, then $g(p)\ge \rho_0$. After that, the proof will be concluded since, by Lebesgue Differentiation theorem, $P$ is a full measure set in $M$.

Let $B_r$ be a closed ball of radius $r>0$ in $M$ centered at $p\in P$. 
Let $f^1=\f\circ f$ and consider the set
\be 
\mathscr{A}_s=\left\{f\in \mathcal{C}^1(M,N)\colon \begin{aligned}
& f(B_s)\subset N_0, \\
&f^1|_{B_s}\colon B_s\hookrightarrow \R^m\text{ is an embedding}
\end{aligned}
\right\}.
\ee
It is straightforward to see that $\mathscr{A}_s$ is an open set in $\mathcal{C}^1(M,N)$. Define a family of continuous functions $\a_{s,\eta}\colon \coo1MN\to [0,1]$ such that $a_{s,\eta}\nearrow\mathbb{1}_{\mathscr{A}_s}$, when $\eta\to +\infty$. 

Let $r\le s$. Now, the random variable $\#_{X\in W_t}(B_r)\a_{s,\eta}(X)$ is bounded by $1$ and converges almost surely as $t\to 0$, because $\E\#_{X\in W}(\de B_r)=0$.
As a consequence, by dominated convergence, we have that
\be
\E\left\{\#_{X\in W_0}(B_r)\a_{s,\eta}(X)\right\}=\lim_{t\to 0}\E\left\{\#_{X\in W_t}(B_r)\a_{s,\eta}(X)\right\}.
\ee
Moreover, arguing as in the proof of point \ref{itm:general:1} of Theorem \ref{thm:general}, we have the following equality for almost every $t\in \R^m$ and every $\eta$:
\be 
\E\left\{\#_{X\in W_t}(B_r)\a_{s,\eta}(X)\right\}=\int_{B_r}\rho_{X\in W_t}^{s,\eta}(u)du,
\ee
where 
\be 
\rho_{X\in W_t}^{s,\eta}(u):=\E\left\{\left|\det\left(\frac{\de X^1(u)}{\de u}\right)\right|\a_{s,\eta}(X)\Bigg|X(u)=q(u,t)\right\}\rho_{X(u)}(q(u,t))\sqrt{\det G(u,t)}.
\ee
The important point here is that $\rho_{X\in W_t}^{s,\eta}(u)$ is continuous at $t=0$\footnote{This is why we defined the continuous functions $\a_{s,\eta}$, instead of simply using the characteristic functions $\mathbb{1}_{{\mathcal{A}_s}}$.} because of the KROK. \ref{itm:krokEcont}, so that we are allowed to do the last step in the following sequence of inequalities. 
\be 
\begin{aligned}
\E\#_{X\in W}(B_r)&\ge 
\E\{\#_{X\in W_0}(B_r)\a_{r,\eta}(X)\} 
\\
&=\lim_{t\to 0}\int_{B_r}\E\{\#_{X\in W_t}\a_{r,\eta}(X)\}\\
&=\lim_{t\to 0}\int_{B_r}\rho_{X\in W_t}^{r,\eta}(u)du\\
&=\int_{B_r}\rho_{X\in W_0}^{r,\eta}(u)du.
\end{aligned}
\ee
Now, because $p\in P$ and $\rho_{X\in W_0}^{s,\eta}$ is continuous, we have that
\be
\begin{aligned}
g(p)&=\lim_{r\to 0}\frac{\E\#_{X\in W}(B_r)}{\mathscr{L}^m(B_r)}
\ge
\lim_{r\to 0}\frac{\int_{B_r}\rho_{X\in W_0}^{s,\eta}(u)du}{\mathscr{L}^m(B_r)}
=\rho_{X\in W_0}^{s,\eta}(p).
\end{aligned}
\ee
Taking the supremum over all $\eta$ we obtain that
\be\label{eq:daje}
g(p)
\ge  \E\left\{\left|\det\left(\frac{\de X^1}{\de u}(p)\right)\right|\mathbb{1}_{\mathscr{A}_s}(X)\Bigg|X(p)=q\right\}\rho_{X(p)}(q)\sqrt{\det G(p,0)},
\ee
where $q=q(p,0)\in W$.
Observe that $\mathscr{A}_s$ is an increasing family (as $s\to 0$) whose union is equal to the set $\mathscr{A}_0$ consisting of all functions $f\in \coo1MN$, such that $f(p)\in N_0$ and such that $d_pf^1$ is invertible. 
Therefore, recalling that $X\transv W$ almost surely, we see that taking the supremum with respect to $s$ in \eqref{eq:daje} we
obtain the thesis, valid for all $p\in P$:
\be 
g(p)\ge \rho_0(p).
\ee 
\end{proof}


%
%
%
%

\subsection{The case of fiber bundles: Proof of Theorem \ref{thm:megafica}}

\begin{proof}[Proof of Theorem \ref{thm:megafica}]
(We refer to Appendix \ref{app:angle} for the notations with frames.)
Let $q\in W$ and $p=\pi(q)$.
Let us take an orthonormal basis $\de_u$ of $T_pM$ and an orthonormal basis $\de_y$ of $T_qS_p\cap T_qW$. Let us complete the latter to an orthonormal basis of $T_qS_p$ by adding an orthonormal frame $\de_t$. Let us denote by $h(de_u)\in T_q S_p^\perp$ the frame such that $d_q\pi(h(\de_u))=\de_u$.
Finally, observe that $T_qW\cap (T_qS_p\cap T_qW)^\perp$ is contained in the space generated by $(h(\de_u),\de_t)$ and it projects surjectively on $T_pM$, hence it has a (not orthonormal) basis of the form $\de_z=h(\de_u)+A\de_t$, for some $m\times m$ matrix $A$.
 (the letters are coherent with the KROK coordinates, see Definition \ref{defi:krokcoords}).
 
 Now, by construction and by Definition \ref{defi:angle}, we have
 \bega
 \frac{d(S_p\cap W)(q) dM(p)}{\sigma(S_p,W)}&=\frac{\vol(\de_z)}{\vol\begin{pmatrix}
 \de_t & \de_z 
 \end{pmatrix}} dy du
 \\
 &= \frac{1}{\vol{\begin{pmatrix}
 \mathbb{1}_m & A \\
 0 & \mathbb{1}_m
 \end{pmatrix}}}\left(\vol\begin{pmatrix}
 \de_y & \de_z
 \end{pmatrix}dydu\right)
\\
&=dW(p).
\eega
\end{proof}
\subsection{Other counting measures: Proof of Theorem \ref{thm:mainalph}}
\begin{proof}[Proof of Theorem \ref{thm:mainalph}]
It is sufficient to prove the case in which $\a$ is continuous, bounded and positive. Indeed, then the result can be extended by monotone convergence to any positive Borel function and finally to any Borel function $\a$ by presenting it as $\a=\a^+-\a^-$.
If $\a$ is continuous, then the hypothesis KROK.\ref{itm:krokEcont} ensures that one can repeat the whole proof of Theorem \ref{thm:maindens} with $\#^\a_{X\in W}$ and $\delta^\a_{X\in W}$. The only thing to check is the very first step, which is provided by the Area and Coarea formulas at equations \eqref{eq:detform} and \eqref{eq:appcoarea} in the proof of the general formula of Theorem \ref{thm:general}. 
For the weighted case, we have to apply the Area formula to the function $\f\circ X$, with
\be 
g(p)=\a(X,p) a\left(\f\left(X(p)\right)\right)1_A(p)1_{N_0}\left(X(p)\right),
\ee
to get a generalization of identity \eqref{eq:detform}:
\bega
\int_{\R^m}a(t)\#^\a_{X\in W_t}(A)dt &=\int_{\R^m}a(t)\left(\sum_{p\in A\cap X^{-1}(W_t)}\a(X,p) \right)dt
\\
&=\int_{A} \a(X,p) a(\f\circ X(p)) 1_{N_0}( X(p)) \delta(\f\circ X)(p).
\eega
On the other hand, via the Coarea formula, the identity \eqref{eq:appcoarea} becomes.
\bega
&\ \E\left\{\a(X,p) a(\f(X(p))1_{N_0}(X(p))\delta_p(\f\circ X)\right\}=
\\
&= \int_{S_p\cap N_0}a(\f(x))\E\{\a(X,p)\delta_p(\f\circ X)|X(p)=x\}\rho_{X(p)}(x)dS_p 
\\
&=\int_{\R^m}a(t)\left(\int_{S_p\cap W_t}\E\left\{\a(X,p)\frac{\delta_p(\f\circ X)}{J_x(\f|_{S_p})}\bigg| X(p)=x\right\}\rho_{X(p)}(x)d(S_p\cap W_t)(x)\right)dt
\\
&=\int_{\R^m}a(t)\delta^\a_{X\in W_t}(p) dt.
\eega
\end{proof}
\section{The Gaussian case: Proofs}
\label{sec:proofgauss}

\subsection{Proof of Theorem \ref{thm:maingau} and Corollaries \ref{cor:maingauW} and \ref{cor:connectedform}}
\begin{proof}[Proof of Theorem \ref{thm:maingau}]
It suffices to show that the couple $(X,W)$ is KROK and use Theorem \ref{thm:main}. The only conditions of Definition \ref{def:krok} that have to be checked are KROK.\ref{itm:kroktransvas} and KROK.\ref{itm:krokEcont}.
 The former is a consequence of the Probabilistic Transversality theorem \ref{thm:transKR} and the latter follows from next Lemma.
\end{proof}

\begin{lemma}\label{lem:gioiellino}
Let $X=(X^0,X^1)\colon \R^m\to\R^k\times\R^h$ be a $\mC^1$ Gaussian random field and assume that $X^0$ is non-degenerate Let $\a\colon \mC^1(\R^m,\R^{k+h})\times\R^m\to \R$ be a continuous function such that
\be\label{eq:subpolynomial}
|\a\left(f,p\right)|\le N\left(\left|f(p)\right|+\left|\de_1 f(p)\right|+\dots +\left|\de_m f(p)\right|\right)^N+N,
\ee
for some constant $N>0$. Then
 the next function is continuous
 \be
 \R^m\times\R^k\ni (p,q)\mapsto \E\left\{\a(X,p)\big|X^0(p)=q\right\}\in \R.
 \ee
\end{lemma}
 \begin{proof}
Let us fix two points $u,p\in\R^m$. By using a standard argument we can find a Gaussian random vector $Y(u,p)\in \R^{k+h}$ that is independent from $X^0(p)$ and such that 
\be\label{eq:gausdegress}
X(u)=A(u,p)X^0(p)+Y(u,p)
\ee
For some $(k+h)\times k$ matrix $A(u,p)$.
Then, defining $K^{*,0}_X(u,p)=\E\{X(u)X^0(p)^T\}$, we have 
\be\label{eq:covagausdegress}
K^{*,0}_X(u,p)\left(K_X^{0,0}(p,p)\right)^{-1}=A(u,p),
\ee
from which we deduce that $A$ is a $\mC^1$ function of $u,p$. 

Consider now the $\mC^1$ Gaussian random field $Y_p=Y(\cdot,p)\colon \R^m\to \R^{k+h}$ defined, for each $p\in \R^m$, by the identity \eqref{eq:gausdegress} above. 
The previous computation ensures that the random vectors $Y_p(u)$ and $X^0(p)$ are independent for every $u\in\R^m$ and therefore $Y_p(\cdot)$ is independent from the vector $X^0(p)$ as a field. It follows that
\be\label{eq:expgausdegress}
\E\{\a(X,p)|X^0(p)=q\}=\E\{\a(Y_p(\cdot)+A(\cdot,p)q,p)\}=\E\{\a(Z_{p,q},p)\},
\ee
for any $q\in\R^k$, where $Z_{p,q}$ is the $\mC^1$ (noncentered) Gaussian random field \bega
Z_{p,q}(\cdot)&=Y_p(\cdot)+A(\cdot,p)q\colon \R^m\to \R^{k+h}\\
&=
X(\cdot)+A(\cdot,p)\left(q-X^0(p)\right)
\eega

To end the proof, let us show that the right hand side of equation \eqref{eq:expgausdegress} depends continuously on $p,q$. Given two converging sequences: $p_n\to p$ in $\R^m$ and $q_n\to q$ in $\R^k$, it is clear that the random map $Z_{p_n,q_n}$ converges almost surely to $Z_{p,q}$ in the space of $\mC^1$ functions, thus
\be\label{eq:fatou1lastlemma}
\E\{\a(Z_{p,q},p)\}\le \liminf_{n\to+\infty} \E\{\a(Z_{p_n,q_n},p_n)\}
\ee
by Fatou's lemma. Now observe that, since $Z_{p,q}(p),\de_1 Z_{p,q}(p),\dots, \de_mZ_{p,q}(p)$ are jointly Gaussian and their covariance matrices depend continuously on $p,q$ because of equation \eqref{eq:covagausdegress}, then the function $\f\colon\R^m\to \R$:
\be
\f(p,q)=\E\{N\left(|Z_{p,q}(p)|+|\de_1Z_{p,q}(p)|+\dots +|\de_m Z_{p,q}(p)|\right)^N+N\}
\ee
is continuous. The hypothesis \eqref{eq:subpolynomial} allows us to apply Fatou's Lemma once again to obtain that $\limsup_{n\to+\infty}\E\{\a(Z_{p_n,q_n},p_n)\}\le \E\{\a(Z_{p,q},p)\}$:
\begin{multline}
\f(p,q)-\E\{\a(Z_{p,q},p)\}
 \\
=\E\left\{\left(N\left(|Z_{p,q}(p)+|\de_1Z_{p,q}(p)|+\dots +|\de_m Z_{p,q}(p)|\right)^N+N\right)-\a(Z_{p,q},p)\right\}\\
\le
\liminf_{n\to+\infty}\left(\f(p_n,q_n)-\E\{\a(Z_{p_n,q_n},p_n)\}\right)\\
=
\f(p,q)-\limsup_{n\to+\infty}\E\{\a(Z_{p_n,q_n},p_n)\}.
\end{multline}
This, together with \eqref{eq:fatou1lastlemma}, implies that $\E\{\a(Z_{p_n,q_n},p_n)\}\to\E\{\a(Z_{p,q},p)\}$.
\end{proof}
\begin{proof}[Proof of Corollary \ref{cor:maingauW}]
Follows directly from the version of the formula for fiber bundles.
The only novelty to prove here is that the density $\delta_{\Gam}\in \mathscr{D}^0(W)$ is smooth in this case. This follows from Lemma \ref{lem:gagagaga} below.
\end{proof}
\begin{proof}[Proof of Corollary \ref{cor:connectedform}]
Since $W\subset E$ is parallel for the structure, we can express the formula of Corollary \ref{cor:maingauW}, in the form described in Remark \ref{rem:connectionForm}, so that $\E\#\left(W \cap V\right)=\int_{W\cap V}\delta_{\Gam}$ where
\be 
\delta_{\Gam}=
\E\left\{|\det\left(\Pi_{T_xW_p^\perp}\circ\left(\nabla^X X\right)_p\right) \Big| X(p)=x\right\}\frac{e^{-\frac{|x|^2}{2}}}{(2\pi)^{\frac s2}}dW(x).
\ee
Here, we already used the fact that the metric structure $g$ is defined by $X$, meaning that $K_p^{-1}\langle x,x\rangle=g(x,x)$ and $\det K_p=1$. Furthermore, from the definition of $\nabla^X$ given by the identity \eqref{eq:nablaX} (compare with \cite{Nicolaescu2016}) we observe that, in a trivialization chart, for any $i,j$ and $v\in T_pM$, we have
\be 
\E\{(\nabla_v^X X)^i_p X^j(p)\}=\E\left\{(D_vX)^i_p X^j(p)-\E\{(D_v X)^i_p |X(p)\}X^j(p)\right\}=0,
\ee
meaning that $(\nabla^X X)_p$ and $X(p)$ are uncorrelated and thus, being Gaussian, independent. Therefore the conditioning can be removed:
\be 
\E\left\{|\det\left(\Pi_{T_xW_p^\perp}\circ\left(\nabla^X X\right)_p\right) \Big| X(p)=x\right\}=\E\left\{|\det\left(\Pi_{T_xW_p^\perp}\circ\left(\nabla^X X\right)_p\right) \right\}
\ee
\end{proof}
\subsection{Proof that the notion of sub-Gaussian concentration is well defined}
\begin{lemma}\label{lem:subGwell}
Let $\pi\colon E=\D^m\times\R^s\to \D^m$ and let $W\subset E$ be a smooth submanifold. If $W\subset E$ has sub-Gaussian concentration for some linearly connected Riemannian metric on $\pi$, then the same holds for any other.
\end{lemma}
\begin{proof}
In this case $E=\D^m\times\R^s\to\D^m$,  $E_u=\{0\}\times \R^s$, a connected Riemannian metric is any Riemannian metric whose matrix is of the form
\be\label{eq:lcRmetric}
g_E(u,x)=
\begin{pmatrix}
g_M(u) & 0
\\ 0 & g_v(u)
\end{pmatrix}+
\begin{pmatrix}
\Gamma(u,x)^Tg_v(u)\Gamma(u,x) & \Gamma(u,x)^Tg_v(u)
\\ g_v(u)\Gamma(u,x) & 0
\end{pmatrix},
\ee
where $g_M(u)$ is the Riemannian metric on $\D^m$, while $\Gamma\colon E\to \R^{m\times s}$ is some smooth function. Moreover, in this case, we have that the horizontal space $H_{(u,x)}=(T_{(u,x)}E_u)^\perp$ is the graph of the linear function $\R^m\ni\dot{u}\mapsto -\Gamma(u,x)\dot{u}\in\R^s$.

The connection is linear precisely when $\Gamma(u,x)$ depends linearly on $x$, so that we can define the Christoffel symbols\footnote{They are the Christoffel symbols of the corresponding covariant derivative: given a smooth section $s(u)=(u,s_v(u))$,  the vertical projection of $d_us$ is given by
\be (\nabla s)_u=d_us_v+\Gamma(u,s_v(u))\in\R^{m\times s}.\ee } $\Gamma_{k,i}^j\in \mC^\infty(\D^m)$ as follows
\be 
\Gamma_{k}^j(u,x)=\sum_{i=1}^s\Gamma_{k,i}^j(u)x^i, \quad \forall k=1,\dots, m \text{ and } j=1,\dots, s.
\ee
Let $g_0$ and $g_1$ be any two linearly connected metrics on $E$ and assume that $W$ has sub-Gaussian concentration for $g_0$. The corresponding vertical norms on $E_u$
\be 
\|x\|_i=\sqrt{g_i\langle\begin{pmatrix}
0 \\ x
\end{pmatrix},\begin{pmatrix}
0 \\ x
\end{pmatrix}\rangle}
\ee
are equivalent and since $u\in \D^m$ there exists a uniform constant $N>0$ such that
\be 
\frac 1 {N}\|\cdot\|_0\le \|\cdot\|_1\le N\|\cdot\|_0.
\ee
Moreover, given the linearity of $\Gamma$  in equation \eqref{eq:lcRmetric}  for both $g_0$ and $g_1$, we can assume that
\be 
dW_1(u,x)\le N(1+R^2)dW_0(u,x),
\ee
for all $x\in B_R^1\cap W$ and every $R>0$. It follows that
\bega 
\vol_{W_1}\left(W\cap B^1_R\right)&=\int_{W\cap B^1_R}dW_1 
\\
&\le \int_{W\cap B_R^1}N(1+R^2)dW_0
\\
&\le \int_{W\cap B_{NR}^0}N(1+R^2)dW_0
\\
&\le N(1+R^2)C_0(\delta)\exp(\delta N^2R^2).
\eega
the last inequality is true for any $\delta>0$ because $W$ has sub-Gaussian concentration with respect to the metric $g_0$. Now, for every fixed $\e>0$,  there exists $C_1(\e)>0$ such that
\be 
N C_0\left(\frac12 \e N^{-2}\right)(1+R^2)\le C_1(\e)\exp\left(\frac12 \e R^2\right),
\ee
simply because $1+R^2=o(e^{\frac12 \e R^2})$. We obtain that $W$ has sub-Gaussian concentration for the metric $g_1$ as well, since for every fixed $\e>0$, setting $\delta=\frac12 \e N^{-2}$, we get
\be 
\vol_{W_1}\left(W\cap B^1_R\right)\le C_1(\e)\exp\left(\frac12 \e R^2\right)\exp(\delta N^2R^2)=C_1(\e)\exp\left( \e R^2\right).
\ee
\end{proof}

\subsection{Proof of theorem \ref{thm:mainEgau}}
Since the statement is local we can assume that $E=\D^m\times\R^s$ and $X\colon \D^m\to\D^m\times \R^m$ is a smooth Gaussian random section of
the trivial bundle, so that $X(u)=(u,V(u))$ for some smooth Gaussian random function $V\colon \D^m\to\R^s$. Let us define the smooth Gaussian random field 
\be 
j^1X\colon \D^m\to \R^{s}\times \R^{m\times s}, \quad u\mapsto j^1_uX=(V(u),d_uV),
\ee
Since $X$ is smooth and nondegenerate, the covariance matrix of $j^1_uV$ 
\bega
K_{j^1X}\colon \D^m \to \mathcal{U}, \quad
K_{j^1_uX}= \begin{pmatrix}
K_0(u) & K_{01}(u)\\ K_{10}(u) & K_1(u)
\end{pmatrix}
\eega
is a smooth function taking values in the open set 
\be 
\mathcal{U}=\left\{K=\begin{pmatrix}
K_0 & K_{01}\\ K_{10} & K_1
\end{pmatrix}\in \R^{(m+ms)\times (m+ms)}\colon \begin{matrix}
\text{$K$ is symmetric, } \\
\text{semipositive}
\\
\text{and $\det(K_0)>0$}
\end{matrix}\right\}.
\ee 
For the proof of the first statement of the theorem, it is sufficient to prove the following Lemma.
\begin{lemma}\label{lem:gagagaga}
Let $E=\D^m\times \R^s$, be endowed with the standard Euclidean metric. Let $E_u:=\{u\}\times\R^s$, 
and let $W\subset E$ be a smooth submanifold of codimension $m$ such that $W\transv E_u$ for every $u\in\D$.\footnote{Notice that $T_{(u,x)}W$ is a well defined $m$ dimensional subspace of $\R^m\times\R^s$ even if $u\in\de\D^m$. In such case the transversality $E_u\transv T_{(u,x)}W$ is still meant in the space $\R^m\times\R^s$.}
There exists a continuous function
\be 
f_W\colon W\times \mathcal{U}\to \R,
\ee
such that for any nondegenerate smooth Gaussian random section $X\colon\D^m\to\D^m\times\R^s$, with $X(u)=(u,V(u))$, we have
\be 
\delta_{\Gam}(u,x)=f_W(u,x,K_{j^1_uX})dW(u,x).
\ee
Moreover, there exists a smooth function $N\in\mathcal{U}\to\R_+$ such that, $\forall (u,x,K)\in W\times\mathcal{U}$,
\be 
f_W(u,x,K)\le N(K)\exp\left(-\frac{1}{N(K)}\|x\|^2\right).
\ee
\end{lemma}
\begin{proof}
By Equation \eqref{eq:gaudelgamma}, we have
\be 
\delta_\Gam=\E\left\{\a(u,x,d_uV)|V(u)=x\right\} \rho_{X(u)}(x)dW(u,x),
\ee
where $\rho_{X(u)}$ is the density of the Gaussian random vector $X(u)$, while where $\a\colon W\times \R^{m\times s}\to \R_+$ is the function defined as follows.
\be \label{eq:alpha}
\a(u,x,J):=\vol\begin{pmatrix}
\mathbb{1}_m \\ J
\end{pmatrix}\sigma\left(\text{Im}\begin{pmatrix}
\mathbb{1}_m \\ J
\end{pmatrix},T_{(u,x)}W\right)=\vol(\Pi_{T_{(u,x)}W^\perp}\begin{pmatrix}
\mathbb{1}_m \\ J
\end{pmatrix}),
\ee
where $\vol$ and $\sigma$ are taken with respect to the Euclidean Riemannian structure on $E=\D^m\times\R^s$ (in the last equality we used Proposition \ref{prop:angleperp} from Appendix \ref{app:angle}). 

From formula \eqref{eq:dgauvec} it is clear that $\rho_{X(u)}(x)=\rho(u,x,K_{j^1_uX})$ with $\rho\colon E\times \mathcal{U}\to\R_+$ being the smooth function
\be \label{eq:rhocazzo}
\rho (x,K)=
\frac{
\exp
\left(
-\frac{1}{2}x^TK_0^{-1}x
\right)
}{
\pi^\frac s2\sqrt{\det\left(K_0\right)}
}
\le N_\rho(K)\exp(-\frac{\|x\|^2}{N_\rho(K)}).
\ee
Clearly there exists a continuous function $N_\rho\colon \mathcal{U}\to\R_+$ such that the inequality \eqref{eq:rhocazzo} holds.

Let us consider the function $\E\{\a(u,x,d_uV)|V(u)=x\}$. We will show that it is a continuous function of $u,x$ and $K_{j^1X}$. 

First, observe that $\a$ is smooth because $W$ is a smooth submanifold. Like we did in the proof of lemma \ref{lem:gioiellino}, let us define a new Gaussian random vector $J(u)$ such that $ 
d_uV=J(u)+A(u)V(u)
$, for some smooth function $A(u)\in\R^{m\times ms}$, linear in $x\in \R^s$, and such that $J(u)$ and $V(u)$ are independent. Then $K_{J(u)}=K_J(K_{j^1_uX})$ and $A(u)=A(K_{j^1_uX})$, where $K_J,A\in\mC^0(\mathcal{U})$ are smooth functions:
\bega\label{eq:dependance}
K_{J}(K)&=K_1-K_{10}K_0^{-1}
K_{01};
\\
A(K)&=K_{10}K_0^{-1}.
\eega
Consider the expression
$E(u,x,K_J,A)=\E\{\a(u,x,J-Ax)\}$. Since $\a$ is of the form \eqref{eq:alpha}, all of its derivatives are bounded by a polynomial in $J$, when $J\to+\infty$ with $u,x$ fixed. From this it follows that
$\E\{\a(u,x,J-Ax)\}$
is finite and depends smoothly on $u,x,A,K_J$, where $K_J$ is the covariance matrix of $J$. Therefore, the function $F(u,x,K)=E\{\a(u,x,K_J(K),A(K)\}$ is smooth in $u,x,K$ and we have
\be 
\E\{\a(u,x,d_uV)|V(u)=x\}=F\left(u,x,K_{j^1_uX}\right).
\ee

Moreover, let $J(K)\in\R^{m\times s}\cong \R^{ms}$ be the Gaussian random matrix corresponding to the Gaussian random vector having covariance matrix $K_J(K)$. Let $\|\cdot\|$ be any norm on $\R^{m\times s}$. Then
\bega
F(u,x,K)&=\E\{\a(u,x,J(K)-A(K)x)\}
\\
&\le \E\left\{\vol\begin{pmatrix}
\mathbb{1}_m \\ J(K)+A(K)x
\end{pmatrix}\right\}
\\
&\le \left(\E\left\{\left\|\vol\begin{pmatrix}
\mathbb{1}_m \\ J(K)
\end{pmatrix}\right\|^m\right\}^\frac1m
+ \| A(K)x\|
\right)^m
\\
&\le N_\a(K)\left(1+\|x\|\right)^m.
\eega
for some big enough continuous function $N_\a\colon \mathcal{U}\to \R_+$. This is because the $m^{th}$ moment of a Gaussian random vector depends continuously on its covariance matrix. 
Combining this with \eqref{eq:rhocazzo}, we see that the function 
\be 
f_W(u,x,K):=F(u,x,K)\rho(x,K).
\ee 
has all the properties that we wanted to show.
\end{proof}

Let us show the continuity statement for a sequence of GRSs $X_d$, with limit $X_\infty$.
 By Theorem \ref{thm:maingau} and Lemma \ref{lem:gagagaga} we have
\be 
\E\#_{X_d\in W}(A)=\int_{\pi^{-1}(A)\cap W}\delta_{\Gamma(X_d,W)}dW
=
\int_{\pi^{-1}(A)\cap W} f_W(u,x,K_d)dW(u,x),
\ee
where $K_d=K_{j^1X_d}$ converges to $K_\infty=K_{j^1X_{\infty}}$ because they depend linearly on the $2^{nd}$ jet of $K_{X_d}$ and $K_{X_\infty}$.
From Lemma \ref{lem:gagagaga} we deduce that $f_W(u,x,K_d)\to f_W(u,x,K_{\infty})$, therefore we could conclude by the dominated convergence theorem if we show that $f_W(u,x,K_d)$  is uniformly bounded by an integrable (on $W$) function. Since $K_d\in\mathcal{U}$ is a convergent sequence, we have that $N(K_d)$ is uniformly bounded by some constant $N\in\N$, thus 
\bega
f_W(u,x,K_d)&\le Ne^{-\frac{|x|^2}{N}}.
\eega
Now, the fact that $W$ has sub-Gaussian concentration implies that the latter function is integrable on $W$:
\bega
\int_W Ne^{-\frac{|x|^2}{N}}dW(u,x) &\le \sum_{R\in\N}\int_{W\cap B_R} Ne^{-\frac{(R-1)^2}{N}}dW(u,x)
\\
&\le 
\sum_{R\in\N}\vol_W\left(W\cap B_R\right) Ne^{-\frac{(R-1)^2}{N}}
\\
&\le 
\sum_{R\in\N}C(\e)e^{\e R^2} Ne^{-\frac{(R-1)^2}{N}}<+\infty,
\eega
for all $\e<\frac{1}{N}$. this concludes the proof of theorem \ref{thm:mainEgau}.

\subsection{Proof of Corollary \ref{cor:mainjgau}}
\begin{proof}[Proof of Corollary \ref{cor:mainjgau}]
Corollary \ref{cor:mainjgau} follows directly from Theorem \ref{thm:mainEgau} and the following observation. The GRS $j^rX$ is a Gaussian random section of the vector bundle $J^rE\to M$. Its first jet $j^1(j^rX)$ is equivalent to $j^{r+1}X$ and thus its covariance tensor is determined by the $(2r+2)^{th}$ jet of $K_X$. 
\end{proof}
%
%
%
%
\begin{subappendices}
\section{Densities}\label{app:densities}
Let $M$ be a smooth manifold. The \emph{density bundle} or, as we call it, the \emph{bundle of density elements} $\Delta M$ is the vector bundle:
\be \Delta M=\wedge^m(T^*M)\otimes L,\ee
where $L$ is the orientation bundle (see \cite[Section 7]{botttu}); $\Delta M$ is a smooth real line bundle and the fiber can be identified canonically with
\bega\label{eq:apdensi}
\Delta_pM 
=\{\delta\colon (T_pM)^m\to \R\colon \delta=\pm |\w|, \text{ for some $\w\in \wedge^mT^*_pM$}\}.
\eega
We call an element $\delta\in \Delta_pM$, a \emph{density element} at $p$.

 Given a set of coordinates $x^1,\dots,x^m$ on $U$, we denote
\be
dx=dx^1\dots dx^m=|dx^1\wedge\dots\wedge dx^m|.
\ee
Using the language of \cite{botttu}, $dx$ is the section $(dx_1\wedge\cdots \wedge dx_m)\otimes e_U$, where $e_U$ is the section of $L|_{U}$ which in the trivialization induced by the chart $(x_1, \ldots, x_m)$ corresponds to the constant function $1$.

We have, directly from the definition, that 
\be 
dx=\left|\det\left(\frac{\de x}{\de y}\right)\right|dy.
\ee
for any other set of coordinates $y^1,\dots,y^m$. It follows that an atlas for $M$ with transition functions $g_{a,b}$ defines a trivializing atlas for the vector bundle $\Delta M$, with transition functions for the fibers given by $|\det(Dg_{a,b})|$. 

Each density element $\delta$, is either positive or negative, respectively if $\delta=|\w|$ or $\delta=-|\w|$, for some skew-symmetric multilinear form $w\in \wedge^mT^*_pM$. In other words, the bundle $\Delta M$ is canonically oriented and thus trivial (but not canonically trivial) and we denote the subbundle of positive elements by $\Delta^+ M\cong M\times [0,+\infty)$ (again, not canonically).

The \emph{modulus} of a density element is defined in coordinates, by the identity $|(\w(x)dx)|=|\w(x)|dx$. This defines a continuous map $|\cdot|\colon \Delta M\to \Delta^+ M$.

The sections of $\Delta M$ are called \emph{densities} and we will usually denote them as maps  $p\mapsto \delta(p)$. We define $\mathscr{D}^r(M)$ to be the space of $\mC^r$ densities and by $\mathscr{D}^r_c(M)$ the subset of the compactly supported ones. For smooth densities, we just write $\mathscr{D}=\mathscr{D}^\infty$.
From the formula for transitition functions it is clear that there is a canonical linear function 
\be 
\int_M \colon \mathscr{D}_c(M)\to \R \qquad \int_M \delta =:\int_{M}\delta(p)dp.
\ee
If $N$ is Riemannian and $f\colon M\to N$ is a $\ci$ map with $\dim M\le \dim N$, we define the \emph{jacobian density} $\delta f\in\mathscr{D}(M)$ by the following formula, in coordinates:
\be\label{eq:indensity}
\delta_u f=
\sqrt{\det\left(\frac{\de f}{\de u}^Tg(f(u))\frac{\de f}{\de u}\right)}du\in \Delta_p^+M,
\ee
In particular $\delta(\text{id}_M)$ is the Riemannian volume density of $M$ and we denote it by $dM$. Similarly, if $f$ is a Riemannian inclusion $M\subset N$, then $\delta f=dM$. 

Given any Riemannian metric on the manifold $M$, one can make the identification $\Delta_p M= \R (dM(p))$ and treat densities as if they were functions. Moreover this identification preserves the sign, since $dM$ is always a positive density: $dM(p)\in\Delta_p^+M$, for every $p\in M$. 
We denote by $B(M)$ the set of all Borel measurable functions $M\to [-\infty,+\infty]$ and by $L(M)=B(M)dM$ the set of all measurable, not necessarily finite, densities (the identification with $B(M)$ depends on the choice of the metric). Let $B^+(M)$ be the set of positive measurable functions $M\to [0,+\infty]$ and $L^+(M)$ be the set of densities of the form $\rho dM$, for some $\rho\in B^+(M)$. In other words 
\be 
L^+(M)=\left\{ \text{measurable functions }M\ni p\mapsto \delta(p)\in\Delta^+_pM \cup \{+\infty\}\right\}
\ee 
is the set of all nonnegative, non necessarily finite measurable densities. 
The integral can be extended in the usual way (by monotone convergence) to a linear function $\int_M\colon L^+(M)\to \R\cup\{\infty\}$. Similarly, we define the spaces $L^1(M),L^\infty(M),L^1_{loc}(M),L^\infty_{loc}(M)$ and their respective topologies by analogy with the standard case $M=\R^m$. 

\begin{defi}\label{def:abscont} We say that a real Radon measure (see \cite{ambrofuscopalla}) $\mu$ on $M$ is \emph{absolutely continuous} if $\mu(A)=0$ on any zero measure subset $A\subset M$. In other words $\mu$ is absolutely continuous if and only if $\f_*\left(\mu|_{U}\right)$ is absolutely continuous with respect to the Lebesgue measure $\mathscr{L}^m$ for any chart $\f\colon U\subset M\to \R^m$.
\end{defi}
In this language, the Radon-Nikodym Theorem takes the following form.
\begin{thm}\label{thm:RadoNiko} Let $\mu$ be an absolutely continuous real Radon measure on $M$. Then there is a density $\delta\in L^1_{loc}(M)$ such that, for all Borel subsets $A\subset M$,
\be 
\mu(A)=\int_A\delta.
\ee
\end{thm}
\section{Angle between subspaces}\label{app:angle}
Let $E, \langle\cdot,\cdot \rangle$ be a euclidean vector space (i.e. a finite dimensional real Hilbert space).
\begin{defi}
Let $f=(f_1,\dots,f_k)$ be a tuple (in row) of vectors $f_i\in E$. We define its volume as
\be 
\vol(f)=\sqrt{\det\langle f^T, f\rangle}.
\ee
If the vectors $f_1,\dots,f_k$ are independent, then $f$ is called a \emph{frame}. The \emph{span} $\text{span}(f)$ of the tuple $f$ is the subspace of $E$ spanned by the vectors $f_1,\dots,f_k$. A tuple $f$ is a \emph{basis} of a subspace $V\subset E$ if and only if it is a frame and its span is $V$. Given two frames $v,w$, we can form the tuple $(v,w)$.
\end{defi}
\begin{defi}\label{defi:angle}
Let $V,W\subset E$ be subspaces. Let $v$ and $w$ be frames in $E$ such that:
\begin{enumerate}[$\bullet$]
    \item $\text{span}(v)=V\cap (V\cap W)^\perp$;
    \item $\text{span}(w)=W\cap (V\cap W)^\perp$.
\end{enumerate}
We define the \emph{angle} 
 between $V$ and $W$ as
\be 
\sigma(V,W)=\begin{cases}
\frac{\vol(v w)}{\vol(v)\vol(w)} \quad &\text{ if $V\not\subset W$ and $W\not\subset V$;}\\
1 &\text{ otherwise.} 
\end{cases}
\ee
\end{defi}
It is easy to see that the definition is well posed, independently from the choices of the frames.  This definition corresponds to Howard's \cite{Howard} in the case when $V\cap W=\{0\}$.
Observe that $\sigma$ is symmetric and that we have $\sigma(V,W)=1$ if and only if $V=A\oplus_\perp B$ and $W=A\oplus_\perp C$ with and $B\perp C$.

When $V$ and $W$ are one dimensional, $\sigma(V,W)=|\sin\theta |$, where $\theta$ is the angle between the two lines. In general $\sigma(V,W)\in [0,1]$ is equal to the product of the sines of the nontrivial principal angles between $V$ and $W$ (see \cite{JordAngles,Zhu_2013_ang,MIAO199281_ang}). In particular, notice that $\sigma(V,W)>0$ always.

It is important to notice that $\sigma\colon \text{Gr}_k(E)\times \text{Gr}_h(E)\to [0,1]$ is not a continuous function. However the restriction to the subset of the pairs of subspaces $(V,W)$ such that $\dim (V+W)=n$ is continuous.
\begin{prop}\label{prop:chi}
Assume that $W\not\subset V$. Let $w$ be a basis for $W\cap (V\cap W)^\perp$, as in Definition \ref{defi:angle}, then
\be 
\sigma(V,W)=\frac{\vol\left(\Pi_{V^\perp}(w)\right)}{\vol(w)}
\ee
\end{prop}
\begin{proof}
First, observe that the projected frame $\Pi_{V^\perp}(w)$ is a basis of the space $(V+W)\cap V^\perp$. Let $\nu$ be an orthonormal basis of the same space, let $v$ be an orthonormal basis for $V\cap (V\cap W)^\perp$ and let $\tau$ be a basis for $V\cap W$. It follows that $(\tau, v,\nu)$ is a basis for $V+W$. Therefore, there is an invertible matrix $B$ and a matrix $A$ such that
\be w=\begin{pmatrix}\tau & v &\nu\end{pmatrix}\begin{pmatrix} 0\\ A\\ B\end{pmatrix}.\ee
Then, by Definition \ref{defi:angle}, we have
\be 
\begin{aligned}
\sigma(V,W)=\frac{\vol(v,w)}{\vol(v)\vol(w)}
=\frac{\left|\det\begin{pmatrix}\mathbb{1} & A \\ 0 & B \end{pmatrix}\right|}{\vol(w)}=\frac{\vol\left(\Pi_{V^\perp}(w)\right)}{\vol(w)}.
\end{aligned}
\ee
\end{proof}
\begin{prop}\label{prop:angleperp}
$\sigma(V^\perp,W^\perp)=\sigma(V,W).$
\end{prop}
\begin{proof}
The statement is trivially true if $V\subset W$ or $W \subset V$, so let us assume that this is not the case. Let $\nu$ be an orthonormal basis of the space $(V+W)\cap V^\perp$ and $v$ be an orthonormal basis of $V\cap (V\cap W)^\perp$. Besides, let $(\tau,w)$ be an orthonormal basis of $W$, such that $\tau$ is a basis for $V\cap W$. We have 
\be 
\begin{aligned}
\sigma(V^\perp,W^\perp)&=\sigma(W^\perp,V^\perp)\\
&=\vol\left(\Pi_{W}(\nu)\right)=\\
&=\det\langle w^T,\nu\rangle=\\
&=\det\langle \nu^T,w\rangle=\\
&= \vol \left(\Pi_{V^\perp}(w)\right)=\sigma(V,W).
\end{aligned}
\ee
\end{proof}
\section{Area and Coarea formula}\label{app:coarea} 
\begin{defi}\label{def:jacob}
Let $\f\colon M\to N$ be a $\ci$ map between $\ci$ Riemannian manifolds. The \emph{Jacobian} (often called \emph{normal Jacobian} when $f$ is a submersion) of $\f$ at $p\in M$ is
\be 
J_p\f:=\begin{cases}
0 & \text{if rank$(d_p\f)$ is not maximal;} \\
\frac{\vol_{N}\left(d_p\f(e)\right)}{\vol_{M}(e)} & \text{otherwise;}
\end{cases}
\ee
where $e=(e_1,\dots, e_k)$ is any basis of $\ker(d_p\f)^\perp\subset T_pM$. 

If $L\colon V_1\to V_2$ is a linear map between metric vector spaces, then we write $J L:=J_0L$ (Clearly $J_p\f=Jd_p\f$). If $V_1,V_2$ have the same dimension, then, to stress this fact, we may write $|\det L|:=JL$, although the sign of $\det L$ is not defined, unless we specify orientations.
\end{defi}
In particolar, let $\f\colon (\R^m;g_1)\to (\R^n;g_2)$ with differential $\frac{\de \f}{\de u}(u)=A$ having maximal rank, then
\be \label{eq:jacobidet}
J_u\f:=\begin{cases}
\sqrt{\frac{\det(A^Tg_2A)}{\det(g_1)}} & \text{if $m\le n$;} \\
\sqrt{\frac{\det(Ag_1^{-1}A^T)}{\det(g_2)^{-1}}} & \text{if $m\ge n$.}
\end{cases}
\ee
\begin{remark} In the case $m\le n$, the density induced on $M$ by a map $\f$, defined in \eqref{eq:indensity} corresponds to
$\delta_p\f=(J_p\f) dM$.
\end{remark}
\begin{thm}[Area formula]\label{thm:area}
Let $f\colon M\to N$ be a Lipschitz map between $\ci$ Riemannian manifolds, with $\dim M=\dim N$. Let $g\colon M\to [0,+\infty]$ be a Borel function, then
\be 
\int_M g(p)(J_pf) dM(p)=\int_N\left[\sum_{p\in f^{-1}(q)}g(p)\right]dN(q).
\ee
\end{thm}
\begin{proof}
See \cite[Theorem 3.2.3]{federer1996}.
\end{proof}
The Area formula is actually much more general than this, in that it holds for $\dim M\le \dim N$ and with the Hausdorff measure instead of $dN$. However, this simplified statement is all that we need in this paper. It also can be thought as a generalization of the following, in the case $\dim M=\dim N$.
\begin{thm}[Coarea formula]\label{thm:coarea}
Let $f\colon M\to N$ be a $\mC^1$ submersion between smooth Riemannian manifolds, with $\dim M\ge \dim N$. Let $g\colon M\to [0,+\infty]$ be a Borel function, then
\be 
\int_M g(p) (J_pf) dM(p)=\int_N\int_{f^{-1}(q)}g(p) d\left(f^{-1}(q)\right)(p) dN(q).
\ee
\end{thm}
\begin{proof}
See \cite{chavel} or deduce it from \cite[Theorem 3.2.12]{federer1996}.
\end{proof}

%
%
%
\end{subappendices}

\chapter{MAXIMAL AND TYPICAL TOPOLOGY OF REAL POLYNOMIAL SINGULARITIES}\label{chap:mttps}
Given a semialgebraic set $W\subseteq J^{r}(S^m, \R^k)$ and a polynomial map $\psi:S^m\to \R^k$ with components of degree $d$, we investigate the structure of the semialgebraic set $j^r\psi^{-1}(W)\subseteq S^m$ (we call such a set a ``singularity''). 

Concerning the upper estimate on the topological complexity of a polynomial singularity, we sharpen the classical bound $b(j^r\psi^{-1}(W))\leq O(d^{m+1})$, proved by Milnor \cite{MilnorBound}, with
\be\label{eq:abstract} b(j^r\psi^{-1}(W))\leq O(d^{m}),\ee
which holds for the generic polynomial map. 

For what concerns the ``lower bound'' on the topology of $j^r\psi^{-1}(W)$, we prove a general semicontinuity result for the Betti numbers of the zero set of $\mathcal{C}^0$ perturbations of smooth maps -- the case of $\mathcal{C}^1$ perturbations is the content of Thom's Isotopy Lemma (essentially the Implicit Function Theorem). This result is of independent interest and it is stated for general maps (not just polynomial); this result implies that small continuous perturbations of $\mathcal{C}^1$ manifolds have a richer topology than the one of the original manifold.

Keeping \eqref{eq:abstract} in mind, we compare the extremal case with a random one and prove that on average the topology of $j^r\psi^{-1}(W)$ behaves as the ``square root'' of its upper bound: for a random Kostlan map $\psi:S^m\to \R^k$ with components of degree $d$ and $W\subset J^{r}(S^m, \R^k)$ semialgebraic, we have:
\be \mathbb{E}b(j^{r}\psi^{-1}(W))=\Theta(d^{\frac{m}{2}}).\ee
This generalizes classical results of Edelman-Kostlan-Shub-Smale from the zero set of a random map, to the structure of its singularities.

\section{Introduction}
In this paper we deal with the problem of understanding the structure of the singularities of polynomial maps
\be \psi:S^m\to \R^k,\ee
where each component of $\psi=(\psi_1, \ldots, \psi_k)$ is the restriction to the sphere of a homogeneous polynomial of degree $d$. For us ``singularity'' means the set of points in the sphere where the $r$-jet extension $j^r\psi:S^n\to J^{r}(S^n, \R^k)$ meets a given semialgebraic set $W\subseteq J^{r}(S^n, \R^k).$ Example of these type of singularities are: zero sets of polynomial functions, critical points of a given Morse index of a real valued function or the set of Whitney cusps of a planar map.

Because we are looking at \emph{polynomial} maps, this problem has two different quantitative faces, which we both investigate in this paper. 

(1) From one hand we are interested in understanding the \emph{extremal} cases, meaning that, for fixed $m, d$ and $k$ we would like to know how complicated can the singularity be, at least in the generic case. 

(2) On the other hand, we can ask what is the \emph{typical} complexity of such a singularity. Here we adopt a measure-theoretic point of view and endow the space of polynomial maps with a natural Gaussian probability measure, for which it makes sense to ask about expected properties of these singularities, such as their Betti numbers.

\subsection{Quantitative bounds, the h-principle and the topology semicontinuity}
Measuring the complexity of  $Z=j^r\psi^{-1}(W)$ with the sum $b(Z)$ of its Betti numbers, problem (1) above means producing a-priori upper bounds for $b(Z)$ (as a function of $m, d, k$) as well as trying to realize given subsets of the sphere as $j^r\psi^{-1}(W)$ for some $W$ and some map $\psi$. 

For the case of the zero set $Z=\psi^{-1}(0)$ of a polynomial function $\psi:S^m\to\R$ of degree $d$, the first problem is answered by a Milnor's type bound\footnote{Milnor's bound \cite{MilnorBound} would give $b(Z)\leq O(d^{m+1})$, whereas \cite[Proposition 14]{LerarioJEMS} gives the improvement $b(Z)\leq O(d^m)$. In the context of this paper the difference between these two bounds is relevant, especially because when switching to the probabilistic setting it will give the so called ``generalized square root law''.} $b(Z)\leq O(d^m)$ and the second problem by Seifert's theorem: every smooth hypersurface in the sphere can be realized (up to ambient diffeomorphisms) as the zero set of a polynomial function. 

In the case of more general singularities, both problems are more subtle. The problem of giving a good upper bound on the complexity of $Z=j^r\psi^{-1}(W)$ will require us to develop a quantitative version of stratified Morse Theory for semialgebraic maps (Theorem \ref{thm:strat}). We use the word ``good'' because there is a vast literature on the subject of quantitative semialgebraic geometry, and it is not difficult to produce a bound of the form $b(Z)\leq O(d^{m+1})$; instead here (Theorem \ref{thm:bound} and Theorem \ref{thm:bound2}) we prove the following result.

\begin{thm}\label{thm:bettibound}For the generic polynomial map $\psi:S^m\to \R^k$ with components of degree $d$, and for $W\subseteq J^{r}(S^m, \R^k)$ semialgebraic, we have: 
\be \label{eq:estib}b(j^{r}\psi^{-1}(W))\leq O(d^m).\ee
(The implied constant depends on $W$.)
\end{thm}

In the case $W$ is algebraic we do not need the genericity assumption on $\psi$ for proving \eqref{eq:estib}, but in the general semialgebraic case some additional complications arise and this assumption allows to avoid them through the use of  Theorem \ref{thm:strat}. We believe, however, that \eqref{eq:estib} is still true even in the general case\footnote{In the algebraic case in fact one can use directly Thom-Milnor bound, but in the general semialgebraic case it is necessary first to ``regularize'' the semialgebraic set, keeping control on its Betti numbers. In the algebraic (or even the basic semialgebraic case) this is the procedure of Milnor \cite{MilnorBound}, in the general semialgebraic case it is not clear what this controlled regularization procedure would be. The nondegeneracy assumption on the jet allows us to avoid this step.}. Moreover, for our scopes the genericity assumption is not restrictive, as it fits in the probabilistic point of view of the second part of the paper, where a generic property is a property holding with probability one.

For what concerns the realizability problem, as simple as it might seem at first glance, given $W\subseteq J^{r}(S^m, \R^k)$ it is not even trivial to find a map $f:S^m\to \R^k$ whose jet is transversal to $W$ and such that $b(j^rf^{-1}(W))>0$ (we prove this in Corollary \ref{cor:topositive}). 

Let us try to explain carefully what is the subtlety here. In order to produce such a map, one can certainly produce a section of the jet bundle $\sigma:S^m\to J^{r}(S^m, \R^k)$ which is transversal to $W$ and such that $b(\sigma^{-1}(W))>0$ (this is easy). However, unless $r=0$, this section needs not to be holonomic, i.e. there might not exist a function $f:S^m\to \R^k$ such that $\sigma=j^{r}f$.  

We fix this first issue using an h-principle argument: the Holonomic Approximation Theorem \cite[p. 22]{eliash} guarantees that, after a small $\mathcal{C}^0$ perturbation of the whole picture, we can assume that there is a map $f:S^m\to \R^k$ whose jet $j^{r}f$ is $\mathcal{C}^0$ close to $\sigma$. 

There is however a second issue that one needs to address. In fact, if the jet perturbation was $\mathcal{C}^1$ small (i.e. if $\sigma$ and $j^rf$ were $\mathcal{C}^1$ close), Thom's Isotopy Lemma would guarantee that $\sigma^{-1}(W)\sim j^rf^{-1}(W)$ (i.e. the two sets are ambient diffeomorphic), but the perturbation that we get from the Holonomic Approximation Theorem is guaranteed to be only $\mathcal{C}^0$ small! To avoid this problem we prove the following general result on the semicontinuity of the topology of small $\mathcal{C}^0$ perturbations (see Theorem \ref{thm:semiconttop} below for a more precise statement).

\begin{thm}\label{thm:nonho}Let $S, J$ be smooth manifolds, $W\subseteq J$ be a closed cooriented submanifold 
 and $\sigma\in \mathcal{C}^{1}(S, J)$ such that $\sigma\pitchfork W$. Then for every $\gamma\in  \mathcal{C}^{1}(S, J) $ which is sufficiently close to $\sigma$ in the $\mathcal{C}^0$-topology and such that $\gamma\pitchfork W$, we have:
\be b(\gamma^{-1}(W))\geq b(\sigma^{-1}(W)).\ee
\end{thm}
In particular we see that if small $\mathcal{C}^1$ perturbations of a regular equation preserve the topology of the zero set, still if we take just small $\mathcal{C}^0$ perturbations the topology of such zero set can only increase.

To apply Theorem \ref{thm:nonho} to our original question we consider $S=S^m$ and $J=J^{r}(S^m, \R^k)$, $W\subseteq J^{r}(S^m, \R^m)$ is the  semialgebraic set defining the singularity and $\sigma:S^m\to J^r(S^m, \R^k)$ is the (possibly non-holonomic) section such that $\sigma\pitchfork W$ and  $b(\sigma^{-1}(W))>0$. Moreover we can construct $\sigma$ in such a way that its image meets only a small (relatively compact and cooriented) subset of the smooth locus of $W$. Then for every $f\in  \mathcal{C}^{r+1}(S^m, \R^k) $ with $\tau=j^rf$ sufficiently close to $\sigma$ and such that $j^rf\pitchfork W$, we have:
\be \label{eq:app}b(j^rf^{-1}(W))\geq b(\sigma^{-1}(W))>0.\ee
(We will use the content of Corollary \ref{cor:topositive} and the existence of a function $f$ such that \eqref{eq:app} holds in the second part of the paper for proving the convergence of the expected Betti numbers of a random singularity.)

\subsection{The random point of view and the generalized square-root law}Switching to the random point of view offers a new perspective on these problems: from Theorem \ref{thm:bettibound} we have an extremal bound \eqref{eq:estib} for the complexity of polynomial singularities, but it is natural to ask how far is this bound from the typical situation. 
Of course, in order to start talking about randomness, we need to choose a probability distribution on the space of (homogeneous) polynomials. It is natural to require that this distribution is gaussian, centered, and that it is invariant under orthogonal changes of variables (in this way there are no preferred points or directions in the sphere). If we further assume that the monomials are independent, this distribution is unique (up to multiples), and called the \emph{Kostlan distribution}. 

To be more precise, this probability distribution is the measure on $\R[x_0,\dots,x_{m}]_{(d)}$ (the space of homogeneous polynomials of degree $d$) induced by the gaussian random polynomial:
\be \label{eq:kosdef}
P(x)=\sum_{ |\a|=d}\xi_\a\cdot \left(\frac{d!}{\alpha_0!\cdots\alpha_m!}\right)^{1/2}  x_0^{\alpha_0}\cdots x_m^{\alpha_m},
\ee
where $\{\xi_\a\}$ is a family of standard independent gaussian variables. A list of $k$ independent Kostlan polynomials $P=(P_1, \ldots, P_k)$ defines a random polynomial map:
\be \psi=P|_{S^m}\to \R^k.\ee
In particular, it is now natural to view such a $\psi$ as a random variable in the space $\mathcal{C}^{\infty}(S^m, \R^k)$ and to study the differential topology of this map, such as the behavior of its singularities, described a preimages of jet submanifolds $W\subseteq J^{r}(S^m, \R^k)$ in the previous section.

In this direction, it has already been observed by several authors, in different contexts, that random real algebraic geometry seems to behave as the ``square root'' of generic complex geometry. Edelman and Kostlan \cite{Kostlan:93, EdelmanKostlan95} were the first to observe this phenomenon: a random Kostlan polynomial of degree $d$ in one variable has $\sqrt{d}$ many real zeroes, on average\footnote{In the notation of the current paper this correspond to the case of $\psi:S^1\to \R$ of degree $d$, whose expected number of zeroes is $2\sqrt{d}$. The multiplicative constant ``$2$'' appears when passing from the projective to the spherical picture}. Shub and Smale \cite{shsm} generalized this result and proved that the expected number of zeroes of a system of $m$ Kostlan equations of degrees $(d_1, \ldots, d_m)$ in $m$ variables is $\sqrt{d_1\cdots d_m}$ (the bound coming from complex algebraic geometry would be $d_1\cdots d_m$).

Moving a bit closer to topology, B\"urgisser \cite{buerg:07} and Podkorytov \cite{Po} proved that the expectation of the Euler characteristic of a random Kostlan algebraic set has the same order of the square-root of the Euler characteristic of its complex part (when the dimension is even, otherwise it is zero).
A similar result for the Betti numbers has also been proved by Gayet and Welschinger \cite{GaWe1, GaWe2, GaWe3}, and by Fyodorov, Lerario and Lundberg \cite{FyLeLu} for invariant distributions. 

Using the language of the current paper, these results correspond to the case of  a polynomial map $\psi:S^m\to \R^k$ and to the ``singularity'' $Z=j^0\psi^{-1}(W)$, where 
\be W=S^m\times \{0\}\subset J^{0}(S^m, \R^k)=S^m\times \R^k\ee
 and $j^0\psi(x)=(x, \psi(x))$ is the section given by the map $\psi$ itself. Here we generalize these results and prove that a similar phenomenon is a very general fact of Kostlan polynomial maps.
 
 \begin{thm}\label{thm:bettiorderintro}
Let $W\subset J^r(S^m,\R^k)$ be a closed intrinsic\footnote{We say that $W\subset J^r(S^m,\R^k)$ is intrinsic if it is invariant under diffeomorphisms of $S^m$, see Definition \ref{def:intrinsic}. This property it is satisfied in all natural examples.} semialgebraic set of positive codimension. If $\psi:S^m\to \R^k$ is a random Kostlan polynomial map, then
\be \label{eq:ineqbetti2}
\E b(j^r\psi^{-1}(W))=\Theta(d^{\frac{m}{2}}).\footnote{We write $f(d)=\Theta(g(d))$ if there exist constants $a_1, a_2>0$ such that $a_1 g(d)\leq f(d)\leq a_2f(d)$ for all $d\geq d_0$ sufficiently large.}
\ee
(The implied constants depend on $W$.)

\end{thm}

We call the previous Theorem \ref{thm:bettiorderintro} the ``generalized square root law'' after comparing it with the extremal inequality $b(j^{r}\psi^{-1}(W))\leq O(d^m)$
 from Theorem \ref{thm:bettibound}, whose proof is ultimately based on bounds coming from complex algebraic geometry\footnote{The reader can now appreciate the estimate $O(d^m)$ instead of $O(d^{m+1})$ from Theorem \ref{thm:bettibound}.}. In the case $W$ has codimension $m$ (i.e. when we expect $j^r\psi^{-1}(W)$ to consist of points), we actually sharpen \eqref{eq:ineqbetti2} and get the explicit asymptotic to the leading order, see Theorem \ref{thm:sqrlaw} below. Moreover, a similar result holds for every fixed Betti number $b_i(j^r\psi^{-1}(W))$ when $i$ is in the range $0\leq i\leq m-\mathrm{codim}(W),$ see Theorem \ref{thm:bettiorder} for a detailed statement.
  
\begin{remark}The ingredients for the proof of Theorem \ref{thm:bettiorderintro} are: Theorem \ref{thm:strat} for the upper bound and Corollary \ref{cor:topositive} for the lower bound. The main property that we use in this context is the fact that a Kostlan map $\psi:S^m\to \R^k$ has a rescaling limit when restricted to a small disk $D_d=D(x, d^{-1/2})$ around any point $x\in S^m$. In other words, one can fix a diffeomorphism $a_d:\mathbb{D}^m\to D_d$ of the standard disk $\mathbb{D}^m$ with the small spherical disk $D(x, d^{-1/2})\subset S^m$ and see that the sequence of random functions:
\be X_d=\psi\circ a_d:\mathbb{D}^m\to \R^k\ee
converges to the Bargmann-Fock field, see Theorem \ref{thm:Kostlan}. 
In Chapter \ref{chap:dtgrf} we introduced a general framework for dealing with random variables in the space of smooth functions and their differential topology -- again we can think of $X_d\in C^{\infty}(\mathbb{D}^m, \R^k)$ as a sequence of random variables of this type. The results from Chapters \ref{chap:dtgrf} and \ref{chap:kr}, applied to the setting of random Kostlan polynomial maps are collected in Theorem \ref{thm:Kostlan} below, which lists the main properties of the rescaled Kostlan polynomial $X_d$. Some of these properties are well-known to experts working on random fields, but some of them seem to have been missed. Moreover, we believe that our language is more flexible and well-suited to the setting of differential topology, whereas classical references look at these random variables from the point of view of functional analysis and stochastic calculus.

Of special interest from Theorem \ref{thm:Kostlan} are properties (2), (5) and (7), which are closely related. In fact (2) and (5) combined together tells that open sets $U\subset \mathcal{C}^{\infty}(\mathbb{D}^m, \R^k)$ which are defined by open conditions on the $r$-jet of $X_d$, have a positive limit probability when $d\to \infty$. Property (7), tells that the law for Betti numbers of a random singularity $Z_d=j^rX_d^{-1}(W)$ has a limit. (Even in the case of zero sets this property was not noticed before, see Example \ref{ex:welsch}.)

We consider Theorem \ref{thm:Kostlan} as a practical tool that people interested in random algebraic geometry can directly use, and we will show how to concretely use this tool in a list of examples that we give in Appendix \ref{sec:examples}.
\end{remark}

\begin{remark}The current paper, and in particular the generalized square-root law Theorem \ref{thm:bettiorderintro}, complement recent work of Diatta and Lerario \cite{DiattaLerario} and Breiding, Keneshlou and Lerario \cite{HaniehPaulAnto}, where tail estimates on the probabilities of the maximal configurations are proved. 
\end{remark}
\subsection{Structure of the paper}In Section \ref{sec:strat} we prove a quantitative semialgebraic version of stratified Morse Theory, which is a technical tool needed in the sequel, and in Section \ref{sec:quantitative} we prove Theorem \ref{thm:bound} and Theorem \ref{thm:bound2} (whose combination give Theorem \ref{thm:bettibound}).  In section \ref{sec:holo} we discuss the semicontinuity of topology under holonomic approximation and prove Theorem \ref{thm:semiconttop} (which is Theorem \ref{thm:nonho} from the Introduction). In Section \ref{sec:RAG} we introduce the random point of view and prove the generalized square-root law. Appendix \ref{sec:examples} contains three short examples of use the random techniques.

 \section{Quantitative bounds, the h-principle and the topology semicontinuity}
 \subsection{Stratified Morse Theory}\label{sec:strat}
 
Let us fix a Whitney stratification $W=\sqcup_{S\in \mathscr{S}}S$ (see \cite[p. 37]{GoreskyMacPherson} for the definition) of the semialgebraic subset $W\subset J^r(S^m,\R^k)=:J$, with each stratum $S\in \mathscr{S}$ being semialgebraic and smooth (such decomposition  exists \cite[p. 43]{GoreskyMacPherson}), so that, by definition a smooth map $f\colon M\to J$, is transverse to $W$ if $f\transv S$ for all strata $S\in \mathscr{S}$. 
 When this is the case, we write $\psi\transv W$ and implicitly consider the subset $\psi^{-1}(W)\subset M$ to be equipped with the Whitney stratification given by $\psi^{-1}\mathscr{S}=\{\psi^{-1}(S)\}_{S\in\mathscr{S}}$.
 \begin{defi}\label{def:mors}
Given a Whitney stratified subset $Z=\cup_{i\in I}S_i$ of a smoooth manifold $M$ (without boundary), we say that a function $g\colon Z\to \R$ is a Morse function if $g$ is the restriction of a smooth function $\tilde{g}\colon M\to \R$ such that
\begin{enumerate}[(a)]
\item $g|_{S_i}$ is a Morse function on $S_i$. 
\item For every critical point $p\in S_i$ and every generalized tangent space $Q\subset T_pM$ (defined as in \cite[p. 44]{GoreskyMacPherson})
we have $d_p\tilde{g}(Q)\neq 0$, except for the case $Q=T_pS_i$.
\end{enumerate}
\end{defi}
Note that the condition of being a Morse function on a stratified space $Z\subset M$ depends on the given stratification of $Z$.
\begin{remark}
The definition above is slightly different than the one given in the book \cite[p. 52]{GoreskyMacPherson} by Goresky and MacPherson, where a Morse function, in addition, must be proper and have distinct critical values.
\end{remark}
%
The following theorem is the quantitative version of stratified Morse theory for semialgebraic maps we need in order to prove Theorem \ref{thm:bettibound}.
\begin{thm}\label{thm:strat}
Let $W\subset J$ be a semialgebraic subset of a real algebraic smooth manifold $J$, with a given semialgebraic Whitney stratification $W=\sqcup_{S\in \mathscr{S}}S$ and let $M$ be a real algebraic smooth manifold. 
There exists a semialgebraic subset $\hat{W}\subset J^{1}(M,J\times \R)$ having codimension larger or equal than $\dim M$, equipped with a semialgebraic Whitney stratification that satisfies the following properties with respect to any couple of smooth maps $\psi\colon M\to J$ and $g\colon M\to \R$.
\begin{enumerate}
\item 
If $\psi\transv W$and $j^1(\psi,g)\transv \hat{W}$, then $g|_{\psi^{-1}(W)}$ is a Morse function with respect to the stratification $\psi^{-1}\mathscr{S}$ and
\be \label{eq:stratcrit}
\text{Crit}(g|_{\psi^{-1}(W)})=\left(j^{1}(\psi,g)\right)^{-1}(\hat{W}).
\ee

\item
There is a constant $N_W>0$ depending only on $W$ and $\mathscr{S}$, such that if $\psi^{-1}(W)$ is compact, $\psi\transv W$ and $j^1(\psi,g)\transv \hat{W}$, then
\be 
b_i(\psi^{-1}(W))\le N_W\#\text{Crit}(g|_{\psi^{-1}(W)}),
\ee
for all $i=0,1,2\dots$
\end{enumerate}
\end{thm}

\begin{proof}
Let $S\in \mathscr{S}$ be a stratum of $W$, hence $S\subset J$ is a smooth submanifold and since $\psi\transv W$ implies that $\psi\transv S$, we also have that $\psi^{-1}(S)$ is a submanifold of $M$ of the same codimension which we denote by $k$. Define
\be 
\begin{aligned}
\hat{S}&=\{j^1_p(F,f)\in J^1(M,J\times\R)\colon F(p)\in S \text{ and }d_pf\in d_pF^*(T_{F(p)}S^\perp)\}\\
&=\{j^1_p(F,f)\in J^1(M,J\times\R)\colon F(p)\in S \text{ and }\exists \lambda\in T_{F(p)}S^\perp \text{ s.t. } d_pf=\lambda\circ d_pF\}.
\end{aligned}
\ee
Orthogonality here is meant in the sense of dual vector spaces: if $Q\subset T$ are vector spaces, then $Q^\perp=\{\xi\in T^*\colon \xi(Q)=0\}$.

It is clear, by this definition, that $\hat{S}$ is semialgebraic and its codimension is equal to the dimension of $M$. 
\begin{claim}\label{claim:stratproof1}
$j^1_{p_0}(\psi,g)\in \hat{S}$ if and only if $p_0$ is a critical point for $g|_{\psi^{-1}(S)}$.
\end{claim}
If $j^1_{p_0}(\psi,g)\in \hat{S}$, then of course $p_0\in \psi^{-1}(S)$ and there exists a (Lagrange multiplier) conormal covector $\lambda\in T_{\psi(p_0)}S^\perp$ such that $d_{p_0}g=\lambda\circ d_{p_0}\psi$. It follows that $d_{p_0}g$ vanishes on $T_{p_0}\psi^{-1}(S)=d_{p_0}\psi^{-1}(T_{\psi(p_0)}S)$. This proves the ``only if'' statement of the Claim as a consequence of the following inclusion
\be 
d_{p_0}\psi^*\left(T_{p_0}S^\perp\right)\subset\left(T_{p_0}\psi^{-1}(S)\right)^\perp.
\ee

To conclude the proof of Claim \ref{claim:stratproof1} we need to show the opposite inclusion. We do this by showing that the dimensions of the two spaces are equal. First observe that, since by hypotheses $\psi\transv S$, the image $d_{p_0}\psi$ is a complement to $T_{\psi(p_0)}S$ in $T_{\psi(p_0)}J$ and this is equivalent (it is the dual statement) to say that the restriction of $d_{p_0}\psi^*$ to $(T_{\psi(p_0)}S)^\perp$ is injective. It follows that
\be 
\begin{aligned}
\dim d_{p_0}\psi^*\left(T_{p_0}S^\perp\right)
&=
\dim \left(T_{p_0}S^\perp\right)=\\
&=
\mathrm{codim}\ S =\\
&=\mathrm{codim}\ \psi^{-1}(S)=\\
&= \dim \left(T_{p_0}\psi^{-1}(S)\right)^\perp.
\end{aligned}
\ee
This concludes the proof of Claim \ref{claim:stratproof1}.

\begin{claim}\label{claim:stratproof2}
Given a Whitney stratification of $\hat{S}$, and a critical point $p_0\in M$ of the map $g|_{\psi^{-1}(S)}$, if $j^1(\psi,g)\transv \hat{S}$ at $p_0$ then this critical point is Morse.
\end{claim}
Let us pass to a coordinate chart $\phi$ defined on a nighborhood $\mathcal{U}\subset J^1(M,J\times \R)$ of $j^1_{p_0}(\psi,g)$:
\be 
\phi=\left(
x=\begin{pmatrix}
x^1 \\ x^2
\end{pmatrix},
y=\begin{pmatrix}
y^1 \\ y^2
\end{pmatrix},
a,
Y=\begin{pmatrix}
Y^1 \\ Y^2
\end{pmatrix},
A
\right)
\colon \mathcal{U}\to \R^m\times \R^{s+k}\times \R\times \R^{(s+k)\times m}\times \R^m
\ee
\be 
j^1_p(F,f)\mapsto \left(x(p),y(F(p)),g(p),\frac{\de (y\circ F)}{\de x},\frac{\de g}{\de x}\right);
\ee
where $y^2=0$ is a local equation for $S$ and $x^2=0$ is a local equation for $\psi^{-1}(S)$. Indeed, by the implicit function theorem (applied to the map $\psi$ in virtue of the transversality assumption $\psi\transv S$) we can assume that $y^2(\psi(x^1,x^2))=x^2.$ In this coordinate chart we have that the restriction of $d_{p}F^*$ to the space ${T_{\psi(p)}S^\perp}$ is represented by the matrix $({Y}^2)^T$, thus 
\be \label{eq:shat}
\hat{S}\cap \mathcal{U}=\left\{y^2=0; A\in \mathrm{Im}\left((Y^2)^T\right)\right\}\cap \phi(\mathcal{U}).
\ee
Let us denote by $x\mapsto(x,\tilde{y}(x),\tilde{a}(x),\tilde{Y}(x),\tilde{A}(x))$ the local expression of the jet map $p\mapsto j^1_p(\psi,g)$ with respect to the above coordinates.
By construction we have that \be \left(\tilde{Y}^2(p_0)\right)^T=\begin{pmatrix} 0 \\ \mathbbm{1}_k\end{pmatrix}.\ee In particular the image of the above matrix is a complement to the subspace spanned by the first $m-k$ coordinates
and we may assume, reducing the size of the neighborhood if needed, that this property holds for every element 
$(x,y,a,Y,A)\in\phi(\mathcal{U})$, so that there exist unique vectors $\lambda\in\R^k$ and $\xi\in\R^{(m-k)}$ such that
\be \label{eq:cheduecoioni}
A=\begin{pmatrix}
A_1 \\ A_2
\end{pmatrix}=\left(Y^2\right)^T \lambda+\begin{pmatrix}
\xi\\ 0
\end{pmatrix}.
\ee 
Now, this defines a smooth function $\xi\colon\mathcal{U}\to \R^k$ such that the equations $y^2=0$; $\xi=0$ are smooth regular equations for $\phi(\hat{S}\cap\mathcal{U})$. 

Notice that this ensures that $\phi(\mathcal{U})$ intersects only the smooth locus of $\hat{S}$. Now, since by hypotheses $j^1(\psi,g)$ is transverse to all the strata of $\hat{S}$ then it must be transverse to the smooth locus in the usual sense, even if the latter is a union of strata (this follows directly from the definition of transversality). Therefore, while proving Claim \ref{claim:stratproof2}, we are allowed to forget about the stratification of $\hat{S}$ and just assume that the map $j^1(\psi,g)$ is transverse to the smooth manifold $\hat{S}\cap \phi(\mathcal{U})$ in the usual sense.

In this setting we can see that if $j^1(\psi,g)\transv \hat{S}$ at $p_0$, then the following matrix has to be surjective:
\be 
\begin{pmatrix}
d y^2 \\ d\xi
\end{pmatrix}\circ d_{p_0}\left(j^1(\psi,g)\right)=
\begin{pmatrix}
0 & \mathbbm{1}_k \\
\frac{\de\tilde{\xi}}{\de x^1}(p_0) &
\frac{\de\tilde{\xi}}{\de x^2}(p_0)
\end{pmatrix}\in \R^{(k+k)\times ((m-k) +k)},
\ee
where $\tilde{\xi}(x)=\xi(x,\tilde{y}(x),\tilde{a}(x),\tilde{Y}(x),\tilde{A}(x))$. 
Therefore the lower left block $\frac{\de\tilde{\xi}}{\de x^1}(p_0)$ is surjective as well and hence invertible. This concludes our proof of Claim \ref{claim:stratproof2} since such matrix is in fact the hessian of the map $g|_{\psi^{-1}(S)}$ at the critical point $p_0$:
\be 
\begin{aligned}
d_{p_0}\left(g|_{\psi^{-1}(S)}\right)&=
\frac{\de}{\de x^1}\Big|_{p_0}\left(\frac{\de g}{\de x^1}\right)
=\frac{\de \tilde{A}_1}{\de x^1}(p_0)= \frac{\de\tilde{\xi}}{\de x^1}(p_0).
\end{aligned}
\ee
The last equality is due to the equation \eqref{eq:cheduecoioni} combined with the observation that $\tilde{Y^2}$ is of the form $\begin{pmatrix}0 &*\end{pmatrix}$ for all $p$ in a neighborhood of $p_0$, since $\frac{\de \tilde{y}^2}{\de x_1}(p)=0$. 

At this point, Claim \ref{claim:stratproof1} and Claim \ref{claim:stratproof2} prove that, for whatever stratification of $\hat{S}$, if $j^1(\psi,g)\transv \hat{S}$ and $\psi\transv S$ then $g|_{\psi^{-1}(S)}$ is a Morse function and that its critical set coincide with the set $\left(j^1(\psi,g)\right)^{-1}(\hat{S})$, so that condition $(a)$ of Definition \ref{def:mors} is satisifed along the stratum $S$. In order to establish when $g|_{\psi^{-1}(W)}$ is a Morse function along the stratum $\psi^{-1}(S)$ on the stratified manifold $W$, in the sense of Definition \ref{def:mors}, we now need to prove a similar statement to ensure condition $(b)$.

Let us consider the set $D_qS$ of degenerate covectors at a point $q\in S$ that are conormal to $S$ (conormal and degenerate covectors are defined as in \cite[p.44]{GoreskyMacPherson}), in other words:
\be 
D_qS=\{\xi\in T^*_q{J}\colon \xi\in T_qS^\perp,\  \xi\in Q^\perp \text{ for some $Q$ generalized tangent space at $q$}\}. 
\ee
It is proved in \cite[p.44]{GoreskyMacPherson} that $DS=\cup_{q\in S}D_qS$ is a semialgebraic subset of codimension greater than $1$ of the conormal bundle $TS^\perp$\footnote{$TS^\perp =T_S^*J$ , in the notation of \cite{GoreskyMacPherson}.} to the stratum S.
We claim that the subset $D\hat{S}\subset \hat{S}$ containing the jets that do not satisfy condition $(b)$ of Definition \ref{def:mors} has the following description:
\be 
D\hat{S}=\{j^1_p(F,f)\in J^1(M,J\times\R)\colon F(p)\in S \text{ and }d_pf\in d_pF^*(D_{F(p)}S)\}.
\ee
In fact, since $\psi\transv W$, then all the generalized tangent spaces of the stratified subset $\psi^{-1}(W)\subset M$ at a point $p\in \psi^{-1}(S)$ are of the form $d_p\psi^{-1}(Q)$. It follows that if a conormal covector $d_pg=\lambda \circ d_p\psi$ is degenerate then $\lambda\in D_{\psi(p)}S$.

Note that $D\hat{S}$ is a subset of $\hat{S}$ of codimension $\ge 1$, thus the codimension of $D\hat{S}$ in $J^1(M,J\times \R)$ is $\ge m+1$. As a consequence we have that $j^1(\psi,g)\transv D\hat{S}$ if and only if $ j^1(\psi,g)\notin D\hat{S}$.
Therefore if $j^1(\psi,g)\transv \hat{S}$ and $j^1(\psi,g)\notin D\hat{S}$ then $\psi\transv S$ and $g|_{\psi^{-1}(W)}$ is a Morse function on $\psi^{-1}(W)$ along the stratum $\psi^{-1}(S)$.

We are now ready to define $\hat{W}=\cup_{S\in\mathscr{S}}\hat{S}$. An immediate consequence of Claim \ref{claim:stratproof1} is that $\hat{W}$ satisfies equation \eqref{eq:stratcrit}. Moreover, since $\hat{S}\supset D\hat{S}$ are semialgebraic, $\hat{W}$ is semialgebraic and admits a semialgebraic Whitney stratification $\hat{\mathscr{S}}$ (refining the one of $\hat{S}$) such that all the subsets $\hat{S}$ and $D\hat{S}$ are unions of strata. With such a stratification, if the jet map $j^1(\psi,g)$ is transverse to $\hat{W}$ then, for each stratum $S\in \mathscr{S}$, it is also transverse to $\hat{S}$ and it avoids the set $D\hat{S}$, so that $g|_{\psi^{-1}(W)}$ is a Morse function, in the sense of Definition \ref{def:mors}.
This proves that $\hat{W}$ satisfies condition $(1)$ of the Theorem.

Let us prove condition (2). Let $Z=\psi^{-1}(W)\subset M$ be compact. Without loss of generality we can assume that each of the critical values $c_1,\dots, c_n$ of $g|_{Z}$ corresponds to only one critical point (this can be obtained by makingcontaining the jets that do not satisfy condition $(b)$ of Definition \ref{def:mors}: a $\mathcal{C}^1$ small perturbation of $g$, which won't affect the number of its critical points). Consider a sequence of real numbers $a_1,\dots a_{n+1}$ such that
\be 
a_1<c_1 <a_2<c_2<\dots <a_{n}<c_n<a_{n+1}.
\ee
By the main Theorem of stratified Morse theory \cite[p. 8, 65]{GoreskyMacPherson}, there is an homeomorphism
\be 
Z\cap \{g\le a_{l+1}\}\cong(Z\cap \{g\le a_l\})\sqcup_B A,
\ee
with
\be 
(A,B)= TMD_p(g)\times NMD_p(g),
\ee
where $TMD_p(g)$ is the tangential Morse data and $NMD_p(g)$ is the normal Morse data. A fundamental result of classical Morse theory is that the tangential Morse data is homeomorphic to a pair 
\be 
TMD_p(g)\cong(\mathbb{D}^\lambda\times \mathbb{D}^{m-\lambda},(\de\mathbb{D}^\lambda)\times \mathbb{D}^{m-\lambda}),
\ee 
while the normal Morse data is defined as the local Morse data of $g|_{N_p}$ for a normal slice (see \cite[p. 65]{GoreskyMacPherson}) at $p$. A consequence of the transversality hypothesis $\psi\transv W$ is that there is a small enough normal slice $N_p$ such that $\psi|_ {N_p}\colon N_p\to J$ is the embedding of a normal slice at $\psi(p)$ for $W$. Therefore the normal data $NMD_p(g)$ belongs to the set $\nu(W)$ of all possible normal Morse data that can be realized (up to homeomeorphisms) by a critical point of a Morse function on $W$. By Corollary $7.5.3$ of \cite[p. 95]{GoreskyMacPherson} it follows that the cardinality of the set $\nu(W)$ is smaller or equal than the number of connected components of the semialgebraic set $\cup_{S\in \mathscr{S}}(TS^\perp\-DS)$, hence finite\footnote{In the book this is proved only for any fixed point $p$, as a corollary of Theorem $7.5.1$ \cite[p.93]{GoreskyMacPherson}. However the same argument generalizes easily to the whole bundle.}.
Let 
\be 
N_W:=\max_{Y \in \nu(W),\ \lambda \in \{0,\dots, m\}} b_i\left(\left(\mathbb{D}^\lambda\times \mathbb{D}^{m-\lambda},(\de\mathbb{D}^\lambda)\times \mathbb{D}^{m-\lambda}\right)\times Y\right)\in \N.
\ee
From the long exact sequence of the pair $\left(Z\cap \{g\ge a_{l+1}\},(Z\cap \{g\ge a_{l}\}\right)$ we deduce that
\be\label{eq:lesbetti}\begin{aligned}
b_i(Z\cap \{g\le a_{l+1}\})-b_i(Z\cap \{g\le a_{l}\})&\le b_i\left(Z\cap \{g\le a_{l+1}\},Z\cap \{g\le a_{l}\}\right)\\
&= b_i\left(A,B\right) \\
&= b_i\left(TMD_p(g)\times NMD_p(g)\right) \\
&\le N_W.
\end{aligned}
\ee
Since $Z$ is compact, the set $Z\cap \{g\le a_1\}$ is empty, hence by repeating the inequality \eqref{eq:lesbetti} for each critical value, we finally get
\be 
b_i(Z)=b_i(Z\cap g\le a_{n+1})\le N_Wn=N_W\#\text{Crit}\left(g|_{\psi^{-1}(W)}\right).
\ee
This concludes the proof of Theorem \ref{thm:strat}.
\end{proof}
Below we will restrict to those semialgebraic sets $W\subset J^{r}(S^m, \R^k)$ that have a differential geometric meaning, as specified in the next definition.
\begin{defi}\label{def:intrinsic}
A submanifold $W\subset J^r(M,\R^k)$ is said to be \emph{intrinsic} if there is a submanifold $W_0\subset J^r(\mathbb{D}^m,\R^k)$, called the \emph{model}, such that for any embedding $\f\colon \mathbb{D}^m\hookrightarrow M$, one has that $j^r\f^*(W)=W_0$, where 
\be 
j^r\f^*\colon J^r\left(\f(\mathbb{D}^m),\R^k\right)\xrightarrow{\cong}J^r\left(\mathbb{D}^m,\R^k\right), \qquad j^r_{\f(p)}f\mapsto j^r_p(f\circ\f ).
\ee
\end{defi}
Intrinsic submanifolds are, in other words, those that have the same shape in every coordinate charts, as in the following examples.
\begin{enumerate}
    \item $W=\{j^r_pf\colon f(p)=0\}$;
    \item $W=\{j^r_pf\colon j^sf(p)=0\}$ for some $s\le r$;
    \item $W=\{j^r_pf\colon \text{rank}(df(p))=s\}$ for some $s\in\N$.
\end{enumerate}
\begin{remark}\label{rem:strat}
In the case when $J=J^{r}(M,\R^k)$ we can consider $\hat{W}$ to be a subset of $ J^{r+1}(M,\R^{k+1})$ taking the preimage via the natural submersion
\be 
J^{r+1}(M,\R^{k+1})\to J^1\left(M,J^{r}(M,\R^k)\times\R \right), \qquad j^{r+1}(f,g)\mapsto j^1(j^rf,g).
\ee
In this setting Theorem \ref{thm:strat} can be translated to a more natural statement by considering $\psi$ of the form $\psi=j^rf$.
Moreover, in this case, observe that if $W$ is intrinsic (in the sense of Definition \ref{def:intrinsic} below), then $\hat{W}$ is intrinsic as well.
\end{remark}

\subsection{Quantitative bounds}\label{sec:quantitative}
In this section we prove Theorem \ref{thm:bettibound}, which actually immediately follows by combining Theorem \ref{thm:bound} and Theorem \ref{thm:bound2}.


Next theorem gives a deterministic bound for on the complexity of $Z=j^r\psi^{-1}(W)$ when the codimension of $W$ is $m$.
\begin{thm}\label{thm:bound}
Let $P\in \R[x_0, \ldots, x_m]_{(d)}^k$ be a polynomial map and consider its restriction $\psi=P|_{S^m}$ to the unit sphere:
\be \psi:S^m\to \R^k.\ee
Let also $j^r\psi:S^m\to J^r(S^m, \R^k)$ be the associated jet map and $W\subset J^{r}(S^m, R^k) $ be a semialgebraic set of codimension $m$. There exists a constant $c>0$ (which only depends on $W$, $m$ and $k$) such that, if $j^r\psi\transv W$,  then:
\be\label{eq:detbound} \#j^r\psi^{-1}(W)\leq c \cdot d^m.\ee
\end{thm}
\begin{proof}
Let us make the identification $J^{r}(\R^{m+1}, \R^k)\simeq \R^{m+1}\times \R^N$, so that the restricted jet bundle $J^{r}(\R^{m+1}, \R^k)|_{S^m}$ corresponds to the semialgebraic subset $S^m\times \R^N$. Observe that the inclusion $S^m\hookrightarrow \R^{m+1}$ induces a semialgebraic map:
\be J^{r}(\R^{m+1}, \R^k)|_{S^m}\stackrel{i^*}{\longrightarrow}  J^r(S^m, \R^k),\ee
that, roughly speaking, forgets the normal derivatives.
Notice that while the map $j^{r}\psi=j^r(P|_{S^m})$ is a section of $J^r(S^m, \R^k)$,  $(j^rP)|_{S^m}$ is a section of $J^{r}(\R^{m+1}, \R^k)|_{S^m}$. These sections are related by the identity \be 
 i^*\circ (j^rP)|_{S^m}=j^{r}\psi.
 \ee
Thus, defining $ \overline{W}=(i^*)^{-1}(W)$, we have
\be j^r\psi^{-1}(W)=\left((j_rP)|_{S^m}\right)^{-1}(\overline W).\ee
Since $\overline{W}$ is a semialgebraic subset of $\R^{m+1}\times \R^N$ , it can be written as:
\be \overline{W}=\bigcup_{j=1}^\ell\left\{f_{j, 1}=0, \ldots, f_{j, \alpha_j}=0, g_{j, 1}>0,\ldots, g_{j, \beta_j}>0\right\},\ee
where the $f_{j,i}$s and the $g_{j, i}$s are polynomials of degree bounded by a constant $b>0.$ For every $j=1, \ldots, \ell$ we can write:
\be \left\{f_{j, 1}=0, \ldots, f_{j, \alpha_j}=0, g_{j, 1}>0,\ldots, g_{j, \beta_j}>0\right\}=Z_j\cap A_j,\ee
where $Z_j$ is algebraic (given by the equations) and $A_j$ is open (given by the inequalities).

Observe also that the map $(j^rP)|_{S^m}$ is the restriction to the sphere $S^m$ of a polynomial map
\be Q:\R^{m+1}\to \R^{m+1}\times \R^N\ee
whose components have degree smaller than $d$.
Therefore for every $j=1\ldots, \ell$ the set $((j^rP)|_{S^m})^{-1}(Z_j)=(Q|_{S^m})^{-1}(Z_j)$ is an algebraic set on the sphere defined by equations of degree less than $b \cdot d$ and, by \cite[Proposition 14]{LerarioJEMS} we have that:
\be\label{eq:ineqb} b_0(Q|_{S^m})^{-1}(Z_j))\leq B d^m\ee
for some constant $B>0$ depending on $b$ and $m$.
The set $(Q|_{S^m})^{-1}(Z_j)$ consists of several components, some of which are zero dimensional (points):
\be (Q|_{S^m})^{-1}(Z_j)=\underbrace{\{p_{j, 1}, \ldots, p_{j, \nu_j}\}}_{P_j}\cup \underbrace{X_{j,1}\cup\cdots \cup X_{j, \mu_j}}_{Y_j}.\ee
The inequality \eqref{eq:ineqb} says in particular that: \be\label{ineqbb}\#P_j\leq Bd^n.\ee
Observe now that if $j^r\psi\transv W$ then, because the codimension of $W$ is $m$, the set $j^{r}\psi^{-1}(W)=(Q|_{S^m})^{-1}(\overline{W})$ consists of finitely many points and therefore, since $(Q|_{S^m})^{-1}(A_j)$ is open, we must have:
\be j^{r}\psi^{-1}(W)\subset \bigcup_{j=1}^\ell P_j.\ee
(Otherwise $ j^{r}\psi^{-1}(W)$ would contain an open, nonempty set of a component of codimension smaller than $m$.)
Inequality \eqref{ineqbb} implies now that:
\be \#j^{r}\psi^{-1}(W)\leq \sum_{j=1}^\ell \#P_j\leq \ell b d^m\leq cd^m.\ee
\end{proof}

Using Theorem \ref{thm:strat} it is now possible to improve Theorem \ref{thm:bound} to the case of any codimension, replacing the cardinality with any Betti number. \begin{thm}\label{thm:bound2}
Let $P\in \R[x_0, \ldots, x_m]_{(d)}^k$ be a polynomial map and consider its restriction $\psi=P|_{S^m}$ to the unit sphere:
\be \psi:S^m\to \R^k.\ee
Let also $j^r\psi:S^m\to J^r(S^m, \R^k)$ be the associated jet map and $W\subset J^{r}(S^m, R^k) $ be a closed semialgebraic set (of arbitrary codimension). There exists a constant $c>0$ (which only depends on $W$, $m$ and $k$) such that, if $j^r\psi\transv W$,  then:
\be b_i\left(j^r\psi^{-1}(W)\right)\leq c \cdot d^m.\ee
\end{thm}
\begin{proof}
Let $J=J^r(S^m,\R^k)$ and let $\hat{W}$ be the (stratified according to a chosen stratification of $W$) subset of $J^{r+1}(S^m,\R^{k+1})$ coming from Theorem \ref{thm:strat} and Remark \ref{rem:strat}.
Let $g$ be a homogeneous polynomial of degree $d$ such that 
\be 
\Psi=(\psi,g)\in \R[x_0, \ldots, x_m]_{(d)}^{k+1}
\ee
satisfies the condition $j^{r+1}\Psi \transv \hat{W}$ (almost every polynomial $g$ has this property by standard arguments) and $(j^{r}\psi)^{-1}(W)$ is closed in $S^m$, hence compact. Then by Theorem \ref{thm:strat}, there is a constant $N_W$, such that
\be 
b_i\left(j^r\psi^{-1}(W)\right)\le N_W \#\{(j^{r+1}\Psi)^{-1}(\hat{W})\}
\ee
and by Theorem \ref{thm:bound}, the right hand side is bounded by $cd^m$.
\end{proof}
Given $P=(P_1, \ldots, P_k)$ with each $P_i$ a homogeneous polynomial of degree $d$ in $m+1$ variables, we denote by
\be \psi_d:S^m\to \R^k
\ee
its restriction to the unit sphere (the subscript keeps track of the dependence on $d$).
\begin{example}[Real algebraic sets] Let us take $W=S^m\times \{0\}\subset J^{0}(S^m, \R^k),$ then $j^0\psi^{-1}(W)$ is the zero set of $\psi_d:S^m\to \R^k$, i.e. the set of solutions of a system of polynomial equations of degree $d$. In this case the inequality \eqref{eq:detbound} follows from \cite{LerarioJEMS}.
\end{example}
\begin{example}[Critical points] If we pick $W=\{j^1f=0\}\subset J^1(S^m, \R),$ then $Z_d=j^1\psi_d^{-1}(W)$ is the set of critical points of $\psi_d:S^m\to \R$. In $2013$ Cartwright and Sturmfels \cite{CS} proved that 
\be 
\#Z_d\le 2(d-1)^m+\dots+(d-1)+1
\ee
(this bounds follows from complex algebraic geometry), and this estimate was recently proved to be sharp by Kozhasov \cite{khazOFRET}.
Of course one can also fix the index of a nondegenerate critical point (in the sense of Morse Theory); for example we can take $W=\{df=0, d^2f>0\}\subset J^2(S^m, \R),$ and $j^2\psi_d^{-1}(W)$ is the set of nondegenerate \emph{minima} of $\psi_d:S^m\to \R$ (similar estimates of the order $d^{m}$ holds for the fixed Morse index, but the problem of finding a sharp bound is very much open).
\end{example}

\begin{example}[Whitney cusps]When $W=\{\textrm{Whitney cusps}\}\subset J^3(S^2, \R^2),$ then $\psi_d^{3}f^{-1}(W)$ consists of the set of points where the polynomial map $\psi_d:S^2\to \R^2$ has a critical point which is a Whitney cusp. In this case \eqref{eq:detbound} controls the number of possible Whitney cusps (the bound is of the order $O(d^2)$).
\end{example}

\subsection{Semicontinuity of topology under holonomic approximation}\label{sec:holo}
Consider the following setting: $M$ and $J$ are smooth manifolds, $M$ is compact, and $W\subset J$ is a smooth cooriented submanifold. Given a smooth map $F\colon M\to J$ which is transversal to $W$, it follows from standard transversality arguments that there exists a small $\mathcal{C}^1$ neighborhood $U_1$ of $F$ such that for every map $\tilde {F}\in U_1$ the pairs $(M, F^{-1}(W))$ and $(M, \tilde{F}^{-1}(W))$ are isotopic (in particular $F^{-1}(W)$ and $\tilde{F}^{-1}(W)$ have the same Betti numbers, this is the so-called ``Thom's isotopy Lemma''). The question that we address is the behavior of the Betti numbers of $\tilde{F}^{-1}(W)$ under small $\mathcal{C}^0$ perturbations, i.e. how the Betti number can change under modifications of the map $F$ \emph{without} controlling its derivative.  
\begin{figure}\begin{center}
\includegraphics[scale=0.11]{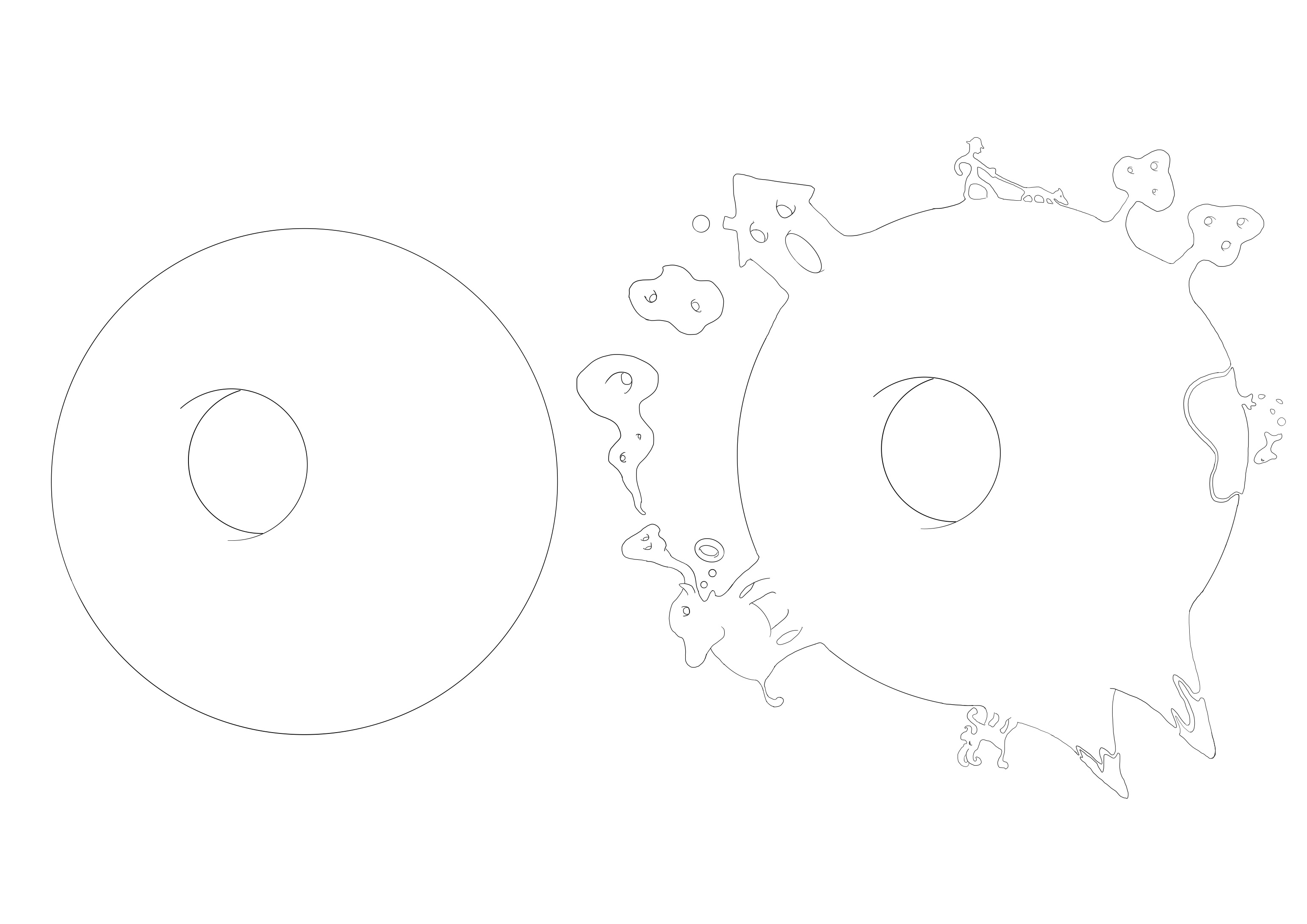}\caption{A small $\mathcal{C}^0$ perturbation of a regular equation can only increase the topology of its zero set.}\label{fig:perturb}
\end{center}\end{figure}
In this direction we prove the following result.
\begin{thm}\label{thm:semiconttop}
Let $M,J$ be smooth manifolds and let $W\subset J$ be a smooth cooriented closed submanifold.
Let $F\colon M\to J$ be a smooth map such that $F\transv W$. If a smooth map $\tilde{F}$ is strongly\footnote{Meaning: in Whitney strong topology. In particular if $C\subset M$ is closed and $U\subset J$ is open, then the set $\{f\in \mathcal{C}^0(M,J)\colon f(C)\subset U\}$ is open, see \cite{Hirsch}.} $\mathcal{C}^0-$close to $F$ such that $\tilde{F}\transv W$, then for all $i\in\N$ there is a group $K^i$ such that
\be \label{eq:semiconttop}
H^i\left(\tilde{F}^{-1}(W)\right)\cong H^i\left(F^{-1}(W)\right)\oplus K^i.
\ee
\end{thm}

\begin{proof}Call $A=F^{-1}(W)$ and $\tilde{A}=\tilde{F}^{-1}(W)$. Let
$E\subset M$ be a closed  tubular neighborhood (it exists because $A$ is closed), meaning that $E=\text{int}(E)\cup\de E$ is diffeomorphic to the unit ball of a metric vector bundle over $A$ (via a diffeomorphism that preserves $A$). Denote by $\pi\colon E\to A$ the retraction map. Since $\tilde{F}$ is $\mathcal{C}^0-$close to $F$ we can assume that there is a homotopy $F_t$ connecting $F=F_0$ and $\tilde{F}=F_1$ such that $F_t(\de E)\subset J\-
W$.
 Define analogously $\tilde{\pi}:\tilde{E}\to \tilde{A}$ in such a way that $\tilde{E}\subset \text{int}(E)$. It follows that there is an inclusion of pairs $u:(E,\de E)\to (E,E\-\tilde{E})$.
By construction, the function $F_t$ induces a well defined mapping of pairs $F_t\colon (E,\de E)\to (J,J\-W)$ for every $t\in[0,1]$, in particular there is a homotopy between $F_0$ and $F_1$ (meant as maps of pairs). Moreover with $t=1$, this map is the composition of $u$ and the map $F_1\colon (E,E\-\tilde{E})\to (J,J\-W)$.

The fact that $W$ is closed and cooriented guarantees the existence of a Thom class $\phi\in H^r(J,J\- W)$, where $r$ is the codimension of $W$. By transversality we have that also $A$ and $\tilde{A}$ are cooriented with Thom classes $F_0^*\phi=\phi_E\in H^r(E,\de E)\cong H^r(E,E\- A)$ and $F_1^*\phi=\phi_{\tilde{E}}\in H^r(\tilde{E},\de\tilde{E})\cong H^r(\tilde{E},\tilde{E}\- \tilde{A})$.
We now claim the commutativity of the diagram below. 
\be 
\begin{tikzcd}
                                                                                               & {H^{*+r}(J, J\backslash W)}\arrow[ldd, "F_1^*"'] \arrow[rdd, "F_1^*=F_0^*"]                         &                                                                  \\
                                                                                               &                                                      &                                                                  \\
{H^{*+r}(\tilde{E}, \de\tilde{E})} \arrow[r, "\eta^{-1}"]  & {H^{*+r}(E, E\backslash \tilde{E})} \arrow[r, "u^*"] & {H^{*+r}(E, \de E)}           \\
H^{*}(\tilde{A}) \arrow[u, "\tilde{\pi}^*(\cdot)\cup \phi_{\tilde{E}}"]                        &                                                      & H^*(A) \arrow[ll, "\pi^*"] \arrow[u, "\pi^*(\cdot)\cup \phi_E"']
\end{tikzcd}\ee
(where $\eta$ is the excision isomorphism).
For what regards the upper triangular diagram, the commutativity simply follows from the fact that all the maps $F_t$ are homotopic and that the excision homomorphism is the inverse of that induced by the inclusion $(E,E\-\tilde{E})\subset (\tilde{E},\de \tilde{E})$.
 To show that the lower rectangle commutes, observe that since $\tilde{\pi}$ is homotopic to the identity of $\tilde{E}$ we have that $\pi\circ\tilde{\pi}$ is homotopic to $\pi|_{\tilde{E}}$. Thus
the commutativity follows from the property of the cup product, saying that for all $\f\in H^*(A)$ we have
\be\begin{aligned}
u^*\circ \eta^{-1}\circ\left(\tilde{\pi}^*\left(\pi|_{\tilde{A}}\right)^*\f\right)\cup \phi_{\tilde{E}}
&=
\left(u^*\circ\eta^{-1}\circ\left(\pi|_{\tilde{E}}\right)^*\f\right)\cup\left( u^*\circ \eta^{-1} \circ F_1^*\phi\right)\\
&=
\pi^*\f\cup \phi_E,
\end{aligned}\ee
where in the last equality we used the identity $u^*\circ \eta^{-1}\circ F_1^*=F_0^*$ implied by the commutativity of the upper triangle.
Since the vertical maps are (Thom) isomorphisms, there exists a homomorphism $U\colon H^*(\tilde{A})\to H^*(A)$ such that $U\circ \pi^*=$id.
\end{proof}
\begin{remark}\label{rem:semiconttop}
The above proof also provides a way to determine how small should the perturbation be. In fact we showed that if $F_t\colon M\to J$ is a homotopy such that $F_1\transv W$ and $F_t(\de E)\subset J\- W$ for all $t\in [0,1]$, where $E$ is a closed tubular neighborhood of $F^{-1}(W)$, then the map $\tilde{F}=F_1$ satisfies \eqref{eq:semiconttop}. Notice that to have such property it is enough that $\tilde{F}\transv W$ and $\tilde{F}|_{\de E}$ is $\mathcal{C}^0-$close to $F|_{\de E}$.  
 This implies that the size of the $\mathcal{C}^0$ neighborhood of $F$ in which the identity \eqref{eq:semiconttop} holds depends only on the restriction of $F$ to a codimension $1$ submanifold.
\end{remark}
\begin{cor}\label{cor:topositive}Let $M$ be a compact manifold of dimension $m$.
Let $W\subset J^r(M,\R^k)$ be a Whitney stratified submanifold of codimension $1\le l\le m$ being transverse to the fibers of the canonical projection $\pi\colon J^r(M,\R^k)\to M$. Then for any number $n \in \N$ there exists a smooth function $\psi\in \mathcal{C}^{\infty}(M,\R^k)$ such that $j^r\psi\transv W$ and 
\be 
b_i\left((j^r\psi)^{-1}(W)\right)\ge n, \quad \forall i=0,\dots, m-l.
\ee 
\end{cor}
\begin{proof}
Let $B\subset J^r(M,\R^k)$ be a small neighbourhood of a regular point $j^r_pf$ of $W$ so that $(B,B\cap W)\cong (\R^{N+l},\R^N\times \{0\})$. Moreover we can assume that there is a neighbourhood $U\cong \R^m$ of $p\in M$ and a commutative diagram of smooth maps
\be\label{eq:diagramabove}
\begin{tikzcd}
                                                         & \R^m\times \R^k\times \{0\} \arrow[rr, hook] &                                                        & \R^m\times\R^k\times\R^l \arrow[d] \\
B\cap W \arrow[rr, hook] \arrow[ru, "\cong" description] &                                              & B \arrow[d, "\pi"'] \arrow[ru, "\cong" description] & \R^m                               \\
                                                         &                                              & U \arrow[ru, "\cong" description]                      &                                   
\end{tikzcd}
\ee
This follows from the fact that $\pi|_{W}$ is a submersion, because of the transversality assumption. 
For any $0\le i\le m-l$ consider the smooth map
\be 
\f_i\colon \R^m\to \R^l, \quad 
u\mapsto \left(\sum_{\ell=1}^{i+1} (u_\ell)^2-1,\sum_{\ell=i+2}^{m} (u_\ell)^2-1 ,u_{m-l+3},\dots,u_m\right)
\ee 
Clearly $0$ is a regular value for $\f_i$, with preimage\footnote{Except for the case  $l=1$. Here one should adjust the definition of $\f_i$ in order to have $b_i(\f_i^{-1}(0))>0$.} $\f_i^{-1}(0)\cong S^i\times S^{m-l-i}$ and it is contained in the unit ball of radius $2$. Let $C\subset\R^m$ be a set of $n(m-l+1)$ points such that $|c-c'|\ge 5$ for all pair of distinct elements $c,c'\in C$. Now choose a partition $C=C_0\sqcup C_1\sqcup\dots C_{m-l}$ in sets of cardinality $n$ and define a smooth map $\f\colon \R^m\to \R^l$ such that $\f(x)=\f_i(x-c)$ whenever $\text{dist}(x,C_i)\le 2$. We may also assume that $0$ is a regular value for $\f$. Notice that $\f^{-1}(0)$ has a connected component 
\be 
S\cong \{1,\dots, n\}\times\left( S^0\times S^{m-l}\sqcup S^1\times S^{m-l-1}\sqcup \dots S^{m-l}\times S^{0}\right).
\ee
Construct a smooth (non necessarily holonomic) section $F\colon U\to J^r(U,\R^k)$ such that $F\transv W$ and such that $F=(u,0,\f)$ on a neighbourhood of $S$, so that $F^{-1}(W)$ still contains $S$ as a connected component, hence $b_i(F^{-1}(W))\ge n$ for all $i=0,\dots, m-l$.

Let $E\subset U$ be a closed tubular neighborhood of $F^{-1}(W)$. 
To conclude we use the holonomic approximation theorem \cite[p. 22]{eliash}, applied to $F\colon U\to J^r(U,\R^k)\cong U\times \R^{k+l}$ near the codimension $1$ submanifold $\de E\subset U$. Such theorem ensures that for any $\e>0$ there exists a diffeomorphism $h\colon U\to U$, an open neighborhood $O_{\de E}\subset U$ of $\de E$ and a smooth function $\psi \colon U\to \R^k$ such that
\be 
\text{dist}_{\mathcal{C}^0}\left((j^rf)|_{h(O_{\de E})},F|_{h(O_{\de E})}\right)<\e,\quad  \text{and }\quad  \text{dist}_{\mathcal{C}^0}\left(h,\mathrm{id}\right)<\e.
\ee
Moreover, we can assume that $j^r\psi\transv W$, by Thom transversality Theorem (see \cite{Hirsch} or \cite{eliash}). 
In particular, it follows that
\be 
\text{dist}_{\mathcal{C}^0}\left((j^rf)\circ h|_{\de E},F|_{\de E}\right)<\left(1+C(F)\right)\cdot \e,
\ee
where $C(F)$ is the lipshitz constant of $F|_{U}$, which can be assumed to be finite (if not, replace $U\cong \R^m$ with an open ball that still contains $F^{-1}(W)$).
Consider the smooth manifold $J=J^r(U,\R^k)$. By the diagram \eqref{eq:diagramabove} it follows that $W\subset J$ is a closed and cooriented smooth submanifold, so that by Theorem \ref{thm:semiconttop} and Remark \ref{rem:semiconttop} we know that if $\e>0$ is small enough, then the map $\tilde{F}=(j^rf)\circ h$ satisfies the identity \eqref{eq:semiconttop}. Therefore for each $i=0,\dots,m-l$, we have
\be 
\begin{aligned}
b_i\left(\left(j^rf\right)^{-1}(W)\right)
&=
b_i\left(\left((j^rf)\circ h\right)^{-1}(W)\right)\\
 &\ge 
b_i\left(\left(F\circ h\right)^{-1}(W)\right) \\
&=
b_i\left(F^{-1}(W)\right)\\
&\ge 
n.
\end{aligned}
\ee
\end{proof}
\section{Random Algebraic Geometry}\label{sec:RAG}
\subsection{Kostlan maps}In this section we give the definition of a random Kostlan polynomial map $P:\R^{m+1}\to \R^{k}$, which is a Gaussian Random Field (GRF) that generalizes the notion of Kostlan polynomial.
\begin{defi}[Kostlan polynomial maps]\label{def:Kostlan}
Let $d,m,k\in \N$. We define the degree $d$ homogeneous Kostlan random map as the measure on $\R[x]_{(d)}^k=\R[x_0,\dots,x_{m}]_{(d)}^k$ induced by the gaussian random polynomial:
\be 
P_d^{m,k}(x)=\sum_{\a\in \N^{m+1},\ |\a|=d} \xi_\a x^\a,
\ee
where $x^\a=x_0^{\a_0}\dots x_m^{\a_m}$ and $\{\xi_\a\}$ is a family of independent gaussian random vectors in $\R^k$ with covariance matrix
\be 
K_{\xi_\a}={d\choose\a}\mathbbm{1}_k=\left(\frac{d!}{\a_0!\dots \a_m!}\right)\mathbbm{1}_k.
\ee
We will call $P_d^{m,k}$ the \emph{Kostlan polynomial} of type $(d,m,k)$ (we will simply write $P_d=P_d^{m,k}$ when the dimensions are understood).
\end{defi}
(In other words, a Kostlan polynomial map $P_d^{m, k}$ is given by a list of $k$ independent Kostlan polynomials of degree $d$ in $m+1$ homogeneous variables.)

There is a non-homogeneous version of the Kostlan polynomial, which we denote as
\be \label{eq:Kd}
p_d(u)=P_d(1,u)=\sum_{\beta\in \N^{m},\ |\beta|\le d} \xi_\beta u^\beta  \in \g \infty{\R^m}k,
\ee
where $u=(u_1,\dots, u_m)\in \R^m$ and $\xi_\beta\sim N\left(0,{d\choose \beta}\mathbbm{1}_k\right)$ are independent. Here we use the notation of Chapter \ref{chap:dtgrf}, where $\g \infty{\R^m}k$ denotes the space of gaussian random field on $\R^m$ with values in $\R^k$ which are $\mathcal{C}^{\infty}$.
Next Proposition collects some well known facts on the Kostlan measure.
\begin{prop}\label{covkostlanprop}Let $P_d$ be the Kostlan polynomial of type $(d,m,k)$ and $p_d$ be its dehomogenized version, as defined in \eqref{eq:Kd}.
\begin{enumerate}
\item
 For every $x,y\in \R^{m+1}$:
\be 
K_{P_d}(x,y)=\left(x^Ty\right)^d\mathbbm{1}_k.
\ee
Moreover, given $R\in O(m+1)$ and $S\in O(k)$ and defined the polynomial $\tilde{P}_d(x)=SP_d(Rx)$, then $P_d$ and $\tilde{P}_d$ are equivalent\footnote{Two random fields are said to be \emph{equivalent} if they induce the same probability measure on $\Cr \infty{\R^m}k$.}.

\item  For every $u, v\in \R^n$
\be K_{p_d}(u,v)=(1+u^Tv)^d\mathbbm{1}_k.\ee
Moreover, if $R\in O(m)$ and $S\in O(k)$ and defined the polynomial $\tilde{p}_d(x)=Sp_d(Rx)$, then $p_d$ and $\tilde{p}_d$ are equivalent.
\end{enumerate}
\end{prop}
\begin{proof}The proof of this proposition simply follows by computing explicitly the covariance functions and observing that they are invariant under orthogonal change of coordinates in the target and the source. For example, in the case of $P_{d}$ we have:
\be
\begin{aligned}
K_{P_d}(x,y)&= \E\{P_d(x)P_d(y)^T\}= \\
&= \sum_{|\a|,|\a '|=d}\E\left\{\xi_\a\xi_{\a '}^T\right\} x^\a y^{\a'}=\\
&=\sum_{|\a|=d}{d\choose\a} (x_0 y_0)^{\a_0}\dots(x_m y_m)^{\a_m}\mathbbm{1}_k=\\
&= (x_0y_0+\dots +x_my_m)^d\mathbbm{1}_k,
\end{aligned}
\ee
from which the orthogonal invariance is clear. The case of $p_d$ follows from the identity:
\be 
K_{p_d}(u,v)=K_{P_d}\left((1,u),(1,v)\right).
\ee .
\end{proof}
\subsection{Properties of the rescaled Kostlan}

The main feature here is the fact that the local model of a Kostlan polynomial has a rescaling limit. The orthogonal invariance is used to prove that the limit does not depend on the point where we center the local model, hence it is enough to work around the point $(1, 0, \ldots, 0)\in S^m$. These considerations lead to introduce the Gaussian Random Field $X_d:\R^m\to\R^{k}$ (we call it the \emph{rescaled Kostlan}) defined by:
\be \label{eq:rescaledKostlan}X_d(u)=P_d^{m,k}\left(1, \frac{u_1}{\sqrt{d}}, \ldots, \frac{u_m}{\sqrt{d}}\right).\ee
Next result gives a description of the properties of the rescaled Kostlan polynomial, in particular its convergence in law as a random element of the space of smooth functions, space which, from now, on we will always assume to be endowed with the weak Whitney's topology as in Chapter \ref{chap:dtgrf}.

\begin{thm}[Properties of the rescaled Kostlan]\label{thm:Kostlan}
Let $X_d:\R^m\to \R^{k}$ be the Gaussian random field defined in \eqref{eq:rescaledKostlan}. 
\begin{enumerate}[1.]
\item \emph{(The limit)} Given a family of independent gaussian random vectors $\xi_\beta \sim N\left(0,\frac{1}{\beta !}\mathbbm{1}_k\right)$, the series 
\be 
X_\infty(u)=\sum_{\beta\in \N^{m}} \xi_\beta u^\beta  ,
\ee
is almost surely convergent in $\Cr \infty {\R^m}k$ to the Gaussian Random Field\footnote{$X_\infty$ is indeed a random analytic function, commonly known as the Bargmann-Fock ensemble.} $X_\infty\in \g \infty{\R^m}k$.
\item \emph{(Convergence)} $X_d\nrw X_\infty$ in $\g \infty {\R^m}k$, that is:
\be 
\lim_{d\to +\infty}\E\{F(X_d)\}=\E\{F(X_\infty)\}
\ee 
for any bounded and continuous function $F\colon \Cr \infty{\R^m}k\to \R$. Equivalently, we have
\be 
\label{eq:inequKostlan}\P\{X_\infty \in \emph{int}(A)\}\le \liminf_{d\to +\infty}\P\{X_d\in A\}\le \limsup_{d\to +\infty}\P\{X_d\in A\}\le \P\{X_\infty\in \overline{A}\}
\ee
for any Borel subset $A\subset \Cr \infty {\R^m}k$.
\item \emph{(Nondegeneracy of the limit)} The support of $X_\infty$ is the whole $\Cr \infty {\R^m}k$. In other words, for any non empty open set $U\subset \Cr\infty{\R^m}k$ we have that $\P\{X_\infty \in U\}>0$.
\item \emph{(Probabilistic Transversality)} For $d\ge r$ and $d=\infty$, we have $\mathrm{supp}(j^r_pX_d)=J_p^r(\R^m,\R^k)$ for every $p\in \R^m$ and consequently for every submanifold $W\subset J^r(\R^m,\R^k)$, we have 
\be 
\P\{j^rX_d\transv W\}=1.
\ee
\item \emph{(Existence of limit probability)} Let $V\subset J^{r}(\R^m, \R^k)$ be an open set whose boundary is a (possibly stratified) submanifold\footnote{For example $V$ could be a semialgebraic set}. Then
\be 
\lim_{d\to +\infty}\P\{j^r_p X_d\in V,\ \forall p\in \R^m\}= \P\{j^r_pX_\infty(\R^m) \in V,\ \forall p\in \R^m\}.
\ee
In other words, we have equality in \eqref{eq:inequKostlan} for sets of the form $U=\{f\colon j^rf\in V\}$.
\item \emph{(Kac-Rice densities)} Let $W\subset J^r(\R^m,\R^k)$ be a semialgebraic subset of codimension $m$, such that\footnote{In this paper the symbol $\transv$ stands for ``it is transverse to''.} $W\transv J_p^r(\R^m,\R^k)$ for all $p\in M$ (i.e. $W$ is transverse to fibers of the projection of the jet space). Then for all $d\ge r$ and for $d=+\infty$ there exists a locally bounded function $\rho_d^W\in L^\infty_{loc} (\R^m)$ such that\footnote{A formula for $\rho_d^W$ is presented in Chapter \ref{chap:kr}, as a generalization of the classical Kac-Rice formula.} 
\be 
\E\#\{u\in A\colon j^r_uX_d\in W\}=\int_A\rho^W_d,
\ee
for any Borel subset $A\subset \R^m$. Moreover $\rho^W_d  \to \rho^W_\infty$ in $L^\infty_{loc}$.
\item \label{thm:Ebetti}\emph{(Limit of Betti numbers)} Let $W\subset J^r(\R^m,\R^k)$ be any closed semialgebraic subset transverse to fibers. Then:
\be \label{eq:localEbetti}
\lim_{d\to +\infty}\E\left\{b_i\left((j^{r}X_d)^{-1}(W)\cap\mathbb{D}^m\right)\right\}=\E\left\{b_i\left((j^{r}X_\infty)^{-1}(W)\cap \mathbb{D}^m\right)\right\},
\ee
where $b_i(Z)=\dim H_i(Z,\R)$. Moreover, if the codimension of $W$ is $l\ge 1$, then the r.h.s. in equation \eqref{eq:localEbetti} is strictly positive for all $i=0,\dots, m-l$.
\end{enumerate}
\end{thm}
\begin{proof} The proof uses a combination of results from Chapters \ref{chap:dtgrf} and \ref{chap:kr}.
\begin{enumerate}[$(1)$]
\item Let $S_d=\sum_{|\beta|\le d}\xi_\beta u^\beta\in \g \infty Mk$. The covariance function of $S_d$ converges in Whitney's weak topology: \be K_{S_d}(u,v)=\sum_{|\beta|\le d}\frac{u^\beta v^\beta}{\beta !}\mathbbm{1}_k\xrightarrow{\mathcal{C}^\infty} \exp(u^Tv)\mathbbm{1}_k.\ee 
It follows by Theorem \ref{thm:1} of Chapter \ref{chap:dtgrf} that $S_d$ converges in $\g \infty Mk$, moreover since all the terms in the series are independent we can conclude with the Ito-Nisio \footnote{It may not be trivial to apply the standard Ito-Nisio theorem, which actually regards convergence of series in a Banach space. See Theorem \ref{converepr} of Chapter \ref{chap:dtgrf} for a statement that is directly applicabile to our situation} Theorem \cite{ItoNisio} that indeed the convergence holds almost surely.
\item By Theorem \ref{thm:2} of Chapter \ref{chap:dtgrf} it follows from convergence of the covariance functions:
\be 
K_{X_d}(u,v)=\left(1+\frac{u^Tv}{d}\right)^d\mathbbm{1}_k \quad \xrightarrow{\mathcal{C^\infty}}\quad  K_{X_\infty}(u,v)=\exp(u^Tv)\mathbbm{1}_k
\ee
\item The support of $X_\infty$ contains the set of polynomial functions $\R[u]^k$, which is dense in $\Cr \infty {\R^m}k$, hence the thesis follows from Theorem \ref{thm:dtgrf:3} of Chapter \ref{chap:dtgrf}.
\item Let $d\ge r$ or $d=+\infty$. We have that
\be 
\begin{aligned}
\spt(j_u^rX_d)&=\{j^r_uf\colon f\in \R[u]^k \text{ of degree $\le d$}\}=\\
&=\textrm{span}\{j^r_uf\colon f(v)=(v-u)^\beta \text{ with $|\beta|\le d$}\}=\\
&=\textrm{span}\{j^r_uf\colon f(v)=(v-u)^\beta \text{ with $|\beta|\le r$}\}=\\
&=J^r_u(\R^m,\R^k).
\end{aligned}
\ee
The fact that $\P\{j^rX_d\transv W\}=1$ follows from Theorem \ref{transthm2} of Chapter \ref{chap:dtgrf}.
\item Let $A=\{f\in \Cr \infty{\R^m}k\colon j^rf\in V\}$. If $f\in \de A$, then $j^rf\in \overline{V}$ and there is a point $u\in \R^m$ such that $j^r_uf\in \de V$. Let $\partial V$ be stratified as $\partial V=\coprod Z_i$ with each $Z_i$ a submanifold. If $j^rf\transv \de V$ then it means that $j^r f$ is transversal to all the $Z_i$ and there exists one of them which contains $j^r_uf$ (i.e. the jet of $f$ intersect $\partial V$). Therefore the intersection would be transversal \emph{and nonempty}, and then there exists a small Whitney-neighborhood of $f$ such that for every $g$ in this neighborhood $j^rg$ still intersects $\partial V.$ This means that there is a neighborhood of $f$ consisting of maps that are not in $A$, which means $f$ has a neighborhood contained in $A^c$. It follows that $f\notin \overline{A}$ and consequently $f\notin \partial A$, which is a contradiction. Therefore we have that
\be 
\de A\subset \{f\in \Cr \infty{R^m}k\colon f \text{ is not transverse to }\de V\}.
\ee
It follows by point $(4)$ that $\P\{X\in \de A\}=0$, so that we can conclude by points $(2)$ and $(3)$.
\item By previous points, we deduce that we can apply the results of Chapter \ref{chap:kr}.
\item This proof is postponed to Section \ref{sec:betti}.
\end{enumerate} 
\end{proof}
Given  a $\mathcal{C}^{\infty}$ Gaussian Random Field $X:\R^m\to \R^k$ , let us denote by $[X]$ the probability measure induced on $\mathcal{C}^{\infty}(\R^m, \R^k)$ and defined by:
\be [X](U)=\mathbb{P}(X\in U),\ee
for every $U$ belonging to the Borel $\sigma-$algebra relative to the weak Whitney topology, see \cite{Hirsch} for details on this topology.
Combining Theorem \ref{thm:Kostlan} with Skorohod Theorem \cite[Theorem 6.7]{Billingsley} one gets that it is possible to represent $[X_d]$ with equivalent fields $\tilde{X}_d$ such that $\tilde{X}_d\to \tilde{X}_\infty$ almost surely in $\Cr\infty{\R^m}k$. This is in fact equivalent to point $(2)$ of Theorem \ref{thm:Kostlan}. In other words there is a (not unique) choice of the gaussian coefficients of the random polynomials in $\eqref{eq:Kd}$, for which the covariances $\E\{\tilde{X}_d\tilde{X}_{d'}^T\}$ are such that the sequence converges almost surely. We leave to the reader to check that a possible choice is the following. 
Let $\{\gamma_\beta\}_{\beta\in\N^m}$ be a family of i.i.d. gaussian random vectors $\sim N(0,\mathbbm{1}_k)$ and define for all $d<\infty$
\be \label{eq:askocond}
\tilde{X}_d=\sum_{|\beta|\le d}{d\choose \beta}^\frac12 \gamma_\beta\left(\frac{u}{\sqrt{d}}\right)^{\beta}
\ee
and
\be\label{eq:askoconf}
\tilde{X}_\infty=\sum_{\beta}{\left(\frac{1}{\beta!}\right)}^\frac12 \gamma_\beta u^{\beta}.
\ee
\begin{prop}
$\tilde{X}_d\to \tilde{X}_{\infty}$ in $\Cr \infty{\R^m}k$ almost surely.
\end{prop}
However, we stress the fact that in most situations: when one is interested in the sequence of probability measures $[X_d]$, it is sufficient to know that such a sequence exists.
\subsection{Limit laws for Betti numbers and the generalized square-root law}\label{sec:betti}
\begin{figure}[h]
\begin{center}
\includegraphics[scale=0.15]{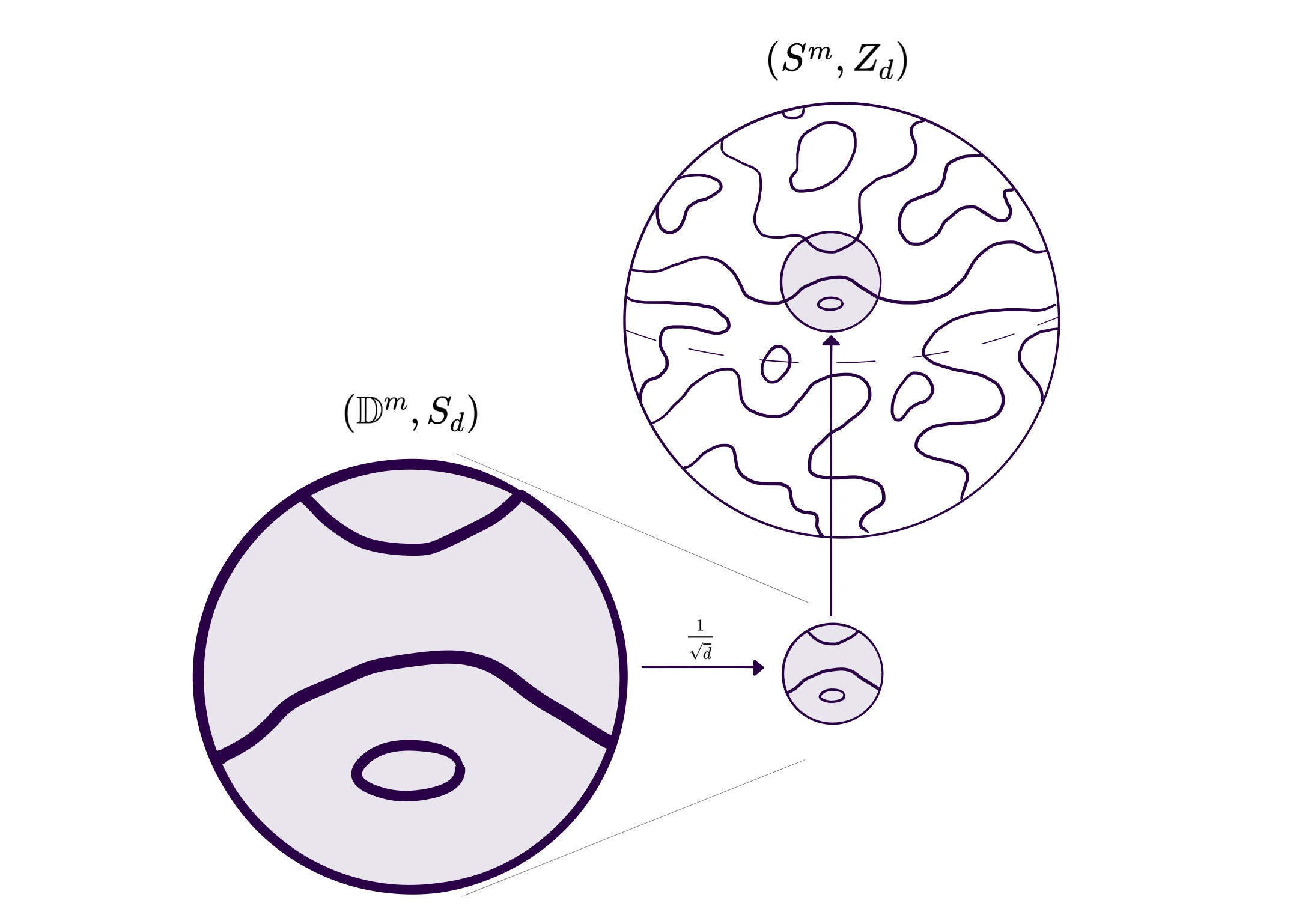}\caption{The random set $S_d=\{X_d=0\}\subset \D^m$ is a rescaled version of $Z_d\cap D(p, d^{-1/2})$, where $Z_d=\{\psi_d=0\}$.}\label{fig:rescale}
\end{center}
\end{figure}

Let $W_0\subset J^r(\R^m,\R^k)$ be a semialgebraic subset.
Consider the random set
\be S_d= \{p\in \mathbb{D}^m\colon j_p^rX_d\in W_0\},\ee
where $X_d\colon \R^m\to \R^k$ is the rescaled Kostlan polynomial from Theorem \ref{thm:Kostlan} (see Figure \ref{fig:rescale}).
We are now in the position of complete the proof of Theorem \ref{thm:Kostlan} by showing point (7).
Let us start by proving the following Lemma.
\begin{lemma}\label{lem:positibetti}
Let $r$ be the codimension of $W_0$ and suppose $0\le i\le m-r\le m-1$. Then \be \E\{b_i(S_\infty)\}>0 .\ee
\end{lemma}
\begin{proof}
From Corollary \ref{cor:topositive} we deduce that there exists a function $f\in \Cr \infty{\mathbb{D}^m} k$ such that $j^rf\transv W_0$ and $b_i\left((j^rf)^{-1}(W_0)\right)\neq 0$.
Since the condition on $f$ is open, there is an open neighbourhood $O$ of $f$ where  $b_i((j^rg)^{-1}(W_0))=c>0$ for all $g\in O$. Thus $\P\{b_i(S_\infty)=c\}>0$ because every open set has positive probability for $X_\infty$, by \ref{thm:Kostlan}.3. therefore $\E\{b_i(S_\infty)\}>0$. 
\end{proof}
We complete the proof of Theorem \ref{thm:Kostlan} with the next Proposition.
\begin{prop}
\be\label{eq:Ebetti}
\lim_{d\to \infty}\E\{b_i(S_d)\}=\E\{b_i(S_\infty)\}.
\ee
\end{prop}
\begin{proof}
Let $b_i(S_d)=b_d$.
Define a random field $Y_d=(X_d,x_d)\colon \R^m\to \R^{k}\times\R$ to be the rescaled Kostlan polynomial of type $(m,k+1)$.
Consider the semialgebraic subset $W'=W\cap J^r(\mathbb{D}^m,\R^k)$ of the real algebraic smooth manifold $J^r(\R^m,\R^k)$ and observe that $S_d=(j^rX_d)^{-1}(W')$ is compact. Now Theorem \ref{thm:strat}, along with Remark \ref{rem:strat}, implies the existence of a semialgebraic submanifold $\hat{W'}\subset J^{r+1}(\R^m,\R^{k+1})$ of codimension $m$ and a constant $C$, such that 
\be 
b_d\le C\#\left\{\left(j^{r+1}(Y_d)\right)^{-1}(\hat{W'})\right\}=:N_d
\ee
whenever $j^rX_d \transv W'$ and $j^{r+1}Y_d\transv \hat{W'}$, hence almost surely, because of Theorem \ref{thm:Kostlan}.4.
Since $Y_d\nrw Y_\infty$ by \ref{thm:Kostlan}.2, we see that $[b_d,N_d]\nrw [b_\infty,N_\infty]$ and it is not restrictive to assume that $(b_i,N_d)\to (b_i,N_\infty)$ almost surely, by Skorokhod's theorem (see \cite[Theorem 6.7]{Billingsley}). Moreover $\E\{N_d\}\to \E\{N_\infty\}$ by Theorem \ref{thm:Kostlan}.6.
Now we can conclude with Fatou's Lemma as follows 
\be 
\begin{aligned}
2\E\{N_\infty\}&=\E\{\liminf_d N_d+N_\infty-|b_d-b_\infty|\}\le \\
&\le \liminf_d\E\{ N_d+N_\infty-|b_d-b_\infty|\}=\\
&=
2\E\{N_\infty\}-\limsup_d\E\{|b_d-b_\infty|\},
\end{aligned}
\ee 
so that
\be 
\limsup_d\E\{|b_d-b_\infty|\}\le 0.
\ee
\end{proof}


In the sequel, with the scope of keeping a light notation, for a given $W\subset J^{r}(S^m, \R^k)$ and $\psi:S^m\to \R^k$ we will denote by $Z_d\subseteq S^m$ the set
\be Z_d=j^r\psi^{-1}(W).\ee
If $W$ is of codimension $m$, then by Theorem \ref{thm:Kostlan}, $Z_d$ is almost surely a finite set of points and the expectation of this number is given by next result.
\begin{thm}[Generalized square-root law for cardinality]\label{thm:sqrlaw}
Let $W\subset J^r(S^m,\R^k)$ be a semialgebraic intrinsic subset of codimension $m$. Then there is a constant $C_W>0$ such that:
\be 
\E\{\#Z_d\}=C_W d^{\frac{m}{2}}+ O(d^{\frac{m}{2}-1}).
\ee
Moreover, the value of  $C_W$ can be computed as follows. Let $Y_\infty=e^{-\frac{|u|^2}{2}}X_\infty\in \g \infty{\mathbb{D}^m}k$ and let $W_0\subset J^r(\mathbb{D}^m,\R^k)$ be the local model for $W$. Then
\be 
C_W=m\frac{\vol (S^m)}{\vol (S^{m-1})} \E\#\{u\in\mathbb{D}^m\colon j^r_uY_{\infty}\in W_0\}.
\ee
\end{thm}

In order to prove Theorem \ref{thm:sqrlaw}, we will need a preliminary Lemma, which ensures that we will be in the position of using the generalized Kac-Rice formula of point (6) from Theorem \ref{thm:Kostlan}.
\begin{lemma}\label{lem:intrilinstab}
If $W\subset J^r(M,\R^k)$ is intrinsic, then $W$ is transverse to fibers.
\end{lemma}
\begin{proof}
Since the result is local it is sufficient to prove it in the case when $M=\R^m$. 
In this case we have a natural identification (see \cite[Chapter 2, Section 4]{Hirsch})

For any point $u\in \R^m$ we consider the embedding $ i_u \colon \mathbb{D}^m\to\R^m$ obtained as the isometric inclusion in the disk with center $u$ and let $\tau_u\colon \R^m\to \R^m$ be the translation map $x\mapsto u+x$.
Let $u,v\in \R^m$ be two points with distance smaller than $1$, he fact that the submanifold $W$ is intrinsic implies that $j^r_vf\in W$ if and only if $(j^ri_u)^*(j^r_vf)\in W_0$, where $W_0\subset J^r(\mathbb{D}^m,\R^k)$ is the model for $W$. From this we deduce that also the jet $j^r_u(f\circ \tau_{v-u})$ is in $W$, since:
\be \begin{aligned}
\left(j^ri_v\right)^*\left(j^r_{v}f\right)&=
j^r(\tau_{v-u}\circ i_u)^*(j^r_vf)
\\
&= 
\left(j^ri_u\right)^*\left(j^r\left(\tau_{v-u}\right)^*\left(j^r_{\tau_{v-u}(u)}f\right)\right)
\\
&= 
\left(j^ri_u\right)^*\left(j^r_u\left(f\circ\tau_{v-u}\right)\right).
\end{aligned}
\ee
By interchanging the role of $u$ and $v$, we conclude that $j^r_u(f\circ\tau_{v-u})\in W$ if and only if $j^r_vf\in W$. Notice that such statement is thus true for any couple of points $u,v\in\R^m$, regardless of their distance.

We thus claim that $T(W)$ is of the form $\R^m\times \bar{W}$, under the natural identification (see \cite[Sec. 2.4]{Hirsch}):
\be 
T\colon J^r(\R^m,\R^k)\cong \R^m\times J^r_0(\R^m,\R^k), \qquad j^r_uf\mapsto (u,j^r_0(f\circ \tau_{u})).
\ee
To see this, observe that if $(v,j^r_0g)\in T(W)$, hence $(v,j^r_0g)= T(j^r_vf)$ for a jet $j^r_vf\in W$ such that $ 
g=f\circ \tau_{v}$,
then $(u,j^r_0g)=T\left(j^r_u(f\circ \tau_{v-u})\right)\in T(W)$.
\end{proof}
\begin{figure}
\begin{center}
\includegraphics[scale=0.11]{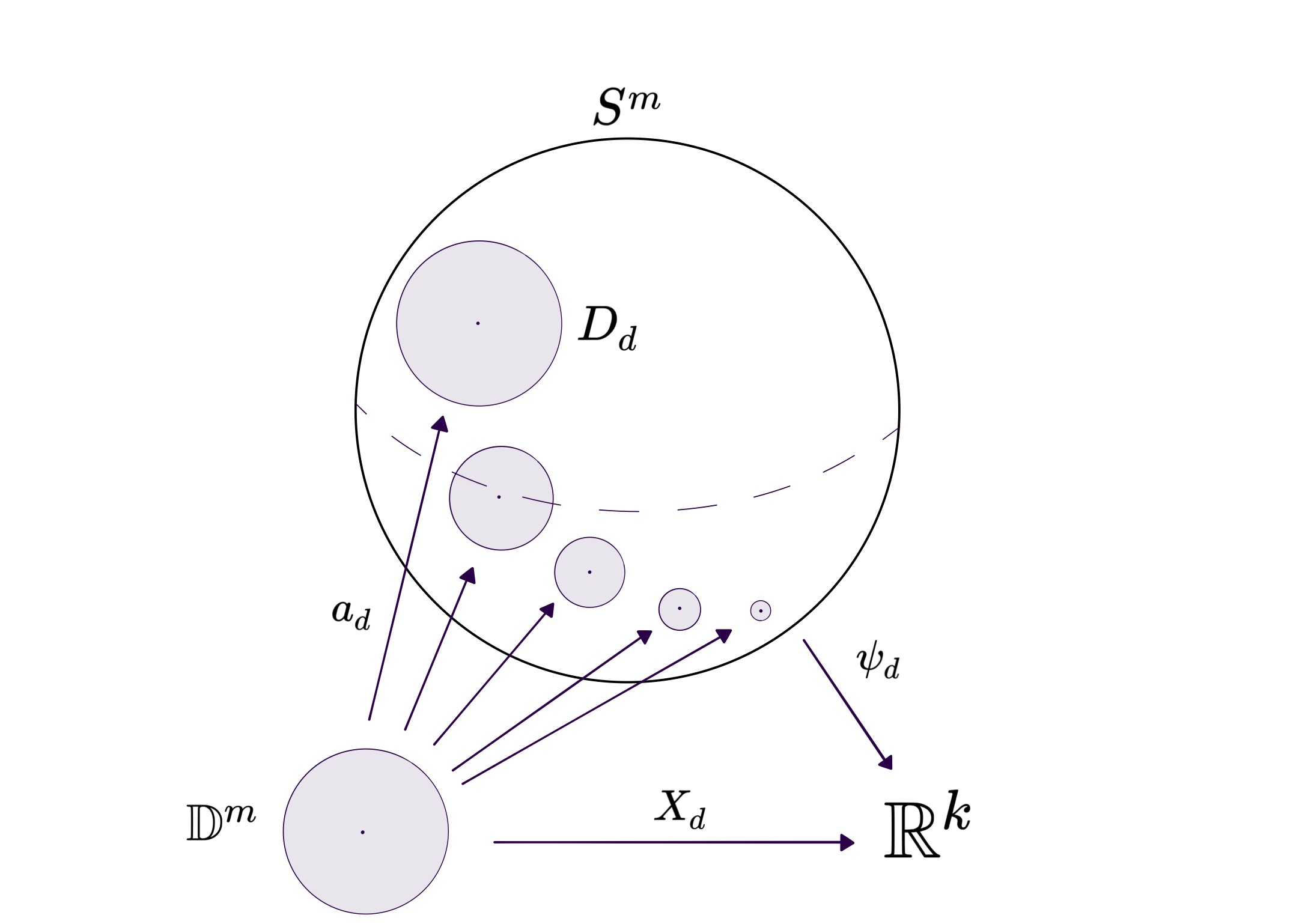}\caption{A family of shrinking embedding of the unit disk.}\label{fig:emb}
\end{center}
\end{figure}
The reason why we consider intrinsic submanifold is to be able to easily pass to the rescaled Kostlan polynomial $X_d\in \g \infty{\mathbb{D}^m}k$ by composing $\psi_d$ with the embedding of the disk $a_d^R$ defined by: 
\be\label{eq:diskemb}
a_d^R\colon \mathbb{D}^m\hookrightarrow S^m,\quad 
u\mapsto \frac{R\begin{pmatrix}
1 \\ \frac{u}{\sqrt{d}}
\end{pmatrix}}{\sqrt{\left(1+\frac{|u|^2}{d}\right)}}
\ee
for any $R\in O(m+1)$ (see Figure \ref{fig:emb}).
\begin{proof}[Proof of Theorem \ref{thm:sqrlaw}]
Let us consider the set function $\mu_d\colon \mathcal{B}(S^m)\mapsto \R$ such that $A\mapsto \E\{\#(j^rX_d)^{-1}(W)\cap A\}$. It is explained in Chapter \ref{chap:dtgrf} that $\mu_d$ is a Radon measure on $S^m$. Because of the invariance under rotation of $P_d$, by Haar's theorem $\mu$ needs to be proportional to the volume measure. Therefore for any Borel subset $A\subset S^m$ we have $\E\{\#Z_d\}=\mu_d(S^m)=\mu_d(A)\vol(A)^{-1}\vol(S^m)$. 
Define $Y_d\in \g \infty {\mathbb{D}^m}k$ as
\be\label{eq:Y}
Y_d= \left(1+\frac{|u|^2}{d}\right)^{-\frac{d}2}X_d.
\ee
Observe that $Y_d\nrw Y_\infty=\exp(-\frac{|u|^2}{2})X_\infty$ and that $Y_d$ is equivalent to the GRF $\psi_d\circ a^R_d$ for any $R\in O(m+1)$.

Now let $W_0\subset J^r(\mathbb{D}^m,\R^k)$ be the (semialgebraic) model of $W$. By the same proof of point (7) from Theorem \ref{thm:Kostlan}, adapted to $Y_d$, there is a convergent sequence of functions $\rho_d\to \rho_{+\infty}\in L^1(\mathbb{D}^m)$ such that 
\be 
\E\{\#(j^rY_d)^{-1}(W_0)\}=\int_{\mathbb{D}^m}\rho_d\to \int_{\mathbb{D}^m}\rho_\infty= \E\{\#(j^rY_\infty)^{-1}(W_0)\}.
\ee
In conclusion we have for $A=a^R_d(\mathbb{D}^m)$, as $d\to +\infty$
\be 
\begin{aligned}
\E\{\#Z_d\}&=\mu_d(A)\vol(A)^{-1}\vol(S^m)\\
&=\E\{\#(j^rY_d)^{-1}(j^r\f^*(W))\}\vol(A)^{-1}\vol(S^m)\\
&=\E\{\#(j^rY_d)^{-1}(W_0)\}\left(\frac{\int_0^{\pi}|\sin\theta|^{m-1}d\theta}{\int_0^{\arctan\left(d^{-\frac12}\right)}|\sin\theta|^{m-1}d\theta} \right)\\
&=\E\{\#(j^rY_\infty)^{-1}(W_0)\}m\frac{\vol(S^m)}{\vol(S^{m-1})}d^{\frac{m}{2}}+O(d^{\frac{m}{2}-1}).
\end{aligned}
\ee

\end{proof}
Building on the previous results, we can now prove the general case for Betti numbers of a random singualrity.
\begin{thm}[Generalized square-root law for Betti numbers]\label{thm:bettiorder}
Let $W\subset J^r(S^m,\R^k)$ be a closed semialgebraic intrinsic (as defined in Definition \ref{def:intrinsic}) of codimension $1\le l\le m$. Then there are constants $b_W, B_W > 0$ depending only on $W$ such that
\be \label{eq:ineqbetti}
b_W d^{\frac{m}{2}}\le\E\{b_i(Z_d)\}\le B_W d^{\frac{m}{2}}\quad \forall i=0,\dots, m-l
\ee
and $\E\{b_i(Z_d)\}=0$ for all other $i$.
\end{thm}

\begin{proof}The proof is divided in two parts, first we prove the upper bound, using the square-root law from Theorem \ref{thm:sqrlaw}, then the we use Theorem \ref{thm:Ebetti} to deduce the lower bound. The globalization step for the lower bound is a generalization of the so-called ``barrier method'' from \cite{NazarovSodin2, GaWe1}.

1. Assume $W$ is smooth with codimension $s$. Let us consider 
\be P^{m, k+1}_d|_{S^m}=\Psi_d=(\psi_d,\psi_d^1)\in \g \infty {S^m}{k+1}
\ee
and Let $\hat{W}\subset J^{r+1}(S^m,\R^{k+1})$ be the intrinsic semialgebraic submanifold coming from Theorem \ref{thm:strat} and Remark \ref{rem:strat}. Thus, using Theorems \ref{thm:strat} and \ref{thm:sqrlaw}, we get
\be 
\E\{b_i(Z_d)\}\le N_W\E\#\{(j^{r+1}\Psi_d)^{-1}(\hat{W})\}\le N_WC_{\hat{W}}d^{\frac{m}{2}}.
\ee

2. Consider the embeddings of the $m$ dimensional disk $a_d^R\colon \mathbb{D}^m\hookrightarrow S^m$ defined in \eqref{eq:diskemb}. For any fixed $d\in\N$, choose a finite subset $F_d\subset O(m+1)$ such that the images of the corresponding embeddings $\{a_d^R(\mathbb{D}^m)\}_{R\in F_d}$ are disjoint. Denoting by $Z_d^R$ the union of all connected components of $Z_d$ that are entirely contained in $a_d^R(\mathbb{D}^m)$, we have
\be 
b_i(Z_d)\ge \sum_{R\in F_d}b_i(Z_d^R).
\ee
Let $W_0\subset J^r(\mathbb{D}^m,\R^k)$ be the model of $W$ as an intrinsic submanifold, it is closed and semialgebraic. By Definition \ref{def:intrinsic}, we have
\be\label{eq:pezzi}
(a_d^R)^{-1}\left((j^r\psi_d)^{-1}(W)\right)= \left(j^r(\psi_d\circ a_d^R)\right)^{-1}(W_0)\subset \mathbb{D}^m.
\ee
Recall that for any $R\in O(m+1)$, the GRF $\psi_d\circ a^R_d$ is equivalent to $Y_d\in \g \infty{\mathbb{D}^m}k$ defined in \ref{eq:Y}, hence taking expectation in Equation \eqref{eq:pezzi} we find
\be 
\E \{b_i(Z_d)\}\ge \#(F_d)\E \{b_i(S_d)\},
\ee
where $S_d=\left(j^r(Y_d)\right)^{-1}(W_0)$.
 is easy to see (repeating the same proof) that Theorem \ref{thm:Kostlan}.7 holds also for the sequence $Y_d\nrw Y_\infty$, so that $\E\{S_d\}\to \E\{S_\infty\}$. We can assume that $\E\{S_\infty\}>0$, because of Lemma \ref{lem:positibetti}, thus for big enough $d$, the numbers $\E\{b_i(S_d)\}$ are bounded below by a constant $C>0$.
Now it remains to estimate the number $\#(F_d)$. Notice that $a_d^R(\mathbb{D}^m)$ is a ball in $S^m$ of a certain radius $\e_d$, hence it is possible to choose $F_d$ to have at least $N_m \e_d^{-1}$ elements, for some dimensional constant $N_m>0$ depending only on $m$. We conclude by observing that 
\be 
\e_d\approx d^{-\frac{m}{2}}.
\ee
\end{proof}
\begin{subappendices}
\section{Examples of applications of Theorem \ref{thm:Kostlan}}\label{sec:examples}
\begin{example}[Zero sets of random polynomials]\label{ex:welsch}Consider the zero set $Z_d\subset \R\P^m$ of a random Kostlan polynomial $P_d=P_{d}^{m+1,1}$. Recently Gayet and Welschinger \cite{GaWe1} have proved that given a compact hypersurface $Y\subset \R^{m}$ there exists a positive constant $c=c(\R^m, Y)>0$ and $d_0=d_0(\R^m,Y)\in \N$ such that for every point $x\in \R\P^m$ and every large enough degree $d\ge d_0$, denoting by $B_d$ any open ball of radius $d^{-1/2}$ in $\R\P^m$, we have:
\be \label{eq:isotopic}\left(B_d, B_d\cap Z_d\right)\cong(\R^m, Y) \ee
(i.e. the two pairs are diffeomorphic) with probability larger than $c$. This result follows from Theorem \ref{thm:Kostlan} as follows. Let $\mathbb{D}^m\subset \R^m$ be the unit disk, and let $U\subset \mathcal{C}^{\infty}(\mathbb{D}^m, \R)$ be the open set consisting of functions $g:\mathbb{D}^m\to \R$ whose zero set is regular (an open $\mathcal{C}^1$ condition satisfied almost surely by $X_d$, because of point (4)), entirely contained in the interior of $\mathbb{D}^m$ (an open $\mathcal{C}^0$ condition) and such that, denoting by $\B\subset \R^m$ the standard unit open ball,  the first two conditions hold and $(\mathbb{B}, \mathbb{B}\cap\{g=0\})$ is diffeomorphic to $(\R^m, Y)$ (this is another open $\mathcal{C}^1$ condition).
Observe that, using the notation above:
\be \left(B_d, B_d\cap Z_d\right)\sim (\mathbb{B}, \mathbb{B}\cap\{X_d=0\}) \ee
(this is simply because $X_d(u)=P_d(1, ud^{-1/2})$). Consequently point (5) of Theorem \ref{thm:Kostlan} implies that:
\begin{align}
\lim_{d\to +\infty}\P\{\eqref{eq:isotopic}\}&=\lim_{d\to \infty} \P\left\{(\mathbb{B}, \mathbb{B}\cap\{X_d=0\})\sim (\R^m, Y)\right\}\\
&=\lim_{d\to \infty} \P\left\{X_d\in U\right\}\\
&=\P\left\{X_\infty \in U \right\}>0.
\end{align}
We stress that, as an extra consequence of Theorem \ref{thm:Kostlan}, compared to \cite{GaWe1} what we get is the existence of the limit of the probability of seeing a given diffeomorphism type.
\end{example}
\begin{example}[Discrete properties of random maps]\label{sec:discrete}
Let $[X_d]\nrw [X_\infty]$ be a converging family of gaussian random fields. In this example we introduce a useful tool for studying the asymptotic probability induced by $X_d$ on discrete sets as $d\to \infty$.
The key example that we have in mind is the case when we consider a codimension-one ``discriminant'' $\Sigma \subset \mathcal{C}^\infty(S^m, \R^k)$  which partitions the set of functions into many connected open sets. For instance $\Sigma$ could be the set of maps for which zero is not a regular value: the complement of $\Sigma$ consists of countably many open connected sets, each one of which corresponds to a rigid isotopy class of embedding of a smooth codimension-$k$ submanifold $Z\subset S^m$. The following Lemma gives a simple technical tool for dealing with these situations.

\begin{lemma}\label{discretelemma}
Let $E$ be a metric space and let $[X_d], [X_\infty]$ be a random fields such that $[X_d]\nrw [X_\infty]$. Let also $Z$ be a discrete space and $\nu \colon U\subset E\to Z$ be a continuous function defined on an open subset $U\subset E$ such that\footnote{Of course, $E\smallsetminus U=\Sigma$ is what we called ``discriminant'' in the previous discussion. Note that we do not require that $\P\{X_d\in U\}=1$, however it will follow that  $\lim_d \P\{X_d\in U\}=1$.} $\P\{X_\infty\in U\}=1$. Then, for any $A\subset Z$ we have:
\be 
\exists \lim_{d\to \infty}\P\left\{X_d\in U,\ \nu(X_d)\in A\right\}=\P\left\{\nu(X_\infty)\in A\right\}.
\ee
\end{lemma}
\begin{proof}
Since $\nu^{-1}(A)$ is closed and open by continuity of $\nu$, it follows that $\de \nu^{-1}(A)\subset E\backslash U$. Therefore $\P\{X_\infty\in \de \nu^{-1}(A)\}=0$ and by Portmanteau's Theorem \cite[Theorem 2.1]{Billingsley}, we conclude that 
\be \label{discreteq}
\P\{X_d\in \nu^{-1}(A)\}\xrightarrow[d\to\infty]{}\P\{X_\infty\in \nu^{-1}(A)\}, \quad \ \forall \ A\subset Z.
\ee
\end{proof}
Equation \eqref{discreteq}, in the case of a discrete topological space such as $Z$, is equivalent to narrow convergence $\nu(X_d)\nrw \nu(X)$, by Portmanteau's Theorem, because $\de A=\emptyset $ for all subsets $A\subset Z$. Note also that to prove narrow convergence of a sequence of measures on $Z$, it is sufficient to show \eqref{discreteq} for all $A=\{z\}$, indeed in that case the inequality 
\be 
\liminf_{d\to \infty}\P\{\nu_d\in A\}=\liminf_{d\to \infty}\sum_{z\in A}\P\{\nu_d=z\}\ge \sum_{z\in A}\P\{\nu=z\}=\P\{\nu\in A\}
\ee
follows automatically from Fatou's lemma.

Following Sarnak and Wigman \cite{SarnakWigman}, let us consider one simple application of this Lemma. Let $H_{m-1}$ be the set of diffeomorphism classes of smooth closed connected hypersurfaces of $\R^{m}$. Consider $U=\{ f\in \Cr \infty {\mathbb{D}^m}{}\,\colon\, f\transv 0\}$ and let $\nu(f)$ be the number of connected components of $f^{-1}(0)$ entirely contained in the interior of  $\mathbb{D}^m$. For $h\in H_{m-1}$ let $\nu_{h}(f)$ be the number of those components which are diffeomorphic to $h\subset \R^{m}$. In the spirit of \cite{SarnakWigman}, we define the probability measure $\mu(f)\in \mathscr{P}(H_{m-1})$ as 
\be 
\mu(f)=\frac{1}{\nu(f)}\sum_{h\in H_{m-1}}\nu_h(f)\delta_h.
\ee 
Let us consider now the rescaled Kostlan polynomial $X_d:\mathbb{D}^m\to \R$ as in Theorem \ref{thm:Kostlan}. The diffeomorphism type of each internal component of $f^{-1}(0)$ remains the same after small perturbations of $f$ inside $U$, hence $\mu\colon U\to \mathscr{P}(H_{m-1})$ is a locally constant map, therefore by Lemma \ref{discretelemma} we obtain that for any subset $A\subset\mathscr{P}(H_{m-1})$,
\be 
\exists \lim_{d\to \infty}\P\{X_d\in U \text{ and }\mu(X_d)\in A\}=\P\{\mu(X_\infty)\in A\}. 
\ee
Moreover since in this case $X_d\in U$ with $\P=1$, for all $d\in \N$ and the support of $X_\infty$ is the whole $\Cr\infty {\mathbb{D}^m}{}$, we have
\be 
\exists \lim_{d\to \infty}\P\{\mu(X_d)\in A\}=\P\{\mu(X_\infty)\in A\}>0. 
\ee

\end{example}
\begin{example}[Random rational maps] \label{ex:rational}
The Kostlan polynomial $P_d^{m, k+1}$ can be used to define random rational maps. In fact, writing $P_{d}^{m, k+1}=(p_0, \ldots, p_k)$, then one can consider the map $\varphi_{d}^{m,k}:\RP^m\dashrightarrow\RP^k$ defined by:
\be\label{eq:rrm} \varphi_d^{m,k}([x_0, \ldots, x_m])=[p_0(x), \ldots, p_m(x)].\ee
(When $m>k$, with positive probability, this map might not be defined on the whole $\RP^m$; when $m\le k$ with probability one we have that the list $(p_0, \ldots, p_k)$ has no common zeroes, and we get a well defined map $\varphi_{d}^{m,k}:\RP^m\to \RP^k.$) Given a point $x\in \RP^m$ and a small disk $D_d=D(x, d^{-1/2})$ centered at this point, the behavior of $\varphi_{d}^{m,k}|_{D_d}$ is captured by the random field $X_d$ defined in \eqref{eq:rescaledKostlan}: one can therefore apply Theorem \ref{thm:Kostlan} and deduce, asymptotic local properties of this map. 

For example, when $m\le k$ for any given embedding of the unit disk $q:\mathbb{D}^m\hookrightarrow \RP^k$ and for every neighborhood $U$ of $q(\partial \mathbb{D}^m)$ there exists a positive constant $c=c(q)>0$ such that for big enough degree $d$ and with probability larger than $c$ the map 
\be X_d=\varphi_{d}^{m,k}\circ a_d:\mathbb{D}^m\to \RP^k\ee (defined by composing $\varphi$ with the rescaling diffeomorphism $a_d:\mathbb{D}^m\to D_d$) is isotopic to $q$ thorugh an isotopy $\{q_t:\mathbb{D}^m\to \RP^k\}_{t\in I}$ such that $q_t(\partial \mathbb{D}^m)\subset U$ for all $t\in I$.

The random map $\f_d^{m,k}$ is strictly related to the random map $\psi^{m,k}_d\colon S^m\to \R^k$:
\be 
\psi^{m,k}_d(x)=P^{m,k}_d(x),
\ee
which is an easier object to work with. For example the random algebraic variety $\{\f_d=0\}$ is the quotient of $\{\psi_d=0\}$ modulo the antipodal map. If we denote by $D_d$ any sequence of disks of radius $d^{-\frac{1}{2}}$ in the sphere, then $\psi_d|_{D_d}\approx X_d$, so that we can understand the local asymptotic behaviour of $\psi_d$ using Theorem \ref{thm:Kostlan} (see Figure \ref{fig:rescale}).
For instance, from point $(7)$ it follows that 
\be 
\E\left\{b_i\left(\left\{\psi_d=0\right\}\cap D_d\right)\right\}\to \E\left\{b_i\left(\left\{X_\infty=0\right\}\cap \mathbb{D}^m\right)\right\}.
\ee
\end{example}

\begin{example}[Random knots]
\label{ex:knots}
Kostlan polynomials offer different possible ways to define a ``random knot''. The first is by considering a random rational map:
\be \varphi_{d}^{1,3}:\RP^1\to \RP^3,\ee
to which the discussion from Example \ref{ex:rational} applies.
(Observe that this discussion has to do with the \emph{local} structure of the knot.)

Another interesting example of random knots, with a more global flavour, can be obtained as follows. Take the random Kostlan map $X_d: \R^2\to \R^3$ (as in \eqref{eq:rescaledKostlan} with $m=2$ and $k=3$) and restrict it to $S^1=\partial \mathbb{D}^m$ to define a random knot:
\be k_d=X_{d}|_{\partial \mathbb{D}^m}:S^1\to \R^3.\ee
The difference between this model and the previous one is that this is global, in the sense that as $d\to \infty$ we get a limit global model $k_\infty=X_\infty|_{\partial D}:S^1\to \R^3$.  
What is interesting for this model is that the Delbruck--Frisch--Wasserman conjecture \cite{delbruck, frischwasserman}, that a typical random knot is non-trivial, does not hold: in fact $k_\infty$ charges every knot (included the unknot) with positive probability.
\begin{prop}The random map:
\be k_d=X_d|_{\partial \D^2}:S^1\to \R^3.\ee
is almost surely a topological embedding (i.e. a knot).
Similarly, the random rational map $\varphi_{d}^{1,3}:\RP^1\to \RP^3$ is almost surely an embedding.
\end{prop}
\begin{proof}We prove the statement for $k_d$, the case of $\varphi_d^{1, 3}$ is similar.
Since $S^1$ is compact, it is enough to prove that $k_d$ is injective with probability one.

Let $F_d=\R[x_0, x_1, x_2]_{(d)}^3$ be the space of triples of homogeeous polynomials of degree $d$ in $3$ variables. Recall that $k_d=X_d|_{\partial \mathbb{D}^2}$, where, if $P\in F_d$, we have set:
\be X_d(u)=P\left(1, \frac{u}{\sqrt{d}}\right),\quad u=(u_1, u_2)\in \R^2.\ee
Let now $S^1=\partial \mathbb{D}^2\subset \R^2$ and $\phi:\left((S^1\times S^1)\backslash \Delta\right)\times F_d\to \R^3$ be the map defined by
\be \phi(x,y, P)=P\left(1, \frac{x}{\sqrt{d}}\right)-P\left(1, \frac{y}{\sqrt{d}}\right).\ee
Observe that $\phi\pitchfork \{0\}.$ By the parametric transversality theorem we conclude that $\phi$ is almost surely transversal to $W=\{0\}$. This imples that, with probability one, the set 
\be \{x\neq y\in S^1\times S^1\,|\, k_d(x)=k_d(y)\}\ee
is a codimension-three submanifold of $S^1\times S^1$ hence it is empty, so that $k_d$ is injective. 
\end{proof}

Theorem \ref{thm:Kostlan} implies now that the random variable $k_d\in C^\infty(S^1, \R^3)$ converges narrowly to $k_\infty\in C^{\infty}(S^1, \R^3)$, which is the restriction to $S^1=\partial \mathbb{D}^2$ of $X_\infty.$ Note that, since the support of $X_\infty$ is all $C^\infty(\mathbb{D}^2, \R^3)$, it follows that the support of $k_\infty$ is all $C^\infty(S^1, \R^3)$ and in particular every knot (i.e. isotopy class of topological embeddings $S^1\to \R^3$, a set with nonempty interior in the $C^\infty$ topology) has positive probability by Theorem \ref{thm:Kostlan}.3. Moreover, denoting by $\gamma_1\sim \gamma_2$ two isotopic knots, we have that
\be 
\P\left(\de\{k_\infty\sim \gamma\}\right)\le\P\{k_\infty \text{ is not an immersion}\}=0
\ee
by Theorem \ref{thm:Kostlan}.4, because the condition of being an immersion is equivalent to that of being transverse to the zero section of $J^1(S^1,\R^3)\to S^1\times \R^3$. 
Theorem \ref{thm:Kostlan}.2, thus implies that for every knot $\gamma:S^1\to \R^3$ we have:
\be \lim_{d\to \infty}\PP \{k_d\sim \gamma\}=\PP\{k_\infty\sim \gamma\}>0.\ee
 
\end{example}
\end{subappendices}

\chapter{SEMICONTINUITY OF TOPOLOGY}\label{chap:semicontop}
\newcommand{\coostr}[3]{\mathcal{C}_{S}^{#1}(#2,#3)}
\newcommand{\dZ}{\textit{d}_Z}
\section{The Semicontinuity Theorem}
A consequence of Thom's Isotopy Lemma is that the set of solutions of a regular smooth equation is stable under $\mathcal{C}^1$-small perturbation, but what happens if the perturbation is just $\mathcal{C}^0$-small? In this section we will show that then the topology cannot decrease.

Instead of a ``regular smooth equation'' we are actually going to consider a more geometric (and more general) situation. Let $f\colon M\to N$ be a $\mC^1$ map and let $Y\subset N$ be a smooth submanifold and assume that $f$ is transverse to $Y$, that is 
\be 
d_pf(T_pM)+T_qY=T_qN \quad \forall p\in f^{-1}(Y),
\ee
and is denoted as $f\transv Y$. It is a standard fact that transversality implies that the set $f^{-1}(Y)$ is a $\mC^1$ submanifold of $M$ having the same codimension as that of $Y$ in $N$.
The main result of this chapter is that the homology groups of $f^{-1}(Y)$ do not decrease when $f$ is perturbed in a $\mC^0$ small way. By this we mean that, given a distance function $\textit{dist}$ on $N$, there exists a continuous function $\e\colon M\to \R_+$ such that if $\tilde{f}$ is another function such that
\be\label{semi:eq:epsilon} 
\sup_{p\in M} \textit{dist}(f(p),\tilde{f}(p)) < \e(p), 
\ee
then, for every $i\in\N$, there is a group $G_i$ such that
\be 
\check{H}^i\left(\tilde{f}^{-1}(W)\right)=H^i\left(f^{-1}(W)\right)\oplus G_i.
\ee
In this sense the homology groups cannot decrease after the perturbation.
By varying the function $\e\colon M\to \R_+$, the family of sets of maps satisfying \eqref{semi:eq:epsilon} form a basis for 
the strong Whitney topology of $\coo 0MN$.

In what follows, we will denote by $\coo rMN$, the space of functions from $M$ to $N$ of class $\mC^r$, for any $r\in[0,+\infty]$. This set has two natural topologies (which coincide when $M$ is compact) called Whitney weak and strong topologies (see \cite{Hirsch}). We will implicitely assume that $\coo rMN$ is endowed with the former, while we will use the notation $\coostr rMN$ to indicate the topological space on the same underlying set, but endowed with Whitney's strong topology.

\begin{defi}
We say that a subset $\mathcal{U} \subset \coostr 0MN$ is \emph{homotopy connected} if for any two elements $f,g\in\mathcal{U}$ there exists a homotopy connecting the two, namely a function $F\in \coo0{M\times [0,1]}N$ such that $F(\cdot,0)=f$ and $F(\cdot,1)=g$.
\end{defi}

\begin{remark}\label{semi:rem:connect}
A set $\mathcal{U}\subset \coo rMN$ is homotopy connected if and only if it is path connected with respect to the weak topology, in fact if $F\in \coo 0{M\times[0,1]}N$ is an homotopy then the function $t\mapsto F(\cdot, t)$ is a continuous path in $\coo 0MN$. 
We are being extra careful here, for if $M$ is not a compact manifold, then the space $\coostr 0MN$ is not locally path connected. Indeed one can show that all the continuous maps $\a\colon [0,1]\to \coostr 0MN$ have the property that 
\be\a(t)|_{M\- K}=\a(0)|_{M\- K}
\ee 
for some compact subset $K\subset M$. 
Nevertheless, the strong topology has a basis consisting of homotopy connected subsets.
\end{remark}

Every open subset $\mathcal{U}\subset \coostr 0MN$ is a disjoint union of homotopy connected open subsets. To see this, we observe that its homotopy connected components $U_i$ are open in $\coostr0MN$. 
Indeed for any $f\in U_i$, there exists a strong open neighborhood $O_f\subset \mathcal{U}$ and, by remark \ref{semi:rem:connect}, we can assume it to be homotopy connected, so that $O_f$ must be contained in the weak connected component of $\mathcal{U}$ containing $f$, which is $U_i$. 
\begin{figure}
    \centering
    \includegraphics[scale=0.2]{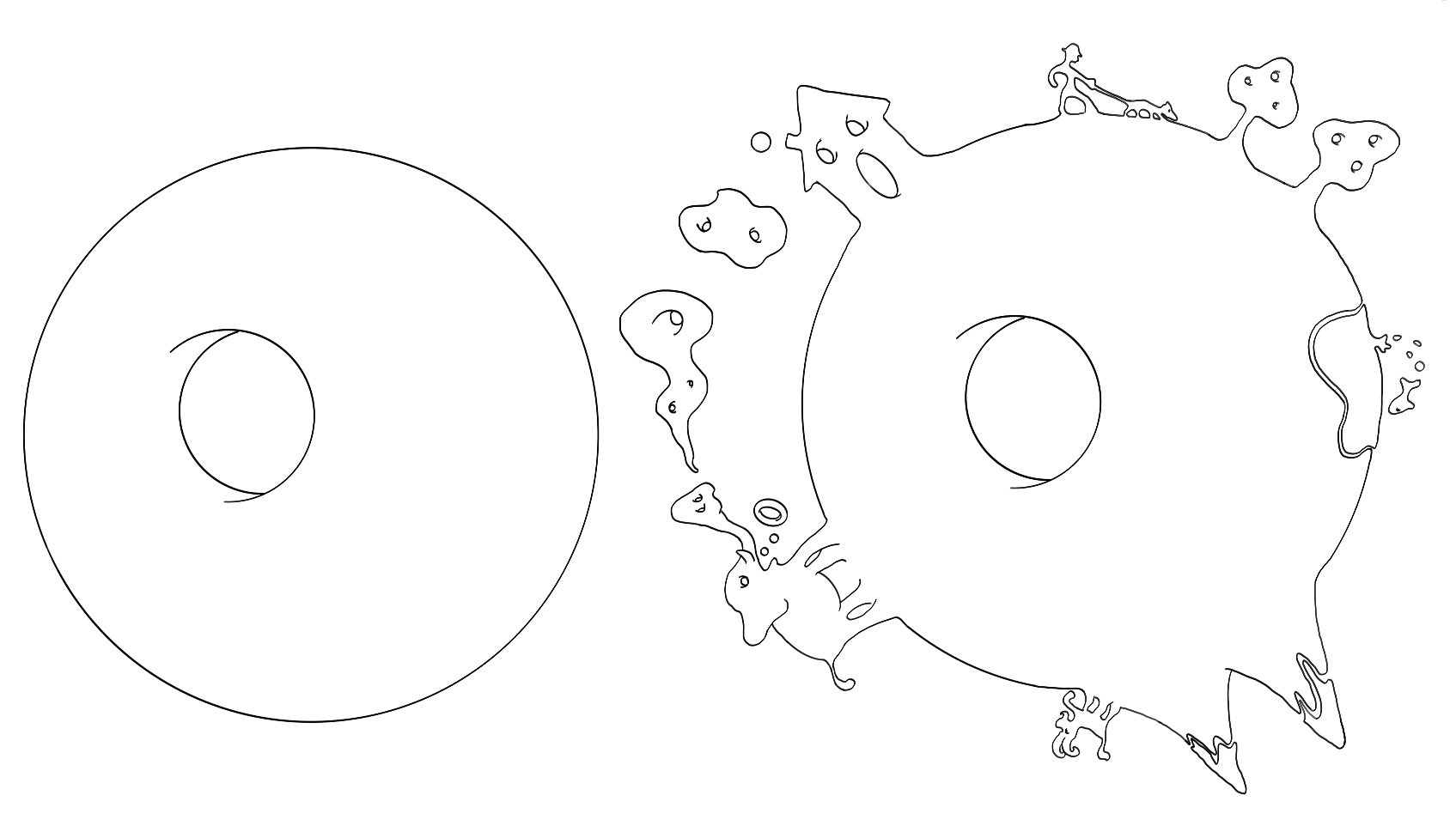}
    \caption{With a $\mC^0$ small perturbation, holes can be \emph{created}, but not  \emph{destroyed}.}
    \label{semi:fig:toricworld}
\end{figure}
\begin{thm}\label{semi:thm:main}
Let $M,N$ be smooth manifolds, let $Y\subset N$ be a closed and cooriented smooth submanifold. Let $f\colon M\to N$ be a $\mC^1$ map such that $f\transv Y$ and $A:=f^{-1}(Y)$. Let $E,E_1\subset M$ be two open tubular neighborhoods\footnote{Equivalently, there exists a $\mC^1$ map $\pi\colon E\to M$ and a $\mC^1$ vector bundle structure on $E$, of which $\pi$ is the projection.} of $A$ such that $\overline{E}\subset E_1$.
\begin{enumerate}
    \item\label{semi:itm:openness} Define the set  $\mathcal{U}_{E,f}$ as the homotopy connected component containing $f$ of the set
    \be\label{semi:eq:U}
    \mathcal{U}_E=\left\{g\in\coo 0MN\colon g(M\- E)\subset N\- Y\right\}.
    \ee
    Then $\mathcal{U}_{E,f}\subset \coostr0MN$ is open with respect to Whitney's strong topology.
    \item If $\Tilde{f}\in\mathcal{U}_{E,f}$ and $\Tilde{A}=\Tilde{f}^{-1}(Y)$, then there exist abelian groups $G_i$, for each $i\in \N$, such that
    \be\label{semi:eq:sumgru}
     \check{H}^i(\Tilde{A})\cong\check{H}^i(A)\oplus G_i.
    \ee
\end{enumerate}
\end{thm}
The groups $\check{H}^i(A)$ we refer to in the above statement are the \v{C}ech cohomology groups (see \cite{bredon}). 

\begin{remark}
When $f\transv Y$ the preimage set $A$ is a $\mC^1$ manifold, therefore $\check{H}^i(A)=H^i(A)$ for every $i\in \N$. In particular, in the right hand side of \eqref{semi:eq:sumgru} we can write the singular cohomology $H^i(A)$, or any other cohomology theory.
\end{remark}

Although in most cases the two objects coincide, the statement would be false if we substitute with the more standard singular cohomology groups $H^i(A)$, as example \ref{semi:exa:counter} below shows. The reason is that, contrary to singular homology, \v{C}ech cohomology theory posesses the following crucial property (see \cite[p. 348]{bredon})
\be\label{semi:eq:crucial}
\check{H}^i(A)=\lim_{\to}\left\{\check{H}^i(B)\colon A\subset B \text{ open}\right\}.
\ee
In fact, we can replace $\check{H}$ in Theorem \ref{semi:thm:main} with any cohomology theory satisfying \eqref{semi:eq:crucial}.

\begin{figure}
    \centering
    \includegraphics[scale=0.15]{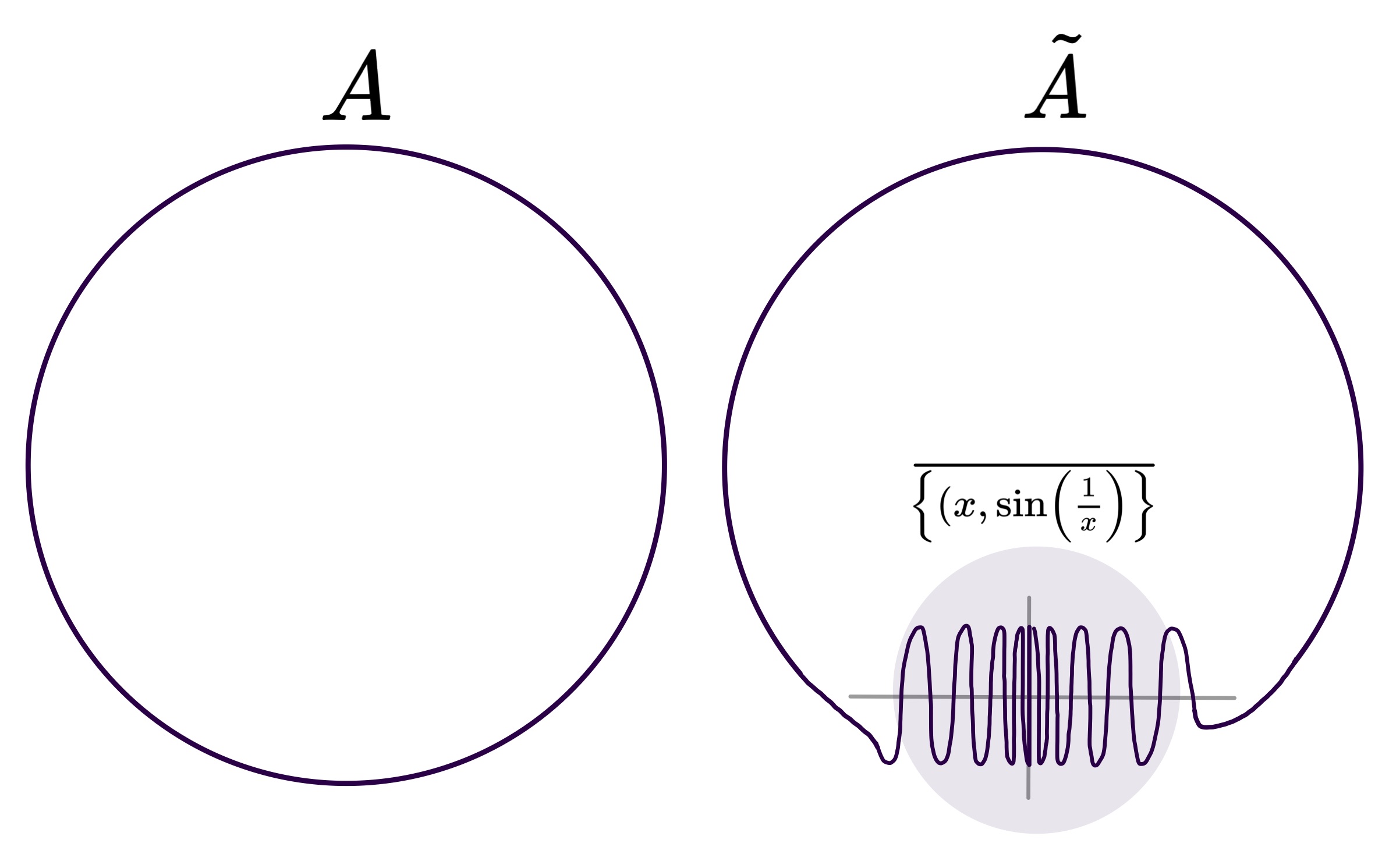}
    \caption{On the left, $A$ is homeomorphic to $S^1$, while $\tilde{A}$, on the right, is made by two contractible path connected components.}
    \label{fig:semi:semicontro}
\end{figure}
\begin{example}\label{semi:exa:counter}
Consider the two topological spaces $A$ and $\tilde{A}$ of figure \ref{fig:semi:semicontro}. Clearly, we can find a smooth function $f\colon M\to \R$ such that $f^{-1}(0)=A$ and such that $0$ is a regular value for $f$, that is $f\transv \{0\}$. Moreover, it is not difficult to show that there are smooth functions $\tilde{f}\colon M\to \R$ that are arbitrarily close to $f$ in $\coostr 0M\R$, whose zero set is homeomorphic to $\tilde{A}$. 
Theorem \ref{semi:thm:main} holds in this case and indeed
\be 
1=b_0(A)\le \check{b}_0(\tilde{A})=1 \quad \text{ and } \quad 1=b_1(A)\le \check{b}_1(\tilde{A})=1,
\ee
where $\check{b}_i$ is the dimension of the $i^{th}$ \v{C}ech cohomology group. 
However, this is false for what concerns the dimension of the singular homology groups, since 
\be 
b_1(\tilde{A})=0.
\ee
\end{example}

The requirement that $Y$ should be coorientable in Theorem $\ref{semi:thm:main}$ can be removed, paying the price of losing the control on the $\mC^0$ neighborhood.
\begin{cor}[Semicontinuity Theorem]\label{semi:cor:mainnotco}
Let $M,N$ be smooth manifolds, let $Y\subset N$ be a closed (not necessarily cooriented) smooth submanifold. Let $f\colon M\to N$ be a $\mC^1$ map such that $f\transv Y$. 
There is a neighborhood $\mathcal{U}\subset \coostr0MN$ of $f$, open with respect to Whitney's strong topology, with the property that for any $\Tilde{f}\in\mathcal{U}$ there exist abelian groups $G_i$, for each $i\in \N$, such that
    \be
     \check{H}^i(\Tilde{f}^{-1}(Y))\cong H^i(f^{-1}(Y))\oplus G_i.
    \ee
\end{cor}
\section{Proofs of Theorem \ref{semi:thm:main} and Corollary \ref{semi:cor:mainnotco}}
\begin{proof}[Proof of Theorem \ref{semi:thm:main}] 
Openness of the set $\mathcal{U}_E$ in $\coostr0MN$ is due to the fact that $M\-E$ is closed and $N\- Y$ is open, thus point $(\ref{semi:itm:openness})$ follows from the discussion preceding the Theorem.

Fix $\tilde{f}\in\mathcal{U}_{E,f}$, so that $\tilde{A}\subset E\subset E_1$ and denote by $\pi\colon E_1\to A$ the retraction map. 
Let $F$ be an homotopy connecting $f$ with $\tilde{f}$ such that $F(\cdot,t)=F_t \in \mathcal{U}_{E,f}$. 

There is an inclusion of pairs.
\be 
u:(E_1,E_1\- E)\to (E_1,E_1\-\tilde{A}).
\ee
Define analogously $\tilde{\pi}:\tilde{E}\to \tilde{A}$ to be a tubular neighborhood of $\tilde{A}$, so small that $\tilde{E}\subset E$.  
Let now $\tilde{B}$ be any open neighborhood of $A$, such that $\tilde{A}\subset \tilde{B}\subset \tilde{E}$.

\begin{figure}[h]
    \centering
    \includegraphics[scale=0.17]{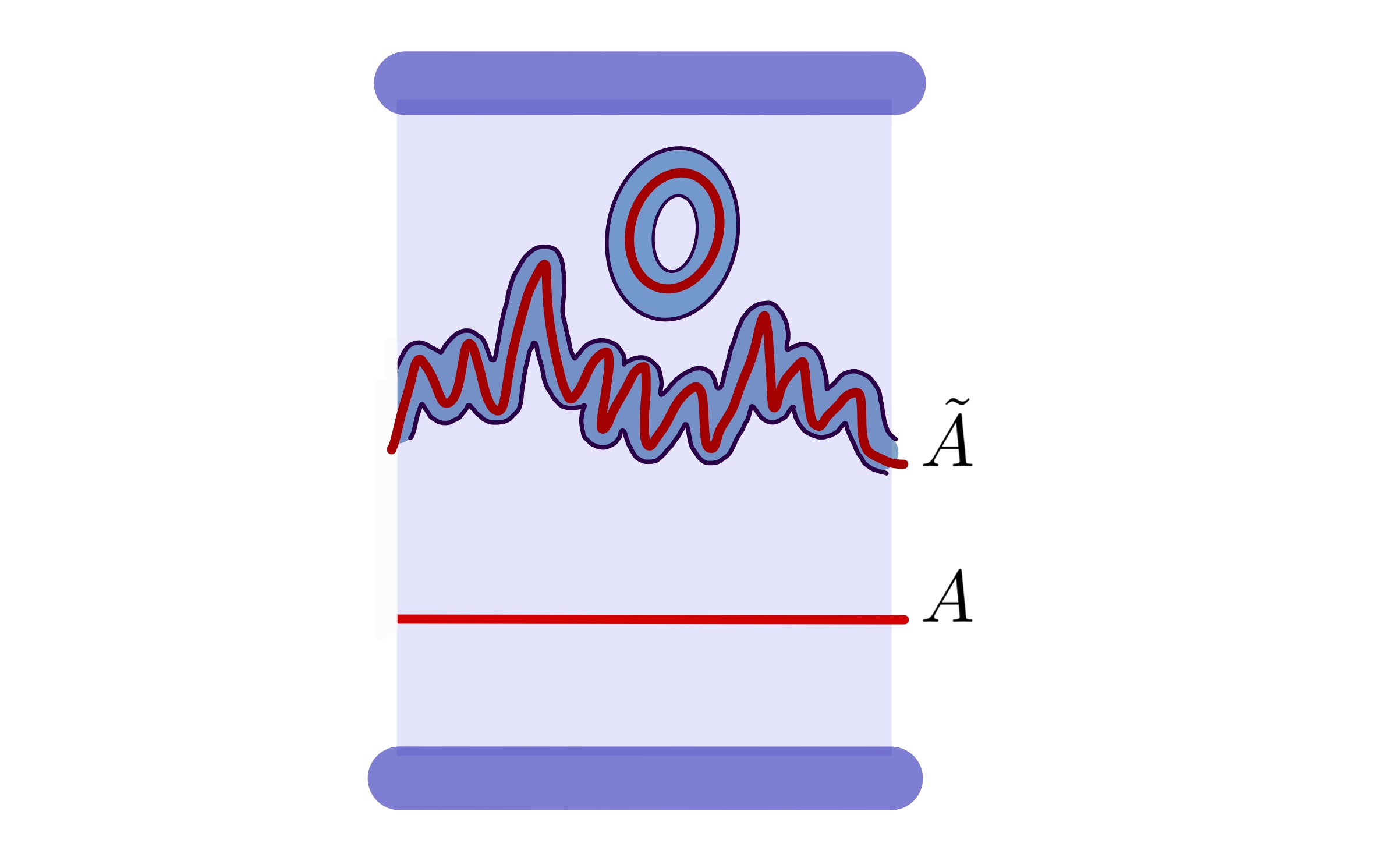}
    \caption{$\tilde{A}$ is contained in a tubular neighborhood of $A$.}
    \label{semi:fig:nbhds}
\end{figure}

The fact that $Y$ is cooriented guarantees the existence of a Thom cohomology class $\phi\in H^r(N,N\- Y)$, where $r$ is the codimension of $Y$. By the transversality condition $\f\transv Y$, we have that also $A$ is cooriented with Thom class $f^*\phi=\phi_E\in H^r(E_1,E_1\- E)\cong H^r(E,E\- A)$ and we have that
\be 
f^*\phi={F_t}^*\phi=\tilde{f}^*\phi,
\ee
because of the homotopy of pairs $F_t\colon (E_1,E_1\- E)\to (N,N\- Y)$.
This and the fact that $\pi\circ\tilde{\pi}$ is homotopic to $\pi|_{\tilde{E}}$, imply the commutativity of the diagram below (where $\eta$ is the excision isomorphism and $\Phi_{\tilde{B}}=\tilde{f}^*(\phi)$). 

\be 
\begin{tikzcd}
                                                                                               & {H^{*+r}(N, N\backslash Y)}\arrow[ldd, "\tilde{f}^*"'] \arrow[rdd, "\tilde{f}^*=f^*"]                         &                                                                  \\
                                                                                               &                                                      &                                                                  \\
{H^{*+r}(\tilde{B}, \tilde{B}\smallsetminus\tilde{A})} \arrow[r, "\eta^{-1}"]  & {H^{*+r}(E_1, E_1\smallsetminus \tilde{A})} \arrow[r, "u^*"] & {H^{*+r}(E_1, E_1\smallsetminus E)}           \\
H^{*}(\tilde{B}) \arrow[u, "\tilde{\pi}^*(\cdot)\cup \phi_{\tilde{B}}"]                        &                                                      & H^*(A) \arrow[ll, "{(\pi|_{\tilde{B}})}^*"] \arrow[u, "\pi^*(\cdot)\cup \phi_E"']
\end{tikzcd}\ee

Since the map $\pi^*(\cdot)\cup \phi_E$ is the Thom isomorphisms, there exists a homomorphism $U_{\tilde{B}}\colon H^*(\tilde{B})\to H^*(A)$ such that $U_{\tilde{B}}\circ (\pi|_{\tilde{B}})^*=\text{id}$. Notice that if $\tilde{B_1}\subset \tilde{B_2}$ are two neighborhoods of $\tilde{A}$, then there is a commutative diagram
\be
\begin{tikzcd}
H^*(A) \arrow[rr, "(\pi|_{\tilde{B_2}})^*"] \arrow[rrd, "(\pi|_{\tilde{B_1}})^*"'] \arrow[rrrr, "\text{id}", bend left] &  & H^*(\tilde{B_2}) \arrow[rr, "U_{\tilde{B_2}}"] \arrow[d] &  & H^*(A) \\
                                                                                                                  &  & H^*(\tilde{B_1}) \arrow[rru, "U_{\tilde{B_1}}"']     &  &       
\end{tikzcd}
\ee
Therefore, by the fundamental property of the direct limit \eqref{semi:eq:crucial} we deduce the existence of a homomorphism $U_{\tilde{A}}$ such that the next diagram is commutative. 
\be
\begin{tikzcd}
H^*(A) \arrow[rr, "(\pi|_{\tilde{B_2}})^*"] \arrow[rrd, "(\pi|_{\tilde{B_1}})^*"'] \arrow[rrrr, "\text{id}", bend left] \arrow[rrddd, "(\pi|_{\tilde{A}})^*"'] &  & H^*(\tilde{B}) \arrow[rr, "U_{\tilde{B_2}}"] \arrow[d]                                                                               &  & H^*(A) \\
                                                                                                                                                             &  & H^*(\tilde{B_1}) \arrow[rru, "U_{\tilde{B_1}}"'] \arrow[d]                                                                         &  &        \\
                                                                                                                                                             &  & \vdots \arrow[d]                                                                                                                   &  &        \\
                                                                                                                                                             &  & \varinjlim_{\tilde{B}\supset \tilde{A}} H^*(\tilde{B})\cong \check{H}^*(\tilde{A}) \arrow[rruuu, "\exists U_{\tilde{A}}"', dashed] &  &       
\end{tikzcd}
\ee
This concludes the proof, since it means that $(\pi|_{\tilde{A}})^*$ has a left inverse $U_{\tilde{A}}$, that is equivalent (in the context of abelian groups) to say that the short exact sequence
\be 
\begin{tikzcd}
0 \arrow[r] & \ker(\pi|_{\tilde{A}})^* \arrow[r] & H^*(A) \arrow[r, "(\pi|_{\tilde{A}})^*"] & \check{H}^*(\tilde{A}) \arrow[r] \arrow[l, "U_{\tilde{A}}"', bend left=49] & 0
\end{tikzcd}
\ee
splits, thus
\be 
\check{H}^*(\tilde{A})\cong H^*(A)\oplus \ker(\pi|_{\tilde{A}})^*.
\ee
\end{proof}
\begin{proof}[Proof of Corollary \ref{semi:cor:mainnotco}]
Let $N_0\subset N$ be a tubular neighborhood of $Y$ in $N$. Define $A=f^{-1}(Y)$ and let $\pi\colon E\to A$ be a tubular neighborhood of $A$ in $M$, small enough that there exists another open tubular neighborhood $E_1$ such that $\overline{E}\subset E_1$ and $f_0(E_1)\subset N_0$. Now define $\mathcal{U}$ to be the homotopy connected component containing $f$ of the open set
\be 
\mathcal{U}_0=\left\{g\in\coo 0MN\colon g(M\- E)\subset N\- Y \text{ and } g\left(\overline{E}\right)\subset N_0\right\}.
\ee

For every $\tilde{f}\in \mathcal{U}$ there is a homotopy $F$ connecting it to $f$ inside $\mathcal{U}$. Let us now focus on its restriction to $E$
\be 
F\colon [0,1]\times E\to N_0.
\ee
Let $\pi'\colon L\to N_0$ be the pull-back over $N_0$ of the determinant (real) line bundle of $Y$ and observe that the natural embedding of $Y$ in the total space of $L$ is canonically cooriented. By standard facts about topological vector bundles, all the pull-back bundles $F_t^*L$ on the manifold $E$ are isomorphic, which implies that there exist a family of linear isomorphisms
\be 
\phi_{(t,p)}\colon L_{F(0,p)}\to L_{F(t,p)},
\ee
depending continuosly on $(t,p)\in [0,1]\times E$ and such that $\phi_{(0,p)}=$id. This allows to define a new map
\be
\hat{F}\colon [0,1]\times F_0^*L\to L
\ee
\be 
\hat{F}_t(p,v)=\hat{F}(t,p,v)=\left(F(t,p), \phi_{(t,p)}(v)\right),
\ee
with the properties that $\hat{F}\transv Y$, and that $\hat{F}_t^{-1}(Y)=F_t^{-1}(Y)$ for every $t\in[0,1]$ and since $F_0^*L$ is already a tubular neighborhood of $A$ we have authomatically that
\be \hat{F}_1\in\mathcal{U}_{F_0^*L,\hat{F}_0},
\ee 
thus we can conclude by applying Theorem $\ref{semi:thm:main}$ to the maps $\hat{F}_0$ and $\hat{F}_1$.
\end{proof}
\section{Application: What is the degree of a smooth hypersurface?}
Assume that $Z$ is the zero set of a polynomial $p\in\R[x_1,\dots,x_n]$ of degree $d$, then the sum of its Betti numbers $b(Z)$ can be estimated by the classical Milnor-Thom bound (Milnor \cite{MilnorBound}, 1964):
\be \label{semi:eq:milthom}
b(Z)\le d(2d-1)^{n-1}\sim d^n.
\ee
In this section we will prove, as an application of Theorem \ref{semi:thm:main}, an analogous result valid for any smooth compact hypersurface\footnote{A closed submanifold of codimension $1$.} in $\R^n$.

Of course, the first question is \emph{What object can play the role of the degree?}
Since we are in $\R^n$, all compact hypersurfaces are the regular zero set of a smooth function, and smooth functions can be approximated by polynomials, thus the idea is to define a quantity that estimate the degree of the best polynomial approximation. 

To this end, we define, for every $\mC^1$ function $f$ defined on a closed embedded disk $D\subset \R^n$, the quantity $\delta(f)$, that measures the amount of regularity of the equation $f=0$.
\be\label{semi:eq:delta} \delta(f, D)=\inf_{x\in D}\left(|f(x)|^2+\|\nabla f (x)\|^2\right)^{1/2}.\ee
Observe that, if the zero set of $f$, denoted by $Z(f)$ is contained in the interior of $D$, then the equation $f=0$ is regular on $D$ if and only if $\delta(f, D)>0$. Therefore we interpret this quantity as a sort of ``distance'' from $f$ to the set of functions with singular zero set in $D$. 
Here the space $\coo \ell D \R$ is endowed with the norm
\be\label{semi:eq:norm} \nu_\ell(f,D)=\|f\|_{C^\ell(D, \R)}=\sup_{x\in D}\left(\sum_{|\alpha|\leq \ell}\left|\frac{\partial ^\alpha f}{\partial x^{\alpha}}(x)\right|^2\right)^{1/2}.\ee
The second one is a measurement of how much the boundary condition $Z(f)\subset \text{int}(D)$ is met,
\be 
m(f,D)=\min_{x\in \de D}|f(x)|.
\ee
We pack together the above data in a single, scaling invariant, parameter which constitutes the piece of information that we require to know on our function $f$. 
\be \kappa_\ell(f, D)=\frac{\nu_\ell(f,D)}{\min\left\{m(f,D),\delta(f,D)\right\}};\ee
this, with $\ell=1$, is in fact the only quantity we use in our estimates. 
\begin{thm}\label{semi:thm:witdash}
Let $f\in \coo 1D\R$ have regular zero set $Z(f)\subset \text{int}(D)$. Then there is a constant $c_1=c_1(D)>0$ such that
\be\label{semi:eq:bobound}
b\left(Z(f)\right)\le {c_1}\cdot\left({\kappa_1(f,D)}\right)^n
\ee
\end{thm}
\begin{remark}It is possible to deduce \eqref{semi:eq:bobound} also from the work of Yomdin \cite{Yomdin}, where bounds on the Betti numbers of $Z(f)$ are stated in terms of the distance from zero to the the set of critical values of $f$.
\end{remark}
The bound is ``asymptotically sharp'', in the sense of the following statement.
\begin{thm}\label{semi:thm:sharpdeg}
There exists a sequence $\{f_k\}_{k\in \mathbb{N}}$ of maps in $\mC^1(D, \R)$, such that the sequence $\kappa_1(f_k,D)$ converges to $+\infty$, and a constant $c_0=c_0(D)>0$ such that for every $k\in \mathbb{N}$ the zero set $Z(f_k)\subset \mathrm{int}(D)$ is regular and
\be b(Z(f_k))\geq c_0\cdot(\kappa_1(f_k, D))^{n}.
\ee 
\end{thm}

\begin{remark}
Actually, what we will prove is that
\be 
b(Z(f))\le \left(a_0(D)\kappa_1(f,D)+1\right)^n,
\ee
where $a_0(D)$ is the constant given by \eqref{eq:w0}. This implies \eqref{semi:eq:bobound}, since $\kappa_1(f)\ge 1$, by definition.
On the other hand, for any given compact hypersurface $Z\subset D$, there exists a sequence of smooth functions $f_k$ such that
\be 
b(Z(f_k))\ge \frac{b(Z)}{\kappa_1(f,D)^n}h(D)\kappa_1(f_k,D)^n,
\ee
where $h(D)$ is the infimum among the numbers $N\e^n$, such that there exists a collection of $N$ disjoint $n$-dimensional 
disks of radius $\e$ contained in $D$.
\end{remark}

\subsection{Every smooth hypersurface in $\R^n$ is the regular zero set of a smooth function}This result is well-known, see for instance \cite[Theorem 7.2.3]{DFN}. The purpose of this section is to produce an estimate on the condition number of a defining equation, in terms of some explicit data of the manifold. 

To this end, given $Z\subset \R^n$ of class $\mC^1$ we define the \emph{reach} of $Z$ as:
\be \rho(Z)=\sup\left\{r>0: \mathrm{dist}(x, Z)<r\implies \exists! z\in Z \,|\,\mathrm{dist}(x, Z)=\mathrm{dist}(x, z)\right\}.
\ee
The reach of a $\mC^1$ manifold doesn't need to be positive, as shown in \cite[Example 4]{KP}, where an example of a $\mC^{2-\epsilon}$ compact curve with zero reach in $\R^2$ is constructed; however $\rho(Z)>0$ if $Z$ is of class $\mC^2.$

We prove the following result.

\begin{thm}\label{semi:thm:reach}
Given $Z\subset \R^n$ a compact hypersurface of class $\mC^1$ with $\rho(Z)>0$, there exists a $\mC^1$ function $f:\R^n\to \R$ whose zero set is $Z$ and such that for every disk $D=D_R$ with $Z\subset \mathrm{int}D_{R-\rho(Z)}$, we have:
\be\label{eq:k1} \kappa_1(f, D)\leq 2\left(1+\frac{1}{\rho(Z)}\right).\ee
\end{thm}
\begin{proof}Denote by $\rho=\rho(Z).$ We consider the function $\dZ:\R^n\to \R$ defined to be the \emph{signed} distance from $Z$. By \cite[Remark 2]{Foote}, if $Z$ is of class $\mC^1$ and with positive reach $\rho(Z)>0$, the function $\dZ$ is $\mC^1$ on the set $\{\dZ<\rho(Z)\}.$ 

We need to consider also an auxiliary function $g:\R\to \R$ of class $\mC^2$ such that $g(t)=-g(-t)$ for all $t\in \R$ and:
\be g(t)=\begin{cases}
                                   t & \text{if $0\leq t\leq \frac{1}{2}$} \\
                                   \textrm{is increasing and concave} & \text{if $\frac{1}{2}\leq t\leq \frac{7}{8}$} \\
 \frac{3}{4} & \text{if $t\geq \frac{7}{8}$}
  \end{cases}\ee
  The existence of such a function is elementary; it can be taken, for instance, to be piecewise polynomial. Denoting by $g_\rho$ the function $t\mapsto \rho\cdot  g(t/\rho)$, we set:
 \be \label{eq:dfnf}f(x)=g_\rho(\dZ(x)).\ee
 
Notice that this condition implies that, by construction, the function $f\equiv \frac{3}{4}\rho$ on $\R^n\backslash \mathrm{int}(D)$ and in particular:
 \be\label{eq:d} \delta(f, D)=\min\left\{\inf_{z\in D}\left(|f(z)|^2+\|\nabla f (z)\|^2\right)^{1/2}, \frac{3}{4}\rho\right\}.\ee
 Observe now that for every $x$ such that $t=\dZ(x)<\rho$ we have:
 \begin{align} |f(x)|^2+\|\nabla f(x)\|^2&=|g_\rho(\dZ(x))|^2+|g'_{\rho}(\dZ(x))|^2\cdot\|\nabla \dZ(x)\|^2\\
 &=|g_\rho(\dZ(x))|^2+|g'_{\rho}(\dZ(x))|^2\\
 \label{eq:red}&=\rho^2|g(t/\rho)|^2+|g'(t/\rho)|^2,
 \end{align}
 where in the second line we have used the fact that $\|\nabla\dZ\|\equiv 1.$
 In particular, partitioning the domain of definition of the function $g_\rho$, it follows that:
 \begin{align} \inf_{z\in D}|f(z)|^2+\|\nabla f (z)\|^2& \ge
 \inf\left\{1+\frac{\rho^2}{4}, \rho^2(3/4)^2, \inf_{\frac{1}{2}\rho\leq t\leq \frac{7}{8}t}\rho^2|g(t/\rho)|^2+|g'(t/\rho)|^2\right\}\\
 &\geq \inf\left\{1+\frac{\rho^2}{4}, \rho^2(3/4)^2, \inf_{\frac{1}{2}\rho\leq t\leq \frac{7}{8}t}\rho^2|g(t/\rho)|^2+ \inf_{\frac{1}{2}\rho\leq t\leq \frac{7}{8}t}|g'(t/\rho)|^2\right\}\\
 &=\inf\left\{1+\frac{\rho^2}{4}, \rho^2(3/4)^2, \frac{\rho^2}{4}\right\}\\
 &=\frac{\rho^2}{4}.
 \end{align}
 Together with \eqref{eq:d}, this gives:
 \be\label{eq:db} \delta(f, D)\geq \frac{\rho}{2}.\ee
 
 Let us now estimate $\nu_1(f , D)$. Again partitioning the domain and using \eqref{eq:red} and the fact that $|g'(t)|\leq 1$ for all $t$, we immediately get:
 \be\label{eq:nb} \nu_1(f, D)\leq 1+\frac{3}{4}\rho.
 \ee
 Combining \eqref{eq:db} with \eqref{eq:nb} gives \eqref{eq:k1}.
\end{proof}
An immediate consequence of Theorem \ref{semi:thm:witdash} and Theorem \ref{semi:thm:reach} is the following.
\begin{cor}\label{semi:cor:reachwitdash}
There is a constant $c_2=c_2(D)>0$ such that if $Z$ is a $\mC^1$ compact hypersurface contained in a disk $D$, then 
\be\label{semi:eq:reachbobound}
b\left(Z\right)\le {c_2}\cdot\left(1+\frac 1{\rho(Z)}\right)^n
\ee
\end{cor}
\begin{proof}
Observe that $\rho(Z)$ cannot be greater than the radius of $D$ unless $Z$ is empty, in which case there is nothing to prove. Define $D'$ to be the disk with a double radius than that of $D$, so that $Z$ and $D'$ satisfy the hypotheses of Theorem \ref{semi:thm:reach}, thus there exists a function $f$ such that
\bega 
b(Z)&=b\left(Z(f)\right)\le {c_1(D')}\cdot\left({\kappa_1(f,D')}\right)^n\le 2^n{c_1(D')}\cdot\left(1+\frac 1{\rho(Z)}\right)^n,
\eega 
where the first inequality is implied by Theorem \ref{semi:thm:witdash}. Taking $c_2(D)=2^nc_1(D')$ we conclude the proof.
\end{proof}
\subsection{Approximation theorems}The second main ingredients, after the semicontinuity Theorem \ref{semi:thm:main}, to prove \ref{semi:thm:witdash} is the following quantitative versions of Weierstrass' Approximation Theorem from \cite{BBL}.

\begin{thm}[Theorem 2 from \cite{BBL}]\label{semi:lem:boh}
Let $f\in C^{r}(D, \R)$. Then, for every $d\geq0$ there is a polynomial $w_d(f)\in \R[x_1, \ldots, x_n]$ of degree at most $d$ such that for every $\ell\leq \min\{r, n\}$:
\be \nu_\ell(f-w_d(f), D)\leq \frac{a}{d^{r-\ell}}\cdot \nu_r(f, D),
\ee
where the constant $a>0$ only depends on $n, r$ and $D$.
\end{thm}
In particular, using the above Theorem, one immediately gets the following.
\begin{enumerate}
\item There exists a constant $a_0>0$ such that for every $f\in C^{1}(D, \R)$ and for every $d\geq 0$ there is a polynomial $w_{0,d}(f)\in \R[x_1, \ldots, x_n]$ of degree at most $d$ such that:
\be\label{eq:w0} \nu_0(f-w_{0,d}(f), D)\leq a_0\cdot \nu_1(f, D)\cdot d^{-1}.
\ee
\item There exists a constant $a_1>0$ such that for every $f\in C^{2}(D, \R)$ and for every $d\geq 0$ there is a polynomial $w_{1,d}(f)\in \R[x_1, \ldots, x_n]$ of degree at most $d$ such that:
\be\label{eq:w1} \nu_1(f-w_{1,d}(f), D)\leq a_1\cdot \nu_2(f, D)\cdot d^{-1}.
\ee
\end{enumerate}
%
%
The estimate \eqref{eq:w1}, has to be considered in relation to the following fact. Given $g,f\in \coo 1D\R$ be such that \be 
\nu_1(g-f,D)< \delta(f,D),
\ee 
one can show that the couples $(D,Z(g))$ and $(D,Z(f))$ are diffeomorphic. This can be proved easily with Thom's isotopy Lemma and it is essentially due to the fact that there is a continuous family of regular equations $f+t(g-f)=0$ that ``connects'' $Z(f)$ with $Z(g)$. Therefore, combining this with Milnor-Thom bound \ref{semi:eq:milthom}, we have
\be 
b\left(Z(f)\right)\le \left({a_1}{\kappa_2(f)}\right)^{n}.
\ee
However, since we are interested just in controlling the topology of $Z(f)$, we can use Theorem \ref{semi:thm:main}, for which 
only $\mC^1$ information on $f$ is needed (in fact we will use only \eqref{eq:w0}).   This results in the better estimate of Theorem \ref{semi:thm:witdash}
\be 
b\left(Z(f)\right)\le c_1\left({\kappa_1(f)}\right)^n\le c_1\left({\kappa_2(f)}\right)^n.
\ee



\subsection{Proofs of Theorem \ref{semi:thm:witdash} and Theorem \ref{semi:thm:sharpdeg}}
\begin{proof}[Proof of Theorem \ref{semi:thm:witdash}]
First, observe that if $\e<m(f,D)$, then the set $E=f^{-1}(-\e,\e)$ is entirely contained in the interior of $D$. Moreover if $\e< \delta(f,D)$, then $f$ has no critical value in the interval $(-\e,\e)$, therefore $E$ is a tubular neighborhood of $Z(f)\subset \R^n$. In fact, by a standard Morse theoretic argument, $E$ is diffeomorphic to $Z(f)\times (-\e,\e)$. From now on we fix $\e>0$ such that
\be\label{semi:eq:arbitra} 
\e<\min\left\{m(f,D),\delta(f,D)\right\}.
\ee
Let $p\in\R[x_1,\dots,x_n]$ be 
a polynomial such that $\nu_0(f-p,D)< \e$. By Theorem \ref{semi:lem:boh} we can assume that its degree $d$ satisfies the bound
\be\label{semi:eq:degreeeeee}
d-1\le a_0\nu_1(f,D)\frac{1}{\e}
\ee
(take $p=w_{0,d}(f)$, where $d$ is the smallest positive integer such that \eqref{semi:eq:degreeeeee} is false).
Let $F_t=f+t(p-f)$ and call $F_t$ its restriction to $M=\text{int}(D)$. Consider the set $\mathcal{U}_E$ defined as in \eqref{semi:eq:U}, where $N=\R$ and $Y=\{0\}$. Suppose that $F_t\in\mathcal{U}_E$ for every $t\in [0,1]$, then we could apply Theorem \ref{semi:thm:main} to deduce that
\be 
b(Z(f))\le b(Z(p)) \le \left(\frac{a_0\nu_1(f,D)}{\e}+1\right)^n
\ee
where the second inequality is due to the Milnor-Thom bound \eqref{semi:eq:milthom} and to \eqref{semi:eq:degreeeeee}. The thesis now would follow by the arbitrariness of $\e$ in \eqref{semi:eq:arbitra}.

Thus to conclude the proof it is suffcient to show that $F_t(D\-E)\subset \R\- \{0\}$, wich is equivalent to say that $F_t^{-1}(0)\subset E$. To see this, let $x\in D$ such that $F_t(x)=0$ and observe that then
\be
|f(x)|=|F_t(x)-t(p(x)-f(x))|\le \nu_0(p-f,D)\le \e.
\ee
\end{proof}
\begin{proof}[Proof of Theorem \ref{semi:thm:sharpdeg}]
Assume that $D=\D^m$ is the standard unit disk in $\R^n$. Let $f\in\coo 1D\R$ be a smooth function such that $Z=Z(f)$. Since $Z\cap \de D=\emptyset$ we can assume that $f>0$ on $\de D$, so that by slightliy modifying $f$ in a neighborhood of $\de D$, we may assume that $f=1$ near $\de D$.  Extend $f$ to the whole space $\R^n$, by setting $f(x)=1$ for all $x\notin D$. Let $k\in\N$ and define $f_{k,z}\in \coo 1{\R^n}\R$ to be the function
\be 
f_{k,z}(x)=f\left(k(x-z)\right),
\ee
so that the submanifold $Z(f_{k,z})$ is contained in the interior of the disk of radius $k^{-1}$ centered in the point $z$ and it is diffeomorphic to $Z$.
Moreover we can observe that, with $k\ge 1$, we have the inequalities
\be 
\delta(f_{k,z},D)=\delta(f_{k,0,D})=\inf_{x\in D}\left(|f(x)|^2+k^{2}\|\nabla f (x)\|^2\right)^{1/2}\ge \delta(f,D);
\ee
\be 
\nu(f_{k,z},D)=\nu(f_{k,0},D)=\sup_{x\in D}\left(|f(x)|^2+k^{2}\|\nabla f (x)\|^2\right)^{1/2}\le k\nu_0(f,D),
\ee
therefore $\kappa(f_{k,z},D)\le k\kappa(f,D)$.

\begin{figure}[h!]
    \centering
    \includegraphics[scale=0.1]{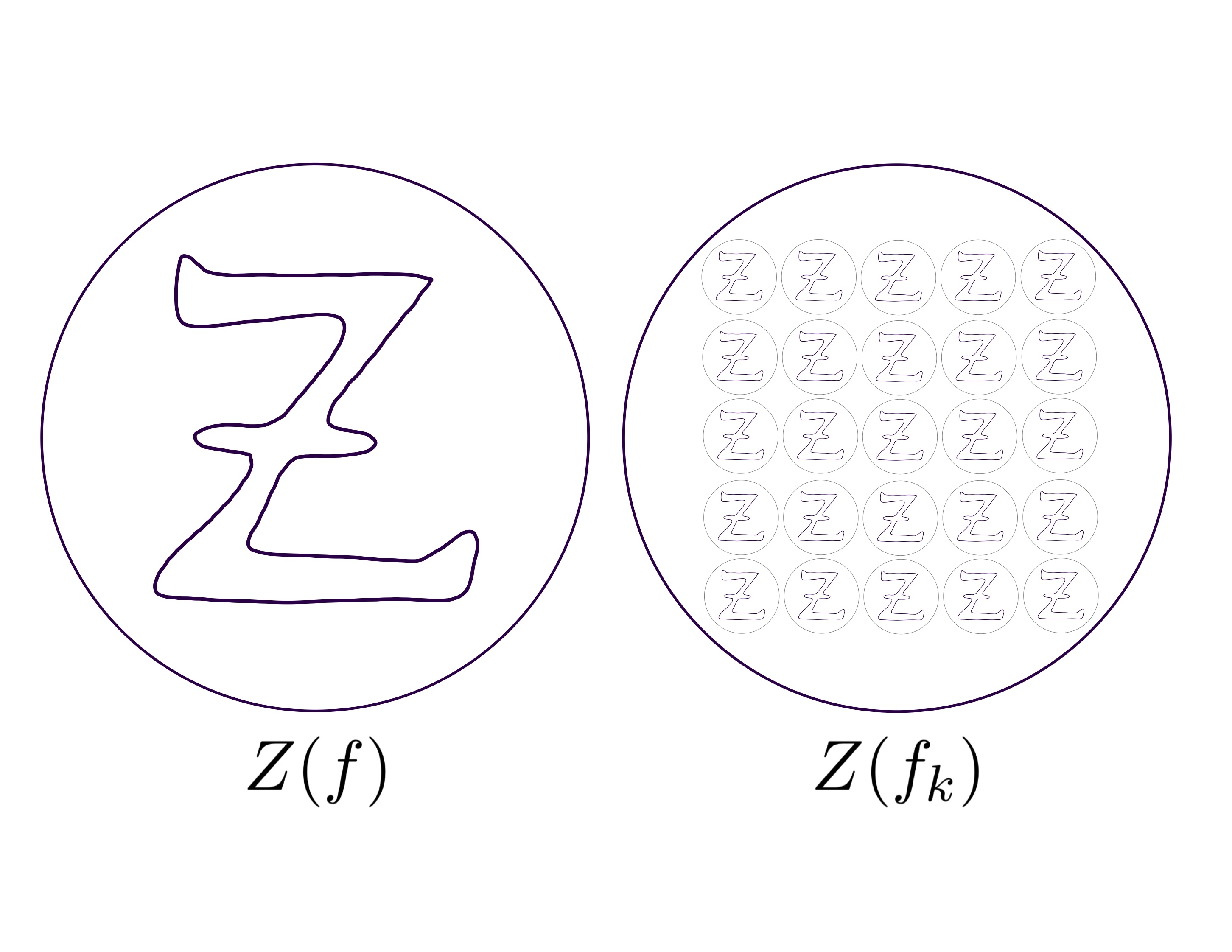}
    \caption{$Z(f_k)$ is a disjoint union of many copies of $Z=Z(f)$.}
    \label{semi:fig:zzz}
\end{figure}

For any $k\in\N$, choose a finite family $I_k$ of 
points $z_{k,i}\in D$, such that the disks $D_{k,i}$, centered in $z_{k,i}$ and with radius $k^{-1}$, are disjoint. Since the (Hausdorff) dimension of $D^n$ is $n$, we can assume that the number of points in such a family is $\#(I_k)\ge hk^n$, for some universal constant $h>0$ (it depends on $n$ only). 
Define the function $f_k\in \coo 1D\R$ as
\be 
f_k=1-\#(I_k)+\sum_{i\in I_k}{f_{k, z_i}},
\ee
so that $f_k$ coincides with $f_{k,i}$ on the disk $D_{k,i}$ and it is constantly equal to $1$ outside the union of all disks. 

Then we have $\kappa_1(f_k,D)\le k\kappa_1(f,D)$ and also $\kappa_1(f_k,D)$ is divergent, when $k\to +\infty$.
Moreover, the zero set of $f_k$ is homeomorphic to a disjoint union of $\#(I_k)$ copies of $Z$, hence
\be 
b(Z(f_k))\ge \#(I_k)b(Z)\ge hb(Z)k^n.
\ee
Putting these two observations together we conclude that 
\be 
b(Z(f_k))\ge hb(Z)\left(\frac{\kappa_1(f_k,D)}{\kappa_1(f,D)}\right)^n.
\ee
\end{proof}
%
%

%

%
%
%
%

%
%

%

%
%
%
%

%
%
%
%
%
%
%

\bibliography{RDT_Arxiv}

\begin{thebibliography}{10}

\bibitem{AdlerTaylor}
R.~J. Adler and J.~E. Taylor.
\newblock {\em Random fields and geometry}.
\newblock Springer Monographs in Mathematics. Springer, New York, 2007.

\bibitem{ambrofuscopalla}
L.~Ambrosio, N.~Fusco, and D.~Pallara.
\newblock {\em Functions of Bounded Variation and Free Discontinuity Problems}.
\newblock Oxford Science Publications. Clarendon Press, 2000.

\bibitem{arnold2012singularities}
V.I. Arnold, S.M. Gusein-Zade, and A.N. Varchenko.
\newblock {\em Singularities of Differentiable Maps, Volume 1: Classification
  of Critical Points, Caustics and Wave Fronts}.
\newblock Modern Birkh{\"a}user Classics. Birkh{\"a}user Boston, 2012.

\bibitem{Wschebor}
Jean-Marc Azais and Mario Wschebor.
\newblock {\em {Level sets and extrema of random processes and fields}}.
\newblock Wiley, Hoboken, NJ, 2009.

\bibitem{BBL}
T.~Bagby, L.~Bos, and N.~Levenberg.
\newblock Multivariate simultaneous approximation.
\newblock {\em Constr. Approx.}, 18(4):569--577, 2002.

\bibitem{Billingsley}
Patrick Billingsley.
\newblock {\em Convergence of probability measures}.
\newblock Wiley Series in Probability and Statistics: Probability and
  Statistics. John Wiley \& Sons, Inc., New York, second edition, 1999.
\newblock A Wiley-Interscience Publication.

\bibitem{bogachev}
V.I. Bogachev and American~Mathematical Society.
\newblock {\em Gaussian Measures}.
\newblock Mathematical surveys and monographs. American Mathematical Society,
  1998.

\bibitem{botttu}
R.~Bott and L.W. Tu.
\newblock {\em Differential Forms in Algebraic Topology}.
\newblock Graduate Texts in Mathematics. Springer New York, 1995.

\bibitem{bredon}
G.E. Bredon.
\newblock {\em Topology and Geometry}.
\newblock Graduate texts in mathematics. Springer-Verlag, 1993.

\bibitem{HaniehPaulAnto}
Paul Breiding, Hanieh Keneshlou, and Antonio Lerario.
\newblock {Quantitative Singularity Theory for Random Polynomials}.
\newblock {\em International Mathematics Research Notices}, 10 2020.
\newblock rnaa274.

\bibitem{buerg:07}
Peter B{\"u}rgisser.
\newblock Average {E}uler characteristic of random real algebraic varieties.
\newblock {\em C. R. Math. Acad. Sci. Paris}, 345(9):507--512, 2007.

\bibitem{CS}
Dustin Cartwright and Bernd Sturmfels.
\newblock The number of eigenvalues of a tensor.
\newblock {\em Linear Algebra Appl.}, 438(2):942--952, 2013.

\bibitem{chavel}
I.~Chavel.
\newblock {\em Riemannian Geometry: A Modern Introduction}.
\newblock Cambridge Studies in Advanced Mathematics. Cambridge University
  Press, 2006.

\bibitem{Erhan}
E.~{\c{C}}{\i}nlar.
\newblock {\em Probability and Stochastics}.
\newblock Graduate Texts in Mathematics. Springer New York, 2011.

\bibitem{delbruck}
M.~Delbruck.
\newblock Knotting problems in biology.
\newblock {\em Plant Genome Data and Information Center collection on
  computational molecular biology and genetics}, 1961.

\bibitem{DiattaLerario}
Daouda~Niang Diatta and Antonio Lerario.
\newblock Low degree approximation of random polynomials, 2018.

\bibitem{DFN}
B.~A. Dubrovin, A.~T. Fomenko, and S.~P. Novikov.
\newblock {\em Modern geometry---methods and applications. {P}art {II}}, volume
  104 of {\em Graduate Texts in Mathematics}.
\newblock Springer-Verlag, New York, 1985.
\newblock The geometry and topology of manifolds, Translated from the Russian
  by Robert G. Burns.

\bibitem{dudley}
R.~M. Dudley.
\newblock {\em Real Analysis and Probability}.
\newblock Cambridge Studies in Advanced Mathematics. Cambridge University
  Press, 2 edition, 2002.

\bibitem{EdelmanKostlan95}
Alan Edelman and Eric Kostlan.
\newblock How many zeros of a random polynomial are real?
\newblock {\em Bull. Amer. Math. Soc. (N.S.)}, 32(1):1--37, 1995.

\bibitem{eliash}
Y.~Eliashberg, N.M. Mishachev, and S.~Ariki.
\newblock {\em Introduction to the $h$-Principle}.
\newblock Graduate studies in mathematics. American Mathematical Society, 2002.

\bibitem{federer1996}
H.~Federer.
\newblock {\em Geometric Measure Theory}.
\newblock Grundlehren der mathematischen Wissenschaften. Springer, 1996.

\bibitem{Foote}
Robert~L. Foote.
\newblock Regularity of the distance function.
\newblock {\em Proc. Amer. Math. Soc.}, 92(1):153--155, 1984.

\bibitem{frischwasserman}
Harry~L. Frisch and Edel Wasserman.
\newblock Chemical topology.
\newblock {\em Journal of the American Chemical Society}, 83 (18):3789--3795,
  1961.

\bibitem{FyLeLu}
Yan~V. Fyodorov, Antonio Lerario, and Erik Lundberg.
\newblock On the number of connected components of random algebraic
  hypersurfaces.
\newblock {\em J. Geom. Phys.}, 95:1--20, 2015.

\bibitem{GaWe1}
Damien Gayet and Jean-Yves Welschinger.
\newblock Lower estimates for the expected {B}etti numbers of random real
  hypersurfaces.
\newblock {\em J. Lond. Math. Soc. (2)}, 90(1):105--120, 2014.

\bibitem{GaWe3}
Damien Gayet and Jean-Yves Welschinger.
\newblock Expected topology of random real algebraic submanifolds.
\newblock {\em J. Inst. Math. Jussieu}, 14(4):673--702, 2015.

\bibitem{GaWe2}
Damien Gayet and Jean-Yves Welschinger.
\newblock Betti numbers of random real hypersurfaces and determinants of random
  symmetric matrices.
\newblock {\em J. Eur. Math. Soc. (JEMS)}, 18(4):733--772, 2016.

\bibitem{GoreskyMacPherson}
M.~Goresky and R.~MacPherson.
\newblock {\em Stratified Morse Theory}.
\newblock Ergebnisse der Mathematik und ihrer Grenzgebiete. Springer-Verlag,
  1988.

\bibitem{Hirsch}
Morris~W. Hirsch.
\newblock {\em Differential topology}, volume~33 of {\em Graduate Texts in
  Mathematics}.
\newblock Springer-Verlag, New York, 1994.
\newblock Corrected reprint of the 1976 original.

\bibitem{Howard}
Ralph Howard.
\newblock The kinematic formula in {R}iemannian homogeneous spaces.
\newblock {\em Mem. Amer. Math. Soc.}, 106(509):vi+69, 1993.

\bibitem{ItoNisio}
Kiyosi Ito and Makiko Nisio.
\newblock On the convergence of sums of independent banach space valued random
  variables.
\newblock {\em Osaka J. Math.}, 5(1):35--48, 1968.

\bibitem{JordAngles}
Camille Jordan.
\newblock Essai sur la g\'eom\'etrie \`a $n$ dimensions.
\newblock {\em Bulletin de la Soci\'et\'e Math\'ematique de France},
  3:103--174, 1875.

\bibitem{kac43}
M.~Kac.
\newblock On the average number of real roots of a random algebraic equation.
\newblock {\em Bull. Amer. Math. Soc.}, 49:314--320, 1943.

\bibitem{Kostlan:93}
Eric Kostlan.
\newblock On the distribution of roots of random polynomials.
\newblock In {\em From {T}opology to {C}omputation: {P}roceedings of the
  {S}malefest ({B}erkeley, {CA}, 1990)}, pages 419--431. Springer, New York,
  1993.

\bibitem{khazOFRET}
K.~Kozhasov.
\newblock On fully real eigenconfigurations of tensors.
\newblock {\em SIAM Journal on Applied Algebra and Geometry}, 2(2):339--347,
  2018.

\bibitem{KP}
Steven~G. Krantz and Harold~R. Parks.
\newblock Distance to {$C^{k}$} hypersurfaces.
\newblock {\em J. Differential Equations}, 40(1):116--120, 1981.

\bibitem{counterSardKupka}
Ivan Kupka.
\newblock Counterexample to the morse-sard theorem in the case of
  infinite-dimensional manifolds.
\newblock {\em Proceedings of the American Mathematical Society},
  16(5):954--957, 1965.

\bibitem{LerarioJEMS}
A.~Lerario.
\newblock Complexity of intersections of real quadrics and topology of
  symmetric determinantal varieties.
\newblock {\em J. Eur. Math. Soc. (JEMS)}, 18(2):353--379, 2016.

\bibitem{Letwo}
Antonio Lerario.
\newblock Random matrices and the average topology of the intersection of two
  quadrics.
\newblock {\em Proc. Amer. Math. Soc.}, 143(8):3239--3251, 2015.

\bibitem{Lerarioshsp}
Antonio Lerario and Erik Lundberg.
\newblock Statistics on {H}ilbert's 16th problem.
\newblock {\em Int. Math. Res. Not. IMRN}, (12):4293--4321, 2015.

\bibitem{Lerariolemniscate}
Antonio Lerario and Erik Lundberg.
\newblock On the geometry of random lemniscates.
\newblock {\em Proc. Lond. Math. Soc. (3)}, 113(5):649--673, 2016.

\bibitem{mttrp}
Antonio Lerario and Michele Stecconi.
\newblock Maximal and typical topology of real polynomial singularities, 2019.

\bibitem{witdoash}
Antonio Lerario and Michele Stecconi.
\newblock What is the degree of a smooth hypersurface?, 2020.

\bibitem{dtgrf}
Antonio Lerario and Michele Stecconi.
\newblock Differential topology of gaussian random fields, 2021.

\bibitem{marinucci_peccati_2011}
Domenico Marinucci and Giovanni Peccati.
\newblock {\em Random Fields on the Sphere: Representation, Limit Theorems and
  Cosmological Applications}.
\newblock London Mathematical Society Lecture Note Series. Cambridge University
  Press, 2011.

\bibitem{MIAO199281_ang}
Jianming Miao and Adi Ben-Israel.
\newblock On principal angles between subspaces in rn.
\newblock {\em Linear Algebra and its Applications}, 171:81 -- 98, 1992.

\bibitem{MilnorBound}
J.~Milnor.
\newblock On the betti numbers of real varieties.
\newblock {\em Proc. Amer. Math. Soc}, pages 275--280, 1964.

\bibitem{milnor-stasheff}
John~W. Milnor and James~D. Stasheff.
\newblock {\em Characteristic classes}.
\newblock Princeton University Press, Princeton, N. J.; University of Tokyo
  Press, Tokyo, 1974.
\newblock Annals of Mathematics Studies, No. 76.

\bibitem{NazarovSodin2}
F.~Nazarov and M.~Sodin.
\newblock Asymptotic laws for the spatial distribution and the number of
  connected components of zero sets of {G}aussian random functions.
\newblock {\em Zh. Mat. Fiz. Anal. Geom.}, 12(3):205--278, 2016.

\bibitem{NazarovSodin1}
Fedor Nazarov and Mikhail Sodin.
\newblock On the number of nodal domains of random spherical harmonics.
\newblock {\em Amer. J. Math.}, 131(5):1337--1357, 2009.

\bibitem{Nicolaescu2016}
Liviu~I. Nicolaescu.
\newblock A stochastic gauss--bonnet--chern formula.
\newblock {\em Probability Theory and Related Fields}, 165(1):235--265, Jun
  2016.

\bibitem{park2013betti}
Changbom Park, Pratyush Pranav, Pravabati Chingangbam, Rien van~de Weygaert,
  Bernard Jones, Gert Vegter, Inkang Kim, Johan Hidding, and Wojciech~A.
  Hellwing.
\newblock Betti numbers of gaussian fields, 2013.

\bibitem{Parth}
K.R. Parthasarathy.
\newblock {\em Probability Measures on Metric Spaces}.
\newblock Ams Chelsea Publishing. Academic Press, 2005.

\bibitem{Po}
S.~S. Podkorytov.
\newblock On the {E}uler characteristic of a random algebraic hypersurface.
\newblock {\em Zap. Nauchn. Sem. S.-Peterburg. Otdel. Mat. Inst. Steklov.
  (POMI)}, 252(Geom. i Topol. 3):224--230, 252--253, 1998.

\bibitem{rice44}
S.~O. Rice.
\newblock Mathematical analysis of random noise.
\newblock {\em Bell System Technical Journal}, 23(3):282--332, 1944.

\bibitem{SarnakWigman}
Peter Sarnak and Igor Wigman.
\newblock Topologies of nodal sets of random band-limited functions.
\newblock {\em Communications on Pure and Applied Mathematics}, 72(2):275--342,
  2019.

\bibitem{schwartz1957}
L.~Schwartz.
\newblock {\em Th{\'e}orie des distributions}.
\newblock Number v. 1-2 in Actualit{\'e}s scientifiques et industrielles.
  Hermann, 1957.

\bibitem{shsm}
M.~Shub and S.~Smale.
\newblock Complexity of {B}ezout's theorem. {II}. {V}olumes and probabilities.
\newblock In {\em Computational algebraic geometry ({N}ice, 1992)}, volume 109
  of {\em Progr. Math.}, pages 267--285. Birkh\"auser Boston, Boston, MA, 1993.

\bibitem{smalesard}
S.~Smale.
\newblock An infinite dimensional version of sard's theorem.
\newblock {\em American Journal of Mathematics}, 87(4):861--866, 1965.

\bibitem{krstec}
Michele Stecconi.
\newblock Kac-rice formula for transverse intersections, 2021.

\bibitem{Wig}
I.~Wigman.
\newblock Fluctuations of the nodal length of random spherical harmonics.
\newblock {\em Commun. Math. Phys.}, 298(787), 2010.

\bibitem{wigEBNRF}
Igor Wigman.
\newblock On the expected betti numbers of the nodal set of random fields.
\newblock {\em Analysis \& pde}, 3 2020.

\bibitem{Yomdin}
Y.~Yomdin.
\newblock Global bounds for the {B}etti numbers of regular fibers of
  differentiable mappings.
\newblock {\em Topology}, 24(2):145--152, 1985.

\bibitem{Zhu_2013_ang}
P.~Zhu and A.V. Knyazev.
\newblock Angles between subspaces and their tangents.
\newblock {\em Journal of Numerical Mathematics}, 21(4), Jan 2013.

\end{thebibliography}

\end{document}